\documentclass[preprint,12pt, nonatbib]{elsarticle}
\usepackage{amssymb}
\usepackage{amsmath}
\usepackage{amsthm}
\usepackage{mathtools}
\usepackage{bm}
\usepackage{mathscinet}
\usepackage[colorlinks]{hyperref}
\AtBeginDocument{
  \hypersetup{
    linkcolor=blue,
    citecolor=green,
  }
}

\makeatletter
\let\c@author\relax
\makeatother
\makeatletter
\def\blx@maxline{77}
\makeatother
\usepackage[backend=bibtex,style=numeric,sorting=nyt,
giveninits=true,isbn=false,doi=false,url=true,maxbibnames=99,giveninits=false]{biblatex}
\renewbibmacro{in:}{\ifentrytype{article}{}{\printtext{\bibstring{in}\intitlepunct}}}
\DeclareFieldFormat*{title}{#1}
\usepackage{geometry}
\usepackage[T1]{fontenc}

\usepackage[mathlines,pagewise]{lineno}
\newcommand*\patchAmsMathEnvironmentForLineno[1]{%
  \expandafter\let\csname old#1\expandafter\endcsname\csname #1\endcsname
  \expandafter\let\csname oldend#1\expandafter\endcsname\csname end#1\endcsname
  \renewenvironment{#1}%
     {\linenomath\csname old#1\endcsname}%
     {\csname oldend#1\endcsname\endlinenomath}}%
\newcommand*\patchBothAmsMathEnvironmentsForLineno[1]{%
  \patchAmsMathEnvironmentForLineno{#1}%
  \patchAmsMathEnvironmentForLineno{#1*}}%
\AtBeginDocument{%
\patchBothAmsMathEnvironmentsForLineno{equation}%
\patchBothAmsMathEnvironmentsForLineno{align}%
\patchBothAmsMathEnvironmentsForLineno{flalign}%
\patchBothAmsMathEnvironmentsForLineno{alignat}%
\patchBothAmsMathEnvironmentsForLineno{gather}%
\patchBothAmsMathEnvironmentsForLineno{multline}%
}

\usepackage{aliascnt}
\usepackage{cleveref}

\newtheorem{theorem}{Theorem}[section]
\Crefname{theorem}{Theorem}{Theorems}

\newaliascnt{lemma}{theorem}      
\newtheorem{lemma}[lemma]{Lemma}  
\aliascntresetthe{lemma}          
\Crefname{lemma}{Lemma}{Lemmas}

\newaliascnt{corollary}{theorem}      
\newtheorem{corollary}[corollary]{Corollary}  
\aliascntresetthe{corollary}          
\Crefname{corollary}{Corollary}{Corollaries}

\newaliascnt{proposition}{theorem}      
\aliascntresetthe{proposition}          
\Crefname{proposition}{Proposition}{Propositions}

\theoremstyle{definition}
\newaliascnt{definition}{theorem}
\newtheorem{definition}[definition]{Definition}
\aliascntresetthe{definition}
\Crefname{definition}{Definition}{Definitions}

\newaliascnt{notation}{theorem}
\newtheorem{notation}[notation]{Notation}
\aliascntresetthe{notation}
\Crefname{notation}{Notation}{Notations}

\theoremstyle{remark}
\newaliascnt{remark}{theorem}
\newtheorem{remark}[remark]{Remark}
\aliascntresetthe{remark}
\Crefname{remark}{Remark}{Remarks}

\numberwithin{equation}{section}

\crefname{equation}{}{}
\crefname{figure}{}{}

\journal{~~}

\begin{document}
\begin{frontmatter}

\title{Pointwise Regularity for Fully Nonlinear Elliptic Equations in General Forms\tnoteref{t1}}
\author[rvt1]{Yuanyuan Lian}
\ead{lianyuanyuan.hthk@gmail.com}
\author[rvt2,rvt3]{Lihe Wang}
\ead{lihe-wang@uiowa.edu}
\author[rvt4]{Kai Zhang\corref{cor1}}
\ead{zhangkaizfz@gmail.com}
\address[rvt1]{Departamento de An\'{a}lisis Matem\'{a}tico, Instituto de Matem\'{a}ticas IMAG, Universidad de Granada}
\address[rvt2]{School of Mathematical Sciences, Shanghai Jiao Tong University, Shanghai, China}
\address[rvt3]{Department of Mathematics, The University of Iowa, Iowa City, IA 52242, USA}
\address[rvt4]{Departamento de Geometr\'{i}a y Topolog\'{i}a, Instituto de Matem\'{a}ticas IMAG, Universidad de Granada}

\cortext[cor1]{Corresponding author. ORCID: \href{https://orcid.org/0000-0002-1896-3206}{0000-0002-1896-3206}}
\tnotetext[t1]{This research has been financially supported by the National Natural Science Foundation of China (Grant No. 12031012, 11831003, 12171299 and 12171313), the Institute of Modern Analysis-A Frontier Research Center of Shanghai, Project PID2020-118137GB-I00 and PID2020-117868GB-I00 funded by MCIN/AEI /10.13039/501100011033.}

\begin{abstract}
In this manuscript, we develop systematically the pointwise regularity for $L^n$-viscosity solutions of fully nonlinear elliptic equations in general forms. In particular, the equations with quadratic growth (called natural growth) in the gradient are considered. We obtain a series of interior and boundary pointwise $C^{k,\alpha}$ regularity ($k\geq 1$ and $0<\alpha<1$). In addition, we also derive the pointwise $C^k$ regularity ($k\geq 1$) and $C^{k,\mathrm{lnL}}$ regularity ($k\geq 0$), which correspond to the endpoints $\alpha=0$ and $\alpha=1$ respectively. Some regularity results are new even for the linear equations. Moreover, the minimum requirements are imposed on the coefficients and the prescribed data to obtain the above regularity and the proofs are relatively simple.
\end{abstract}

\end{frontmatter}

\tableofcontents
\section{Introduction}\label{S1}
In this paper, we investigate the interior and boundary pointwise regularity for $L^n$-viscosity solutions (viscosity solutions for short) of fully nonlinear uniformly elliptic equations
\begin{equation}\label{FNE}
  F(D^2u,Du,u,x)=f\quad\mbox{in}~~\Omega
\end{equation}
and the corresponding Dirichlet problems
\begin{equation}\label{FNE2}
\left\{\begin{aligned}
&F(D^2u,Du,u,x)=f&& \quad\mbox{in}~~\Omega;\\
&u=g&& \quad\mbox{on}~~\Gamma\subset\partial \Omega.
\end{aligned}\right.
\end{equation}
respectively. Here, $\Omega\subset \mathbb{R}^n$ ($n\geq 2$) is a bounded domain and $F$ is a real fully nonlinear operator defined in
$\mathcal{S}^n\times \mathbb{R}^n\times \mathbb{R}\times \bar\Omega$, where $\mathcal{S}^n$ denotes the set of $n\times n$ symmetric matrices (see \Cref{no1.1}).

When studying fully nonlinear equations, some structure condition is necessary. We always assume the following structure condition for any operator $F$ throughout this paper: for any $M,N\in \mathcal{S}^n$,
$p,q\in \mathbb{R}^n$ and $s,t\in \mathbb{R}$ with $|s|,|t|\leq K$, we have for $a.e.~x\in \bar \Omega$,
\begin{equation}\label{SC2}
  \begin{aligned}
    &\mathcal{M}_{\lambda,\Lambda}^{-}(M-N)-\mu|p-q|(|p|+|q|)-b(x)|p-q|-c(x)\omega_0(K,|s-t|)\\
    &\leq F(M,p,s,x)-F(N,q,t,x)\leq \\
    &\mathcal{M}_{\lambda,\Lambda}^{+}(M-N)+\mu|p-q|(|p|+|q|)+b(x)|p-q|+c(x)\omega_0(K,|s-t|),\\
  \end{aligned}
\end{equation}
where $0<\lambda\leq \Lambda$, $\mu$ are nonnegative constants; $b,c$ are nonnegative functions; $\mathcal{M}_{\lambda,\Lambda}^{-},\mathcal{M}_{\lambda,\Lambda}^{+}$ denote the Pucci's operators (see \Cref{d-Sf}), and $\omega_0$ is a modulus of continuity depending on $K$, i.e., for any $K>0$, $\omega_0(K,\cdot)$ is a nonnegative non-decreasing function and $\omega_0(K,s)\rightarrow 0$ as $s\rightarrow 0$. In this paper, we assume that all functions are measurable in $x$. We will drop ``$a.e.$'' in the following arguments and use ``for any'' instead.

This structure condition allows equations to have quadratic growth in the gradient. The following are two typical examples:
\begin{equation}\label{e.linear-1}
  a^{ij}(x)u_{ij}(x)+\mu^{ij}(x)u_i(x)u_j(x)+b^i(x)u_i(x)+c(x)h(u)=f(x)
\end{equation}
and
\begin{equation}\label{e.nonlinear-1}
\mathcal{M}_{\lambda,\Lambda}^{+}(D^2u)+\mu|Du|^2+b(x)|Du|+c(x)h(u)=f(x),
\end{equation}
where the Einstein summation convention is used (similarly hereinafter), i.e., repeated indices mean summation. In addition, the structure condition \cref{SC2} allows $b,c\in L^p$ for some $p>n$ and is more general compared with previous structure conditions (see \cite{MR1376656,MR3980853,MR1606359,MR2486925} etc.).

Quadratic growth in the gradient is also called natural growth, which means that the equation is invariant under nonlinear transformation. For instance, if $u$ is a solution of \cref{FNE} and $v=T(u)$ where
$T\in C^2$ and $T'>0$, then $v$ is a solution of some equation satisfying the structure condition \cref{SC2}. Equations with quadratic growth in the gradient have been studied to some extent. Ladyzhenskaya and Ural'tseva obtained a priori estimates and proved the existence of smooth solutions for quasilinear elliptic equations (see \cite[Theorem III]{MR0149075} and \cite[Chapter 4 and Chapter 6]{MR0244627}). Trudinger \cite{MR701522} extended these results to fully nonlinear elliptic equations. For regularity theory, Sirakov \cite{MR2592289} proved the interior and boundary $C^{\alpha}$ regularity for viscosity solutions of fully nonlinear elliptic equations. Recently, Nornberg \cite{MR3980853} obtained the $C^{1,\alpha}$ regularity. For more work on the equations with quadratic growth in the gradient, we refer to the references in \cite{MR3980853} and \cite{MR2592289}.

We also consider equations under the following special structure condition:
\begin{equation}\label{SC1}
  \begin{aligned}
    &\mathcal{M}_{\lambda,\Lambda}^{-}(M-N)-b(x)|p-q|-c(x)|s-t|\\
    &\leq F(M,p,s,x)-F(N,q,t,x)\leq\\
    &\mathcal{M}_{\lambda,\Lambda}^{+}(M-N)+b(x)|p-q|+c(x)|s-t|.
  \end{aligned}
\end{equation}
In this case, \cref{FNE} is the natural generalization of the linear uniformly elliptic equations in nondivergence form
\begin{equation}\label{e.linear}
a^{ij}(x)u_{ij}(x)+b^{i}(x)u_i(x)+c(x)u(x)=f(x)\quad\mbox{in}~~ \Omega.
\end{equation}
With the special structure condition \cref{SC1}, we can obtain explicit estimates for solutions.

In this paper, we aim to develop the interior and boundary pointwise regularity systematically for viscosity solutions under structure conditions \cref{SC2} (or \cref{SC1}). We take the $C^{2,\alpha}$ regularity for example to clarify here several different types of regularity:
\begin{itemize}
  \item Interior pointwise regularity: $ x_0\in \Omega, \quad f\in C^{\alpha}(x_0)\Rightarrow u\in C^{2,\alpha}(x_0)$;
  \item Interior local regularity: $f\in C^{\alpha}(\bar{\Omega})\Rightarrow u\in C^{2,\alpha}(\bar{\Omega}'),~\forall \Omega'\subset\subset \Omega$;
  \item Boundary pointwise regularity:
  \begin{equation*}
   x_0\in \partial\Omega, \quad f\in C^{\alpha}(x_0),~g,\Gamma\in C^{2,\alpha}(x_0)\Rightarrow u\in C^{2,\alpha}(x_0);
  \end{equation*}
  \item Boundary local regularity:
  \begin{equation*}
  f\in C^{\alpha}(\bar\Omega), g\in C^{2,\alpha}(\bar{\Gamma}), \Gamma\in C^{2,\alpha}\Rightarrow u\in C^{2,\alpha}(\bar{\Omega}'),~\forall \Omega'\subset\subset \Omega\cup \Gamma;
  \end{equation*}
  \item Global regularity:
  \begin{equation*}
  f\in C^{\alpha}(\bar\Omega), g\in C^{2,\alpha}(\partial \Omega),
  \partial \Omega\in C^{2,\alpha} \Rightarrow u\in C^{2,\alpha}(\bar{\Omega}).
  \end{equation*}
\end{itemize}
It is well-known that the local and global regularity can be obtained directly from the pointwise regularity based on the equivalence between the classical definition of $C^{k,\alpha}$ space and the pointwise characterization (see \Cref{r-1.1}).

In this paper, we impose minimal requirements on the coefficients and the prescribed data to obtain various pointwise regularity. In other words, our results are optimal. In particular, part of results are new even for the linear equations. In addition, our proofs are relatively simple compared with previous results.

Essentially, the behavior of a solution near some point is determined by the coefficients and the prescribed data near the same point. The pointwise regularity shows clearly how these data influence the behavior of the solution.  Moreover, the assumptions for pointwise regularity could be weaker than that for local or global regularity. For example, for the boundary pointwise $C^{k,\alpha}$ regularity, the boundary may be not a graph of some function locally (see \Cref{r-12}), which is necessary for usual boundary regularity since one need to flatten the boundary by a transformation.

To study fully nonlinear elliptic equations, one may first assume that the solution is smooth and obtain a priori estimates. Based on a priori estimates, the existence of smooth solutions can be proved by the method of continuity. The benefit is that we can differentiate the equation directly. This method has its limitations. It usually relies on higher smoothness assumptions on the operator, the solution and the domain etc. (e.g. \cite[Chapter 9]{MR1351007} and \cite{MR701522}). In addition, it often brings the global (or local) estimates rather than pointwise estimates. Moreover, the proofs are relatively complicated compared with the proofs in this paper. In fact, this method is more appropriate for non-uniformly elliptic equations (see \cite{MR739925,MR806416,MR780073,MR1368245} etc.).

Another way to approach an equation is proving the existence of a solution in some weak sense first and then obtaining the regularity later. Viscosity solution is a kind of weak solution which is introduced by Crandall and Lions \cite{MR690039} (see also \cite{MR503721,MR597451}) and suitable for elliptic equations in nondivergence form, especially for fully nonlinear elliptic equations. The related theories of existence, uniqueness and regularity have been studied extensively (see \cite{MR1351007,MR1118699,MR2084272} and references therein). We also refer to \cite{MR3457602,MR4186265} for an exposition of the regularity theory.

Among various regularity results for fully nonlinear elliptic equations, pointwise regularity occupies an important position. Caffarelli \cite{MR1005611} (see also \cite{MR1351007}) proved the interior pointwise $C^{1,\alpha}$ and $C^{2,\alpha}$ regularity. Kovats \cite{MR1713596} obtained the pointwise $C^{2}$ regularity. Teixeira \cite{MR3158810} derived the pointwise $C^{0,\mathrm{lnL}}$ and $C^{1,\mathrm{lnL}}$ regularity. For equations with lower terms, Savin \cite{MR2334822} proved the interior pointwise $C^{2,\alpha}$ regularity for small solutions without the usual assumption that $F$ is convex or concave in $M$. For boundary pointwise regularity, Silvestre and Sirakov \cite{MR3246039} proved the $C^{1,\alpha}$ and $C^{2,\alpha}$ regularity on flat boundaries for equations depending the gradient. Lian and Zhang \cite{MR4088470} obtained the pointwise $C^{1,\alpha}$ and $C^{2,\alpha}$ regularity on general boundaries. We point out that the pointwise regularity also attract a lot of attention for other types of equation, such as dengenerate equations (see \cite{MR4523464}), parabolic equations (see \cite{MR1135923,MR1139064,MR1151267}) and the Monge-Amp\`{e}re equation (see \cite{MR2983006}) etc.

For equations with quadratic growth in the gradient, Sirakov \cite{MR2592289} proved the interior and boundary pointwise $C^{\alpha}$ regularity. The interior pointwise $C^{1,\alpha}$ regularity and boundary pointwise $C^{1,\alpha}$ regularity on flat domains were obtained by Nornberg \cite{MR3980853}. Recently, da Silva and Nornberg \cite{MR4304555} considered equations with more general nonlinear growth in the gradient. They obtained local $C^{k,\alpha}$($k=0,1,2$) and $C^{k,\mathrm{lnL}}$ ($k=0,1$) regularity, as well as some regularity in Sobolev and $BMO$ spaces.

The perturbation and compactness techniques are used in this paper. The perturbation technique is motivated originally by \cite{MR1351007} and the application to boundary regularity is inspired by \cite{MR3780142} and \cite{MR4088470}.  The compactness technique has been inspired by \cite{MR3246039} and \cite{Wang_Regularity}. As stated in \cite[P. 17]{Wang_Regularity}, the advantage of compactness technique is that we do not need to solve an equation and use its solution to approximate the original solution. In fact, our proofs in this paper do not rely on any solvability.

Next, we explain briefly the key idea used in this paper. Consider the following linear equations for example
\begin{equation}\label{e.linear-2}
\left\{\begin{aligned}
&a^{ij}u_{ij}+b^{i}u_i+cu=f&&\quad\mbox{in}~~ \Omega\cap B_1;\\
&u=g&& \quad\mbox{on}~~\partial \Omega\cap B_1,
\end{aligned}\right.
\end{equation}
where $0\in \partial \Omega$ and we study the regularity at $0$.

The main idea is perturbation, which can be tracked at least to \cite{MR1005611}. Roughly speaking, if we have enough regularity for harmonic functions, the regularity for \cref{e.linear-2} can be obtained by a perturbation argument. In the usual perturbation technique, the coefficients $a^{ij}$ and the right-hand term $f$ are regarded as a perturbation of a constant matrix and $0$. In this paper, we move one step forward and regard the coefficients $b^i,c_i$, the boundary value $g$ and the curved boundary $\partial \Omega\cap B_1$ as the perturbation of $0,0,0$ and a hyperplane respectively.

More precisely, take the boundary $C^{1,\alpha}$ regularity for instance. The proof contains mainly two steps. First, if \cref{e.linear-2} is quite close to
\begin{equation}\label{e.2.2}
\left\{\begin{aligned}
&\Delta u=0&&\quad\mbox{in}~~B_1^+;\\
&u=0&& \quad\mbox{on}~~T_1,
\end{aligned}\right.
\end{equation}
then the solution can be approximated by a linear polynomial in $\Omega\cap B_{\eta}$ for some $0<\eta<1$. For example, the closedness can be measured by
\begin{equation*}
\max\left(\|a^{ij}-\delta_{ij}\|_{L^{\infty}}, \|b\|_{L^{\infty}}, \|c\|_{L^{\infty}}, \|f\|_{L^{\infty}}, \|g\|_{L^{\infty}},\mathrm{osc}_{B_1}\partial \Omega\right) \leq \delta,
\end{equation*}
where $0<\delta<1$ is a small constant.

This step can be proved by the method of compactness. Indeed, if the conclusion is false, we will have a sequence of solutions to the problems in the form of \cref{e.linear-2} whose coefficients and prescribed data converge to that of \cref{e.2.2}. If this sequence of solutions are compact (e.g. by the uniform H\"{o}lder continuity \Cref{l-3Ho}), there exists a subsequence of solutions converging to some function $\bar u$. Combining with the closedness result (e.g. by \Cref{l-35}), $\bar u$ is a solution of \cref{e.2.2}. Then $\bar u$ can be approximated by a linear polynomial, which will lead to a contradiction.

The second step is a scaling argument, i.e. a sequence of repetitions of the first step. By a scaling argument, we have a sequence of estimates in $\Omega\cap B_{\eta^m} (m\geq 1)$, which implies the boundary $C^{1,\alpha}$ regularity. The scaling invariance of equations is the key to this step.

We point out that Nornberg \cite{MR3980853} also regarded $b^i,c$ as perturbation for $C^{1,\alpha}$ regularity. We deal with this in a more delicate way especially for higher regularity. Silvestre and Sirakov \cite{MR3246039} also obtained the boundary pointwise $C^{2,\alpha}$ regularity for flat boundaries. However, they first established the regularity for the equation $\Delta u+b_0^iu_i=0$ where $b_0$ is constant vector. Then they proved $C^{2,\alpha}$ regularity for general equations $a^{ij} u_{ij}+b^i(x)u_i=f$ by a perturbation argument. On the contrast, we obtain the pointwise $C^{2,\alpha}$ regularity directly based on the regularity of $\Delta u=0$. Finally, we remark here that some proofs presented in this paper could possibly be extended to models governed by degenerate/singular operators.

This paper is organized as follows. In \Cref{No}, we introduce some basic notions, including the pointwise characterization of smoothness of functions and domains. \Cref{P1} is devoted to prepare some preliminary results, such as the regularity for model equations (see \Crefrange{l-3modin1}{l-32}), the compactness and the closedness for a family of viscosity solutions.

The pointwise regularity will be stated and proved in the subsequent sections, along with some comments. In particular, we will compare our results with related previous ones. We give the interior $C^{1,\alpha}$ regularity, $C^{2,\alpha}$ regularity and $C^{k,\alpha}$ regularity ($k\geq 3$) in \Crefrange{In-C1a-mu}{In-Cka-mu} respectively. The corresponding boundary regularity are proved in \Cref{C1a-mu} to \Cref{Cka-mu}. Since the proofs of $C^k$ regularity and $C^{k,\mathrm{lnL}}$ regularity are similar to that of the $C^{k,\alpha}$ regularity, we only give proofs of several results. Precisely, in \Cref{C1C2Ck} we provide proofs of the interior $C^{2}$ regularity and boundary $C^1$ regularity. In \Cref{CLlnL}, we derive the interior $C^{k,\mathrm{lnL}}$ ($k\geq 2$) and boundary $C^{1,\mathrm{lnL}}$ regularity.

Symbols frequently used in this paper are listed below.
\begin{notation}\label{no1.1}
\begin{enumerate}~~\\
\item $\{e_i\}^{n}_{i=1}\colon$ the standard basis of $\mathbb{R}^n$, i.e., $e_i=(0,...0,\underset{i^{th}}{1},0,...0)$.
\item Given $x\in \mathbb{R}^n$, we write
\begin{equation*}
x=(x_1,...,x_n)=(x',x_n), \quad x'=(x_1,x_2,...,x_{n-1}).
\end{equation*}
\item $|x|\coloneqq\left(\sum_{i=1}^{n} x_i^2\right)^{1/2}$ for $x\in \mathbb{R}^n$.
\item $\mathbb{R}^n_+\coloneqq\{x\in \mathbb{R}^n\big|x_n>0\}$.
\item $B_r(x_0)\coloneqq B(x_0,r)\coloneqq\{x\in \mathbb{R}^{n}\big| |x-x_0|<r\}$, $B_r\coloneqq B_r(0)$, $B_r^+(x_0)\coloneqq B_r(x_0)\cap \mathbb{R}^n_+$ and $B_r^+\coloneqq B^+_r(0)$.
\item $T_r(x_0)\ \coloneqq \{(x',0)\in \mathbb{R}^{n}\big| |x'-x_0'|<r\}$, $T_r\coloneqq T_r(0)$.
\item $\mathcal{S}^{n}\colon$ the set of $n\times n$ symmetric matrices. For any $M\in \mathcal{S}^{n}$, we use $M_{ij}$ or $M^{ij}$ to denote the entry in the $i$-th row and $j$-th column, $1\leq i,j\leq n$.
\item $|M|\coloneqq$ the spectral radius of $M$; $tr(M)\coloneqq\sum_{i=1}^{n}M_{ii}$, the trace of $M$ for any $M\in \mathcal{S}^{n}$.
\item $I\coloneqq \delta_{ij}\colon$ the unit matrix in $\mathcal{S}^n$; $\tilde I\colon$ the matrix whose entries are all $0$ except $\tilde{I}_{nn}=1$.
\item Let $\Omega,\Omega'\subset \mathbb{R}^n$. Define $\Omega^c\colon$ the complement of $\Omega,\quad$  $\bar \Omega \colon$ the closure of $\Omega$ and we call $\Omega'\subset\subset \Omega$ if $\bar{\Omega}'\subset \Omega$.
\item $\mathrm{diam}(\Omega)\colon$ the diameter of $\Omega$; $\mathrm{dist}(\Omega_1,\Omega_2)\colon$ the distance between $\Omega_1$ and $\Omega_2$, where $ \Omega_1,\Omega_2\subset \mathbb{R}^n$.
\item $\Omega_r\coloneqq \Omega\cap B_r$, $(\partial\Omega)_r\coloneqq \partial\Omega\cap B_r$.
\item $a^+\coloneqq \max(a,0)$, the positive part of $a$; $a^-\coloneqq \max(-a,0)$, the negative part of $a$ for $a\in \mathbb{R}$.
\item Given a function $\varphi\colon \mathbb{R}^n\to \mathbb{R}$, define $\varphi _i\coloneqq \partial \varphi/\partial x _{i}$, $\varphi _{ij}\coloneqq \partial ^{2}\varphi/\partial x_{i}\partial x_{j}$ etc.
\item $D^0\varphi\coloneqq \varphi$, $D \varphi\coloneqq (\varphi_1 ,...,\varphi_{n} )$ and $D^2 \varphi \coloneqq \left(\varphi _{ij}\right)_{n\times n}$ etc. In addition, $D_{x'} \varphi\coloneqq (\varphi_1 ,...,\varphi_{n-1} )$ and $D_{x'}^2 \varphi\coloneqq \left(\varphi _{ij}\right)_{(n-1)\times (n-1)}$ etc.
\item We also use the standard multi-index notation. Let $\sigma=(\sigma_1,...,\sigma_n)\in \mathbb{N}^n$, i.e., each component $\sigma_i$ is a nonnegative integer. Define
\begin{equation*}
|\sigma|\coloneqq \sum_{i=1}^{n}\sigma_i,\quad\sigma!\coloneqq \prod_{i=1}^{n}(\sigma_i!),\quad
x^{\sigma}\coloneqq \prod_{i=1}^{n} x_i^{\sigma_i},\quad
D^{\sigma}\varphi \coloneqq \frac{\partial^{|\sigma|} \varphi }{\partial x_1^{\sigma_1}\cdots \partial x_n^{\sigma_n}}.
\end{equation*}
\item $|D^k\varphi |\coloneqq \left(\sum_{|\sigma|=k}|D^{\sigma}\varphi| ^2\right)^{1/2}$ for $k\geq 0$.
\item Given $F\colon \mathcal{S}^n\times \mathbb{R}^n\times \mathbb{R}\times \Omega\to \mathbb{R}$, define
\begin{equation*}
F_{M_{ij}}\coloneqq \frac{\partial F}{\partial M_{ij}},\quad F_{p_i}\coloneqq \frac{\partial F}{\partial p_{i}},\quad F_{s}\coloneqq \frac{\partial F}{\partial s},\quad F_{x_i}\coloneqq \frac{\partial F}{\partial x_{i}},
\end{equation*}
where $1\leq i,j\leq n$. Moreover, let $\xi\in \mathbb{N}^{n\times n}$ denote the matrix-valued multi-index. Then define
\begin{equation*}
D^{\xi}_MF\coloneqq \frac{\partial^{|\xi|} F}{\partial M_{ij}^{\xi_{ij}}},\quad
D^k_MF\coloneqq\left\{\frac{\partial^k F}{\partial M^{\xi}}\colon  |\xi|=k\right\},\quad
|D^k_MF|\coloneqq \left(\sum_{|\xi|=k}\left|\frac{\partial^k F}{\partial M^{\xi}}\right| ^2\right)^{1/2}.
\end{equation*}
Similarly, we can define $D^k_{p}F$, $D^k_{s}F$ and $D^k_{x}F$ etc.
\item $\mathcal{P}_k (k\geq 0)\colon$ the set of polynomials of degree less than or equal to $k$. Any $P\in \mathcal{P}_k$ can be written as
\begin{equation*}
P(x)=\sum_{|\sigma|\leq k}\frac{a_{\sigma}}{\sigma!}x^{\sigma},
\end{equation*}
where $a_{\sigma}$ are constants. Define
\begin{equation*}
\|P\|\coloneqq \sum_{|\sigma|\leq k}|a_{\sigma}|.
\end{equation*}
\item $\mathcal{HP}_k (k\geq 0)\colon$ the set of homogeneous polynomials of degree $k$. Any $P\in \mathcal{HP}_k$ can be written as
\begin{equation*}
P(x)=\sum_{|\sigma|= k}\frac{a_{\sigma}}{\sigma!}x^{\sigma}.
\end{equation*}

\item $\mathcal{SP}_k (k\geq 1)\colon$ the set of homogeneous polynomials of degree $k$ in a special form, i.e., $P\in \mathcal{SP}_k$ if and only if $P\in \mathcal{HP}_k$ can be written as
\begin{equation*}
P(x)=\sum_{|\sigma|=k,\sigma_n\geq 1}\frac{a_{\sigma}}{\sigma!}x^{\sigma}.
\end{equation*}

\end{enumerate}
\end{notation}
~\\

\section{Notions and terminology}\label{No}
In this section, we introduce some notions and terminology. They will make the statements of our results and the proofs concise and readable, which is one of our goals. Our paper treats the regularity of solutions in H\"{o}lder spaces and the following is the classical definition of H\"{o}lder spaces. Recall that a function $\omega:\mathbb{R}^+\to \mathbb{R}^+$ is called a modulus of continuity if $\omega$ is non-decreasing and $\omega(r)\to 0$ as $r\to 0$.
\begin{definition}\label{de2.1}
Let $k\geq 0$, $\Omega\subset \mathbb{R}^n$ be a bounded domain and $\omega$ be a modulus of continuity. We say that $f\in C^{k,\omega}(\bar\Omega)$ if $f$ has continuous derivatives up to order $k$ and for some $K>0$,
\begin{equation}\label{e2.2}
|D^kf(x)-D^kf(y)|\leq K\omega(|x-y|),~\forall ~x,y\in \Omega.
\end{equation}
Moreover, $C^{k,\omega}(\bar\Omega)$ is endowed with the semi-norm
\begin{equation*}
[f]_{C^{k,\omega}(\bar{\Omega})}\coloneqq \min \left\{K\omega(r_0) \big | \cref{e2.2} ~\mbox{holds with}~K\right\}, \quad r_0\coloneqq\mathrm{diam}(\Omega)
\end{equation*}
and norm
\begin{equation*}
\|f\|_{C^{k,\omega}(\bar{\Omega})}
\coloneqq \|f\|_{C^{k}(\bar{\Omega})}+[f]_{C^{k,\omega}(\bar{\Omega})}
\coloneqq \sum_{i=0}^{k} \|D^{i}f\|_{L^{\infty}(\bar\Omega)}+[f]_{C^{k,\omega}(\bar{\Omega})}.
\end{equation*}

If $\omega$ is a Dini function, i.e.,
\begin{equation}\label{e.1.Dini}
I_{\omega}\coloneqq \int_{0}^{r_0}\frac{\omega(r)}{r} dr<\infty,
\end{equation}
we say that $f\in C^{k,\mathrm{Dini}}(\bar{\Omega})$. Then we set
\begin{equation*}
[f]_{C^{k,\mathrm{Dini}}(\bar{\Omega})}\coloneqq \min \left\{KJ_{\omega}\big | \cref{e2.2} ~\mbox{holds with}~K\right\}, \quad J_{\omega}\coloneqq I_{\omega}+\omega(r_0)
\end{equation*}
and
\begin{equation}\label{e2.3}
\|f\|_{C^{k, \mathrm{Dini}}(\bar{\Omega})}\coloneqq \|f\|_{C^{k}(\bar{\Omega})}+[f]_{C^{k, \mathrm{Dini}}(\bar{\Omega})}.
\end{equation}

Furthermore, if $\omega(r)=r^{\alpha}$ for some $0<\alpha\leq 1$, we say $f\in C^{k,\alpha}(\bar{\Omega})$ and use
$ \|f\|_{C^{k,\alpha}(\bar{\Omega})}$ to denote the corresponding norm in \cref{e2.3} with
\begin{equation*}
[f]_{C^{k,\alpha}(\bar{\Omega})}\coloneqq \min \left\{K \big | \cref{e2.2} ~\mbox{holds with}~K\right\}.
\end{equation*}

Similarly, if $\omega(r)=r\big|\ln\min(r,1/2)\big|$, we call $f\in C^{k,\mathrm{lnL}}(\bar{\Omega})$ and denote the norm in \cref{e2.3} by
$\|f\|_{C^{k,\mathrm{lnL}}(\bar{\Omega})}$ along with
\begin{equation*}
[f]_{C^{k,\mathrm{lnL}}(\bar{\Omega})}\coloneqq \min \left\{K \big | \cref{e2.2} ~\mbox{holds with}~K\right\}.
\end{equation*}
\end{definition}
~\\

In this paper, we mainly consider the pointwise regularity and use the following definition of pointwise $C^{k,\omega}$ for a function, which is first introduced by Campanato \cite{MR156188, MR167862}.
\begin{definition}\label{d-f}
Let $k\geq 0$, $\Omega\subset \mathbb{R}^n$ be a bounded set (may be not a domain) and $\omega$ be a modulus of continuity. We say that $f$ is $C^{k,\omega}$ at $x_0\in \Omega$ or $f\in C^{k, \omega}(x_0)$ if there exist $P\in \mathcal{P}_k$ (space of polynomials of degree $k$, see \Cref{no1.1}) and constants $K,r_0>0$ such that
\begin{equation}\label{m-holder}
  |f(x)-P(x)|\leq K|x-x_0|^{k}\omega(|x-x_0|),~\forall~x\in \Omega\cap B_{r_0}(x_0).
\end{equation}
Then we call $P$ the Taylor polynomial of $f$ at $x_0$ and define
\begin{equation*}
D^mf(x_0)\coloneqq D^mP(x_0),\quad\|f\|_{C^{k}(x_0)}\coloneqq \sum_{m=0}^{k}|D^m P(x_0)|
\end{equation*}
and
\begin{equation*}
[f]_{C^{k,\omega}(x_0)}\coloneqq \min \left\{K\omega(r_0)\big | \cref{m-holder} ~\mbox{holds with}~K \right\},~~
\|f\|_{C^{k,\omega}(x_0)}\coloneqq \|f\|_{C^{k}(x_0)}+[f]_{C^{k,\omega}(x_0)}.
\end{equation*}
If $f\in C^{k,\omega}(x)$ for any $x\in \Omega$ with the same $\omega, r_0$ and
\begin{equation*}
  \|f\|_{C^{k,\omega}(\bar{\Omega})}\coloneqq  \sup_{x\in \Omega} \|f\|_{C^{k,\omega}(x)}<+\infty,
\end{equation*}
we say that $f\in C^{k,\omega}(\bar{\Omega})$.

If $\omega$ is a Dini function, (i.e. \cref{e.1.Dini} holds), we say that $f\in C^{k,\mathrm{Dini}}(x_0)$. Then we define
\begin{equation*}
[f]_{C^{k,\mathrm{Dini}}(x_0)}=\min \left\{KJ_{\omega} \big | \cref{m-holder} ~\mbox{holds with}~K\right\}
\end{equation*}
and
\begin{equation*}
\|f\|_{C^{k, \mathrm{Dini}}(x_0)}=\|f\|_{C^{k}(x_0)}+[f]_{C^{k, \mathrm{Dini}}(x_0)}.
\end{equation*}
If $f\in C^{k,\mathrm{Dini}}(x)$ for any $x\in \Omega$ with the same $\omega,r_0$ and
\begin{equation*}
  \|f\|_{C^{k,\mathrm{Dini}}(\bar{\Omega})}\coloneqq  \sup_{x\in \Omega} \|f\|_{C^{k}(x)}+\sup_{x\in \Omega} [f]_{C^{k, \mathrm{Dini}}(x)}<+\infty,
\end{equation*}
we say that $f\in C^{k,\mathrm{Dini}}(\bar{\Omega})$.

Similarly, if $\omega(r)=r^{\alpha}$ ($0<\alpha\leq 1$) (resp. $\omega(r)=r|\ln \min(r,1/2)|$), we can define
$f\in C^{k, \alpha}(x_0)$ (resp. $f\in C^{k,\mathrm{lnL}}(x_0)$) with
\begin{equation*}
[f]_{C^{k,\alpha}(x_0)}=\min \left\{K\big | \cref{m-holder} ~\mbox{holds with}~K\right\},\quad
\|f\|_{C^{k, \alpha}(x_0)}=\|f\|_{C^{k}(x_0)}+[f]_{C^{k, \alpha}(x_0)}
\end{equation*}
(resp.
\begin{equation*}
[f]_{C^{k,\mathrm{lnL}}(x_0)}=\min \left\{K \big | \cref{m-holder} ~\mbox{holds with}~K\right\},\quad
\|f\|_{C^{k, \mathrm{lnL}}(x_0)}=\|f\|_{C^{k}(x_0)}+[f]_{C^{k, \mathrm{lnL}}(x_0)}).
\end{equation*}
Furthermore, we define $f\in C^{k,\alpha}(\bar{\Omega})$ (resp. $f\in C^{k,\mathrm{lnL}}(\bar{\Omega})$) in a similar way.

If $\omega$ is only a modulus of continuity rather than a Dini function,
we may simply say that $f\in C^{k}(x_0)$ instead of $f\in C^{k,\omega}(x_0)$.
\end{definition}

\begin{remark}\label{r-1.1}
If $\Omega$ is a bounded Lipschitz domain, the definitions of $C^{k,\alpha}(\bar{\Omega})$ etc. in \Cref{de2.1} and \Cref{d-f} are equivalent (see \cite{MR167862}, \cite[Step 4 in the proof Theorem 2.1]{MR4156124} and \cite{MR1713596}). In this paper, we mainly treat the cases that $\Omega$ is a bounded domain or $\Omega$ is part of the boundary of a domain.
\end{remark}
~\\

Next, we define some other types of continuity.
\begin{definition}\label{d-2}
Let $\Omega$, $f$ and $\omega$ be as in \Cref{d-f}. We say that $f$ is $C^{-1,\omega}$ at $x_0$ or $f\in C^{-1,\omega}(x_0)$ if there exist $K,r_0>0$ such that
\begin{equation}\label{e.c-1}
\|f\|_{L^n(\bar{\Omega}\cap B_r(x_0) )}\leq K\omega(r), ~\forall ~0<r<r_0.
\end{equation}

If $\omega$ is a Dini function, we say that $f\in C^{-1,\mathrm{Dini}}(x_0)$ and define
\begin{equation*}
\|f\|_{C^{-1,\mathrm{Dini}}(x_0)}=\min \left\{KJ_{\omega} \big | \cref{e.c-1} ~\mbox{holds with}~\omega\right\}.
\end{equation*}
If $f\in C^{-1,\mathrm{Dini}}(x)$ for any $x\in \Omega$ with the same $r_0$ and
\begin{equation*}
\|f\|_{C^{-1,\mathrm{Dini}}(\bar{\Omega})}\coloneqq  \sup_{x\in \Omega} \|f\|_{C^{-1,\mathrm{Dini}}(x)}<+\infty,
\end{equation*}
we say that $f\in C^{-1,\mathrm{Dini}}(\bar{\Omega})$.

Finally, we can define $f\in C^{-1, \alpha}(x_0)$ and $f\in C^{-1,\alpha}(\bar{\Omega})$ ($0<\alpha\leq 1$) similarly to the previous.
\end{definition}

Since the boundary pointwise regularity is also considered, we give the definitions of the pointwise geometric conditions on the domain. Usually, if we say that $\partial \Omega$ is $C^{k,\alpha}$ near $x_0\in \partial \Omega$, it means that $\partial \Omega\cap B_{r_0}(x_0)$ (for some $r_0>0$) can be represented as a graph of a $C^{k,\alpha}$ function. Here, we use a more general pointwise definition for the smoothness of the boundary; a notion borrowed from \cite{MR4088470}.

\begin{definition}\label{d-re} Let $\Omega$ be a bounded domain, $x_0\in \Gamma\subset \partial \Omega$ and $\omega$ be a modulus of continuity. We say that $\Gamma$ is $C^{k,\omega}$ ($k\geq 1$) at $x_0$ or $\Gamma\in C^{k,\omega}(x_0)$ if there exist constants $K,r_0>0$, a coordinate system $\{x_1,...,x_n \}$ (isometric to the original coordinate system) and $P(x')\in \mathcal{P}_k$ with $P(0)=0$ and $DP(0)=0$ such that $x_0=0$ in this coordinate system,
\begin{equation}\label{e-re}
B_{r_0} \cap \{(x',x_n)\big |x_n>P(x')+|x'|^k\omega(|x'|)\} \subset B_{r_0}\cap \Omega
\end{equation}
and
\begin{equation}\label{e-re2}
B_{r_0} \cap \{(x',x_n)\big |x_n<P(x')-|x'|^k\omega(|x'|)\} \subset B_{r_0}\cap \Omega^c.
\end{equation}
Then we call $P$ the Taylor polynomial of $\Gamma$ at $x_0$ and define
\begin{equation*}
\|\Gamma\|_{C^{k}(x_0)}=\sum_{m=0}^{k}|D^m P(0)|.
\end{equation*}

If $\omega$ is a Dini function, we say that $\Gamma\in C^{k,\mathrm{Dini}}(x_0)$ and define
\begin{equation*}
[\Gamma]_{C^{k,\mathrm{Dini}}(x_0)}=\min \left\{KJ_{\omega}\big | \cref{e-re} \mbox{ and } \cref{e-re2} ~\mbox{hold with}~\omega\right\}
\end{equation*}
and
\begin{equation*}
\|\Gamma\|_{C^{k, \mathrm{Dini}}(x_0)}=\|\Gamma\|_{C^{k}(x_0)}+[\Gamma]_{C^{k, \mathrm{Dini}}(x_0)}.
\end{equation*}
If for any $\Gamma'\subset\subset \Gamma$, there exist $r'>0$ such that $\Gamma\in C^{k, \mathrm{Dini}}(x)$ with $r'$ for any $x\in \Gamma'$ and
\begin{equation*}
  \|\Gamma'\|_{C^{k,\mathrm{Dini}}}\coloneqq  \sup_{x\in \Gamma'} \|\Gamma\|_{C^{k}(x)}+\sup_{x\in \Gamma'}~[\Gamma]_{C^{k, \mathrm{Dini}}(x)}<+\infty,
\end{equation*}
we say that $\Gamma\in C^{k,\mathrm{Dini}}$.

As before, we can define $\Gamma\in C^{k,\alpha}(x_0)$ and $\Gamma\in C^{k,\alpha}$ ($0<\alpha\leq 1$) similarly. If $\omega$ is only a modulus of continuity rather than a Dini function,
we may simply say that $\Gamma\in C^{k}(x_0)$ and $\Gamma\in C^{k}$.
\end{definition}

\begin{remark}\label{r-12}
One feature of this definition is that $\partial \Omega$ may not be represented as a graph of a function near $x_0$. For example, let
\begin{equation*}
  \Omega=B_1\cap \left\{(x',x_n)\big | x_n>\frac{1}{2}|x'|^2\right\}\backslash \left\{(x',x_n)\big | x_n=\frac{1}{2}|x'|^2+|x'|^4, |x|\leq \frac{1}{2} \right\}.
\end{equation*}
Then $\partial \Omega$ is $C^{3,\alpha}$ at $0$ for any $0<\alpha\leq 1$ by the definition.
\end{remark}

\begin{remark}\label{re1.5}
Throughout this paper, if we assume that $f\in C^{k,\omega}(x_0)$ ($\Gamma\in C^{k,\omega}(x_0)$), we will use $P_f$ ($P_{\Omega}$) by default to denote its Taylor polynomial in \Cref{d-f} (\Cref{d-re}).

In addition, if we assume that $f\in C^{k,\mathrm{Dini}}$ ($k\geq -1$) ($\Gamma\in C^{k,\mathrm{Dini}}(x_0)$ ($k\geq 1$)), we always use $\omega_f$ ($\omega_{\Omega}$) to denote its corresponding Dini function.
\end{remark}
~\\

%
%

In the following, we introduce some notions with respect to $L^p$-viscosity solutions, which are standard (see \cite{MR1376656}, \cite{MR1351007} and \cite{MR1118699}).
\begin{definition}\label{d-viscoF}
We say that $u\in C(\Omega)$ is an $L^p$-viscosity ($p>n/2$) subsolution (resp., supersolution) of \cref{FNE} if
\begin{equation*}
  \begin{aligned}
    ess~\underset{y\to x}{\lim\sup}\left(F(D^2\varphi(y),D\varphi(y),u(y),y)-f(y)\right)\geq 0\\
    \left(\mathrm{resp.},\quad ess~\underset{y\to x}{\lim\inf}\left(F(D^2\varphi(y),D\varphi(y),u(y),y)-f(y)\right)\leq 0\right)
  \end{aligned}
\end{equation*}
provided that for $\varphi\in W^{2,p}(\Omega)$, $u-\varphi$ attains its local maximum (resp., minimum) at $x\in\Omega$.

We call $u\in C(\Omega)$ an $L^p$-viscosity solution of \cref{FNE} if it is both an $L^p$-viscosity subsolution and supersolution of \cref{FNE}.
\end{definition}

\begin{remark}\label{r-14}
If all functions are continuous and $\varphi\in W^{2,p}(\Omega)$ is replaced by $\varphi\in C^2(\Omega)$, we arrive at the definition of $C$-viscosity solution (equivalent to the $L^p$-viscosity solution \cite[Proposition 2.9]{MR1376656}). It has been adopted by Caffarelli \cite{MR1005611} to study the interior regularity for fully nonlinear elliptic equations. In 1996, Caffarelli, Crandall, Kocan and \'{S}wi\k{e}ch \cite{MR1376656} studied $L^p$-viscosity solutions to equations with measurable ingredients. In this paper, we deal with the $L^n$-viscosity solution and write ``viscosity solution'' for short.
\end{remark}
~\\

We introduce the Pucci's class, which is frequently used for studying fully nonlinear elliptic equations.
\begin{definition}\label{d-Sf}
For $M\in \mathcal{S}^n$, denote its eigenvalues by $\lambda_i$ ($1\leq i\leq n$) and set
\begin{equation*}
\mathcal{M}_{\lambda,\Lambda}^{+}(M)=\Lambda\left(\sum_{\lambda_i>0}\lambda_i\right)
+\lambda\left(\sum_{\lambda_i<0}\lambda_i\right),~~
\mathcal{M}_{\lambda,\Lambda}^{-}(M)=\lambda\left(\sum_{\lambda_i>0}\lambda_i\right)
+\Lambda\left(\sum_{\lambda_i<0}\lambda_i\right).
\end{equation*}
Then we can define the Pucci's class as follows. Let $b\in L^p(\Omega) (p>n)$ and $f\in L^n(\Omega)$. We say that $u\in \underline{S}(\lambda,\Lambda,\mu,b,f)$ if $u$ is an $L^n$-viscosity subsolution of
\begin{equation}\label{SC2Sf}
  \mathcal{M}_{\lambda,\Lambda}^{+}(D^2u)+\mu|Du|^2+b|Du|= f.
\end{equation}
Similarly, we denote $u\in \bar{S}(\lambda,\Lambda,\mu,b,f)$ if $u$ is an $L^n$-viscosity supersolution of
\begin{equation}\label{SC2Sf-}
  \mathcal{M}_{\lambda,\Lambda}^{-}(D^2u)-\mu|Du|^2-b|Du|= f.
\end{equation}

We also define
\begin{equation*}
  S^*(\lambda,\Lambda,\mu,b,f)=\underline{S}(\lambda,\Lambda,\mu,b,-|f|)\cap \bar{S}(\lambda,\Lambda,\mu,b,|f|).
\end{equation*}
 We will denote $\underline{S}(\lambda,\Lambda,\mu,b,f)$ ($\bar{S}(\lambda,\Lambda,\mu,b,f)$, $S^*(\lambda,\Lambda,\mu,b,f)$) by $\underline{S}(\mu,b,f)$ ($\bar{S}(\mu,b,f)$, $S^*(\mu,b,f)$) and $\mathcal{M}_{\lambda,\Lambda}^{+}(M)$ ($\mathcal{M}_{\lambda,\Lambda}^{-}(M)$) by $\mathcal{M}^{+}(M)$ ($\mathcal{M}^{-}(M)$) for short since $\lambda,\Lambda$ are fixed constants throughout this paper.
\end{definition}

\section{Preliminary results}\label{P1}

In this section, we gather some preliminary results which will be used for proving our main regularity in following sections. First, we introduce the well-known regularity for some basic problems, which are called model problems in this paper. The various general problems are regarded as perturbations of them. Precisely, there exist two constants
\begin{equation}\label{e2.17}
0<\bar{\alpha}<1, \quad \bar{C}>0\quad (\mbox{fixed throughout this paper})
\end{equation}
depending only on $n,\lambda$ and $\Lambda$ such that the following four lemmas hold.

The first result is the interior  $C^{1,\alpha}$ regularity, which is almost a direct result of the interior H\"{o}lder regularity for the Pucci's class and Jensen's uniqueness theorem (see \cite[Corollary 5.7]{MR1351007}).
\begin{lemma}\label{l-3modin1}
Let $u$ be a viscosity solution of
\begin{equation}\label{e2.4}
F(D^2u)=0\quad\mbox{ in}~B_1.
\end{equation}
Then $u\in C^{1,\bar{\alpha}}(\bar{B}_{1/2})$ and
\begin{equation}\label{e.l3modin1-1}
\|u\|_{C^{1,\bar{\alpha}}(\bar{B}_{1/2})}\leq \bar{C}\|u\|_{L^{\infty}(B_1)}.
\end{equation}

In particular, $u\in C^{1,\bar{\alpha}}(0)$, i.e., there exists $P\in \mathcal{P}_1$ such that
\begin{equation}\label{e.3.1}
  |u(x)-P(x)|\leq \bar{C} |x|^{1+\bar{\alpha}}\|u\|_{L^{\infty }(B_1)}, ~~\forall ~x\in B_{1}
\end{equation}
and
\begin{equation}\label{e.3.2}
  |Du(0)|\leq \bar{C}\|u\|_{L^{\infty }(B_1)}.
\end{equation}
\end{lemma}
\begin{remark}\label{re2.3}
Silvestre and Teixeira \cite[Corollary 1.2]{MR3494911} obtained $C^{1,\alpha}$ regularity for \cref{e2.4} with every $0<\alpha<1$, provided that $F$ is ``asymptotically convex at infinity''.
\end{remark}
~\\

Next lemma concerns the interior $C^{2,\alpha}$ regularity, which was first proved independently by Evans \cite{MR649348} and Krylov \cite{MR661144} for smooth solutions. The proof for viscosity solutions was given by Cabr\'{e} and Caffarelli (see \cite{MR1391943} and \cite[Theorem 6.6]{MR1351007}).
\begin{lemma}\label{l-3modin2}
Let $F$ be convex (or concave) in $M$ and $u$ be a viscosity solution of
\begin{equation*}
F(D^2u)=0\quad\mbox{ in}~B_1.
\end{equation*}
Then $u\in C^{2,\bar{\alpha}}(\bar{B}_{1/2})$ and
\begin{equation}\label{e.l3modin2-1}
\|u\|_{C^{2,\bar{\alpha}}(\bar{B}_{1/2})}\leq \bar{C}\|u\|_{L^{\infty}(B_1)}.
\end{equation}

In particular, $u\in C^{2,\bar{\alpha}}(0)$, i.e., there exists $P\in \mathcal{P}_2$ such that
\begin{equation}\label{e.3.3}
  |u(x)-P(x)|\leq \bar{C} |x|^{2+\bar{\alpha}}\|u\|_{L^{\infty }(B_1)}, ~~\forall ~x\in B_{1},
\end{equation}
\begin{equation}\label{e.3.5}
F(D^2u(0))=0
\end{equation}
and
\begin{equation}\label{e.3.4}
  |Du(0)|+|D^2u(0)|\leq \bar{C}\|u\|_{L^{\infty }(B_1)}.
\end{equation}
\end{lemma}

\begin{remark}\label{re2.1}
In general, the condition ``$F$ is convex (or concave) in $M$ '' cannot be removed for interior $C^{2,\alpha}$ regularity. Indeed, Nadirashvili and Vl\u{a}du\c{t} \cite{MR3125267} proved that for any $0<\varepsilon<1$, there exists a smooth fully nonlinear uniformly elliptic operator $F$ such that $u\notin C^{1,\varepsilon}(\bar B_{1/2})$, where $u$ is a viscosity solution of
\begin{equation*}
F(D^2u)=0\quad\mbox{ in}~B_1\subset \mathbb{R}^n,~n\geq 5.
\end{equation*}
We also refer to \cite{MR2964616,MR2373018,MR2391642,MR2683755,MR2824567,MR3084701} for more results in this direction.

On the other hand, to obtain interior $C^{2,\alpha}$ regularity or a priori estimate, the convexity of $F$ in $M$ can be replaced by one of the following conditions:\\
\begin{description}
  \item (i) $n=2$ (see \cite[Theorem I]{MR0064986}, \cite[Theorem 4.9]{MR4560756} or \cite[Remark 2]{MR1793687});
  \item (ii) $F\in C^{1,1}, u\in C^{1,1}$ and $\left\{M\colon F(M)=0\right\}\cap \left\{M\colon tr(M)=t\right\}$ is strictly convex for any $t\in \mathbb{R}$ (see \cite{MR1793687});
  \item (iii) $F(M)=\min\left(F_1(M),F_2(M)\right)$, where $F_1$ is concave and $F_2$ is convex (see \cite{MR1995493}).
  \item (iv) $F\in C^{1}$ and $\|u\|_{L^{\infty}}\leq \delta\ll 1$ (see \cite{MR2334822} or \cite[Section 4]{MR2928094}).
  \item (v) $F\in C^{1,\alpha}$ for some $0<\alpha\leq 1$ and $u\in C^{2}$ (see \cite{MR2775423}).
  \item (vi) $\Lambda/\lambda-1\leq \delta\ll 1$ (see \cite{MR4503817} and a special case in \cite{MR4251799}).
\end{description}

\end{remark}
~\\

The following is the boundary $C^{1,\alpha}$ regularity, which was first proved by Krylov \cite{MR688919} for smooth solutions. Later, it was simplified by Caffarelli (see \cite[Theorem 9.31]{MR1814364} and \cite[Theorem 4.28]{MR787227}), which is applicable to the Pucci's class.

\begin{lemma}\label{l-31}
Let $u$ satisfy
\begin{equation*}
\left\{\begin{aligned}
&u\in S(\lambda,\Lambda,0)&& \quad\mbox{in}~~B_1^+;\\
&u=0&& \quad\mbox{on}~~T_1.
\end{aligned}\right.
\end{equation*}
Then $u\in C^{1,\bar{\alpha}}(0)$ and for some constant $a$,
\begin{equation}\label{e.l31-1}
  |u(x)-ax_n|\leq \bar{C} |x|^{1+\bar{\alpha}}\|u\|_{L^{\infty }(B_1^+)}, ~~\forall ~x\in B_{1}^+
\end{equation}
and
\begin{equation}\label{e.l31-2}
  |a|\leq \bar{C}\|u\|_{L^{\infty }(B_1^+)}.
\end{equation}
\end{lemma}
~\\

Finally, we present the boundary $C^{2,\alpha}$ regularity. The a priori estimates was first proved for smooth solutions by Krylov \cite{MR688919}. For viscosity solutions, it was proved by Silvestre and Sirakov \cite[Lemma 4.1]{MR3246039} by combining the $C^{1,\alpha}$ regularity and the technique of difference quotient.
\begin{lemma}\label{l-32}
Let $u$ be a viscosity solution of
\begin{equation*}
\left\{\begin{aligned}
&F(D^2u)=0&& \quad\mbox{in}~~B_1^+;\\
&u=0&& \quad\mbox{on}~~T_1.
\end{aligned}\right.
\end{equation*}
Then $u\in C^{2,\bar{\alpha}}(0)$, i.e., there exists $P\in \mathcal{P}_2$ such that
\begin{equation}\label{e.l32-1}
  |u(x)-P(x)|\leq \bar{C} |x|^{2+\bar{\alpha}}\|u\|_{L^{\infty }(B_1^+)}, ~~\forall ~x\in B_{1}^+,
\end{equation}
\begin{equation}\label{e.l32-2}
  F(D^2P)=0
\end{equation}
and
\begin{equation}\label{e.l32-3}
|Du(0)|+|D^2u(0)|\leq \bar{C}\|u\|_{L^{\infty }(B_1^+)}.
\end{equation}
Moreover, $P$ can be written as
\begin{equation*}
P(x)=\sum_{|\sigma|\leq 2, \sigma_n\geq 1} \frac{1}{\sigma !} a_{\sigma} x^{\sigma}.
\end{equation*}
\end{lemma}
~\\

We remark here that the Harnack inequality is a common key ingredient for all above regularity results.\\
(i) From the Harnack inequality, the interior H\"{o}lder regularity for the Pucci's class can be easily derived.\\
(ii) Since the difference quotient of a solution belongs to the Pucci's class, we have the interior $C^{1,\alpha}$ regularity (\Cref{l-3modin1}). \\
(iii) The key to the interior $C^{2,\alpha}$ regularity (\Cref{l-3modin2}) is that $u_{ee}$ (the second derivative of $u$ along some unit vector $e\in \mathbb{R}^n$) is a subsolution and then the weak Harnack inequality is applicable to $\sup u_{ee} -u_{ee}$. \\
(iv) By combining the Harnack inequality and a proper barrier, the boundary $C^{1,\alpha}$ regularity (\Cref{l-31}) follows. \\
(v) With regards to the boundary $C^{2,\alpha}$ regularity (\Cref{l-32}), note that $u_i$ ($1\leq i\leq n-1$) belongs to the Pucci's class and then apply \Cref{l-31} to $u_i$ to conclude that $u_i\in C^{1,\alpha}$. Then the boundary $C^{2,\alpha}$ regularity can be obtained and the convexity of $F$ in $M$ is not needed.

If we study linear uniformly elliptic equations in nondivergence form, equations in \Cref{l-3modin1} and \Cref{l-3modin2} will be reduced to
\begin{equation*}
\Delta u=0\quad \mbox{ in}~B_1
\end{equation*}
and equations in \Cref{l-31} and \Cref{l-32} will be reduced to
\begin{equation*}
\left\{\begin{aligned}
&\Delta u=0&& \quad\mbox{in}~~B_1^+;\\
&u=0&& \quad\mbox{on}~~T_1.
\end{aligned}\right.
\end{equation*}
The regularity for above two problems are well-known, i.e., harmonic function theory. Then through a perturbation argument as in this paper, we can obtain the interior and boundary pointwise $C^{k,\alpha}$ regularity for any $0<\alpha<1$ and $k\geq 1$. Indeed, we have not seen any general theory of pointwise $C^{k,\alpha}$ ($k\geq 3$) regularity even for linear equations. Above argument shows that we can deduce systematically pointwise $C^{k,\alpha}$ regularity for linear equations in a simple way.

As explained above, the proofs of \Crefrange{l-3modin1}{l-32} involve the Harnack inequality. Hence, we have only ``certain'' small decay for the oscillation of the solution, which implies the H\"{o}lder regularity with some small exponent $\bar{\alpha}$. In addition, the constants ``$\bar{\alpha}$'' in above lemmas could be different. The precise statement should be ``there exist four constants $\bar{\alpha}_1,\bar{\alpha}_2,\bar{\alpha}_3$ and $\bar{\alpha}_4$ such that \Crefrange{l-3modin1}{l-32} hold respectively''. For simplicity, we use one symbol $\bar{\alpha}$ to denote these constants.

Since we use the compactness method, we introduce the following closedness result for viscosity solutions (see \cite[Theorem 3.8]{MR1376656}, \cite[Proposition 9.4]{MR2552914} and \cite[Proposition 2.3]{MR3980853}). It states that the limit of a sequence of solutions is again a solution of some equation.
\begin{lemma}\label{l-35}
Let $F,\{F_m\}_{m\geq 1}$ be fully nonlinear operators and $u_m\in C(\Omega)$ be $L^n$-viscosity subsolutions (or supersolutions) of
\begin{equation*}
F_m(D^2u_m,Du_m,u_m,x)= f_m\quad\mbox{in}~~\Omega.
\end{equation*}
Suppose that $b\in L^p (\Omega) (p>n)$, $f,f_m\in L^n(\Omega)$,
\begin{equation*}
u_m\rightarrow u~~\mbox{ in}~~L^{\infty}_{\mathrm{loc}}(\Omega)~~\mbox{ as}~~m\rightarrow \infty
\end{equation*}
and for any ball $B\subset \subset \Omega$, $\varphi \in W^{2,n}(B)$,
\begin{equation*}
\|(\psi_m-\psi)^+\|_{L^n(B)}~\left(\mbox{or}~\|(\psi_m-\psi)^-\|_{L^n(B)}\right)\rightarrow 0~~\mbox{ as}~~ m\rightarrow \infty,
\end{equation*}
where
\begin{equation*}\label{FkCg}
\psi_m(x)\coloneqq F_m(D^2\varphi,D\varphi,u_m,x)-f_m(x),\quad\psi(x)\coloneqq F(D^2\varphi,D\varphi,u,x)-f(x).
\end{equation*}
Then $u$ is an $L^n$-viscosity subsolution (or supersolution) of
\begin{equation}\label{FkC}
F(D^2u,Du,u,x)= f\quad\mbox{in}~~\Omega.
\end{equation}

If $F$ and $f$ are continuous in $x$, it is enough to take $\varphi \in C^2(\bar B)$, in which case $u$ is a $C$-viscosity subsolution (or supersolution) of \cref{FkC}.
\end{lemma}

\begin{remark}\label{r-3.2}
In this paper, we always use this lemma in the case that $F$ is continuous in $x$ and $f\equiv 0$. Hence, we always take $\varphi \in C^2(\bar B)$ in the proofs (see the proof of \Cref{In-l-C1a-mu} etc.).

If all functions concerned in \Cref{l-35} are continuous in $x$ and we consider $C$-viscosity solutions, the closedness result is easily to establish by the definition of viscosity solution and assuming $\|(\psi_m-\psi)^+\|_{L^{\infty}(B)}~(\mbox{or}~\|(\psi_m-\psi)^-\|_{L^{\infty}(B)})\rightarrow 0$ (see \cite[Proposition 2.9]{MR1351007}).
\end{remark}
~\\

\Cref{l-35} requires that $u_m$ converges uniformly. The following lemma (see \cite[Theorem 2]{MR2592289}) states the interior H\"{o}lder regularity, which provides necessary compactness for solutions and then guarantee the uniform convergence. Moreover, it also contains the boundary pointwise H\"{o}lder regularity, which will be used in the normalization procedure for boundary pointwise $C^{1,\alpha}$ regularity (see the proof after \cref{e.c1a-v}).

\begin{lemma}\label{l-3Ho}
Let $u$ satisfy
\begin{equation*}
u\in S^*(\lambda,\Lambda,\mu,b,f) \quad\mbox{in}~~\Omega.
\end{equation*}
Then there exists $0<\alpha<1$ depending only on $n,\lambda,\Lambda,p$ and $\|b\|_{L^p(\Omega)}$ such that $u\in C^{\alpha}(\bar\Omega')$ for any $\Omega'\subset\subset\Omega$ and
\begin{equation*}
\|u\|_{C^{\alpha}(\bar\Omega')}\leq C,
\end{equation*}
where $C$ depends only on $n,\lambda,\Lambda,p,\mu,\|b\|_{L^p(\Omega)},\|f\|_{L^n(\Omega)}$, $\Omega'$, $\mathrm{dist}(\Omega',\partial \Omega)$ and $\|u\|_{L^{\infty}(\Omega)}$.

In addition, suppose that $u$ satisfies
\begin{equation*}
\left\{\begin{aligned}
&u\in S^*(\lambda,\Lambda,\mu,b,f)&& \quad\mbox{in}~~\Omega_1;\\
&u=g&& \quad\mbox{on}~~(\partial \Omega)_1,
\end{aligned}\right.
\end{equation*}
where $\Omega$ satisfies the exterior cone condition at $0\in \partial \Omega$ and $g\in C^{\alpha_1}(0)$ for some $0<\alpha_1<1$. Then there exists $0<\alpha\leq \alpha_1$ depending only on $n,\lambda,\Lambda,p,\|b\|_{L^p(\Omega)}$ and the cone such that $u\in C^{\alpha}(0)$ and
\begin{equation*}
  |u(x)-u(0)|\leq C|x|^{\alpha}, ~\forall ~x\in \Omega_1,
\end{equation*}
where $C$ depends only on $n,\lambda,\Lambda,p,\mu,\|b\|_{L^p(\Omega_1)},\|f\|_{L^n(\Omega_1)}$, $\|u\|_{L^{\infty}(\Omega_1)}$ and the cone.
\end{lemma}
~\\

\Cref{l-3Ho} treats the Pucci's class. In fact, viscosity solutions of any equation belongs to the Pucci's class.
\begin{lemma}\label{le2.5}
Let $u$ be a viscosity solution of
\begin{equation*}
  F(D^2u,Du,u,x)=f\quad\mbox{ in}~~\Omega.
\end{equation*}
Then
\begin{equation*}
u\in S^*(\lambda,\Lambda,\mu,b,|f|+|F(0,0,0,\cdot)|
+c\omega_0(\|u\|_{L^{\infty}(\Omega)},\|u\|_{L^{\infty}(\Omega)})).
\end{equation*}
\end{lemma}
\proof By \Cref{d-viscoF}, for any $x\in \Omega$ and $\varphi\in W^{2,p}(\Omega)$ with $u-\varphi$ attaining its local maximum at $x$, we have
\begin{equation*}
ess~\underset{y\to x}{\lim\sup}\left(F(D^2\varphi(y),D\varphi(y),u(y),y)-f(y)\right)\geq 0.
\end{equation*}
By the structure condition \Cref{SC2},
\begin{equation*}
  \begin{aligned}
F&(D^2\varphi(y),D\varphi(y),u(y),y)-F(0,0,0,y)\\
&\leq\mathcal{M}^{+}(D^2\varphi(y))+\mu|D\varphi(y)|^2+b(y)|D\varphi(y)|
+c(y)\omega_0(\|u\|_{L^{\infty}(\Omega)},\|u\|_{L^{\infty}(\Omega)}).\\
  \end{aligned}
\end{equation*}
Let
\begin{equation*}
\tilde{f}(y)=f(y)-F(0,0,0,y)-c(y)\omega_0(\|u\|_{L^{\infty}(\Omega)},\|u\|_{L^{\infty}(\Omega)}).
\end{equation*}
Then
\begin{equation*}
  \begin{aligned}
ess~&\underset{y\to x}{\lim\sup}
\left(\mathcal{M}^{+}(D^2\varphi(y))+\mu|D\varphi(y)|^2+b(y)|D\varphi(y)|-\tilde{f}(y)\right)\\
&\geq ess~\underset{y\to x}{\lim\sup}\left(F(D^2\varphi(y),D\varphi(y),u(y),y)-f(y)\right)\geq 0.
  \end{aligned}
\end{equation*}
Hence, $u$ is a subsolution. Similarly, we can prove that $u$ is a supersolution and then the conclusion follows.~\qed~\\

We also need the compactness of solutions up to the boundary for the boundary regularity. In the following, we prepare the compactness up to the boundary following the idea of \cite{MR4088470}. First, we introduce the Alexandrov-Bakel'man-Pucci maximum principle (see \cite[Proposition 2.5]{MR3980853}).

\begin{lemma}\label{l-3ABP}
Let $u$ be a viscosity solution of
\begin{equation*}
  \mathcal{M}_{\lambda,\Lambda}^{+}(D^2u)+\mu|Du|^2+b|Du|\geq  f\quad\mbox{in}~~\Omega,
\end{equation*}
where $b\in L^p(\Omega)$ $(p>n)$ and $f\in L^n(\Omega)$. Suppose that
\begin{equation}\label{e2.28}
  \mu \|f^-\|_{L^n(\Omega)}\cdot\mathrm{diam}(\Omega)\leq \delta,
\end{equation}
where $\delta>0$ depends only on
$n,\lambda,\Lambda,p,\mathrm{diam}(\Omega)$ and $\|b\|_{L^p(\Omega)}$. Then
\begin{equation}\label{e2.29}
  \max_{\bar{\Omega}}u\leq \max_{\partial \Omega}u+C \|f^-\|_{L^n(\Omega)},
\end{equation}
where $C$ depends only on $n,\lambda,\Lambda,p,\mathrm{diam}(\Omega)$ and $\|b\|_{L^p(\Omega)}$.
\end{lemma}
\begin{remark}\label{re2.4}
There is another important type of the Alexandrov-Bakel'man-Pucci maximum principle. Precisely, the term $\|f^-\|_{L^n(\Omega)}$ in \cref{e2.29} is replaced by $\|f^-\|_{L^n(\Gamma_u)}$, where $\Gamma_u$ is the contact set (see \cite[Definition 3.1]{MR1351007}). This type of maximum principle is crucial for proving the Harnack inequality (see \cite[Lemma 4.5]{MR1351007}). For this type of maximum principle , we refer to \cite{MR4515258} for the latest development and the references therein for more results.
\end{remark}
~\\

Next, we consider the ``equicontinuity'' of solutions up to the boundary and use the following notation to describe the oscillation of $\partial \Omega$ when $\partial \Omega$ is not smooth. For $r>0$, define
\begin{equation}\label{e1.1}
\underset{B_{r}}{\mathrm{osc}}~\partial\Omega= \underset{x\in \partial \Omega\cap B_r}{\sup} x_n -\underset{x\in \partial \Omega\cap B_r}{\inf} x_n.
\end{equation}
The following lemma provides a uniform estimate up to the boundary.
\begin{lemma}\label{l-33}
Let $0<\delta<1$ be as in \Cref{l-3ABP} and $u$ satisfy
\begin{equation*}
\left\{\begin{aligned}
&u\in S^*(\lambda,\Lambda,\mu,b,f)&& \quad\mbox{in}~~\Omega_1;\\
&u=g&& \quad\mbox{on}~~(\partial \Omega)_1.
\end{aligned}\right.
\end{equation*}
Suppose that
\begin{equation*}
\|u\|_{L^{\infty}(\Omega_{1})}\leq 1,\quad\mu\leq \frac{\delta}{C_0},\quad
\max(\|b\|_{L^{p}(\Omega_1)},\|f\|_{L^{n}(\Omega_1)},\|g\|_{L^{\infty}((\partial \Omega)_1)},
\underset{B_1}{\mathrm{osc}}~\partial\Omega) \leq \delta,
\end{equation*}
where $C_0$ (to be specified later) depends only on $n,\lambda$ and $\Lambda$.

Then
\begin{equation*}
  |u(x)|\leq Cx_n+C\delta,~\forall ~x\in \Omega_{1/2},
\end{equation*}
where $C$ depends only on $n,\lambda,\Lambda$ and $p$.
\end{lemma}

\proof Set $B^{+}=B^{+}_{1}-\delta e_n $, $T=T_1-\delta e_n$ and then $\partial \Omega\cap B_{1/4}\subset B^+$. Take
\begin{equation*}
  v(x)=C\big((1-\delta)^{-\alpha}-|x+e_n|^{-\alpha}\big).
\end{equation*}
Then by choosing proper constants $C$ and $\alpha$ (depending only on $n,\lambda$ and $\Lambda$), we have $v(-\delta e_n)=0$ and
\begin{equation*}
\left\{\begin{aligned}
 &\mathcal{M}_{\lambda,\Lambda}^{+}(D^2v)\leq 0 &&\mbox{in}~~B^+; \\
 &v\geq 0 && \mbox{in}~~ B^+;\\
 &v\geq 1 &&\mbox{on}~~\partial B^+\backslash T.
\end{aligned}
\right.
\end{equation*}
Let $w=u-v$ and then $w$ satisfies
\begin{equation*}
    \left\{
    \begin{aligned}
      &w\in \underline{S}(\lambda,\Lambda,2\mu,b,\tilde{f}) &&\mbox{in}~~ \Omega \cap B^+; \\
      &w\leq g &&\mbox{on}~~\partial \Omega \cap \bar{B}^+;\\
      &w\leq 0 &&\mbox{on}~~\partial B^+\cap \Omega,
    \end{aligned}
    \right.
\end{equation*}
where $\tilde{f}=-|f|-2\mu|Dv|^2-b|Dv|$.

From the definition of $v$,
\begin{equation*}
  |v(x)|\leq C(x_n+\delta)\quad\mbox{ on}~~\left\{(x',x_n)\in \Omega\colon x'=0,-\delta\leq x_n\leq 1/2\right\}.
\end{equation*}
where $C$ depends only on $n,\lambda$ and $\Lambda$. For $w$, note that
\begin{equation*}
\|\tilde{f}\|_{L^n(\Omega_1)}\leq \|f\|_{L^n(\Omega_1)}+\|2\mu |Dv|^2\|_{L^n(\Omega_1)}
+\|b|Dv|\|_{L^n(\Omega_1)}\leq C_0/2,
\end{equation*}
where $C_0$ depends only on $n,\lambda$ and $\Lambda$. Then
\begin{equation*}
2\mu \|\tilde{f}\|_{L^n(\Omega_1)}\leq \delta.
\end{equation*}
By applying the Alexandrov-Bakel'man-Pucci maximum principle to $w$ (\Cref{l-3ABP}),
\begin{equation*}
  \begin{aligned}
\sup_{\Omega \cap B^+} w&
\leq\|g\|_{L^{\infty }(\partial \Omega \cap B^+)}+C\left(\|f\|_{L^n}+
\mu\|Dv\|^2_{L^{\infty}}+\|b\|_{L^n} \|Dv\|_{L^{\infty}}\right)
\leq C\delta,
  \end{aligned}
\end{equation*}
where $C$ depends only on $n,\lambda,\Lambda$ and $p$. Thus,
\begin{equation*}
u=v+w\leq Cx_n+C\delta\quad\mbox{ on}~\{(x',x_n)\in \Omega\colon x'=0,-\delta\leq x_n\leq 1/2\}.
\end{equation*}
By considering $v(x'-x'_0,x_n)$ for any $x'_0\in T_{1/2}$ and similar arguments, we obtain
\begin{equation}\label{e2.1}
u\leq Cx_n+C\delta\quad\mbox{ in}~~\Omega_{1/2}.
\end{equation}

The proof for
\begin{equation*}
u \geq -Cx_n-C\delta \quad\mbox{ in}~~\Omega_{1/2}
\end{equation*}
is similar and we omit it here. Therefore, the proof is completed.\qed~\\

\begin{remark}\label{r-3.1}
If $\Omega$ satisfies the exterior cone condition at every boundary point (or $\partial \Omega\cap B_1\in C^{0,1}$), we have the equicontinuity up to the boundary by \Cref{l-3Ho}. However, it is not assumed in this paper since we study the pointwise regularity.
\end{remark}
~\\

By \Cref{d-f}, to obtain the pointwise $C^{k,\alpha}(0)$ ($k\geq 1$) regularity, we need to find $P\in \mathcal{P}_k$ such that
\begin{equation*}
  |u(x)-P(x)|\leq C|x|^{k+\alpha},~\forall~x\in B_1.
\end{equation*}
Usually, it is not easy to find $P$ directly. Instead, we prove a sequence of estimates with different polynomials in different scales, which is a standard argument in the regularity theory (e.g. \cite[Chapter 8]{MR1351007}). For the reader's convenience, we give a proof here.
\begin{lemma}\label{le2.1}
Suppose that there exist a sequence of $P_m\in \mathcal{P}_k$ ($m\geq -1$, $P_{-1}\equiv 0$) such that for all $m\geq 0$,
\begin{equation}\label{e2.6}
\|u-P_m\|_{L^{\infty }(B _{\eta^{m}})}\leq K\eta ^{m(k+\alpha )}
\end{equation}
and
\begin{equation}\label{e2.7}
\sum_{i=0}^{k}\eta^{mi}|D^iP_m(0)-D^iP_{m-1}(0)|\leq K\eta^{m(k+\alpha)},
\end{equation}
where $0<\eta<1$ and $K>0$ are constants. Then $u\in C^{k,\alpha}(0)$, i.e., there exists $P\in \mathcal{P}_k$ such that
\begin{equation*}
  |u(x)-P(x)|\leq CK|x|^{k+\alpha}, ~\forall ~x\in B_{1},
\end{equation*}
and
\begin{equation}\label{e2.5}
\|P\|\leq CK,
\end{equation}
where $C$ depends only on $k,n,\alpha$ and $\eta$.
\end{lemma}
\proof By \cref{e2.7}, for any $0\leq i\leq k$ and $m_2>m_1$,
\begin{equation}\label{e2.8}
  \begin{aligned}
|D^iP_{m_2}(0)-D^iP_{m_1}(0)|\leq& \sum_{m=m_1+1}^{m_2}|D^iP_m(0)-D^iP_{m-1}(0)|\\
\leq& K \sum_{m=m_1+1}^{m_2}\eta^{m(k-i+\alpha)}\leq CK\eta^{(m_1+1)(k-i+\alpha)}.
  \end{aligned}
\end{equation}
Hence, $\{D^iP_m(0)\}_{m\geq -1}$ are Cauchy sequences for all $0\leq i\leq k$. Then there exist constants $a_{\sigma}$ ($\sigma\in \mathbb{N}^n,|\sigma|\leq k$) such that $D^{\sigma}P_m(0)\to a_{\sigma}$ as $m\to \infty$. In addition, by fixing $m_1$ and letting $m_2\to \infty$ in \cref{e2.8}, we have
\begin{equation}\label{e2.9}
  |D^{\sigma}P_m(0)-a_{\sigma}|\leq  CK\eta ^{(m+1)(k-|\sigma|+\alpha)},~\forall ~m\geq 0,~|\sigma|\leq k.
\end{equation}
Therefore,
\begin{equation*}
  \begin{aligned}
|a_{\sigma}|\leq &CK\eta ^{(m+1)(k-|\sigma|+\alpha)}+|D^{\sigma}P_m(0)|\\
    \leq &CK+\sum_{l=0}^{m}|D^{\sigma}P_l(0)-D^{\sigma}P_{l-1}(0)|\\
    \leq &CK+K\sum_{l=0}^{m}\eta ^{l(k-|\sigma|+\alpha)}
    \leq CK.
  \end{aligned}
\end{equation*}

Let
\begin{equation*}
P(x)=\sum_{|\sigma|\leq k} \frac{a_{\sigma}}{\sigma!} x^{\sigma}
\end{equation*}
and then
\begin{equation}\label{e2.18}
\|P\|\leq CK.
\end{equation}
Moreover, for any $x\in B_1$, there exists $m\geq 0$ such that $\eta^{m+1}\leq |x|< \eta^{m}$. Hence, by \cref{e2.6} and \cref{e2.9},
\begin{equation}\label{e2.19}
  \begin{aligned}
    |u(x)-P(x)|\leq & |u(x)-P_m(x)|+|P_m(x)-P(x)|\\
    \leq & K\eta^{m(k+\alpha )}+\sum_{|\sigma|\leq k} \frac{|D^{\sigma}P_m(0)-a_{\sigma}|}{\sigma!} |x|^{|\sigma|}\\
    \leq &\frac{K}{\eta^{(k+\alpha)}}\cdot \eta^{(m+1)(k+\alpha)}+CK\eta^{(m+1)(k-|\sigma|+\alpha )}|x|^{|\sigma|}\\
\leq &CK|x|^{k+\alpha}.
  \end{aligned}
\end{equation}
Therefore, the proof is completed.~\qed~\\

Similarly, to prove the pointwise $C^{k}(0)$ regularity, we only need to prove a sequence of discrete estimates.
\begin{lemma}\label{le2.3}
Let $\omega$ be a Dini function (see \cref{e.1.Dini}) with
\begin{equation*}
J_{\omega}=\int_{0}^{1}\frac{\omega(r)}{r} dr+\omega(1)\leq 1.
\end{equation*}
Define
\begin{equation}\label{e2.14}
A_{-1}=A_0=1,~A_m=\max(\omega(\eta^{m}),\eta^{\alpha} A_{m-1})~(m\geq 1),
\end{equation}
where $0<\alpha\leq 1$ and $0<\eta<1$. Suppose that there exist a sequence of $P_m\in \mathcal{P}_k$ ($m\geq -1$, $P_{-1}\equiv 0$) such that for all $m\geq 0$,
\begin{equation}\label{e2.10}
\|u-P_m\|_{L^{\infty }(B _{\eta^{m}})}\leq K\eta ^{mk}A_m
\end{equation}
and
\begin{equation}\label{e2.13}
\sum_{i=0}^{k}\eta^{mi}|D^iP_m(0)-D^iP_{m-1}(0)|\leq K\eta^{mk}A_{m-1}.
\end{equation}

Then $u\in C^{k,\omega_u}(0)$, i.e., there exist $P\in \mathcal{P}_k$ and $\omega_u$ such that
\begin{equation*}
  |u(x)-P(x)|\leq CK|x|^k\omega_u(|x|), ~\forall ~x\in B_{\eta^2},
\end{equation*}
and
\begin{equation}\label{e2.24}
\|P\|\leq CK,
\end{equation}
where
\begin{equation*}
\omega_u(r)\coloneqq r^{\alpha}+r^{\alpha}\int_{r}^{1} \frac{\omega(\tau)}{\tau^{1+\alpha}}d\tau
+\int_{0}^{r/\eta^2}\frac{\omega(\tau)}{\tau}d\tau,~\forall ~0<r<\eta^2
\end{equation*}
and $C$ depends only on $k,n,\alpha$ and $\eta$.
\end{lemma}
\proof Throughout the proof, $C$ always denotes a constant depending only on $k,n,\alpha$ and $\eta$. Since $\omega$ is a Dini function, for any $m_2\geq m_1\geq 1$,
\begin{equation}\label{e2.20}
  \begin{aligned}
\sum_{m=m_1}^{m_2}\omega(\eta^m)
&=\sum_{m=m_1}^{m_2}\frac{\omega(\eta^m)\left(\eta^{m-1}-\eta^m\right)}{\eta^{m-1}-\eta^m}
=\frac{1}{1-\eta}
\sum_{m=m_1}^{m_2}\frac{\omega(\eta^m)}{\eta^{m-1}}\left(\eta^{m-1}-\eta^m\right)\\
&\leq \frac{1}{1-\eta}
\sum_{m=m_1}^{m_2}\int_{\eta^{m}}^{\eta^{m-1}}\frac{\omega(\tau)}{\tau}d\tau
= \frac{1}{1-\eta}\int_{\eta^{m_2}}^{\eta^{m_1-1}}\frac{\omega(\tau)}{\tau}d\tau.
  \end{aligned}
\end{equation}

By \cref{e2.14}, $A_m\to 0$ decreasingly as $m\to \infty$. Since
\begin{equation*}
A_m\leq \omega(\eta^m)+\eta^{\alpha} A_{m-1},~\forall ~m\geq 1,
\end{equation*}
we have
\begin{equation*}
\sum_{m=m_1}^{m_2} A_m\leq \sum_{m=m_1}^{m_2}\omega(\eta^m)+\eta^{\alpha}\sum_{m=m_1}^{m_2} A_m+\eta^{\alpha}A_{m_1-1},~\forall ~m_2\geq m_1\geq 1.
\end{equation*}
Thus, by noting \cref{e2.20}, for any $m_2\geq m_1\geq 1$,
\begin{equation}\label{e2.22}
\begin{aligned}
\sum_{m=m_1}^{m_2} A_m
\leq\frac{1}{1-\eta^{\alpha}}\left(\sum_{m=m_1}^{m_2}\omega(\eta^m)+A_{m_1-1}\right)
\leq C\left(\int_{\eta^{m_2}}^{\eta^{m_1-1}}\frac{\omega(\tau)}{\tau}d\tau
+A_{m_1-1}\right).
\end{aligned}
\end{equation}
That is, the series $\sum A_m$ converges and
\begin{equation*}
\sum_{m=1}^{\infty} A_m\leq  C\left(\int_{0}^{1}\frac{\omega(\tau)}{\tau}d\tau+1\right).
\end{equation*}
Furthermore, for any $m\geq 1$,
\begin{equation}\label{e2.23}
  \begin{aligned}
A_m\leq &\omega(\eta^{m})+\eta^{\alpha}A_{m-1}
\leq \sum_{l=1}^{m}\eta^{(m-l)\alpha}\omega(\eta^{l})+\eta^{m\alpha}
=\eta^{m\alpha} \sum_{l=1}^{m}\frac{\omega(\eta^{l})}{\eta^{l\alpha}}+\eta^{m\alpha}\\
\leq & C\eta^{m\alpha}\sum_{l=1}^{m}
\frac{\omega(\eta^{l})}{\eta^{(l-1)\alpha}\eta^{l-1}}(\eta^{l-1}-\eta^{l})
+\eta^{m\alpha}
\leq  C\eta^{m\alpha}\int_{\eta^{m}}^{1} \frac{\omega(\tau)}{\tau^{1+\alpha}}d\tau+\eta^{m\alpha}.
  \end{aligned}
\end{equation}

Then with the aid of above analysis, we can prove the lemma similar to that of \Cref{le2.1}. By \cref{e2.13}, for any $0\leq i\leq k$ and $m_2>m_1\geq 0$,
\begin{equation*}
  \begin{aligned}
|D^iP_{m_2}(0)-D^iP_{m_1}(0)|\leq \sum_{m=m_1+1}^{m_2}|D^iP_m(0)-D^iP_{m-1}(0)|
\leq K \sum_{m=m_1+1}^{m_2}\eta^{m(k-i)}A_{m-1}.
  \end{aligned}
\end{equation*}
Since $\sum A_m$ converges, $\{D^iP_m(0)\}_{m\geq -1}$ are Cauchy sequences for all $0\leq i\leq k$ (in fact, we only use the convergence of $\sum A_m$ for showing that $\{D^kP_m(0)\}_{m\geq -1}$ is a Cauchy sequence). Then there exist constants $a_{\sigma}$ ($\sigma\in \mathbb{N}^n,|\sigma|\leq k$) such that
\begin{equation}\label{e2.21}
  |D^{\sigma}P_m(0)-a_{\sigma}|\leq  K\sum_{l=m+1}^{\infty}\eta^{l(k-|\sigma|)}A_{l-1}
\leq K\eta^{(m+1)(k-|\sigma|)}\sum_{l=m}^{\infty}A_{l},~\forall ~m\geq 0,~|\sigma|\leq k
\end{equation}
and
\begin{equation*}
  \begin{aligned}
|a_{\sigma}|\leq&\sum_{m=0}^{\infty}|D^{\sigma}P_m(0)-D^{\sigma}P_{m-1}(0)|
    \leq K\sum_{m=0}^{\infty}\eta ^{m(k-|\sigma|)}A_m,~\forall ~|\sigma|\leq k.
  \end{aligned}
\end{equation*}
Hence,
\begin{equation}\label{e2.11}
|a_{\sigma}|\leq\left\{  \begin{aligned}
&K\sum_{m=0}^{\infty}\eta ^{m(k-|\sigma|)}\leq CK,&&~\forall ~|\sigma|\leq k-1,\\
&K\sum_{m=0}^{\infty}A_m\leq CK\left(\int_{0}^{1}\frac{\omega(\tau)}{\tau}d\tau+1\right),&&~\forall ~|\sigma|=k.
  \end{aligned}\right.
\end{equation}

As in the proof of \Cref{le2.1}, let
\begin{equation*}
P(x)=\sum_{|\sigma|\leq k} \frac{a_{\sigma}}{\sigma!} x^{\sigma}.
\end{equation*}
For any $x\in B_{\eta^2}$, there exists $m\geq 2$ such that $\eta^{m+1}\leq |x|< \eta^{m}$. Hence,
\begin{equation*}
  \begin{aligned}
|u(x)-P(x)|\leq& |u(x)-P_m(x)|+|P_m(x)-P(x)|\\
\leq& |u(x)-P_m(x)|+\sum_{|\sigma|\leq k} \frac{|D^{\sigma}P_m(0)-a_{\sigma}|}{\sigma!} |x|^{|\sigma|}\\
\leq & K\eta^{mk}A_m+CK|x|^{k}\sum_{l=m}^{\infty}A_l ~\quad(\mbox{by}~\cref{e2.10}~\mbox{and}~\cref{e2.21})\\
\leq & CK\eta^{mk}\left(A_m+A_{m-1}
+\int_{0}^{\eta^{m-1}}\frac{\omega(\tau)}{\tau}d\tau\right)
~\quad(\mbox{by}~\cref{e2.22})\\
\leq & CK\eta^{mk}\left(A_{m}
+\int_{0}^{\eta^{m-1}}\frac{\omega(\tau)}{\tau}d\tau\right)\\
\leq &CK\eta^{mk}\left(\eta^{m\alpha}\int_{\eta^{m}}^{1} \frac{\omega(\tau)}{\tau^{1+\alpha}}d\tau+\eta^{m\alpha}
+\int_{0}^{\eta^{m-1}}\frac{\omega(\tau)}{\tau}d\tau\right) ~\quad(\mbox{by}~\cref{e2.23})\\
\leq &CK|x|^k\left(|x|^{\alpha}+|x|^{\alpha}\int_{|x|}^{1} \frac{\omega(\tau)}{\tau^{1+\alpha}}d\tau
+\int_{0}^{|x|/\eta^2}\frac{\omega(\tau)}{\tau}d\tau\right).
  \end{aligned}
\end{equation*}
For $0<r<\eta^2$, define
\begin{equation*}
  \omega_u(r)=r^{\alpha}+r^{\alpha}\int_{r}^{1} \frac{\omega(\tau)}{\tau^{1+\alpha}}d\tau
+\int_{0}^{r/\eta^2}\frac{\omega(\tau)}{\tau}d\tau.
\end{equation*}
Then $\omega_u$ is a modulus of continuity. Moreover,
\begin{equation}\label{e2.12}
\omega_u(\eta^2)\leq \eta^{2\alpha}+2\int_{0}^{1} \frac{\omega(\tau)}{\tau}d\tau\leq 3.
\end{equation}
Therefore, $u\in C^{k,\omega_u}(0)$ and \cref{e2.24} holds. ~\qed~\\

\begin{remark}\label{re2.2}
From above proof we know that only $|D^ku(0)|$ and $[u]_{C^{k,\omega_u}(0)}$ depend on $I_{\omega}\coloneqq\int_{0}^1 \frac{\omega(\tau)}{\tau}d\tau$ (see \cref{e2.11} and \cref{e2.12}). Other quantities in $\|u\|_{C^{k,\omega_u}(0)}$ do not depends on $I_{\omega}$.
\end{remark}
~\\

For the pointwise $C^{k,\mathrm{lnL}}(0)$ regularity, we need the following lemma.
\begin{lemma}\label{le2.4}
Suppose that there exist a sequence of $P_m\in \mathcal{P}_{k+1}$ ($m\geq -1$, $P_{-1}\equiv 0$) such that for all $m\geq 0$,
\begin{equation}\label{e2.15}
\|u-P_m\|_{L^{\infty }(B _{\eta^{m}})}\leq K\eta ^{m(k+1)}
\end{equation}
and
\begin{equation}\label{e2.16}
\sum_{i=0}^{k+1}\eta^{mi}|D^iP_m(0)-D^iP_{m-1}(0)|\leq K\eta^{m(k+1)}.
\end{equation}

Then $u\in C^{k,\mathrm{lnL}}(0)$, i.e., there exists $P\in \mathcal{P}_k$ such that
\begin{equation*}
  |u(x)-P(x)|\leq CK|x|^{k+1}\big|\ln|x|\big|, ~\forall ~x\in B_{1/2},
\end{equation*}
and
\begin{equation}\label{e2.27}
\|P\|\leq CK,
\end{equation}
where $C$ depends only on $k,n$ and $\eta$.
\end{lemma}
\proof The proof is quite similar to that of \Cref{le2.1}. The main difference is that $\{D^{k+1}P_m(0)\}_{m\geq -1}$ is no longer a Cauchy sequence. Indeed, similar to the proof of \Cref{le2.1}, by \cref{e2.16}, $\{D^iP_m(0)\}_{m\geq -1}$ are Cauchy sequences for all $0\leq i\leq k$. Moreover, there exist constants $a_{\sigma}$ ($\sigma\in \mathbb{N}^n,|\sigma|\leq k$) such that
\begin{equation}\label{e2.25}
  |D^{\sigma}P_m(0)-a_{\sigma}|\leq  CK\eta ^{(m+1)(k+1-|\sigma|)},~\forall ~m\geq 0,~|\sigma|\leq k
\end{equation}
and $|a_{\sigma}|\leq CK$. In addition, by \cref{e2.16},
\begin{equation}\label{e2.26}
  |D^{k+1}P_m|\leq \sum_{l=0}^{m}|D^{k+1}P_l(0)-D^{k+1}P_{l-1}(0)|
\leq (m+1)K=\frac{K}{|\ln \eta|}|\ln \eta^{m+1}|.
\end{equation}

Let
\begin{equation*}
P(x)=\sum_{|\sigma|\leq k} \frac{a_{\sigma}}{\sigma!} x^{\sigma}
\end{equation*}
and then $\|P\|\leq CK$. Moreover, for any $x\in B_{1/2}$, there exists $m\geq 0$ such that $\eta^{m+1}\leq |x|< \eta^{m}$. Hence, by \cref{e2.15}, \cref{e2.25} and \cref{e2.26},
\begin{equation*}
  \begin{aligned}
    |u(x)-P(x)|\leq & |u(x)-P_m(x)|+|P_m(x)-P(x)|\\
    \leq & |u(x)-P_m(x)|
    +\sum_{|\sigma|\leq k} \frac{|D^{\sigma}P_m(0)-a_{\sigma}|}{\sigma!} |x|^{|\sigma|}
    +\sum_{|\sigma|=k+1} \frac{|D^{\sigma}P_m(0)|}{\sigma!} |x|^{|\sigma|}\\
    \leq & K\eta^{m(k+1)}+CK|x|^{k+1}+CK|x|^{k+1}|\ln \eta^{m+1}|\\
    \leq &CK|x|^{k+1}+CK|x|^{k+1}\big|\ln |x|\big|\\
\leq &CK|x|^{k+1}\big|\ln |x|\big|.
  \end{aligned}
\end{equation*}
Therefore, $u\in C^{k,\mathrm{lnL}}(0)$. ~\qed~\\

We also need the following two facts. The first will be used in the proof of the interior $C^{k,\alpha}$ regularity to construct a sequence of appropriate polynomials (see \cref{e6.8}). The second is for the boundary $C^{k,\alpha}$ regularity.
\begin{lemma}\label{le2.2}
Let $a^{ij}$ be a symmetric matrix whose eigenvalues lie in $[\lambda,\Lambda]$. Define the linear operator $L\colon \mathcal{HP}_{k+2}\to \mathcal{HP}_{k}$ as
\begin{equation*}
L(P)\coloneqq a^{ij}P_{ij}.
\end{equation*}
Then $L$ is surjective for any $k\geq 0$. That is, for any $P\in \mathcal{HP}_{k}$, there exists $Q\in \mathcal{HP}_{k+2}$ such that
\begin{equation*}
  a^{ij} Q_{ij}=P
\end{equation*}
and
\begin{equation*}
\|Q\|\leq C\|P\|,
\end{equation*}
where $C$ depends only on $n,\lambda,\Lambda$ and $k$.
\end{lemma}
\proof The proof is taken from \cite{Misc_1}. Up to a coordinate system transformation, we can assume $a^{ij}$ is the unit matrix without loss of generality. Then $L(P)=\Delta P$.


For $1\leq l\leq k$, denote
\begin{equation*}
\mathcal{HP}_k^l\coloneqq \left\{P\in \mathcal{HP}_k \colon P(x)=\sum_{\sigma} \frac{a_{\sigma}}{\sigma!} x^{\sigma} \mbox{ where } |\sigma|=k\mbox{ and }\max(\sigma_1,...,\sigma_n)=l\right\}.
\end{equation*}
Clearly, any $P\in \mathcal{HP}_k$ can be expressed as a sum of elements from $\mathcal{HP}_k^l$ ($1\leq l\leq k$). Since $\Delta$ is linear, we only need to prove the lemma for each $\mathcal{HP}_k^l$ ($1\leq l\leq k$). We prove this by induction on $l$. If $l=k$, any $P\in \mathcal{HP}_k^k$ can be written as
\begin{equation*}
P(x)=\sum_{i=1}^{n}a_ix_i^k.
\end{equation*}
where $a_i$ are constants. Hence, we can choose
\begin{equation*}
Q(x)=\sum_{i=1}^{n}\frac{a_i}{(k+2)(k+1)}x_i^{k+2}
\end{equation*}
and the conclusion holds.

Suppose that the conclusion holds for all $l\geq l_0+1$ and we need to prove the conclusion for $l_0$. Since $\Delta$ is a linear operator, we only need to consider that $P\in \mathcal{HP}_k^{l_0}$ is a monomial. Without loss of generality, we assume
\begin{equation*}
P(x)=ax^{\sigma},~|\sigma|=k,\sigma_1=l_0.
\end{equation*}
Let
\begin{equation*}
Q_1(x)=\frac{a}{(l_0+2)(l_0+1)} x^{\sigma+2e_1}
=\frac{a}{(l_0+2)(l_0+1)} x_1^{l_0+2}x_2^{\sigma_2}\cdots x_n^{\sigma_n}.
\end{equation*}
Then
\begin{equation*}
\Delta Q_1=P+\frac{a}{(l_0+2)(l_0+1)} x_1^{l_0+2}\Delta (x_2^{\sigma_2}\cdots x_n^{\sigma_n})\coloneqq P+P_2.
\end{equation*}
Note that $P_2\in \mathcal{HP}_{k}^{l_0+2}$. By induction, there exists $Q_2\in \mathcal{HP}_{k+2}$ such that $\Delta Q_2=P_2$. Let $Q=Q_1+Q_2$. Then
\begin{equation*}
\Delta Q=P
\end{equation*}
and
\begin{equation*}
\|Q\|\leq C\|P\|.
\end{equation*}
By induction, the proof is completed.~\qed~\\

\begin{lemma}\label{le2.6}
Let $a^{ij}$ be a symmetric matrix whose eigenvalues lie in $[\lambda,\Lambda]$. Then for any $P\in \mathcal{SP}_{k}$ ($k\geq 1$), there exists $Q\in \mathcal{SP}_{k+2}$ such that
\begin{equation*}
L(Q)\coloneqq  a^{ij} Q_{ij}=P
\end{equation*}
and
\begin{equation*}
\|Q\|\leq C\|P\|,
\end{equation*}
where $C$ depends only on $n,\lambda,\Lambda$ and $k$.
\end{lemma}
\proof The proof is similar to that of \Cref{le2.2}. We assume $L(P)=\Delta P$. For $1\leq l\leq k$, denote
\begin{equation*}
\mathcal{SP}_k^l\coloneqq \left\{P\in \mathcal{SP}_k \colon \sigma_n=l\right\}.
\end{equation*}
Similar to \Cref{le2.2}, we prove the lemma for each $\mathcal{SP}_k^l$ ($1\leq l\leq k$) by induction on $l$.

If $l=1$, any $P\in \mathcal{SP}_k^1$ can be written as
\begin{equation*}
P(x)=x_n\tilde{P}(x'),
\end{equation*}
where $\tilde{P}\in \mathcal{HP}_{k-1}$. By \Cref{le2.2}, there exists $\tilde{Q}\in \mathcal{HP}_{k+1}$ such that
\begin{equation*}
\Delta \tilde{Q}=\tilde{P}.
\end{equation*}
Hence, we can choose
\begin{equation*}
Q=x_n\tilde{Q}
\end{equation*}
and then
\begin{equation*}
\Delta Q=P.
\end{equation*}

Suppose that the conclusion holds for all $l\leq l_0-1$ and we need to prove the conclusion for $l_0$. Note that any $P\in \mathcal{SP}_k^{l_0}$ can be written as
\begin{equation*}
P(x)=x_n^{l_0}\tilde{P}(x'),
\end{equation*}
where $\tilde{P}\in \mathcal{HP}_{k-l_0}$. By \Cref{le2.2} again, there exists $\tilde{Q}\in \mathcal{HP}_{k-l_0+2}$ such that
\begin{equation*}
\Delta \tilde{Q}=\tilde{P}.
\end{equation*}
Let
\begin{equation*}
Q_1=x_n^{l_0}\tilde{Q}(x')
\end{equation*}
and then
\begin{equation*}
\Delta Q_1=l_0(l_0+1)x_n^{l_0-2}\tilde{Q}+x_n^{l_0}\tilde{P}.
\end{equation*}
By induction, there exists $Q_2\in \mathcal{SP}_{k+2}$ such that
\begin{equation*}
\Delta Q_2=l_0(l_0+1)x_n^{l_0-2}\tilde{Q}.
\end{equation*}
Hence, we can choose $Q=Q_1-Q_2$ and then
\begin{equation*}
\Delta Q=P
\end{equation*}
and
\begin{equation*}
\|Q\|\leq C\|P\|.
\end{equation*}
By induction, the proof is completed.~\qed~\\

\section{Interior \texorpdfstring{$C^{1,\alpha}$}{C1,a} regularity}\label{In-C1a-mu}

From now on, we state our main results and prove them respectively in the following sections. Before proceeding forward, we set some conventions.

\begin{itemize}
\item For interior pointwise regularity, we always consider the equation in $B_1$ and concerns the regularity for the solutions at the center $0$;
\item For boundary pointwise regularity, we always consider the equation in $\Omega\cap B_1$ and investigate the regularity for the solutions at $0\in \partial \Omega$;
\item In both cases, we always assume that $r_0=1$ in \Cref{d-f}, \Cref{d-2} and \Cref{d-re};
\item If we say that $\partial \Omega\in C^{k,\alpha}(0)$, it always indicates that \cref{e-re} and \cref{e-re2} hold with $P(0)=0$ and $DP(0)=0$.
\end{itemize}

To prove regularity results, it is necessary to make assumptions on the oscillation of $F$ in $x$, besides the structure condition \cref{SC2}. Suppose that we consider an equation in a domain $\Omega$. Let $x_0\in \bar{\Omega}$ be the point where we intend to study the regularity of the solution. For the $C^{1,\alpha}$ regularity, we introduce the following condition to control the oscillation of $F$ in $x$ near $x_0$:~\\
There exist a fully nonlinear operator $G$ and $r_0>0$ such that for any $M\in \mathcal{S}^n$ and $x\in \bar\Omega\cap B_{r_0}(x_0)$,
\begin{equation}\label{e.C1a.beta}
|F(M,0,0,x)-G(M)|\leq \beta_1(x,x_0)|M|+\gamma_1(x,x_0) ,
\end{equation}
where $G(0)=0$ and $\beta_1,\gamma_1\geq 0$. We also require that $G$ satisfies the same structure condition as $F$ (i.e. \cref{SC2} or \cref{SC1}) (similarly hereinafter).

If we study pointwise regularity, we always assume (similarly for subsequent higher regularity) that $x_0=0$ and $r_0=1$ in \cref{e.C1a.beta}. For simplicity, we denote $\beta_1(x)=\beta_1(x,0)$ and $\gamma_1(x)=\gamma_1(x,0)$. We remark here that if $F$ is continuous in $x$, we may use $G(M)\coloneqq F(M,0,0,x_0)$ to measure the oscillation of $F$ in $x$ in \cref{e.C1a.beta}.

Clearly, linear equation \cref{e.linear} satisfies \cref{e.C1a.beta} for $\beta_1(x,x_0)=|a^{ij}(x)-a^{ij}_{x_0}|$ and $\gamma_1\equiv 0$, where $a^{ij}_{x_0}$ is a constant matrix with eigenvalues lying in $[\lambda,\Lambda]$.

Now, we state the interior pointwise $C^{1,\alpha}$ regularity.
\begin{theorem}\label{t-C1a-i}
Let $0<\alpha<\bar{\alpha}$ and $u$ be a viscosity solution of
\begin{equation}\label{e.2.1}
F(D^2u, Du, u,x)=f \quad\mbox{in}~~ B_1.
\end{equation}
Suppose that $F$ satisfies \cref{SC2} and \cref{e.C1a.beta}. Assume that
\begin{equation*}
  \begin{aligned}
    &b\in L^{p}(B_1),\quad c\in C^{-1,\alpha}(0),\quad
    \|\beta_1\|_{C^{-1,1}(0)}\leq \delta_0,\quad\gamma_1\in C^{-1,\alpha}(0),\quad f\in C^{-1,\alpha}(0),
  \end{aligned}
\end{equation*}
where $p=n/(1-\alpha)$ and $0<\delta_0<1$ (small) depends only on $n,\lambda,\Lambda$ and $\alpha$.

Then $u\in C^{1,\alpha}(0)$, i.e., there exists $P\in \mathcal{P}_1$ such that
\begin{equation}\label{e.C1a-1-i-mu}
  |u(x)-P(x)|\leq C |x|^{1+\alpha}, ~\forall ~x\in B_{1}
\end{equation}
and
\begin{equation}\label{e.C1a-2-i-mu}
|Du(0)|\leq C,
\end{equation}
where $C$ depends only on $n,\lambda,\Lambda,\alpha,\mu,\|b\|_{L^p(B_1)}, \|c\|_{C^{-1,\alpha}(0)},\omega_0,\|\gamma_1\|_{C^{-1,\alpha}(0)}$, $\|f\|_{C^{-1,\alpha}(0)}$ and $\|u\|_{L^{\infty }(B_1)}$.

In particular, if $F$ satisfies \cref{SC1}, we have the following explicit estimates
\begin{equation}\label{e.C1a-1-i}
  |u(x)-P(x)|\leq C |x|^{1+\alpha}\left(\|u\|_{L^{\infty }(B_1)}+\|f\|_{C^{-1,\alpha}(0)}+\|\gamma_1\|_{C^{-1,\alpha}(0)}\right), ~\forall ~x\in B_{1},
\end{equation}
and
\begin{equation}\label{e.C1a-2-i}
|Du(0)|\leq C \left(\|u\|_{L^{\infty}(B_1)}
+\|f\|_{C^{-1,\alpha}(0)}+\|\gamma_1\|_{C^{-1,\alpha}(0)}\right),
\end{equation}
where $C$ depends only on $n, \lambda, \Lambda,\alpha,\|b\|_{L^p(B_1)}$ and $\|c\|_{C^{-1,\alpha}(0)}$.
\end{theorem}

\begin{remark}\label{r-1}
The $C^{1,\alpha}$ estimate is also called Cordes-Nirenberg type estimate (see \cite[P. 2]{MR1351007}). The name originates from the work of Cordes \cite{MR91400} and Nirenberg \cite{MR0064986,MR0066532}.
For fully nonlinear elliptic equations, Trudinger \cite{MR931007} derived interior $C^{1,\alpha}$ regularity for Lipschitz continuous viscosity solutions. Caffarelli \cite{MR1005611,MR1351007} proved the interior pointwise $C^{1,\alpha}$ regularity for $C$-viscosity solutions of equations without lower order terms (i.e., $\mu= b= c= 0$). \'{S}wi\k{e}ch \cite{MR1606359} obtained interior $C^{1,\alpha}$ regularity for $\mu=0$, $b,c\in L^{\infty}$ and less general $\omega_0$.

For equations with quadratic growth in the gradient, Trudinger \cite{MR701522} obtained $C^{1,\alpha}$ a priori estimate. Nornberg \cite{MR3980853} proved $C^{1,\alpha}$ (for \emph{some} $0<\alpha<\bar{\alpha}$) regularity with $c\in L^{\infty}$. \Cref{t-C1a-i} is an extension of the result in \cite{MR3980853}.

Recently, da Silva and Nornberg \cite{MR4304555} derived interior $C^{1,\alpha}$ regularity for equations with more general nonlinear growth in the gradient. In particular, the following equation is a model:
\begin{equation*}
\mathcal{M}_{\lambda,\Lambda}^{+}(D^2u)+\mu(x)|Du|^m+b(x)|Du|=f(x),~m\in (0,2].
\end{equation*}

Silvestre and Teixeira \cite{MR3494911} studied fully nonlinear equations with the operator $F$ being asymptotic convex. They obtained interior $C^{1,\alpha}$ regularity for any $0<\alpha<1$ for this kind of operator. We point out that although the results in \cite{MR1606359,MR3980853,MR4304555,MR3494911} are stated in the form of local regularity, their proofs are applicable to derive pointwise regularity.
\end{remark}

\begin{remark}\label{r-2}
Note that the constant $\bar{\alpha}$ comes from \Cref{l-3modin1}. The condition $b\in L^p$ can be replaced by a more general one:
\begin{equation*}
\|b\|_{L^{p_0}(B_r)}\leq K r^{\alpha-1+n/p_0},~\forall ~1<r<1,
\end{equation*}
where $p_0>n$. The reason for requiring $b\in L^p$ for some $p>n$ is the dependence on the closedness (\Cref{l-35}) and the H\"{o}lder regularity (\Cref{l-3Ho}) for viscosity solutions.

If $\bar{\alpha}\leq \alpha<1$, \Cref{t-C1a-i} can also be proved by a similar proof upon an additional assumption that $F$ and $G$ are convex in $M$. Similarly, interior and boundary pointwise $C^{k,\alpha}$ ($\bar{\alpha}\leq \alpha<1,k\geq 1$) can be obtained by strengthening the assumption on $F$ correspondingly.
\end{remark}

\begin{remark}\label{r-2.4}
Since we study fully nonlinear elliptic equations in general form, the statement of the theorem seems a little complicated. In fact, the assumptions can be divided into three groups (similarly hereinafter): the structure condition \cref{SC2}, the oscillation of $F$ in $x$ \cref{e.C1a.beta} and the prescribed data appearing in the equation \cref{e.2.1} (i.e., $f$). The \cref{SC2} describes the dependence of $F$ on $M,p$ and $s$, and \cref{e.C1a.beta} describes the dependence of $F$ on $x$. Both of them are necessary for regularity of fully nonlinear elliptic equations. Hence, we need to make necessary assumptions on the parameters appearing in \cref{SC2} (i.e., $\mu,b,c$ and $\omega_0$) and  \cref{e.C1a.beta} (i.e., $\beta_1$ and $\gamma_1$). For boundary regularity, the prescribed data also contain $g$ and $\partial \Omega$ (see boundary regularity in later sections).
\end{remark}
~\\

In the following, we give the proof of \Cref{t-C1a-i}. First, we show that the solution can be approximated by a linear function in certain scale provided that the coefficients and the prescribed data are small enough.

\begin{lemma}\label{In-l-C1a-mu}
Suppose that $F$ satisfies \cref{e.C1a.beta}. For any $0<\alpha<\bar{\alpha}$, there exists $\delta>0$ depending only on $n,\lambda,\Lambda$ and $\alpha$ such that if $u$ satisfies
\begin{equation*}
F(D^2u, Du, u,x)=f \quad\mbox{in}~~ B_1
\end{equation*}
with
\begin{equation*}
  \begin{aligned}
&\max\left(\|u\|_{L^{\infty}(B_1)},\omega_0(1,1)\right)\leq 1,\\
&\max\left(\mu,\|b\|_{L^{p}(B_1)},\|c\|_{L^{n}(B_1)},\|\beta_1\|_{L^n(B_1)},
\|\gamma_1\|_{L^n(B_1)},\|f\|_{L^{n}(B_1)}\right)\leq \delta,\\
  \end{aligned}
\end{equation*}
then there exists $P\in \mathcal{P}_1$ such that
\begin{equation}\label{In-e.l4.0-mu}
  \|u-P\|_{L^{\infty}(B_{\eta})}\leq \eta^{1+\alpha}
\end{equation}
and
\begin{equation}\label{In-e.14.2-mu}
|DP(0)|\leq \bar{C},
\end{equation}
where $\bar{C}$ comes from \Cref{l-3modin1} and $0<\eta<1$ depends only on $n,\lambda,\Lambda$ and $\alpha$.
\end{lemma}


\begin{remark}\label{In-r-41-mu}
\Cref{In-l-C1a-mu} is sometimes called the ``key lemma'' (or ``key step'')  in the proof of regularity (see \cite[Lemma 1.3]{Wang_Regularity}). One may use a solution of a homogenous equation to approximate $u$ (e.g. \cite[Lemma 8.2]{MR1351007} and \cite[Lemma 3.7]{MR3980853}). Although the compactness method is applied there, they need the solvability of some equation, which is avoided in our proof.
\end{remark}
~\\

\noindent\textbf{Proof of \Cref{In-l-C1a-mu}.} We prove the lemma by contradiction. Suppose that the lemma is false. Then there exist $0<\alpha<\bar{\alpha}$ and a sequence of $(F_m,u_m, f_m)_{m\in \mathbb{N}}$ satisfying
\begin{equation*}
F_m(D^2u_m, Du_m, u_m,x)=f_m \quad\mbox{in}~~ B_1.
\end{equation*}
In addition, $F_m$ satisfy the structure condition \cref{SC2} (with $\lambda,\Lambda,\mu_m,b_m,c_m,\omega_m$) and \cref{e.C1a.beta} (with $G_{m},\beta_m,\gamma_m$). Furthermore,
\begin{equation}\label{e4.1}
  \begin{aligned}
&\max\left(\|u_m\|_{L^{\infty}(B_1)},\omega_m(1,1)\right)\leq 1,\\
&\max\left(\mu_m,\|b_m\|_{L^{p}(B_1)},\|c_m\|_{L^{n}(B_1)},
\|\beta_m\|_{L^n(B_1)},\|\gamma_m\|_{L^n(B_1)},\|f_m\|_{L^{n}(B_1)}\right)\leq \frac{1}{m},\\
  \end{aligned}
\end{equation}
Finally, we have
\begin{equation}\label{In-e.lC1a.1-mu}
  \|u_m-P\|_{L^{\infty}(B_{\eta})}> \eta^{1+\alpha},~ \forall~P~~\mbox{with}~~|DP(0)|\leq \bar{C},
\end{equation}
where $0<\eta<1$ is taken small such that
\begin{equation}\label{In-e.lC1a.2-mu}
\bar{C}\eta^{\bar{\alpha}-\alpha}<1/2.
\end{equation}

Clearly, $u_m$ are uniformly bounded (i.e., $\|u_m\|_{L^{\infty}(B_1)}\leq 1$). By \Cref{le2.5},
\begin{equation*}
u_m\in S^*(\lambda,\Lambda,\mu_m,b_m,|\tilde{f}_m|),\quad\tilde{f}_m\coloneqq |f_m|+|F_m(0,0,0,\cdot)|
+c_m\omega_m(1,1).
\end{equation*}
From \cref{e.C1a.beta} and $G_m(0)=0$,
\begin{equation*}
|F_m(0,0,0,x)|\leq \gamma_m(x),~\forall ~x\in B_1.
\end{equation*}
Hence,
\begin{equation*}
|\tilde{f}_m|\leq |f_m|+\gamma_m+c_m.
\end{equation*}
Therefore, by \Cref{l-3Ho}, there exists $\alpha_0>0$ such that for any $\Omega'\subset \subset B_1$,
\begin{equation*}
\|u_m\|_{C^{\alpha_0}(\bar{\Omega}')}\leq C,
\end{equation*}
where $C$ is independent of $m$.

Since $C^{\alpha_0}(\bar{\Omega}')$ is compactly embedded into $C(\bar{\Omega}')$, there exist a subsequence (denoted by $u_m$ again) and $\bar u\colon B_1\rightarrow \mathbb{R}$ such that
\begin{equation}\label{e3.1}
u_m\rightarrow \bar u\quad \mbox{ in }~ L^{\infty}_{\mathrm{loc}}(B_1).
\end{equation}
The structure condition \cref{SC2} implies that $G_{m}$ are Lipschitz continuous in $M$ with a uniform Lipschitz constant depending only on $n,\lambda$ and $\Lambda$. Note that $G_{m}(0)=0$. Then there exist a subsequence (denoted by $G_{m}$ again) and $\bar{G}\colon \mathcal{S}^n\rightarrow \mathbb{R}$ such that
\begin{equation}\label{e4.2}
G_{m}\rightarrow  \bar{G}~\quad\mbox{ in }~L^{\infty}_{\mathrm{loc}}(\mathcal{S}^n).
\end{equation}

Next, for any ball $B\subset\subset B_1$ and $\varphi\in C^{2}(\bar{B})$, let $\psi_m=F_m(D^2\varphi,D\varphi,u_m,x)-f_m(x)$ and $\psi=\bar{G}(D^2\varphi)$. Then
\begin{equation}\label{In-gkC1a-mu}
  \begin{aligned}
\psi_m-\psi=&F_m(D^2\varphi,D\varphi,u_m,x)-\bar{G}(D^2\varphi)-f_m\\
=&F_m(D^2\varphi,D\varphi,u_m,x)-F_m(D^2\varphi,0,0,x)\\
    &+F_m(D^2\varphi,0,0,x)-G_{m}(D^2\varphi)+G_{m}(D^2\varphi)-\bar{G}(D^2\varphi)-f_m.
  \end{aligned}
\end{equation}
By the structure condition \Cref{SC2},
\begin{equation*}
  \begin{aligned}
|F_m(D^2\varphi,D\varphi,u_m,x)-F_m(D^2\varphi,0,0,x)|\leq & \mu_m|D\varphi|^2+b_m|D\varphi|
     +c_m\omega_m(1,|u_m|).
  \end{aligned}
\end{equation*}
From \cref{e.C1a.beta},
\begin{equation*}
  \begin{aligned}
|F_m(D^2\varphi,0,0,x)-G_{m}(D^2\varphi)|\leq\beta_m|D^2\varphi|+\gamma_m.
  \end{aligned}
\end{equation*}
By combining above inequalities together and H\"{o}lder inequality,
\begin{equation*}
  \begin{aligned}
 \|\psi_m-&\psi\|_{L^n(B)}\\
 \leq&\mu_m \||D\varphi|^2\|_{L^n(B)}
+\|b_m\|_{L^p(B)}\|D\varphi\|_{L^{\frac{np}{p-n}}(B)}+\|c_m\|_{L^n(B)}\omega_m(1,1)\\
&+\|\beta_m\|_{L^n(B)} \|D^2\varphi\|_{L^{\infty}(B)}+\|\gamma_m\|_{L^n(B)} \\
&+\|G_{m}(D^2\varphi)-\bar{G}(D^2\varphi)\|_{L^n(B)}+\|f_m\|_{L^n(B)}\\
\leq& C(\mu_m+\|b_m\|_{L^p(B)}+\|c_m\|_{L^n(B)}+\|\beta_m\|_{L^n(B)}+\|\gamma_m\|_{L^n(B)}
+\|f_m\|_{L^n(B)})\\
&+\|G_{m}(D^2\varphi)-\bar{G}(D^2\varphi)\|_{L^n(B)}.
  \end{aligned}
\end{equation*}
Thus, with the aid of \cref{e4.1} and \cref{e4.2}, $\|\psi_m-\psi\|_{L^n(B)}\rightarrow 0$ as $m\rightarrow \infty$. Therefore, by \Cref{l-35}, $\bar u$ is a viscosity solution of
\begin{equation*}
\bar{G}(D^2\bar u)=0 \quad\mbox{in}~~ B_1.
\end{equation*}

By \Cref{l-3modin1}, there exists $\bar P\in \mathcal{P}_1$ such that
\begin{equation*}
  |\bar u(x)-\bar P(x)|\leq \bar{C} |x|^{1+\bar{\alpha}}, ~~\forall ~x\in B_{1}
\end{equation*}
and
\begin{equation*}
  |D\bar P(0)|\leq \bar{C}.
\end{equation*}
Hence, by noting \cref{In-e.lC1a.2-mu}, we have
\begin{equation}\label{In-e.lC1a.3-mu}
  \|\bar u-\bar P\|_{L^{\infty}(B_{\eta})}\leq \frac{1}{2}\eta^{1+\alpha}.
\end{equation}
However, from \cref{In-e.lC1a.1-mu},
\begin{equation*}
  \|u_m-\bar P\|_{L^{\infty}(B_{\eta})}> \eta^{1+\alpha}.
\end{equation*}
Let $m\rightarrow \infty$, we have
\begin{equation*}
    \|\bar u-\bar P\|_{L^{\infty}(B_{\eta})}\geq \eta^{1+\alpha},
\end{equation*}
which contradicts with \cref{In-e.lC1a.3-mu}.  ~\qed~\\

\begin{remark}\label{re3.1}
\Cref{In-l-C1a-mu} is the most important step towards the interior $C^{1,\alpha}$ regularity. The compactness of $u_m$ is vital to conclude that $u_m$ converges (see \cref{e3.1}). Hence, this method is called the method of compactness.
\end{remark}
~\\

Now, we can prove the interior $C^{1,\alpha}$ regularity.~\\
\noindent\textbf{Proof of \Cref{t-C1a-i}.} We make some normalization first such that \Cref{In-l-C1a-mu} can be applied. Let $\delta$ be the constant as in \Cref{In-l-C1a-mu} and we take $\delta_0=\delta$ ($\delta_0$ is from the assumptions of \Cref{t-C1a-i}). We can assume that
\begin{equation}\label{In-e.tC1a-ass-mu}
\begin{aligned}
&\|u\|_{L^{\infty}(B_1)}\leq 1,\\
&\mu\leq \frac{\delta}{8C_0^2},\quad\|b\|_{L^p(B_1)}\leq \frac{\delta}{8C_0}, \quad\|c\|_{C^{-1,\alpha}(0)}\leq \frac{\delta}{4},\quad\omega_0(1+C_0,C_0)\leq1,\\
&\|\beta_1\|_{C^{-1,1}(0)}\leq \delta,\quad\|\gamma_1\|_{C^{-1,\alpha}(0)}\leq \frac{\delta}{4},\quad\|f\|_{C^{-1,\alpha}(0)}\leq \delta,
\end{aligned}
\end{equation}
where $C_0$ is a constant (depending only on $n,\lambda,\Lambda$ and $\alpha$) to be specified later.

Otherwise, we consider for some $0<\rho<1$,
\begin{equation}\label{e.4.1}
  \bar{u}(y)=\frac{u(x)-u(0)}{\rho^{\alpha_0}},
\end{equation}
where $y=x/\rho$ and $0<\alpha_0<1$ is a H\"{o}lder exponent (depending only on $n,\lambda,\Lambda$ and $\|b\|_{L^p(B_1)}$) such that $u\in C^{2\alpha_0}(0)$ (by \Cref{l-3Ho}). Then $\bar{u}$ satisfies
\begin{equation*}
\bar{F}(D^2\bar{u}, D\bar{u}, \bar{u},y)=\bar{f} \quad\mbox{in}~~ B_1,
\end{equation*}
where  for $(M,p,s,y)\in \mathcal{S}^n\times \mathbb{R}^n\times \mathbb{R}\times \bar B_1$,
\begin{equation*}\label{In-e.tC1a-n1-mu}
\begin{aligned}
  &\bar{F}(M,p,s,y)
  =\rho^{2-\alpha_0}F(\rho^{\alpha_0-2}M,\rho^{\alpha_0-1}p,\rho^{\alpha_0}s+u(0),x),
  \quad\bar{f}(y)=\rho^{2-\alpha_0}f(x).
  \end{aligned}
\end{equation*}

It is easy to check that $\bar F$ satisfies the structure condition \cref{SC2} with $\lambda$, $\Lambda$, $\bar{\mu}$, $\bar{b}$, $\bar{c}$ and $\bar{\omega}_0$, where
\begin{equation*}
\begin{aligned}
&\bar{\mu}=\rho^{\alpha_0}\mu,\quad\bar{b}(y)=\rho b(x),\quad\bar{c}(y)=\rho c(x), \quad
  \bar{\omega}_0(\cdot,\cdot)=\rho^{1-\alpha_0}\omega_0(\cdot+|u(0)|,\cdot).
  \end{aligned}
\end{equation*}
In addition, we take $\bar{G}(M)=\rho^{2-\alpha_0}G(\rho^{\alpha_0-2}M)$ for any $M\in \mathcal{S}^n$, where $G$ is from \cref{e.C1a.beta}. Then $\bar{G}$ is uniformly elliptic with constants $\lambda$ and $\Lambda$, and
\begin{equation*}
  \begin{aligned}
|\bar{F}(M,0,0,y)-\bar{G}(M)|=&\rho^{2-\alpha_0}|F(\rho^{\alpha_0-2}M,0,u(0),x)
-G(\rho^{\alpha_0-2}M)|\\
\leq &\rho^{2-\alpha_0}|F(\rho^{\alpha_0-2}M,0,u(0),x)
-F(\rho^{\alpha_0-2}M,0,0,x)|\\
&+\rho^{2-\alpha_0}|F(\rho^{\alpha_0-2}M,0,0,x)-G(\rho^{\alpha_0-2}M)|\\
\leq & \rho^{2-\alpha_0}c(x)\omega_0(|u(0)|,|u(0)|)
+\beta_1(x)|M|+\rho^{2-\alpha_0}\gamma_1(x)\\
\coloneqq &\bar\beta_1(y)|M|+\bar\gamma_1(y),
  \end{aligned}
\end{equation*}
where
\begin{equation*}
\bar\beta_1(y)=\beta_1(x),\quad
\bar\gamma_1(y)=\rho^{2-\alpha_0}\gamma_1(x)+\rho^{2-\alpha_0}c(x)\omega_0(|u(0)|,|u(0)|).
\end{equation*}

From above arguments,
\begin{equation*}
  \begin{aligned}
\|\bar{u}\|_{L^{\infty}(B_1)}\leq& \rho^{-\alpha_0}\rho^{2\alpha_0}[u]_{C^{2\alpha_0}(0)}
    =\rho^{\alpha_0}[u]_{C^{2\alpha_0}(0)},\quad \bar\mu=\rho^{\alpha_0}\mu,\\
\|\bar{b}\|_{L^{p}(B_1)}=&\rho^{1-\frac{n}{p}}\|b\|_{L^{p}(B_{\rho})}\leq \rho^{\alpha}\|b\|_{L^{p}(B_1)},\qquad  ~~\|\bar{c}\|_{C^{-1,\alpha}(0)}\leq \rho^{\alpha}\|c\|_{C^{-1,\alpha}(0)},\\
\bar{\omega}_0(1+C_0,1)=&\rho^{1-\alpha_0}\omega_0(1+C_0+|u(0)|,1),~~ \qquad~~
\|\bar\beta_1\|_{C^{-1,1}(0)}\leq \|\beta_1\|_{C^{-1,1}(0)}\leq \delta,\\
\|\bar{\gamma}_1\|_{C^{-1,\alpha}(0)}\leq&\rho^{1-\alpha_0+\alpha}(\|\gamma_1\|_{C^{-1,\alpha}(0)}
    +\omega_0(|u(0)|,|u(0)|)\|c\|_{C^{-1,\alpha}(0)}),\\
\|\bar{f}\|_{C^{-1,\alpha}(0)}\leq&\rho^{1-\alpha_0+\alpha}\|f\|_{C^{-1,\alpha}(0)}.
     \end{aligned}
\end{equation*}
Therefore, by taking $\rho$ small enough (depending only on $n$, $\lambda$, $\Lambda$, $\alpha$, $\mu$, $\|b\|_{L^{p}(B_1)}$, $\|c\|_{C^{-1,\alpha}(0)}$, $\omega_0$,
$\|\gamma_1\|_{C^{-1,\alpha}(0)}$, $\|f\|_{C^{-1,\alpha}(0)}$ and $\|u\|_{L^{\infty}(B_1)}$), the assumptions in \cref{In-e.tC1a-ass-mu} for $\bar{u}$ can be guaranteed. Then the regularity of $u$ can be derived from that of $\bar{u}$. Hence, without loss of generality, we assume that \cref{In-e.tC1a-ass-mu} holds for $u$.

Now, we prove that $u$ is $C^{1,\alpha}$ at $0$. By \Cref{le2.1}, we only need to prove the following. There exist a sequence of $P_m\in \mathcal{P}_1$ ($m\geq -1$, $P_{-1}\equiv 0$) such that for all $m\geq 0$,
\begin{equation}\label{In-e.tC1a-3-mu}
\|u-P_m\|_{L^{\infty }(B _{\eta^{m}})}\leq \eta ^{m(1+\alpha )}
\end{equation}
and
\begin{equation}\label{In-e.tC1a-4-mu}
|P_m(0)-P_{m-1}(0)|+\eta^m|DP_m(0)-DP_{m-1}(0)|\leq \bar{C}\eta^{m(1+\alpha)},
\end{equation}
where $\eta$ is as in \Cref{In-l-C1a-mu}.

We prove \cref{In-e.tC1a-3-mu} and \cref{In-e.tC1a-4-mu} by induction. For $m=0$, by setting $P_0\equiv 0$, the conclusion holds clearly. Suppose that the conclusion holds for $m\leq m_0$. We need to prove that the conclusion holds for $m=m_0+1$.

Let $r=\eta ^{m_{0}}$, $y=x/r$ and
\begin{equation}\label{In-e.tC1a-v1-mu}
  v(y)=\frac{u(x)-P_{m_0}(x)}{r^{1+\alpha}}.
\end{equation}
Then $v$ satisfies
\begin{equation}\label{In-e.C1as-F-mu}
\tilde{F}(D^2v, Dv, v,y)=\tilde{f} \quad\mbox{in}~~ B_1,
\end{equation}
where for $(M,p,s,y)\in \mathcal{S}^n\times \mathbb{R}^n\times \mathbb{R}\times \bar B_1$,
\begin{equation*}
  \begin{aligned}
    &\tilde{F}(M,p,s,y)=r^{1-\alpha}F(r^{\alpha-1}M,r^{\alpha}p+DP_{m_0},r^{1+\alpha}s
    +P_{m_0}(x),x),\quad \tilde{f}(y)=r^{1-\alpha}f(x).
  \end{aligned}
\end{equation*}
In addition, we define $\tilde{G}(M)=r^{1-\alpha}G(r^{\alpha-1}M)$ for $M\in \mathcal{S}^n$.

In the following, we show that \cref{In-e.C1as-F-mu} satisfies the assumptions of \Cref{In-l-C1a-mu}. First, it is easy to verify that
\begin{equation*}
  \begin{aligned}
\|v\|_{L^{\infty}(B_1)}\leq& 1, ~\quad(\mathrm{by}~ \cref{In-e.tC1a-3-mu} ~\mathrm{and}~ \cref{In-e.tC1a-v1-mu})\\
 \|\tilde{f}\|_{L^{n}(B_1)}=& r^{-\alpha}\|f\|_{L^{n}(B_r)}\leq \|f\|_{C^{-1,\alpha}(0)}\leq \delta, ~\quad(\mathrm{by}~ \cref{In-e.tC1a-ass-mu})\\
 \tilde{G}(0)=&r^{1-\alpha}G(0)=0.
  \end{aligned}
\end{equation*}

By \cref{In-e.tC1a-4-mu}, there exists a constant $C_0$ depending only on $n,\lambda,\Lambda$ and $\alpha$ such that
\begin{equation*}
|P_{m}(0)|+|DP_m(0)|\leq C_0,~\forall~0\leq m\leq m_0.
\end{equation*}
Then it is easy to verify that $\tilde{F}$ and $\tilde{G}$ satisfy the structure condition \cref{SC2} with $\lambda,\Lambda,\tilde{\mu}$, $\tilde{b}$, $\tilde{c}$ and $\tilde{\omega}_0$, where
\begin{equation*}\label{In-e.tC1a-n2-mu}
  \begin{aligned}
    \tilde{\mu}=r^{1+\alpha}\mu,\quad\tilde{b}(y)=rb(x)+2C_0r\mu,\quad
    \tilde{c}(y)=r^{1-\alpha}c(x),\quad
    \tilde{\omega}_0(\cdot,\cdot)=\omega_0(\cdot+C_0,\cdot).
  \end{aligned}
\end{equation*}
Hence, by combining with \cref{In-e.tC1a-ass-mu},
\begin{equation*}
  \begin{aligned}
\tilde{\mu}=& r^{1+\alpha}\mu\leq \delta,\\
\|\tilde{b}\|_{L^{p}(B_1)}\leq& r^{1-\frac{n}{p}}\|b\|_{L^{p}(B_r)}+2C_0|B_1|^{1/p}r\mu\leq r^{\alpha}\|b\|_{L^{p}(B_1)}+4C_0\mu \leq \delta,\\
\|\tilde{c}\|_{L^n(B_1)}=&r^{-\alpha}\|c\|_{L^n(B_r)}\leq  \|c\|_{C^{-1,\alpha}(0)}\leq \delta,\\
\tilde{\omega}_0(1,1)=&\omega_0(1+C_0,1)\leq \delta.\\
  \end{aligned}
\end{equation*}

Finally, we compute the oscillation of $\tilde{F}$ in $y$:
\begin{equation*}\label{In-e.tC1a-n2-mu-2}
  \begin{aligned}
|\tilde{F}(M&,0,0,y)-\tilde{G}(M)|\\
=&r^{1-\alpha}|F(r^{\alpha-1}M,DP_{m_0},P_{m_0}(x),x)
-G(r^{\alpha-1}M)|\\
=&r^{1-\alpha}|F(r^{\alpha-1}M,DP_{m_0},P_{m_0}(x),x)
-F(r^{\alpha-1}M,0,0,x)\\
&+F(r^{\alpha-1}M,0,0,x)-G(r^{\alpha-1}M)|\\
\leq & r^{1-\alpha} (C_0^2\mu+C_0 b(x)+c(x)\omega_0(C_0,C_0))
+\beta_1(x)|M|+r^{1-\alpha}\gamma_1(x)\\
\coloneqq &\tilde{\beta}_1(y)|M|+\tilde{\gamma}_1(y),
\end{aligned}
\end{equation*}
where
\begin{equation*}
\tilde{\beta}_1(y)=\beta_1(x),\quad
\tilde{\gamma}_1(y)=r^{1-\alpha}\gamma_1(x)+ r^{1-\alpha} \big(C_0^2\mu+C_0 b(x)+c(x)\omega_0(C_0,C_0)\big).\\
\end{equation*}
Then with the aid of \cref{In-e.tC1a-ass-mu},
\begin{equation*}\label{In-e.tC1a-n22-mu}
  \begin{aligned}
\|\tilde{\beta}_1\|_{L^n(B_1)}=&r^{-1}\|\beta_1\|_{L^n(B_r)}
    \leq \|\beta_1\|_{C^{-1,1}(0)}\leq \delta,\\
\|\tilde{\gamma}_1\|_{L^n(B_1)}\leq& r^{-\alpha}\|\gamma_1\|_{L^n(B_r)}+
    C_0^2|B_1|^{1/n}r^{1-\alpha}\mu+C_0|B_1|^{\alpha/n}\|b\|_{L^{p}(B_r)}\\
    &+\omega_0(C_0,C_0)r^{-\alpha}\|c\|_{L^n(B_r)}\\
\leq & \|\gamma_1\|_{C^{-1,\alpha}(0)}+2C_0^2\mu
    +2C_0\|b\|_{L^{p}(B_1)}+\omega_0(C_0,C_0)\|c\|_{C^{-1,\alpha}(0)} \leq \delta.\\
  \end{aligned}
\end{equation*}

Therefore, \cref{In-e.C1as-F-mu} satisfies the assumptions of \Cref{In-l-C1a-mu} and there exists $\tilde{P}\in \mathcal{P}_1$ such that
\begin{equation*}
\begin{aligned}
    \|v-\tilde{P}\|_{L^{\infty }(B _{\eta})}&\leq \eta ^{1+\alpha}
\end{aligned}
\end{equation*}
and
\begin{equation*}
|D\tilde{P}(0)|\leq \bar{C}.
\end{equation*}

Let $P_{m_0+1}(x)=P_{m_0}(x)+r^{1+\alpha}\tilde{P}(y)$. Then \cref{In-e.tC1a-4-mu} holds for $m_0+1$. By recalling \cref{In-e.tC1a-v1-mu}, we have
\begin{equation*}
  \begin{aligned}
\|u-P_{m_0+1}\|_{L^{\infty}(B_{\eta^{m_0+1}})}&= \|u-P_{m_0}-r^{1+\alpha}\tilde{P}(x/r)\|_{L^{\infty}(B_{\eta r})}\\
&= \|r^{1+\alpha}v-r^{1+\alpha}\tilde{P}\|_{L^{\infty}(B_{\eta})}\\
&\leq r^{1+\alpha}\eta^{1+\alpha}=\eta^{(m_0+1)(1+\alpha)}.
  \end{aligned}
\end{equation*}
Hence, \cref{In-e.tC1a-3-mu} holds for $m=m_0+1$. By induction, the proof is completed.

Finally, we consider the special case, i.e., $F$ satisfies \cref{SC1}. Set
\begin{equation*}
K=\|u\|_{L^{\infty}(B_1)}
+\delta^{-1}\left(\|f\|_{C^{-1,\alpha}(0)}+4\|\gamma\|_{C^{-1,\alpha}(0)}\right)
\end{equation*}
and define for $0<\rho<1$
\begin{equation*}
  \bar{u}(y)=\frac{u(x)-u(0)}{K},
\end{equation*}
where $y=x/\rho$. By taking $\rho$ small enough (depending only on $n,\lambda,\Lambda,\alpha, \|b\|_{L^{p}(B_1)}$ and $\|c\|_{C^{-1,\alpha}(0)}$), the assumptions \cref{In-e.tC1a-ass-mu} can be guaranteed. Then, by the same proof, we have the explicit estimates \cref{e.C1a-1-i} and \cref{e.C1a-2-i}. \qed~\\

\begin{remark}\label{r-4.3}
The proof is also called ``scaling argument'' in the regularity theory (see \cite[P. 9 and P. 10]{Wang_Regularity}. Essentially, it is just a sequence of repetitions of \Cref{In-l-C1a-mu} in different scales.

\Cref{t-C1a-i} is optimal in the sense that the conditions imposed on the coefficients and the prescribed data are minimal. The reason is that these conditions are only used to ensure that the equation is scaling invariant. Precisely, after scaling, the equation (see \cref{In-e.C1as-F-mu}) and the solution (see \cref{In-e.tC1a-v1-mu}) satisfy the assumptions of \Cref{In-l-C1a-mu} as well. This scaling invariance is the most basic requirement in the regularity theory for uniformly elliptic equations.
\end{remark}

\begin{remark}\label{r-4.7}
Taking the transformation form \cref{e.4.1} is motivated by \cite{MR3980853} (see the definition of $\tilde{u}$ in Claim 3.2), which has also been used in \cite{MR1139064} (see Section 1.3).
\end{remark}
~\\

\section{Interior \texorpdfstring{$C^{2,\alpha}$}{C2,a} regularity}\label{In-C2a-mu}

In this section, we prove the interior pointwise $C^{2,\alpha}$ regularity. For the $C^{1,\alpha}$ ($0<\alpha<1$) regularity, assuming $b\in L^p$ and $c\in L^n$ is appropriate in \cref{SC2}. For $C^{1,\mathrm{lnL}}$ regularity (see \Cref{t-C1ln-i}) and higher regularity, $b,c\in L^p$ is not enough. Instead, we always assume that \cref{SC2} holds with $b,c\in L^{\infty}$ without loss of generality, i.e.,
\begin{equation}\label{e5.6}
b\equiv b_0, \quad c\equiv c_0 \quad\mbox{ in } \cref{SC2}
\end{equation}
for some positive constants $b_0,c_0$.

In addition, for the $C^{2,\alpha}$ regularity, we require that $\omega_0$ satisfies the following homogeneous condition: there exists a constant $K_0>0$ such that
\begin{equation}\label{e.omega0}
\omega_0(\cdot,r s)\leq K_0r^{\alpha}\omega_0(\cdot,s), ~\forall ~0<r<1,s>0.
\end{equation}
The above condition is necessary. For example, if we consider \cref{e.linear-1} with $h(u)=u^q$, we must require that $q\geq \alpha$ for the $C^{2,\alpha}$ regularity.

Similar to the $C^{1,\alpha}$ regularity, we use the following to estimate the oscillation of $F$ in $x$ near $x_0\in\bar{\Omega}$:~\\
There exist a fully nonlinear operator $G$ and $r_0>0$ such that for any
$(M,p,s)\in \mathcal{S}^n\times \mathbb{R}^n\times \mathbb{R}$ and
$x\in \bar\Omega\cap B_{r_0}(x_0)$,
\begin{equation}\label{e.C2a-KF}
|F(M,p,s,x)-G(M,p,s)|\leq \beta_2(x,x_0)(|M|+1)\omega_2(|p|,|s|),
\end{equation}
where $G(0,0,0)=0$, $\beta_2\geq 0$ and $\beta_2(x,x)\equiv 0$. Here, $\omega_2\geq 0$ is non-decreasing in each variable.

Note that if $F$ satisfies \cref{e.C2a-KF}, it satisfies \cref{e.C1a.beta} with
\begin{equation*}
G(M)\coloneqq G(M,0,0), \quad  \beta_1(x,x_0)\coloneqq \gamma_1(x,x_0)\coloneqq \omega_2(0,0)\beta_2(x,x_0).
\end{equation*}
This is important since we will use the $C^{1,\alpha}$ regularity in the proof of $C^{2,\alpha}$ regularity. For the $C^{1,\alpha}$ regularity, the assumption on $\beta_1$ is different from that on $\gamma_1$ and hence we adopt \cref{e.C1a.beta}. For the $C^{2,\alpha}$ regularity, \cref{e.C2a-KF} is enough.

This paper also considers the following special case for which we can obtain explicit estimates:
\begin{equation}\label{e.C2a-KF-0}
|F(M,p,s,x)-G(M,p,s)|\leq \beta_2(x,x_0) \left(|M|+|p|+|s|\right)+\gamma_2(x,x_0),
\end{equation}
where $\gamma_2\geq 0$ with $\gamma_2(x,x)\equiv 0$.

Take the linear equation \cref{e.linear} for example. The corresponding operator is
\begin{equation*}
F(M,p,s,x)\coloneqq a^{ij}(x)M_{ij}+b^{i}(x)p_i+c(x)s.
\end{equation*}
Assume that $a^{ij},b^{i},c$ are continuous at $x_0$. Let
\begin{equation*}
G(M,p,s)\coloneqq a^{ij}(x_0)M_{ij}+b^{i}(x_0)p_i+c(x_0)s.
\end{equation*}
Then
\begin{equation*}
|F(M,p,s,x)-G(M,p,s)|\leq \beta_2(x,x_0)(|M|+|p|+|s|),
\end{equation*}
where
\begin{equation*}
\beta_2(x,x_0)=\max \left(|a^{ij}(x)-a^{ij}(x_0)|,|b^{i}(x)-b^{i}(x_0)|,|c(x)-c(x_0)|\right).
\end{equation*}

Now, we state the interior pointwise $C^{2,\alpha}$ regularity.
\begin{theorem}\label{t-C2a-i}
Let $0<\alpha<\bar{\alpha}$ and $u$ be a viscosity solution of
\begin{equation*}
F(D^2u, Du, u,x)=f \quad\mbox{in}~~ B_1.
\end{equation*}
Suppose that $F$ satisfies \Cref{SC2} and \cref{e.C2a-KF} with $G$ being convex in $M$. Assume that $\omega_0$ satisfies \cref{e.omega0}, $\beta_2\in C^{\alpha}(0)$ and $f\in C^{\alpha}(0)$.

Then $u\in C^{2,\alpha}(0)$, i.e., there exists $P\in \mathcal{P}_2$ such that
\begin{equation}\label{e.C2a-1-i-mu}
  |u(x)-P(x)|\leq C |x|^{2+\alpha}, ~~\forall ~x\in B_{1},
\end{equation}
\begin{equation}\label{e.C2a-3-i-mu}
|F(D^2P,DP(x),P(x),x)-f(0)|\leq C|x|^{\alpha}, ~~\forall ~x\in B_{1}
\end{equation}
and
\begin{equation}\label{e.C2a-2-i-mu}
|Du(0)|+|D^2u(0)|\leq C,
\end{equation}
where $C$ depends only on $n,\lambda, \Lambda,\alpha,\mu,b_0,c_0,\omega_0,
\|\beta_2\|_{C^{\alpha}(0)},\omega_2,\|f\|_{C^{\alpha}(0)}$ and $\|u\|_{L^{\infty }(B_1)}$.

In particular, if $F$ satisfies \cref{SC1} and \cref{e.C2a-KF-0} with $\gamma_2\in C^{\alpha}(0)$,
\begin{equation}\label{e.C2a-1-i}
|u(x)-P(x)|\leq \tilde{C} |x|^{2+\alpha}, ~~\forall ~x\in B_{1},
\end{equation}
\begin{equation}\label{e.C2a-3-i}
|F(D^2P,DP(x),P(x),x)-f(0)|\leq \tilde{C}|x|^{\alpha}, ~~\forall ~x\in B_{1},
\end{equation}
\begin{equation}\label{e.C2a-2-i}
|Du(0)|+|D^2u(0)|\leq \tilde{C}
\end{equation}
and
\begin{equation*}
\tilde{C}=C\left(\|u\|_{L^{\infty }(B_1)}+\|f\|_{C^{\alpha}(0)}
+\|\gamma_2\|_{C^{\alpha}(0)}\right),
\end{equation*}
where $C$ depends only on $n,\lambda, \Lambda,\alpha, b_0,c_0$ and $\|\beta_2\|_{C^{\alpha}(0)}$.
\end{theorem}

\begin{remark}\label{r-2.1}
This is the classical $C^{2,\alpha}$ regularity (Schauder estimate). The first $C^{2,\alpha}$ a priori estimate for fully nonlinear elliptic equations was obtained by Evans \cite{MR649348} and Krylov \cite{MR661144}. Safonov \cite{MR765302} proved interior $C^{2,\alpha}$ estimate for the Bellman's equation. Trudinger \cite{MR769173} derived $C^{2,\alpha}$ regularity for $W^{2,n}$ strong solutions. The pointwise $C^{2,\alpha}$ regularity for viscosity solutions was developed by Caffarelli \cite{MR1005611,MR1351007}.

The convexity of $G$ is only used to apply \Cref{l-3modin2}, which is not necessary in dimension $2$ (see (i) in \Cref{re2.1}). The convexity assumption was relaxed by Caffarelli and Yuan \cite{MR1793687}. There are several other extensions. For a special fully nonlinear operator
\begin{equation*}
F(M)=\min(F_1(M),F_2(M)), ~\forall ~M\in \mathcal{S}^n,
\end{equation*}
where $F_1$ ($F_2$) is concave (convex) in $M$, Cabr\'{e} and Caffarelli \cite{MR1995493} proved interior $C^{2,\alpha}$ regularity. Savin \cite{MR2334822} (see also \cite[Section 4]{MR2928094}) obtained the interior $C^{2,\alpha}$ regularity if $F\in C^{1}$ and $\|u\|_{L^{\infty}}$ is small, which is called small perturbation regularity. This idea has been used to estimate the Hausdorff dimension of singular set for fully nonlinear elliptic equations (e.g. \cite{MR4556790} for free boundary problem and \cite{MR2928094} for partial regularity result). Cao, Li and Wang \cite{MR2775423} gave the interior $C^{2,\alpha}$ regularity under the assumption that $F\in C^{1,\tilde{\alpha}}$ ($0<\tilde\alpha\leq 1$) and $u\in C^2$. Niu and Wu \cite{MR4503817} derived interior $C^{2,\alpha}$ regularity if $\Lambda$ is close to $\lambda$, which extends the result of Bhattacharya and Warren \cite{MR4251799}.

In regards to equations with quadratic growth in the gradient, Trudinger \cite{MR701522} obtained $C^{2,\alpha}$ a priori estimate. Recently, da Silva and Nornberg \cite{MR4304555} obtained interior local $C^{2,\alpha}$ regularity for equations with more general nonlinear growth in the gradient.

As far as we know, \Cref{t-C2a-i} is the first interior pointwise $C^{2,\alpha}$ regularity for equations with quadratic growth in the gradient.
\end{remark}

\begin{remark}\label{re1.3}
In above theorem, the assumptions on $\beta_2$ and $f$ can be replaced with the following weaker ones:
\begin{equation*}
\|\beta_2\|_{L^n(B_r)}\leq Cr^{1+\alpha},\quad
\|f-f_0\|_{L^n(B_r)}\leq Cr^{1+\alpha},~\forall ~0<r<1,
\end{equation*}
where $f_0$ is a constant. To make the statement clear, we work under the assumptions in \Cref{t-C2a-i}.

If $F(0,0,0,x)\equiv 0$, we can take $\gamma_2\equiv 0$ in \cref{e.C2a-KF-0} to estimate the oscillation of $F$ in $x$. Then we obtain estimates similar to \crefrange{e.C2a-1-i}{e.C2a-2-i} without ``$\gamma_2$'' involved in the right-hand side.

As in \cite{MR1005611,MR1351007}, we can assume the following instead of the condition ``$G$ is convex in $M$'':\\
Any viscosity solution $v$ of
\begin{equation*}
G(D^2v)=0\quad\mbox{ in}~B_1
\end{equation*}
belongs to $C^{2,\bar{\alpha}}(\bar{B}_{1/2})$ and
\begin{equation*}
\|v\|_{C^{2,\bar{\alpha}}(\bar{B}_{1/2})}\leq \bar{C}\|v\|_{L^{\infty}(B_1)}.
\end{equation*}
\end{remark}

\begin{remark}\label{r-2.12}
For the linear equation \cref{e.linear}, if the coefficients $a^{ij},b^{i}$ and $c$ are $C^{\alpha}$ at $0$, then $\beta_2\in C^{\alpha}(0)$. Hence, the classical Schauder estimate for linear equations is a special case of \Cref{t-C2a-i}.
\end{remark}
~\\

The following is the ``key step'' for the interior pointwise $C^{2,\alpha}$ regularity, which is similar to \Cref{In-l-C1a-mu}.
\begin{lemma}\label{In-l-C2a-mu}
Suppose that $F$ satisfies \cref{e.C2a-KF} with $G$ being convex in $M$. For any $0<\alpha<\bar{\alpha}$, there exists $\delta>0$ depending only on $n,\lambda,\Lambda,\alpha$ and $\omega_2$ such that if $u$ satisfies
\begin{equation*}
F(D^2u, Du, u,x)=f \quad\mbox{in}~~ B_1
\end{equation*}
with
\begin{equation*}
  \begin{aligned}
&u(0)=|Du(0)|=0,\quad\max\left(\|u\|_{L^{\infty}(B_1)},\omega_0(1,1)\right)\leq 1,~\\
&\max\left(\mu,b_0,c_0,\|\beta_2\|_{L^{\infty}(B_1)},\|f\|_{L^{\infty}(B_1)}\right)\leq \delta, \\
  \end{aligned}
\end{equation*}
then there exists $P\in \mathcal{HP}_2$ such that
\begin{equation*}
  \|u-P\|_{L^{\infty}(B_{\eta})}\leq \eta^{2+\alpha},
\end{equation*}
\begin{equation*}
G(D^2P,0,0)=0
\end{equation*}
and
\begin{equation*}
\|P\|\leq \bar{C}+1,
\end{equation*}
where $0<\eta<1$ depends only on $n,\lambda,\Lambda$ and $\alpha$.
\end{lemma}

\begin{remark}\label{r-4.5}
Note that by \Cref{t-C1a-i}, $u\in C^{1,\alpha}(0)$. Hence, $Du(0)$ is well defined.

If we only intend to derive that there exists $P\in \mathcal{P}_2$ such that
\begin{equation*}
  \|u-P\|_{L^{\infty}(B_{\eta})}\leq \eta^{2+\alpha},
\end{equation*}
it is enough to make the same assumptions as that of \Cref{In-l-C1a-mu} and the proof is almost exact as that of \Cref{In-l-C1a-mu}. This indicates that we can approximate the solution in some scale only if the coefficients and the prescribed data are small in some ``soft'' norms. The requirements that $\beta_2\in C^{\alpha}(0)$ and $f\in C^{\alpha}(0)$ are only used in the following ``scaling argument'' (see \Cref{In-t-C2as-mu} below).
\end{remark}
~\\

\noindent\textbf{Proof of \Cref{In-l-C2a-mu}.} Similar to \Cref{In-l-C1a-mu}, we prove the lemma by contradiction. Suppose that the lemma is false. Then there exist $0<\alpha<\bar{\alpha},\omega_2$ and a sequence of $(F_m,u_m,f_m)_{m\in \mathbb{N}}$ satisfying
\begin{equation*}
F_m(D^2u_m, Du_m, u_m,x)=f_m \quad\mbox{in}~~ B_1.
\end{equation*}
In addition, $F_m$ satisfy the structure condition \cref{SC2} (with $\lambda,\Lambda,\mu_m,b_m,c_m,\omega_m$) and \cref{e.C2a-KF} (with $G_{m},\beta_m,\omega_2$). Furthermore, $G_{m}$ are convex in $M$ and
\begin{equation*}
  \begin{aligned}
&u_m(0)=|Du_m(0)|=0,\quad\max\left(\|u_m\|_{L^{\infty}(B_1)},\omega_m(1,1)\right)\leq 1,\\
&\max\left(\mu_m,b_m,c_m,\|\beta_{m}\|_{L^{\infty}(B_1)},\|f_m\|_{L^{\infty}(B_1)}\right)\leq \frac{1}{m}.
  \end{aligned}
\end{equation*}
Finally, for any $P\in\mathcal{HP}_2$ satisfying $\|P\|\leq \bar{C}+1$ and
\begin{equation}\label{In-e.lC2a-3-mu}
G_{m}(D^2P,0,0)=0,
\end{equation}
we have
\begin{equation}\label{In-e.lC2a-1-mu}
  \|u_m-P\|_{L^{\infty}(B_{\eta})}> \eta^{2+\alpha},
\end{equation}
where $0<\eta<1$ is taken small such that
\begin{equation}\label{In-e.lC2a-2-mu}
\bar{C}\eta^{\bar{\alpha}-\alpha}<1/2.
\end{equation}

Similar to the $C^{1,\alpha}$ regularity (see the proof of \Cref{In-l-C1a-mu}), $u_m$ are uniformly bounded and equicontinuous. Then there exist a subsequence (denoted by $u_m$ again) and $\bar u\colon B_1\rightarrow \mathbb{R}$ such that
\begin{equation*}
u_m\rightarrow \bar u \quad\mbox{ in } L^{\infty}_{\mathrm{loc}}(B_1).
\end{equation*}
In addition, there exist a subsequence of $G_{m}$ (denoted by $G_{m}$ again) and $\bar G\colon \mathcal{S}^n\rightarrow \mathbb{R}$ such that
\begin{equation*}
G_{m}(\cdot,0,0)\rightarrow \bar{G}\quad \mbox{ in } L^{\infty}_{\mathrm{loc}}(\mathcal{S}^n).
\end{equation*}
Since $G_{m}$ are convex in $M$, $\bar G$ is convex in $M$ as well.

Next, for any ball $B\subset\subset B_1$ and $\varphi\in C^2(\bar{B})$, let $\psi_m=F_m(D^2\varphi,D\varphi,u_m,x)-f_m$, $\psi=\bar{G}(D^2\varphi)$ and $r=\|D^2\varphi\|_{L^{\infty}(B)}+\|D\varphi\|_{L^{\infty}(B)}+1$. Similar to the proof of \Cref{In-l-C1a-mu},
\begin{equation*}
  \begin{aligned}
\|\psi_m-\psi\|_{L^n(B)}=&\|F_m(D^2\varphi,D\varphi,u_m,x)-\bar{G}(D^2\varphi)-f_m\|_{L^n(B)}\\
=& \|F_m(D^2\varphi,D\varphi,u_m,x)-G_{m}(D^2\varphi,D\varphi,u_m)
+G_{m}(D^2\varphi,D\varphi,u_m)\\
&-G_{m}(D^2\varphi,0,0)+G_{m}(D^2\varphi,0,0)-\bar{G}(D^2\varphi)-f_m\|_{L^n(B)}\\
\leq& \|\beta_{m}(|D^2\varphi|+1)\omega_2(r,r)+r^2\mu_m+rb_m+c_m\omega_m(1,1)\|_{L^n(B)}\\
&+\|G_{m}(D^2\varphi,0,0)-\bar{G}(D^2\varphi)\|_{L^{n}(B)}+\|f_m\|_{L^n(B)}\rightarrow 0.
  \end{aligned}
\end{equation*}
By \Cref{l-35}, $\bar u$ is a viscosity solution of
\begin{equation}\label{e5.1}
\bar G(D^2\bar u)=0 \quad\mbox{in}~~ B_1.
\end{equation}

From the $C^{1,\alpha}$ estimate for $u_m$ (see \Cref{t-C1a-i}) and noting  $u_m(0)=0$ and $Du_m(0)=0$, we have
\begin{equation*}
|u_m(x)|\leq C|x|^{1+\alpha}, ~~\forall~x\in B_{1},
\end{equation*}
where $C$ depends only on $n,\lambda,\Lambda,\alpha$ and $\omega_2$. Since $u_m$ converges to $\bar u$,
\begin{equation*}
|\bar u(x)|\leq C|x|^{1+\alpha}, ~~\forall~x\in B_{1},
\end{equation*}
which implies
\begin{equation}\label{e5.5}
\bar u(0)=|D\bar u(0)|=0.
\end{equation}

By applying \Cref{l-3modin2} to \cref{e5.1} and noting \cref{e5.5}, there exists $\bar{P}\in\mathcal{HP}_2$ such that
\begin{equation}\label{In-e.c2a-3-mu}
  |\bar u(x)-\bar{P}(x)|\leq \bar{C} |x|^{2+\bar{\alpha}}, ~~\forall ~x\in B_{1},
\end{equation}
\begin{equation*}
  \bar G(D^2\bar{P})=0
\end{equation*}
and
\begin{equation*}
 \|\bar{P}\|\leq \bar{C}.
\end{equation*}
By Combining \cref{In-e.lC2a-2-mu} with \cref{In-e.c2a-3-mu}, we have
\begin{equation}\label{In-e.lC2a-4-mu}
  \|\bar u-\bar {P}\|_{L^{\infty}(B_{\eta})}\leq \frac{1}{2}\eta^{2+\alpha}.
\end{equation}

Since $G_{m}(D^2\bar{P},0,0)\rightarrow \bar G(D^2\bar{P})=0$, there exist $t_m\rightarrow 0$ and $|t_m|\leq 1$ (for $m$ large) such that
\begin{equation*}
  G_m(D^2P_m,0,0)=G_{m}(D^2\bar{P}+t_m\tilde{I},0,0)=0,
\end{equation*}
where $\tilde{I}$ denotes the matrix whose entries are all $0$ except $\tilde{I}_{nn}=1$ (see \Cref{no1.1}) and
\begin{equation*}
P_m(x)\coloneqq \bar {P}(x)+\frac{t_m}{2}x_n^2.
\end{equation*}
Moreover, $P_m\in\mathcal{HP}_2$ and $\|P_m\|\leq \bar{C}+1$.

Hence, \cref{In-e.lC2a-1-mu} holds for $P_m$, i.e.,
\begin{equation*}
  \|u_m-P_m\|_{L^{\infty}(B_{\eta})}> \eta^{2+\alpha}.
\end{equation*}
Let $m\rightarrow \infty$, we have
\begin{equation*}
    \|\bar u-\bar{P}\|_{L^{\infty}(B_{\eta})}\geq \eta^{2+\alpha},
\end{equation*}
which contradicts with \cref{In-e.lC2a-4-mu}.  ~\qed~\\

Next, we show the interior pointwise $C^{2,\alpha}$ regularity in a special case.
\begin{lemma}\label{In-t-C2as-mu}
Suppose that $F$ satisfies \cref{e.C2a-KF} with $G$ being convex in $M$. Let $0<\alpha <\bar{\alpha}$, $\omega_0$ satisfy \cref{e.omega0} and $u$ satisfy
\begin{equation*}
F(D^2u, Du, u,x)=f \quad\mbox{in}~~ B_1.
\end{equation*}
Assume that
\begin{equation}\label{In-e.C2as-be-mu}
\begin{aligned}
&\|u\|_{L^{\infty}(B_1)}\leq 1,\quad u(0)=|Du(0)|=0, \\
&\mu\leq \frac{\delta_1}{4C_0},\quad b_0\leq \frac{\delta_1}{2},\quad c_0\leq \frac{\delta_1}{K_0},\quad \omega_0(1+C_0,1)\leq 1,\\
&|\beta_2(x)|\leq \delta_1|x|^{\alpha},\quad|f(x)|\leq \delta_1|x|^{\alpha}, ~\forall ~x\in B_1,\\
\end{aligned}
\end{equation}
where $\delta_1$ depends only on $n,\lambda,\Lambda,\alpha,\omega_0$ and $\omega_2$, and $C_0$ depends only on $n,\lambda,\Lambda$ and $\alpha$.

Then $u\in C^{2,\alpha}(0)$, i.e., there exists $P\in\mathcal{HP}_2$ such that
\begin{equation}\label{In-e.C2as-1-mu}
  |u(x)-P(x)|\leq C |x|^{2+\alpha}, ~~\forall ~x\in B_{1},
\end{equation}
\begin{equation}\label{In-e.C2as-3-mu}
G(D^2P,0,0)=0
\end{equation}
and
\begin{equation}\label{In-e.C2as-2-mu}
\|P\| \leq C,
\end{equation}
where $C$ depends only on $n, \lambda, \Lambda$ and $\alpha$.
\end{lemma}
\proof As before, we only need to prove the following. There exist a sequence of $P_m\in\mathcal{HP}_2$ ($m\geq -1$) such that for all $m\geq 0$,

\begin{equation}\label{In-e.C2as-4-mu}
\|u-P_m\|_{L^{\infty }(B _{\eta^{m}})}\leq \eta ^{m(2+\alpha )},
\end{equation}
\begin{equation}\label{In-e.C2as-5-mu}
G(D^2P_m,0,0)=0
\end{equation}
and
\begin{equation}\label{In-e.C2as-6-mu}
\|P_m-P_{m-1}\|\leq (\bar{C}+1)\eta ^{(m-1)\alpha},
\end{equation}
where $0<\eta<1$ is as in \Cref{In-l-C2a-mu}.

We prove above conclusion by induction. Clearly, \crefrange{In-e.C2as-4-mu}{In-e.C2as-6-mu} hold for $m=0$ by setting $P_0\equiv P_{-1}\equiv 0$. Suppose that they hold for $m\leq m_0$. Let $r=\eta ^{m_{0}}$, $y=x/r$ and
\begin{equation}\label{In-e.C2as-v-mu}
  v(y)=\frac{u(x)-P_{m_0}(x)}{r^{2+\alpha}}.
\end{equation}
Then $v$ satisfies
\begin{equation}\label{In-e.C2as-F-mu}
 \tilde{F}(D^2v,Dv,v,y)=\tilde{f} \quad\mbox{in}~~ B_1,
\end{equation}
where for $(M,p,s,y)\in \mathcal{S}^n\times \mathbb{R}^n\times \mathbb{R}\times \bar B_1$,
\begin{equation*}
  \begin{aligned}
&\tilde{F}(M,p,s,y)=r^{-\alpha}F(r^{\alpha}M+D^2P_{m_0},r^{1+\alpha}p
+DP_{m_0}(x),r^{2+\alpha}s+P_{m_0}(x),x),~\\
&\tilde{f}(y)=r^{-\alpha}f(x).\\
  \end{aligned}
\end{equation*}
In addition, define
\begin{equation*}
\tilde{G}(M,p,s)=r^{-\alpha}G(r^{\alpha}M+D^2P_{m_0},r^{1+\alpha}p
,r^{2+\alpha}s).
\end{equation*}

In the following, we show that \cref{In-e.C2as-F-mu} satisfies the assumptions of \Cref{In-l-C2a-mu}. First, it is easy to verify that $\tilde{G}$ is convex in $M$ and
\begin{equation*}
\begin{aligned}
\|v\|_{L^{\infty}(B_1)}\leq& 1, \quad  v(0)=|Dv(0)|=0,
~\quad(\mathrm{by}~ \cref{In-e.C2as-be-mu},~ \cref{In-e.C2as-4-mu} ~\mbox{and}~ \cref{In-e.C2as-v-mu})\\
\|\tilde{f}\|_{L^{\infty}(B_1)}=&r^{-\alpha}\|f\|_{L^{\infty}(B_r)}\leq \delta_1,
 ~\quad(\mathrm{by}~ \cref{In-e.C2as-be-mu})\\
\tilde{G}(0,0,0)=&r^{-\alpha}G(D^2P_{m_0},0,0)=0.~\quad(\mathrm{by}~ \cref{In-e.C2as-5-mu})\\
\end{aligned}
\end{equation*}

By \cref{In-e.C2as-6-mu},
\begin{equation*}
\|P_m\|\leq C_0,~\forall~0\leq m\leq m_0,
\end{equation*}
where $C_0$ depends only on $n,\lambda,\Lambda$ and $\alpha$. Then $\tilde{F}$ and $\tilde{G}$ satisfy the structure condition \cref{SC2} with $\lambda,\Lambda,\tilde{\mu},\tilde{b},\tilde{c}$ and $\tilde\omega_0$, where (note that $\omega_0$ satisfies \cref{e.omega0})
\begin{equation*}
\tilde{\mu}=r^{2+\alpha}\mu,\quad\tilde{b}=rb_0+2C_0r\mu,\quad\tilde{c}=K_0 r^{\alpha+\alpha^2}c_0, \quad\tilde{\omega}_0(\cdot,\cdot)=\omega_0(\cdot+C_0,\cdot).
\end{equation*}
Hence, from \cref{In-e.C2as-be-mu},
\begin{equation*}
\begin{aligned}
 \tilde{\mu}\leq \mu\leq \delta_1,\quad\tilde{b}\leq b_0+2C_0\mu\leq \delta_1,\quad
 \tilde{c}\leq K_0 c_0\leq \delta_1, \quad\tilde{\omega}_0(1,1)\leq 1.
 \end{aligned}
\end{equation*}

Finally, we consider
\begin{equation}\label{e5.7}
\begin{aligned}
|\tilde{F}&(M,p,s,y)-\tilde{G}(M,p,s)|\\
=&r^{-\alpha}\big|F(r^{\alpha}M+D^2P_{m_0},r^{1+\alpha}p+DP_{m_0}(x),
r^{2+\alpha}s+P_{m_0}(x),x)\\
&-G(r^{\alpha}M+D^2P_{m_0},r^{1+\alpha}p,r^{2+\alpha}s)\big|\\
\leq& r^{-\alpha}\big|F(r^{\alpha}M+D^2P_{m_0},r^{1+\alpha}p
+DP_{m_0}(x),r^{2+\alpha}s+P_{m_0}(x),x)\\
&-F(r^{\alpha}M+D^2P_{m_0},r^{1+\alpha}p,r^{2+\alpha}s,x)\big|\\
&+r^{-\alpha}\big|F(r^{\alpha}M+D^2P_{m_0},r^{1+\alpha}p,r^{2+\alpha}s,x)
-G(r^{\alpha}M+D^2P_{m_0},r^{1+\alpha}p,r^{2+\alpha}s)\big|\\
\coloneqq& A_1+A_2.
\end{aligned}
\end{equation}
By the structure condition \cref{SC2} and \cref{In-e.C2as-be-mu},
\begin{equation*}
  \begin{aligned}
A_1\leq &r^{-\alpha}\left(2C_0r^{1+\alpha}\mu |p||x|+C_0^2\mu|x|^2+C_0 b_0|x|+K_0 c_0\omega_0(|s|+C_0,C_0)|x|^{2\alpha}\right)\\
\leq&\delta_1|p|+\delta_1C_0+\delta_1\omega_0(|s|+C_0,C_0)|y|^{\alpha}.
  \end{aligned}
\end{equation*}
Next, by \cref{e.C2a-KF} and \cref{In-e.C2as-be-mu},
\begin{equation*}
  \begin{aligned}
A_2\leq &r^{-\alpha}\beta_2(x)(|M|+C_0+1)\omega_2(|p|,|s|)
\leq\delta_1(|M|+C_0+1)\omega_2(|p|,|s|)|y|^{\alpha}.
  \end{aligned}
\end{equation*}
Hence,
\begin{equation*}
\begin{aligned}
|\tilde{F}&(M,p,s,y)-\tilde{G}(M,p,s)|\\
\leq&\big((|M|+C_0+1)\omega_2(|p|,|s|)+|p|+C_0+\omega_0(|s|+C_0,C_0)\big)\delta_1|y|^{\alpha}\\
\coloneqq & \tilde{\beta}_2(y)(|M|+1)\tilde{\omega}_2(|p|,|s|),
\end{aligned}
\end{equation*}
where
\begin{equation*}
 \tilde{\beta}_2(y)=\delta_1 |y|^{\alpha},\quad \tilde{\omega}_2(|p|,|s|)=(C_0+1)
 \omega_2(|p|,|s|)+|p|+C_0+\omega_0(|s|+C_0,C_0).
\end{equation*}
Then $\|\tilde{\beta}_2\|_{L^{\infty}(B_1)}\leq \delta_1$.


Choose $\delta_1$ small enough (depending only on $n,\lambda,\Lambda,\alpha,\omega_0$ and $\omega_2$) such that \Cref{In-l-C2a-mu} holds for $\tilde{\omega}_0,\tilde{\omega}_2$ and $\delta_1$. Since \cref{In-e.C2as-F-mu} satisfies the assumptions of \Cref{In-l-C2a-mu}, there exists $\tilde{P}(y)\in\mathcal{HP}_2$ such that
\begin{equation*}
\begin{aligned}
    \|v-\tilde{P}\|_{L^{\infty }(B_{\eta})}&\leq \eta ^{2+\alpha},
\end{aligned}
\end{equation*}
\begin{equation*}
\tilde{G}(D^2\tilde{P},0,0)=0
\end{equation*}
and
\begin{equation*}
\|\tilde{P}\|\leq \bar{C}+1.
\end{equation*}

Let
\begin{equation*}
P_{m_0+1}(x)=P_{m_0}(x)+r^{2+\alpha}\tilde{P}(y)=P_{m_0}(x)+r^{\alpha}\tilde{P}(x).
\end{equation*}
Then \cref{In-e.C2as-5-mu} and \cref{In-e.C2as-6-mu} hold for $m_0+1$. By recalling \cref{In-e.C2as-v-mu}, we have
\begin{equation*}
  \begin{aligned}
\|u-P_{m_0+1}\|_{L^{\infty}(B_{\eta^{m_0+1}})}&= \|u-P_{m_0}-r^{\alpha}\tilde{P}\|_{L^{\infty}(B_{\eta r})}\\
&= \|r^{2+\alpha}v-r^{2+\alpha}\tilde{P}\|_{L^{\infty}(B_{\eta})}\\
&\leq r^{2+\alpha}\eta^{2+\alpha}=\eta^{(m_0+1)(2+\alpha)}.
  \end{aligned}
\end{equation*}
Hence, \cref{In-e.C2as-4-mu} holds for $m=m_0+1$. By induction, the proof is completed.\qed~\\

Now, we give the~\\
\noindent\textbf{Proof of \Cref{t-C2a-i}.} In fact, \Cref{In-t-C2as-mu} has contained the essential ingredients for the $C^{2,\alpha}$ regularity. In some sense, the following proof is just a normalization procedure that makes the assumptions in \Cref{In-t-C2as-mu} satisfied. We prove \Cref{t-C2a-i} in two cases.

\textbf{Case 1:} the general case, i.e., $F$ satisfies \cref{SC2} and \cref{e.C2a-KF}. Throughout the proof for this case, $C$ always denotes a constant depending only on $n, \lambda,\Lambda,\alpha,\mu,b_0,c_0,\omega_0$, $\|\beta_2\|_{C^{\alpha}(0)}$, $\omega_2$, $\|f\|_{C^{\alpha}(0)}$ and $\|u\|_{L^{\infty }(B_1)}$. Let
\begin{equation*}
F_1(M,p,s,x)=F(M,p,s,x)-f(0),~\forall ~(M,p,s,x)\in \mathcal{S}^n\times \mathbb{R}^n\times \mathbb{R}\times \bar B_1.
\end{equation*}
Then $u_1$ satisfies
\begin{equation*}
F_1(D^2u,Du,u,x)=f_1\quad\mbox{in}~~B_1,
\end{equation*}
where $f_1(x)=f(x)-f(0)$. Hence,
\begin{equation*}
  |f_1(x)|\leq C|x|^{\alpha}, ~~\forall ~x\in B_1.
\end{equation*}

Since $u\in C^{1,\alpha}(0)$, set
\begin{equation*}
u_1=u-P_u, \quad F_2(M,p,s,x)=F_1(M,p+DP_u,s+P_u(x),x).
\end{equation*}
Recall that $P_u\in \mathcal{P}_1$ denotes the Taylor polynomial of $u$ at $0$ (see \Cref{d-f} and \Cref{re1.5}). Then $u_1$ satisfies
\begin{equation*}
F_2(D^2u_1,Du_1,u_1,x)=f_1\quad\mbox{in}~~B_1
\end{equation*}
and $u_1(0)=|Du_1(0)|=0$.

Next, for $\tau\in \mathbb{R}$ (to be specified later, see \cref{In-e.tC2a-2-mu}), take
\begin{equation*}
u_2=u_1-\tau x_n^2, \quad F_3(M,p,s,x)=F_2(M+2\tau\tilde{I},p+2\tau x_n,s+\tau x_n^2,x).
\end{equation*}
Then $u_2$ satisfies
\begin{equation*}
F_3(D^2u_2,Du_2,u_2,x)=f_1\quad\mbox{in}~~B_1
\end{equation*}
and $u_2(0)=|Du_2(0)|=0$.

We introduce fully nonlinear operators $G_{1},G_{2}$ and $G_{3}$ in a similar way as $F_1,F_2$ and $F_3$. To be more precise,
\begin{equation*}
  \begin{aligned}
&G_1(M,p,s)\coloneqq G(M,p,s)-f(0), \quad G_2(M,p,s)\coloneqq G_1(M,p+DP_u,s+P_u(0)),\\
&G_3(M,p,s)\coloneqq G_2(M+2\tau\tilde{I},p,s).
  \end{aligned}
\end{equation*}

By the structure condition \cref{SC2}, there exists $\tau\in \mathbb{R}$ such that $G_{3}(0,0,0)=0$ and
\begin{equation}\label{In-e.tC2a-2-mu}
  \begin{aligned}
|\tau|&\leq |G_{2}(0,0,0)|/\lambda\leq |G(0,Du(0),u(0))-f(0)|/\lambda\leq C.
  \end{aligned}
\end{equation}

For $0<\rho<1$, define $y=x/\rho$,
\begin{equation*}
  \begin{aligned}
&u_3(y)=\rho^{-1}u_2(x), \quad
F_4(M,p,s,y)=\rho F_3(\rho^{-1}M, p, \rho s,x),\\
&G_{4}(M,p,s)=\rho G_{3}(\rho^{-1}M, p,\rho s).\\
  \end{aligned}
\end{equation*}
Then $u_3$ satisfies
\begin{equation}\label{In-F5-mu}
F_4(D^2u_3,Du_3,u_3,y)=f_2\quad\mbox{in}~~B_1,
\end{equation}
where $f_2(y)=\rho f_1(x)$.

Now, we can check that \cref{In-F5-mu} satisfies the conditions of \Cref{In-t-C2as-mu} by choosing a proper $\rho$. First, it can be verified easily that
\begin{equation*}
  \begin{aligned}
    u_3(0)=|Du_3(0)|=0, \quad |f_2(y)|&= \rho|f_1(x)|\leq C\rho^{1+\alpha}|y|^{\alpha}~\quad\mbox{in}~B_1.
  \end{aligned}
\end{equation*}
Next, by the interior $C^{1,\alpha}$ regularity for $u$,
\begin{equation*}
  \|u_3\|_{L^{\infty}(B_1)}\leq \rho^{-1} \|u_2\|_{L^{\infty}(B_{\rho})}
  \leq  \rho^{-1}\left(\|u_1\|_{L^{\infty}(B_{\rho})}+C\rho^{2}\right)
  \leq \rho^{-1}\left(C\rho^{1+\alpha}+C\rho^{2}\right)\leq C\rho^{\alpha}.
\end{equation*}
Furthermore, $G_{4}(0,0,0)=\rho G_{3}(0,0,0)=0$ and $G_{4}$ satisfies the structure condition \cref{SC2} with $\lambda,\Lambda,\tilde{\mu},\tilde{b},\tilde{c}$ and $\tilde\omega_0$, where
\begin{equation*}
\tilde{\mu}= \rho\mu,\quad\tilde{b}= \rho b_0+C\rho\mu,\quad\tilde{c}= \rho^{1/2} c_0, \quad\tilde{\omega}_0(\cdot,\cdot)=\rho^{1/2} \omega_0(\cdot+C,\cdot).
\end{equation*}
Then $\tilde{\omega}_0$ satisfies \cref{e.omega0}.

Finally, we check the oscillation of $F_4$ in $y$.
\begin{equation*}
  \begin{aligned}
F_4(M&,p,s,y)-G_{4}(M,p,s)\\
=&\rho\Big(F(\rho^{-1}M
+2\tau\tilde{I},p+2\tau x_n+DP_u,\rho s+\tau x_n^2+P_u(x),x)\\
&-G(\rho^{-1}M+2\tau\tilde{I},p+DP_u,\rho s+P_u(0))\Big)\\
= &\rho\Big(F(\rho^{-1}M+2\tau\tilde{I},p+2\tau x_n+DP_u,
\rho s+\tau x_n^2+P_u(x),x)\\
&-F(\rho^{-1}M+2\tau\tilde{I},p+DP_u,\rho s+P_u(0),x)\Big)\\
&+\rho\Big(F(\rho^{-1}M+2\tau\tilde{I},p+DP_u,\rho s+P_u(0),x)\\
&-G(\rho^{-1}M+2\tau\tilde{I},p+DP_u,\rho s+P_u(0))\Big)\\
\coloneqq& A_1+A_2.
  \end{aligned}
\end{equation*}
By the structure condition \cref{SC2},
\begin{equation*}
  \begin{aligned}
|A_1|\leq& C\rho\mu(2|p|+C|x|+C)|x|+C\rho b_0|x|+K_0\rho c_0\omega_0(|s|+C,C)|x|^{\alpha}\\
\leq& C\rho|x|^{\alpha}\big(|p|+1+\omega_0(|s|+C,C)\big).
  \end{aligned}
\end{equation*}
By \cref{e.C2a-KF},
\begin{equation*}
  \begin{aligned}
|A_2|\leq& \rho\beta_2(x)(\rho^{-1}|M|+C)\omega_2(|p|+C,|s|+C)\\
\leq& C\rho^{\alpha}|y|^{\alpha}(|M|+1)\omega_2(|p|+C,|s|+C).
  \end{aligned}
\end{equation*}
Hence,
\begin{equation}\label{e5.8}
  \begin{aligned}
|F_4(M&,p,s,y)-G_{4}(M,p,s)|\\
\leq& C\rho|x|^{\alpha}\big(|p|+1+\omega_0(|s|+C,C)\big)
+C\rho^{\alpha}|y|^{\alpha}(|M|+1)\omega_2(|p|+C,|s|+C)\\
\leq& C\rho^{\alpha}|y|^{\alpha}(|M|+1)\Big(\omega_2(|p|+C,|s|+C)+|p|+1+\omega_0(|s|+C,C)\Big)\\
\coloneqq  &\tilde{\beta}_2(y)(|M|+1)\tilde{\omega}_2(|p|,|s|),
  \end{aligned}
\end{equation}
where
\begin{equation}\label{e5.9}
\tilde{\beta}_2(y)=C\rho^{\alpha}|y|^{\alpha},\quad
\tilde{\omega}_2(|p|,|s|)=\omega_2(|p|+C,|s|+C)+|p|+1+\omega_0(|s|+C,C).
\end{equation}

Take $\delta_1$ small enough such that \Cref{In-t-C2as-mu} holds with $\tilde{\omega}_0,\tilde{\omega}_2$ and $\delta_1$. From above arguments, we can choose $\rho$ small enough (depending only on $n, \lambda,\Lambda,\alpha,\mu,b_0,c_0$, $\omega_0,\|\beta_2\|_{C^{\alpha}(0)}$, $\omega_2$, $\|f\|_{C^{\alpha}(0)}$ and $\|u\|_{L^{\infty }(B_1)}$) such that
\begin{equation*}
\begin{aligned}
&\|u_3\|_{L^{\infty }(B_1)}\leq 1,\quad\tilde\mu \leq \frac{\delta_1}{4C_0},\quad
\tilde{b}\leq \frac{\delta_1}{2},\quad \tilde c\leq \frac{\delta_1}{K_0},\quad
\tilde{\omega}_0(1+C_0,1)\leq 1,\\
&|\tilde\beta_2(y)|\leq \delta_1|y|^{\alpha},\quad|f_2(y)|\leq \delta_1|y|^{\alpha}, ~\forall ~y\in B_1,
\end{aligned}
\end{equation*}
where $C_0$ depending only on $n,\lambda,\Lambda$ and $\alpha$, is as in \Cref{In-t-C2as-mu}. Therefore, the assumptions in \Cref{In-t-C2as-mu} are satisfied for \cref{In-F5-mu}. By \Cref{In-t-C2as-mu}, $u_3$ and hence $u$ is $C^{2,\alpha}$ at $0$, and the estimates \crefrange{e.C2a-1-i-mu}{e.C2a-2-i-mu} hold. In fact, from above proof, we obtain that there exists $P\in \mathcal{P}_2$ such that
\begin{equation*}
  G(D^2P,0,0)=f(0).
\end{equation*}
Since $F$ satisfies \cref{e.C2a-KF} and $G$ satisfies \cref{SC2},
\begin{equation*}
\begin{aligned}
|F(D^2&P,DP(x),P(x),x)-f(0)|\\
\leq& |F(D^2P,DP(x),P(x),x)-G(D^2P,DP(x),P(x))|\\
&+|G(D^2P,DP(x),P(x))-G(D^2P,0,0)|+|G(D^2P,0,0)-f(0)|\\
\leq &C |x|^{\alpha}, ~~\forall ~x\in B_{1}.
\end{aligned}
\end{equation*}
That is, \cref{e.C2a-3-i-mu} holds.

\textbf{Case 2:} $F$ satisfies \cref{SC1} and \cref{e.C2a-KF-0}. Let
\begin{equation*}
K=\|u\|_{L^{\infty }(B_1)}+\|f\|_{C^{\alpha}(0)}+\|\gamma_2\|_{C^{\alpha}(0)}, \quad u_1=u/K.
\end{equation*}
Then $u_1$ satisfies
\begin{equation}\label{e.5.1}
F_1(D^2u_1,Du_1,u_1,x)=f_1\quad\mbox{in}~~B_1,
\end{equation}
where $F_1(M,p,s,x)=F(KM,Kp,Ks,x)/K$ and $f_1=f/K$.

Obviously,
\begin{equation*}
\|u_1\|_{L^{\infty }(B_1)}\leq 1\quad\mbox{and}\quad\|f_1\|_{C^{\alpha}(0)}\leq 1.
\end{equation*}
In addition, $F_1$ satisfies the structure condition \cref{SC1} with the same $\lambda,\Lambda,b_0$ and $c_0$. Define $G_{1}(M,p,s)=G(KM,Kp,Ks)/K$. Then
\begin{equation*}
  \begin{aligned}
    |F_1(M,p,s,x)-G_{1}(M,p,s)|= & K^{-1}|F(KM,Kp,Ks,x)-G(KM,Kp,Ks)|\\
    \leq & \beta_2(x)\left(|M|+|p|+|s|\right)+\tilde\gamma_2(x),
  \end{aligned}
\end{equation*}
where $\tilde\gamma_2(x)=K^{-1}\gamma_2(x)$ and hence $\|\tilde{\gamma}_2\|_{C^{\alpha}(0)}\leq 1$.

By applying \textbf{Case 1} to \cref{e.5.1}, we obtain that $u_1$ and hence $u$ is $C^{2,\alpha}$ at $0$ and the estimates \crefrange{e.C2a-1-i}{e.C2a-2-i} hold. \qed~\\

\section{Interior \texorpdfstring{$C^{k,\alpha}$}{Ck,a} regularity}\label{In-Cka-mu}

In this section, we prove the interior pointwise $C^{k,\alpha}$ ($k\geq 3$) regularity. We impose the following necessary condition on $\omega_0$: for some constant $K_0>0$,
\begin{equation}\label{e.omega0-2}
\omega_0(\cdot,rs)\leq K_0r\omega_0(\cdot,s), ~\forall ~r>0,s>0,
\end{equation}
which is necessary for higher regularity. For the equation \cref{e.linear-1}, if $h\in C^{1}(\mathbb{R})$, \cref{e.omega0-2} holds.

Before stating the regularity result, we introduce the following definition to estimate the oscillation of $F$ in $x$.
\begin{definition}\label{d-FP}
Let $k\geq 1$, $\Omega\subset \mathbb{R}^n$ be a bounded domain, $\omega$ be a modulus of continuity and $F\colon \mathcal{S}^n\times \mathbb{R}^n\times \mathbb{R}\times \bar{\Omega}\rightarrow \mathbb{R}$. We say that $F$ is $C^{k,\omega}$ at $x_0$ or $F\in C^{k,\omega}(x_0)$ if the following holds:

There exist a fully nonlinear operator $G$ and constants $K,r_0>0$ such that for any
$(M,p,s)\in \mathcal{S}^n\times \mathbb{R}^n\times \mathbb{R}$ and
$x\in \bar\Omega\cap B_{r_0}(x_0)$,
\begin{equation}\label{holder}
\begin{aligned}
 |F(M,p,s,x)-G(M,p,s,x)|\leq&  K|x-x_0|^k\omega(|x-x_0|)(|M|+1)\omega_3(|p|,|s|),\\
\end{aligned}
\end{equation}
where~\\
(i) $G$ is convex in $M$;~\\
(ii) $G \in C^{k,\bar{\alpha}}(\mathcal{S}^n\times \mathbb{R}^n\times \mathbb{R}\times \bar{\Omega})$ (recall that $\bar{\alpha}$ originates from \Crefrange{l-3modin1}{l-32} and is fixed throughout this paper);~\\
(iii) there exists $K_1>0$ such that
\begin{equation*}
  \begin{aligned}
&\left|D_MG(\xi)-D_MG(\zeta)\right|+\left|D_pG(\xi)-D_pG(\zeta)\right|+\left|G_s(\xi)-G_s(\zeta)\right|\\
&\quad\leq K_1\left(|M-N|+|p-q|+|s-t|+|x-y|\right),\\
&~~\quad\forall ~\xi=(M,p,s,x),\zeta=(N,q,t,y)\in \mathcal{S}^n\times \mathbb{R}^n\times\mathbb{R}\times \bar{\Omega};
  \end{aligned}
\end{equation*}
(iv) $\omega_3\geq 0$ is non-decreasing in each variable.

If $\omega(r)=r^{\alpha}$ ($0<\alpha< \bar{\alpha}$), we call $F\in C^{k,\alpha}(x_0)$ and define
\begin{equation*}
  \|F\|_{C^{k,\alpha}(x_0)}=\min \left\{K\big | \cref{holder} ~\mbox{holds with}~K\right\}
\end{equation*}
and
\begin{equation*}
\omega_4(r)=\|G\|_{C^{k,\bar{\alpha}}(\bar{\mathbf{B}}_r\times \bar{\Omega})},~\forall ~r>0,
\end{equation*}
where
\begin{equation}\label{e6.2}
\mathbf{B}_r\coloneqq \left\{(M,p,s)\big| |M|+|p|+|s|<r\right\}.
\end{equation}
If $\bar{\alpha}\leq \alpha\leq 1$, we require that $G\in C^{k+1,\bar{\alpha}}$ and set
\begin{equation*}\label{e.cka.F}
\omega_4(r)=\|G\|_{C^{k+1,\bar{\alpha}}(\bar{\mathbf{B}}_r\times \bar{\Omega})},~\forall ~r>0.
\end{equation*}

If $F\in C^{k, \alpha}(x)$ for any $x\in \bar{\Omega}$ with the same $r_0,K_1,\omega_3$ and $\omega_4$, and
\begin{equation*}
\|F\|_{C^{k,\alpha}(\bar{\Omega})}\coloneqq\sup_{x_0\in \Omega} \|F\|_{C^{k,\alpha}(x_0)}<+\infty,
\end{equation*}
we say that $F\in C^{k, \alpha}(\bar{\Omega})$.

Similarly, we can define $F\in C^{k,\mathrm{Dini}}$ and $F\in C^{k,\mathrm{lnL}}$. If $\omega$ is only a modulus of continuity rather than a Dini function, we may say that $F\in C^{k}$.
\end{definition}

\begin{remark}\label{re1.4}
The condition (iii) looks weird. To prove $C^{k,\alpha}$ ($k\geq 3$) regularity, it is necessary that the solution of the equation with the ``good'' operator $G$ possesses the $C^{2,\alpha}$ regularity, and hence higher regularity since $G$ is smooth (see \Cref{In-l-62}). Thus, we need to propose a condition to guarantee this. It is more natural to propose the following condition:
\begin{equation*}
|G(M, p,s,x)-G(M,p,s,y)| \leq K_1(|M|+|p|+|s|+1)|x-y|.
\end{equation*}
However, we could not prove the pointwise higher regularity (\Cref{t-Cka-i}) under this condition. The reason is that the condition cannot be satisfied in the scaling argument (see the proof of \Cref{In-t-Ckas-mu}).

On the other hand, if $u\in C^2(\bar{B}_1)$ is a classical solution of
\begin{equation*}
F(D^2u,Du,u,x)=f(x)\quad\mbox{ in}~B_1,
\end{equation*}
and $F$ is a smooth operator, then we can redefine $F$ in $\mathbf{B}_{R}^c\times \bar{B}_1$ such that the new $F$ satisfies (iii) ($K_1$ depends on $\|u\|_{C^2(\bar{B}_1)}$) and $u$ is still a solution (see also \Cref{re8.1}), where $\mathbf{B}$ is defined in \cref{e6.2} and $R\coloneqq \|u\|_{C^2(\bar{B}_1)}$.
\end{remark}
~\\

Take the linear equation \cref{e.linear} for example again. If $a^{ij},b^{i},c\in C^{k,\alpha}(x_0)$, we can take $G$ as
\begin{equation*}
  G(M,p,s,x)=P_{a^{ij}}(x)M_{ij}+P_{b^{i}}(x)p_i+P_{c}(x)s,
\end{equation*}
where $P_{a^{ij}},P_{b^{i}}$ and $P_c$ are the Taylor polynomials of $a^{ij},b^{i},c$ at $x_0$ respectively (see \Cref{d-f}). Then $G \in C^{\infty}(\mathcal{S}^n\times \mathbb{R}^n\times \mathbb{R}\times \bar{\Omega})$ and $F\in C^{k,\alpha}(x_0)$. Roughly speaking, if the coefficients are $C^{k,\alpha}$ at $x_0$, the operator is $C^{k,\alpha}$ at $x_0$.

The following is the interior pointwise $C^{k,\alpha}$ ($k\geq 3$) regularity.
\begin{theorem}\label{t-Cka-i}
Let $k\geq 3$, $0<\alpha<\bar{\alpha}$ and $u$ be a viscosity solution of
\begin{equation*}
F(D^2u, Du, u,x)=f \quad\mbox{in}~~ B_1.
\end{equation*}
Suppose that $F\in C^{k-2,\alpha}(0)$ satisfies \cref{SC2}, $\omega_0$ satisfies \cref{e.omega0-2} and $f\in C^{k-2,\alpha}(0)$.

Then $u\in C^{k,\alpha}(0)$, i.e., there exists $P\in \mathcal{P}_k$ such that
\begin{equation}\label{e.Cka-1-i}
  |u(x)-P(x)|\leq C |x|^{k+\alpha}, ~~\forall ~x\in B_{1},
\end{equation}
\begin{equation}\label{e.Cka-3-i}
|F(D^2P(x),DP(x),P(x),x)-P_f(x)|\leq C|x|^{k-2+\alpha}, ~~\forall ~x\in B_{1}
\end{equation}
and
\begin{equation}\label{e.Cka-2-i}
|Du(0)|+\cdots+|D^ku(0)|\leq C,
\end{equation}
where $P_f$ is the Taylor polynomial of $f$ at $0$ and $C$ depends only on $k,n,\lambda, \Lambda,\alpha,\mu, b_0$, $c_0,\omega_0,\|F\|_{C^{k-2,\alpha}(0)},K_1,\omega_3,\omega_4$, $\|f\|_{C^{k-2,\alpha}(0)}$ and $\|u\|_{L^{\infty}(B_1)}$.
\end{theorem}

\begin{remark}\label{r-2.10}
For this \emph{higher pointwise regularity}, our results are new even for fully nonlinear equations without lower order terms. For interior local regularity, once we have $C^{2,\alpha}$ regularity, the higher $C^{k,\alpha}$ ($k\geq 3$) regularity can be derived easily by taking derivatives on the equation. However, for pointwise regularity, this technique is not applicable.

 Note that we do not have explicit estimates since the coefficients of the linearized equation depend on the solution itself (see \Cref{In-l-62}). But for linear equation \cref{e.linear}, we have the following explicit estimates:
\begin{equation*}
  |u(x)-P(x)|\leq C |x|^{k+\alpha}\left(\|u\|_{L^{\infty }(B_1)}
  +\|f\|_{C^{k-2,\alpha}(0)}\right), ~~\forall ~x\in B_{1},
\end{equation*}
\begin{equation*}
|a^{ij}P_{ij}+b^iP_i+cP-P_f|\leq C|x|^{k-2+\alpha}\left(\|u\|_{L^{\infty }(B_1)}
+\|f\|_{C^{k-2,\alpha}(0)}\right), ~~\forall ~x\in B_{1}
\end{equation*}
and
\begin{equation*}
|Du(0)|+\cdots+|D^ku(0)|\leq C\left(\|u\|_{L^{\infty }(B_1)}
+\|f\|_{C^{k-2,\alpha}(0)}\right),
\end{equation*}
where $C$ depends only on $k,n,\lambda,\Lambda,\alpha,\|a^{ij}\|_{C^{k-2,\alpha}(0)},\|b^i\|_{C^{k-2,\alpha}(0)}$ and $\|c\|_{C^{k-2,\alpha}(0)}$.
\end{remark}

\begin{remark}\label{r-2.8}
The observation \cref{e.Cka-3-i} is necessary for higher regularity and we refer to the proof of \cref{In-e.cka-6} for details.
\end{remark}
~\\

Before proving the interior pointwise $C^{k,\alpha}$ regularity, we first prove the interior local $C^{k,\alpha}$ regularity (similar to \Crefrange{l-3modin1}{l-32}). Recall that $0<\bar{\alpha}<1$ originates from \Crefrange{l-3modin1}{l-32} and is fixed throughout this paper.
\begin{lemma}\label{In-l-62}
Let $G$ be a fully nonlinear operator satisfying (i)-(iii) in \Cref{d-FP} (with $G \in C^{k-2,\bar{\alpha}}$ instead of $G\in C^{k,\bar{\alpha}}$) and $\omega_0$ satisfy \cref{e.omega0-2}. Suppose that $u$ is a viscosity solution of
\begin{equation*}
G(D^2u,Du,u,x)=0 \quad\mbox{in}~~B_1.
\end{equation*}
Then $u\in C^{k,\alpha}(\bar{B}_{1/2})$ for any $0<\alpha<\bar{\alpha}$ and
\begin{equation*}
\|u\|_{C^{k,\alpha}(\bar{B}_{1/2})}\leq C_k,
\end{equation*}
where $C_k$ depends only on $k,n,\lambda, \Lambda,\alpha,\mu,b_0,c_0,\omega_0,K_1,\omega_4$ and $\|u\|_{L^{\infty }(B_1)}$.

In particular, $u\in C^{k,\alpha}(0)$ and there exists $P\in \mathcal{P}_k$ such that
\begin{equation}\label{In-e.l62-1}
  |u(x)-P(x)|\leq C_k |x|^{k+\alpha}, ~~\forall ~x\in B_{1},
\end{equation}
\begin{equation}\label{In-e.l62-2}
  |G(D^2P(x),DP(x),P(x),x)|\leq C_k|x|^{k-2+\bar{\alpha}}, ~~\forall ~x\in B_{1}
\end{equation}
and
\begin{equation}\label{In-e.l62-3}
  \|P\|\leq C_k.
\end{equation}
\end{lemma}

\proof Since $G$ is smooth, for any $x_0,x\in B_1$,
$(M,p,s)\in \mathcal{S}^n\times \mathbb{R}^n\times \mathbb{R}$,
\begin{equation}\label{e5.2}
  \begin{aligned}
G(M&,p,s,x)-G(M,p,s,x_0)\\
= &G(M,p,s,x)-G(0,0,0,x)+G(0,0,0,x)-G(0,0,0,x_0)\\
    &+G(0,0,0,x_0)-G(M,p,s,x_0)\\
=&\int_{0}^{1}\big(G_{M_{ij}}(\xi)M_{ij}+G_{p_i}(\xi)p_i+G_{s}(\xi)s\big) d\tau+G(0,0,0,x)-G(0,0,0,x_0)\\
&-\int_{0}^{1}\big(G_{M_{ij}}(\zeta)M_{ij}+G_{p_i}(\zeta)p_i+G_{s}(\zeta)s\big)d\tau\\
=&\int_{0}^{1}\Big((G_{M_{ij}}(\xi)-G_{M_{ij}}(\zeta))M_{ij}+(G_{p_i}(\xi)-G_{p_i}(\zeta))p_i
+(G_{s}(\xi)-G_{s}(\zeta))s\Big)d\tau\\
&+G(0,0,0,x)-G(0,0,0,x_0),
  \end{aligned}
\end{equation}
where
\begin{equation*}
\xi=(\tau M,\tau p,\tau s,x),\quad\zeta=(\tau M,\tau p,\tau s,x_0).
\end{equation*}
By combining with $G \in C^{k-2,\bar{\alpha}}$ and (iii) in \Cref{d-FP}, we have
\begin{equation*}
  \begin{aligned}
|G&(M,p,s,x)-G(M,p,s,x_0)|\\
\leq&\int |D_MG(\xi)-D_MG(\zeta)||M|+|D_pG(\xi)-D_pG(\zeta)||p|
+|G_s(\xi)-G_s(\zeta)||s|d\tau\\
&+|G(0,0,0,x)-G(0,0,0,x_0)|\\
\leq& C\left(|M|+|p|+|s|+1\right)|x-x_0|.
  \end{aligned}
\end{equation*}
Hence, from \Cref{t-C2a-i}, $u\in C^{2,\alpha}(x_0)$ for any $0<\alpha<\bar{\alpha}$ and $x_0\in \bar{B}_{3/4}$. Hence, $u\in C^{2,\alpha}(\bar B_{3/4})$ and
\begin{equation*}
\|u\|_{C^{2,\alpha}(\bar{B}_{3/4})}\leq C,
\end{equation*}
where $C$ depends only on $k,n,\lambda, \Lambda,\alpha,\mu,b_0,c_0,\omega_0,K_1,\omega_4$ and $\|u\|_{L^{\infty }(B_1)}$. Unless stated otherwise, $C$ always has the same dependence in the following proof.

Now, we show that $u\in C^{3,\alpha}$, which can be proved by the standard technique of difference quotient. Let $h>0$ be small and $1\leq l\leq n$. Take the difference quotient on both sides of the equation and we have
\begin{equation*}
  a^{ij}(\Delta_l^h u)_{ij}=f \quad\mbox{in}~B_{5/8},
\end{equation*}
where
\begin{equation*}
  \begin{aligned}
&\Delta_l^h u(x)=\frac{1}{h}\left(u(x+he_l)-u(x)\right),\quad a^{ij}(x)=\int_{0}^{1}G_{M_{ij}}(\xi)d\tau,\\
&f(x)=-\int_{0}^{1}\Big(G_{p_i}(\xi)\Delta_l^h u_i+G_{s}(\xi)\Delta_l^h u+G_{x_l}(\xi)\Big)d\tau\\
  \end{aligned}
\end{equation*}
and
\begin{equation*}
  \begin{aligned}
\xi=\tau\left(D^2u(y),Du(y),u(y),y\right)+(1-\tau)\left(D^2u(x),Du(x),u(x),x\right),\quad y=x+he_l.
  \end{aligned}
\end{equation*}
Note that $a^{ij}$ is uniformly elliptic with ellipticity constants depending only on $n,\lambda$ and $\Lambda$. Indeed, the structure condition \cref{SC2} implies that $G_{M_{ij}}$ is uniformly elliptic (see the argument around \cite[(2.6)]{MR1351007}). Moreover, it is easy to check that
\begin{equation*}
\|a^{ij}\|_{C^{\alpha^2}(\bar{B}_{5/8})}\leq C, \quad \|f\|_{C^{\alpha^2}(\bar{B}_{5/8})}\leq C.
\end{equation*}

By the Schauder's estimate for linear equations, we have $\Delta_l^h u \in C^{2,\alpha^2}(\bar{B}_{1/2})$. Hence, $u_l \in C^{2,\alpha^2}(\bar{B}_{1/2})$ and
\begin{equation*}
\|u_l\|_{C^{2,\alpha^2}(\bar{B}_{1/2})}\leq C,
\end{equation*}
which implies $u\in C^{3,\alpha^2}$ and
\begin{equation*}
\|u\|_{C^{3,\alpha^2}(\bar{B}_{1/2})}\leq C.
\end{equation*}

Since $u\in C^{3,\alpha^2}$, for $1\leq l \leq n$, $u_l$ satisfies
\begin{equation}\label{In-uksat}
  a^{ij}(u_l)_{ij}=-f,
\end{equation}
where
\begin{equation}\label{In-aijuk}
  a^{ij}(x)=G_{M_{ij}}(D^2u,Du,u,x)
\end{equation}
and
\begin{equation}\label{In-Guk}
  f(x)=G_{p_i}(D^2u,Du,u,x)u_{il}+G_{s}(D^2u,Du,u,x)u_l+G_{x_l}(D^2u,Du,u,x).
\end{equation}
Then it can be seen that $a^{ij}\in C^{\alpha}(\bar{B}_{1/2})$ and $f\in C^{\alpha}(\bar{B}_{1/2})$. By the Schauder's estimate for linear equations again, $u_l\in C^{2,\alpha}(\bar{B}_{1/4})$. Then $u\in C^{3,\alpha}(\bar{B}_{1/4})$ and
\begin{equation*}
\|u\|_{C^{3,\alpha}(\bar{B}_{1/4})}\leq C.
\end{equation*}

Since $G\in C^{k-2,\alpha}$, by considering \crefrange{In-uksat}{In-Guk} iteratively, we obtain
\begin{equation*}
\|u\|_{C^{k,\alpha}(\bar{B}_{2^{-k}})}\leq C.
\end{equation*}
By a standard covering argument, $u\in C^{k,\alpha}(\bar{B}_{1/2})$ and
\begin{equation*}
\|u\|_{C^{k,\alpha}(\bar{B}_{1/2})}\leq C.
\end{equation*}

In particular, $u\in C^{k,\alpha}(0)$ and there exists $P\in \mathcal{P}_k$ such that \cref{In-e.l62-1} and \cref{In-e.l62-3} hold. From \cref{In-e.l62-1}, for any $0\leq l\leq k$,
\begin{equation*}\label{In-e.Q-mu}
  |D^lu(x)-D^lP(x)|\leq C |x|^{k-l+\alpha}, ~~\forall ~x\in B_{1}.
\end{equation*}
Thus,
\begin{equation*}
  \begin{aligned}
    |G&(D^2P(x),DP(x),P(x),x)|\\
    &=|G(D^2P(x),DP(x),P(x),x)-G(D^2u(x),Du(x),u(x),x)|\\
    &=\bigg|\int_{0}^{1} G_{M_{ij}}(\xi)(P_{ij}-u_{ij})+G_{p_i}(\xi)(P_{i}-u_{i})+G_{s}(\xi)(P-u)d\tau \bigg|\\
    &\leq C|x|^{k-2+\alpha},
  \end{aligned}
\end{equation*}
where $\xi=\tau\left(D^2P(x),DP(x),P(x),x\right)+(1-\tau)\left(D^2u(x),Du(x),u(x),x\right)$. Then
\begin{equation*}
  D^{l}_x \Big(G(D^2P(x),DP(x),P(x),x)\Big)\Big|_{x=0}=0,~\forall ~0\leq l\leq k-2.
\end{equation*}
Since $G\in C^{k-2,\bar{\alpha}}$, \cref{In-e.l62-2} holds. \qed~\\

\begin{remark}\label{r-6.1}
In above proof, we do not assume that $G(0,0,0,x_0)=0$. Thus, \Cref{t-C2a-i} cannot be applied directly. However, this assumption is not essential and we can consider $G(M,p,s,x)-G(0,0,0,x_0)$.
\end{remark}

\begin{remark}\label{r-6.2}
If we consider the linear equation \cref{e.linear}, the coefficients are independent of $u$. Then we obtain
\begin{equation*}
\|u\|_{C^{k,\alpha}(\bar{B}_{1/2})}\leq C_k\|u\|_{L^{\infty }(B_1)},
\end{equation*}
where $C_k$ depends only on $k,n,\lambda,\Lambda,\alpha,
\|a^{ij}\|_{C^{k-2,\alpha}(\bar{B}_1)},\|b^{i}\|_{C^{k-2,\alpha}(\bar{B}_1)}$ and $\|c\|_{C^{k-2,\alpha}(\bar{B}_1)}$.

Based on this estimate, we can obtain the higher order pointwise regularity with explicit estimates as pointed out in \Cref{r-2.10}.
\end{remark}
~\\

In the following, we prove the interior pointwise $C^{k,\alpha}(k\geq 3)$ regularity \Cref{t-Cka-i} by induction. For $k=3$, $F\in C^{1,\alpha}(0)$ and $F(M,p,s,0)\equiv G(M,p,s,0)$ for any
$(M,p,s)\in \mathcal{S}^n\times \mathbb{R}^n\times \mathbb{R}$ (by \cref{holder}). Hence, by an argument similar to \cref{e5.2},
\begin{equation}\label{e5.3}
  \begin{aligned}
|F&(M,p,s,x)-G(M,p,s,0)|\\
&\leq |F(M,p,s,x)-G(M,p,s,x)|+|G(M,p,s,x)-G(M,p,s,0)|\\
&\leq \|F\|_{C^{1,\alpha}(0)}(|M|+1)\omega_3(|p|,|s|)|x|^{1+\alpha}
    +C(|M|+|p|+|s|+1)|x|\\
&\leq C(|M|+1)\left(\omega_3(|p|,|s|)+|p|+|s|+1\right)|x|.
  \end{aligned}
\end{equation}
Then from the interior pointwise $C^{2,\alpha}$ regularity \Cref{t-C2a-i}, $u\in C^{2,\alpha}(0)$. Hence, we can assume (and do assume throughout this section) that the $C^{k-1,\alpha}(0)$ regularity holds if $F\in C^{k-2,\alpha}(0)$ and we need to derive the $C^{k,\alpha}(0)$ regularity.

The following lemma is a higher order counterpart of \Cref{In-l-C1a-mu} and \Cref{In-l-C2a-mu}.
\begin{lemma}\label{In-l-Cka-mu}
Let $0<\alpha<\bar{\alpha}$, $F\in C^{k-2,\alpha}(0)$ and $\omega_0$ satisfy \cref{e.omega0-2}. Then there exists $\delta>0$ depending only on $k,n,\lambda,\Lambda,\alpha,\omega_0,K_1,\omega_3$ and $\omega_4$ such that if $u$ satisfies
\begin{equation*}
F(D^2u,Du,u,x)=f \quad\mbox{in}~~B_1
\end{equation*}
with
\begin{equation*}
  \begin{aligned}
&u(0)=|Du(0)|=\cdots=|D^{k-1}u(0)|=0, \quad \max(\|u\|_{L^{\infty}(B_1)},\mu,b_0,c_0)\leq 1,\\
&\max\left(\|F\|_{C^{k-2,\alpha}(0)},\|f\|_{C^{k-3,\alpha}(0)}\right)\leq \delta,
  \end{aligned}
\end{equation*}
then there exists $P\in\mathcal{HP}_k$ such that
\begin{equation*}
  \|u-P\|_{L^{\infty}(B_{\eta})}\leq \eta^{k+\alpha},
\end{equation*}
\begin{equation*}
  |G(D^2P(x),DP(x),P(x),x)|\leq C |x|^{k-2+\bar{\alpha}}, ~\forall ~x\in B_1,
\end{equation*}
and
\begin{equation*}
\|P\|\leq C,
\end{equation*}
where $C$ and $0<\eta<1$ depend only on $k,n,\lambda, \Lambda,\alpha,\omega_0,K_1,\omega_3$ and $\omega_4$.
\end{lemma}
\begin{remark}\label{r-6-1}
Since we have assumed that the interior pointwise $C^{k-1,\alpha}$ regularity holds by induction, $Du(0),\cdots,D^{k-1}u(0)$ are well defined.
\end{remark}
~\\
\noindent\textbf{Proof of \Cref{In-l-Cka-mu}.} Throughout this proof, $C$ always denotes a constant depending only on $k,n,\lambda, \Lambda,\alpha,\omega_0,K_1,\omega_3$ and $\omega_4$. As before, we prove the lemma by contradiction. Suppose that the lemma is false. Then there exist $0<\alpha<\bar{\alpha},\omega_0,K_1,\omega_3,\omega_4$ and a sequence of $(F_m,u_m,f_m)_{m\in \mathbb{N}}$ satisfying
\begin{equation*}
F_m(D^2u_m,Du_m,u_m,x)=f_m \quad\mbox{in}~~B_1.
\end{equation*}
In addition, $F_m\in C^{k-2,\alpha}(0)$ (with $G_{m},K_1,\omega_3$) satisfy the structure condition \cref{SC2} (with $\lambda,\Lambda,\mu\equiv 1,b\equiv 1,c\equiv 1,\omega_0$) and
\begin{equation}\label{e6.1}
  \begin{aligned}
&\|u_m\|_{L^{\infty}(B_1)}\leq 1,\quad u_m(0)=|Du_m(0)|=\cdots=|D^{k-1}u_m(0)|=0,\\
&\|F_{m}\|_{C^{k-2,\alpha}(0)}\leq \frac{1}{m},\quad \|f_m\|_{C^{k-3,\alpha}(0)}\leq\frac{1}{m},\\
&\|G_{m}\|_{C^{k-2,\bar{\alpha}}(\bar{\mathbf{B}}_r\times \bar{B}_1)}\leq \omega_4(r),~\forall ~r>0.
  \end{aligned}
\end{equation}
Finally, for any $P\in\mathcal{HP}_k$ satisfying $\|P\|\leq C_k$ and
\begin{equation}\label{In-e.cka.Q-mu}
|G_{m}(D^2P(x),DP(x),P(x),x)|\leq C_k |x|^{k-2+\bar{\alpha}}, ~\forall ~x\in B_1,
\end{equation}
we have
\begin{equation}\label{In-e.lCka.1-mu}
  \|u_m-P\|_{L^{\infty}(B_{\eta})}> \eta^{k+\alpha},
\end{equation}
where $C_k$ is to be specified later and $0<\eta<1$ is taken small such that
\begin{equation}\label{In-e.lCka.2-mu}
C_k\eta^{(\bar{\alpha}-\alpha)/2}<\frac{1}{2}.
\end{equation}

As before, $u_m$ are uniformly bounded and equicontinuous. Indeed, by \cref{holder},
\begin{equation*}
|F_m(0,0,0,x)|\leq |G_m(0,0,0,x)|
+\|F_{m}\|_{C^{k-2,\alpha}(0)}|x|^{k-2+\alpha}\omega_3(0,0),~\forall ~x\in B_1.
\end{equation*}
In addition, from \cref{e6.1},
\begin{equation*}
|G_m(0,0,0,x)|\leq C,~\forall ~m\geq 1,~x\in B_1.
\end{equation*}
Then \Cref{l-3Ho} implies that $u_m$ are equicontinuous (cf. the proof of \Cref{In-l-C1a-mu}). Hence, there exists $\bar u$ such that
\begin{equation*}
u_m\rightarrow \bar u \quad\mbox{ in } L^{\infty}_{\mathrm{loc}}(B_1).
\end{equation*}
Since $\|G_{m}\|_{C^{k-2,\bar{\alpha}}(\bar U)}$ are uniformly bounded for any
$U\subset\subset \mathcal{S}^n\times \mathbb{R}^n\times \mathbb{R}\times \bar{B}_1$, there exists $\bar G\in C^{k-2,\bar{\alpha}}(\mathcal{S}^n\times \mathbb{R}^n\times \mathbb{R}\times \bar{B}_1)$ such that
\begin{equation*}
G_{m}\rightarrow \bar{G}\quad\mbox{ in}~C^{k-2,(\bar{\alpha}+\alpha)/2}(\bar U),~\forall ~U\subset\subset
\mathcal{S}^n\times \mathbb{R}^n\times \mathbb{R}\times \bar{B}_1
\end{equation*}
and
\begin{equation*}
\|\bar{G}\|_{C^{k-2,\bar{\alpha}}(\bar{\mathbf{B}}_r\times \bar{B}_1)}\leq \omega_4(r),~\forall ~r>0.
\end{equation*}
Moreover, it is easy to verify that $\bar{G}$ satisfies (i)-(iii) in \Cref{d-FP} (with $K_1$) and the structure condition \cref{SC2} (with $\lambda,\Lambda,\mu\equiv 1,b\equiv  1,c\equiv 1,\omega_0$).

Next, for any ball $B\subset\subset B_1$ and $\varphi\in C^2(\bar{B})$, let $r=\|D^2\varphi\|_{L^{\infty}(B)}+\|D\varphi\|_{L^{\infty}(B)}+1$, $\psi_m=F_m(D^2\varphi,D\varphi,u_m,x)-f_m$ and $\psi(x)=\bar{G}(D^2\varphi,D\varphi,\bar u,x)$. By noting
\begin{equation*}
  \begin{aligned}
|F_m(&D^2\varphi,D\varphi,u_m,x)-\bar{G}(D^2\varphi,D\varphi,\bar u,x)|\\
\leq&|F_m(D^2\varphi,D\varphi,u_m,x)-G_{m}(D^2\varphi,D\varphi,u_m,x)|\\
    &+|G_{m}(D^2\varphi,D\varphi,u_m,x)-\bar{G}(D^2\varphi,D\varphi,u_m,x)|\\
    &+|\bar{G}(D^2\varphi,D\varphi,u_m,x)-\bar{G}(D^2\varphi,D\varphi,\bar{u},x)|\\
\leq& \|F_{m}\|_{C^{k-2,\alpha}(0)}(|D^2\varphi|+1)\omega_3(r,r)\\
    &+|G_{m}(D^2\varphi,D\varphi,u_m,x)-\bar{G}(D^2\varphi,D\varphi,u_m,x)|\\
    &+|\bar{G}(D^2\varphi,D\varphi,u_m,x)-\bar{G}(D^2\varphi,D\varphi,\bar{u},x)|,
\end{aligned}
\end{equation*}
we have $\|\psi_m-\psi\|_{L^n(B)}\rightarrow 0$ as before. From \Cref{l-35}, $\bar{u}$ is a viscosity solution of
\begin{equation*}
\bar{G}(D^2\bar{u},D\bar{u},\bar{u},x)=0 \quad\mbox{in}~~B_{1}.
\end{equation*}

By the $C^{k-1,\alpha}$ estimate for $u_m$ (by induction) and noting
\begin{equation*}
u_m(0)=\cdots =|D^{k-1}u_m(0)|=0,
\end{equation*}
we have
\begin{equation*}
|u_m(x)|\leq C|x|^{k-1+\alpha} ,~\forall~x\in B_1.
\end{equation*}
Since $u_m$ converges to $\bar{u}$,
\begin{equation*}
|\bar{u}(x)|\leq C|x|^{k-1+\alpha} ,~\forall~x\in B_1,
\end{equation*}
which implies
\begin{equation*}
\bar{u}(0)=\cdots=|D^{k-1}\bar{u}(0)|=0.
\end{equation*}
By combining with \Cref{In-l-62}, there exists $\bar{P}\in\mathcal{HP}_k$ such that
\begin{equation}\label{In-e.cka-5}
  |\bar{u}(x)-\bar{P}(x)|\leq C_1 |x|^{k+(\bar{\alpha}+\alpha)/2}, ~~\forall ~x\in B_{1},
\end{equation}
\begin{equation}\label{In-e.cka-4}
  |\bar{G}(D^2\bar{P}(x),D\bar{P}(x),\bar{P}(x),x)|\leq C_1|x|^{k-2+\bar{\alpha}}, ~~\forall ~x\in B_{1}
\end{equation}
and
\begin{equation*}
  \|\bar{P}\|\leq C_1,
\end{equation*}
where $C_1$ depends only on $k,n,\lambda, \Lambda,\alpha,\omega_0,K_1,\omega_3$ and $\omega_4$.

Now, we construct a sequence of $\bar{P}_m\in\mathcal{HP}_k$ such that \cref{In-e.cka.Q-mu} holds for $\bar{P}_m$ and $\bar{P}_m\rightarrow \bar{P}$ as $m\rightarrow \infty$. Let $P_m\in\mathcal{HP}_k$ ($m\geq 1$) with $\|P_m\|\leq 1$ to be chosen later and
\begin{equation*}
\bar{P}_m=P_m+\bar{P}.
\end{equation*}
Denote
\begin{equation*}
h_m(x)=G_{m}(D^2\bar{P}_m(x),D\bar{P}_m(x),\bar{P}_m(x),x).
\end{equation*}
Since $h_m\in C^{k-2,\bar{\alpha}}(\bar{B}_{1})$, \cref{In-e.cka.Q-mu} holds for $\bar{P}_m$ if we can show
\begin{equation*}
D^ih_m(0)=0,~\forall ~1\leq i\leq k-2.
\end{equation*}

Note that $u_m(0)=\cdots=|D^{k-1}u_m(0)|=0$ and from the interior pointwise $C^{k-1,\alpha}$ regularity \Cref{t-Cka-i} (i.e., $P\equiv P_f\equiv 0$ in \cref{e.Cka-3-i}) and \cref{holder},
\begin{equation*}
|G_{m}(0,0,0,x)|\leq C|x|^{k-3+\bar{\alpha}}.
\end{equation*}
(We remark here that if $k=3$, we should use the interior pointwise $C^{2,\alpha}$ regularity \Cref{t-C2a-i} with $P\equiv f(0)\equiv 0$ in \cref{e.C2a-3-i-mu}.) Thus,
\begin{equation*}
  \begin{aligned}
|h&_m(x)|\\
&\leq |G_{m}(D^2\bar{P}_m(x),D\bar{P}_m(x),\bar{P}_m(x),x)-G_{m}(0,0,0,x)|
    +|G_{m}(0,0,0,x)|\\
&=|\int_{0}^{1} G_{m,M_{ij}}(\xi)\bar{P}_{m,ij}+G_{m,p_i}(\xi)\bar{P}_{m,i}
+G_{m,s}(\xi)\bar{P}_md\tau|+|G_{m}(0,0,0,x)|\\
&\leq C|x|^{k-3+\bar{\alpha}},
  \end{aligned}
\end{equation*}
where $\xi=\tau\left(D^2\bar{P}_m(x),D\bar{P}_m(x),\bar{P}_m(x),x\right)
+(1-\tau)\left(0,0,0,x\right)$.

Hence, to verify \cref{In-e.cka.Q-mu} for $\bar{P}_m$, we only need to prove
\begin{equation*}
D^{k-2}h_m(0)=D^{k-2}\left( G_{m}(D^2\bar{P}_m(x),D\bar{P}_m(x),\bar{P}_m(x),x)\right)\bigg|_{x=0}=0.
\end{equation*}
Indeed, since $\bar{P}_m\in \mathcal{HP}_k$,
\begin{equation*}
  \begin{aligned}
D^{k-2}h_m(0)&=G_{m,M_{ij}}(0)D^{k-2}\bar{P}_{m,ij}+(D^{k-2}_{x}G_{m})(0)\\
&=G_{m,M_{ij}}(0)D^{k-2}\bar{P}_{ij}+(D^{k-2}_{x}G_{m})(0)+G_{m,M_{ij}}(0)D^{k-2}P_{m,ij},
  \end{aligned}
\end{equation*}
where $G_{m,M_{ij}}(0)$ is short for $G_{m,M_{ij}}(0,0,0,0)$ (similarly in the following proof). By the structure condition \cref{SC2},
\begin{equation*}
\lambda I\leq D_MG_m(0)\leq \Lambda I.
\end{equation*}
Hence, by \Cref{le2.2}, we can choose $P_m\in \mathcal{HP}_k$ such that
\begin{equation}\label{e6.8}
G_{m,M_{ij}}(0)P_{m,ij}=\sum_{|\sigma|=k-2} \frac{a^{m}_{\sigma}}{\sigma !}x^{\sigma}
\end{equation}
and
\begin{equation*}
\|P_m\|\leq C\sum_{|\sigma|=k-2}|a^m_{\sigma}|,
\end{equation*}
where
\begin{equation*}
a^m_{\sigma}\coloneqq G_{m,M_{ij}}(0)D^{\sigma}\bar{P}_{ij}+(D_x^{\sigma}G_m)(0),~\forall ~|\sigma|=k-2.
\end{equation*}

With this choice of $P_m$,
\begin{equation*}
D^{k-2}h_m(0)=0,~\forall ~m\geq 1.
\end{equation*}
Moreover, note that
\begin{equation*}
  \begin{aligned}
a^m_{\sigma}&=G_{m,M_{ij}}(0)D^{\sigma}\bar{P}_{ij}+(D^{\sigma}_{x}G_{m})(0)\longrightarrow \bar G_{M_{ij}}(0)D^{\sigma}\bar{P}_{ij}+(D^{\sigma}_{x}\bar G)(0)\\
&=D^{\sigma}\left(\bar G(D^2\bar{P}(x),D\bar{P}(x),\bar{P}(x),x)\right)\bigg |_{x=0}=0~~(\mbox{by}~ \cref{In-e.cka-4}).
  \end{aligned}
\end{equation*}
Hence, $\|P_m\|\to 0$.

Therefore, there exists a constant $C_k(\geq C_1+1)$ depending only on $k,n,\lambda, \Lambda,\alpha,\omega_0$, $K_1,\omega_3$ and $\omega_4$ such that $\|\bar{P}_m\|\leq C_k$ for any $m\geq 1$ and
\begin{equation*}
  \begin{aligned}
|h_m(x)|=|G_{m}(D^2\bar{P}_m(x),D\bar{P}_m(x),\bar{P}_m(x),x)|\leq C_k|x|^{k-2+\bar{\alpha}}.
  \end{aligned}
\end{equation*}
Thus, \cref{In-e.lCka.1-mu} holds for $\bar{P}_m$. Let $m\rightarrow \infty$ and we have
\begin{equation*}
    \|\bar{u}-\bar{P}\|_{L^{\infty}(B_{\eta})}\geq \eta^{k+\alpha}.
\end{equation*}
However, by \cref{In-e.lCka.2-mu} and \cref{In-e.cka-5}, we have
\begin{equation*}
  \|\bar{u}-\bar{P}\|_{L^{\infty}(B_{\eta})}\leq \frac{1}{2}\eta^{k+\alpha},
\end{equation*}
which is a contradiction.  ~\qed~\\

Now, similar to the $C^{2,\alpha}$ regularity, we prove the interior pointwise $C^{k,\alpha}$ regularity in a special case.
\begin{lemma}\label{In-t-Ckas-mu}
Let $0<\alpha <\bar{\alpha}$, $F\in C^{k-2,\alpha}(0)$ and $\omega_0$ satisfy \cref{e.omega0-2}. Let $u$ be a viscosity solution of
\begin{equation*}
F(D^2u,Du,u,x)=f \quad\mbox{in}~~B_1.
\end{equation*}
Assume that
\begin{equation}\label{In-e.tCkas-be-mu}
\begin{aligned}
&\|u\|_{L^{\infty}(B_1)}\leq 1,\quad u(0)=|Du(0)|=\cdots=|D^{k-1}u(0)|=0,\\
&\mu\leq \frac{1}{4C_0},\quad b_0\leq \frac{1}{2}, \quad c_0\leq \frac{1}{K_0},\\
&\|F\|_{C^{k-2,\alpha}(0)}\leq \frac{\delta_1}{C_0},\quad|f(x)|\leq \delta_1|x|^{k-2+\alpha}, ~~\forall ~x\in B_1,\\
\end{aligned}
\end{equation}
where $\delta_1\leq \delta$ ($\delta$ is as in \Cref{In-l-Cka-mu}) and $C_0$ depend only on $k,n,\lambda, \Lambda,\alpha,\omega_0,K_1,\omega_3$ and $\omega_4$.

Then $u\in C^{k,\alpha}(0)$ and there exists $P\in\mathcal{HP}_k$ such that
\begin{equation*}\label{In-e.tCkas-1-mu}
  |u(x)-P(x)|\leq C |x|^{k+\alpha}, ~~\forall ~x\in B_{1},
\end{equation*}
\begin{equation*}
  |G(D^2P(x),DP(x),P(x),x)|\leq C |x|^{k-2+\bar{\alpha}}, ~~\forall ~x\in B_{1}
\end{equation*}
and
\begin{equation*}\label{In-e.tCkas-2-mu}
\|P\|\leq C,
\end{equation*}
where $C$ depends only on $k,n,\lambda, \Lambda,\alpha,\omega_0,K_1,\omega_3$ and $\omega_4$.
\end{lemma}

\proof As before, to prove that $u$ is $C^{k,\alpha}$ at $0$, we only need to prove the following. There exist a sequence of $P_m\in\mathcal{HP}_k$ ($m\geq 0$) such that for all $m\geq 1$,

\begin{equation}\label{In-e.tCkas-6-mu}
\|u-P_m\|_{L^{\infty }(B _{\eta^{m}})}\leq \eta ^{m(k+\alpha )},
\end{equation}
\begin{equation}\label{In-e.tCkas-9-mu}
|G(D^2P_m(x),DP_m(x),P_m(x),x)|\leq \tilde{C}|x|^{k-2+\bar{\alpha}}, ~~\forall ~x\in B_{1}
\end{equation}
and
\begin{equation}\label{In-e.tCkas-7-mu}
\|P_m-P_{m-1}\|\leq \tilde{C}\eta ^{(m-1)\alpha},
\end{equation}
where $\tilde{C}$ and $\eta$ depends only on $k,n,\lambda, \Lambda,\alpha,\omega_0,K_1,\omega_3$ and $\omega_4$.

We prove the above by induction. For $m=1$, by \Cref{In-l-Cka-mu}, there exists $P_1\in\mathcal{HP}_k$ such that \crefrange{In-e.tCkas-6-mu}{In-e.tCkas-7-mu} hold for some $C_1$ and $\eta_1$ depending only on $k,n,\lambda, \Lambda,\alpha,\omega_0,K_1,\omega_3$ and $\omega_4$ where $P_0\equiv 0$. Take $\tilde{C}\geq C_1,\eta\leq \eta_1$ (to be specified later) and then the conclusion holds for $m=1$. Suppose that the conclusion holds for $m\leq m_0$. We need to prove that the conclusion holds for $m=m_0+1$.

Let $r=\eta ^{m_{0}}$, $y=x/r$ and
\begin{equation}\label{In-e.tCkas-v-mu}
  v(y)=\frac{u(x)-P_{m_0}(x)}{r^{k+\alpha}}.
\end{equation}
Then $v$ satisfies
\begin{equation}\label{In-e.Ckas-F-mu}
\tilde{F}(D^2v,Dv,v,y)=\tilde{f} \quad\mbox{in}~~ B_1,
\end{equation}
where for $(M,p,s,y)\in \mathcal{S}^n\times \mathbb{R}^n\times \mathbb{R}\times \bar B_1$,
\begin{equation*}
  \begin{aligned}
&\tilde{F}(M,p,s,y)=r^{-l}F(r^{l}M+D^2P_{m_0}(x),r^{l+1}p+DP_{m_0}(x),r^{l+2}s+P_{m_0}(x),x),\\
&\tilde{f}(y)=r^{-l}f(x), \quad l=k-2+\alpha.
  \end{aligned}
\end{equation*}
In addition, define $\tilde{G}$ similarly, i.e.,
\begin{equation*}
  \begin{aligned}
&\tilde{G}(M,p,s,y)=r^{-l}G(r^{l}M+D^2P_{m_0}(x),r^{l+1}p+DP_{m_0}(x),r^{l+2}s+P_{m_0}(x),x).
  \end{aligned}
\end{equation*}

In the following, we show that \cref{In-e.Ckas-F-mu} satisfies the assumptions of \Cref{In-l-Cka-mu}. First, it is easy to verify that
\begin{equation*}
\begin{aligned}
\|v\|_{L^{\infty}(B_1)}\leq& 1, \quad v(0)=\cdots=|D^{k-1}v(0)|=0,~\quad(\mathrm{by}~ \cref{In-e.tCkas-be-mu},~\cref{In-e.tCkas-6-mu} ~\mbox{and}~ \cref{In-e.tCkas-v-mu})\\
|\tilde{f}(y)|=& r^{-(k-2+\alpha)}|f(x)|\leq \delta_1|y|^{k-2+\alpha}, ~\forall ~y\in B_1. ~\quad(\mathrm{by}~\cref{In-e.tCkas-be-mu})\\
  \end{aligned}
\end{equation*}

By \cref{In-e.tCkas-7-mu},
\begin{equation*}
\|P_m\|\leq C_{0},~\forall ~0\leq m\leq m_0,
\end{equation*}
where $C_{0}$ depends only on $k,n,\lambda, \Lambda,\alpha,\omega_0,K_1,\omega_3$ and $\omega_4$. It is easy to check that $\tilde{F}$ and $\tilde{G}$ satisfy the structure condition \cref{SC2} with $\lambda,\Lambda,\tilde{\mu},\tilde{b},\tilde{c}$ and $\tilde{\omega}_0$, where (note that $\omega_0$ satisfies \cref{e.omega0-2})
\begin{equation*}
\tilde{\mu}=r^{k+\alpha}\mu,\quad\tilde{b}= rb_0+2 C_0 r\mu ,\quad\tilde{c}= K_0 r^{2} c_0, \quad\tilde{\omega}_0(\cdot,\cdot)=\omega_0(\cdot+C_0,\cdot).
\end{equation*}
Hence, $\tilde{\omega}_0$ satisfies \cref{e.omega0-2} and from \cref{In-e.tCkas-be-mu},
\begin{equation*}
\begin{aligned}
&\tilde{\mu}\leq \mu\leq 1,\quad\tilde{b}\leq b_0+2C_0\mu\leq 1,\quad\tilde{c}\leq c_0\leq 1.
\end{aligned}
\end{equation*}

In addition, for $(M,p,s,y)\in \mathcal{S}^n\times \mathbb{R}^n\times \mathbb{R}\times \bar B_1$,
\begin{equation}\label{e.6.1}
\begin{aligned}
|\tilde{F}(M&,p,s,y)-\tilde{G}(M,p,s,y)|\\
=&r^{-(k-2+\alpha)}\bigg(F(r^{k-2+\alpha}M+D^2P_{m_0}(x),r^{k-1+\alpha}p+DP_{m_0}(x),
    r^{k+\alpha}s+P_{m_0}(x),x)\\
&-G(r^{k-2+\alpha}M+D^2P_{m_0}(x),r^{k-1+\alpha}p+DP_{m_0}(x),
    r^{k+\alpha}s+P_{m_0}(x),x)\bigg)\\
\leq& r^{-(k-2+\alpha)}
\|F\|_{C^{k-2,\alpha}(0)}|x|^{k-2+\alpha}(|M|+C_0)\omega_3(|p|+C_0,|s|+C_0)\\
\leq& \delta_1(|M|+1)\omega_3(|p|+C_0,|s|+C_0)|y|^{k-2+\alpha}\\
\coloneqq &\delta_1|y|^{k-2+\alpha}(|M|+1)\tilde{\omega}_3(|p|,|s|).
\end{aligned}
\end{equation}
Hence,
\begin{equation*}
\|F\|_{C^{k-2,\alpha}(0)}\leq \delta_1.
\end{equation*}

Obviously, $\tilde{G}$ is convex in $M$. Moreover,
\begin{equation*}
  \begin{aligned}
|\tilde G_{M_{ij}}&(M,p,s,y_1)-\tilde G_{M_{ij}}(N,q,t,y_2)|\\
=&|G_{M_{ij}}(r^{k-2+\alpha}M+D^2P_{m_0}(x_1),r^{k-1+\alpha}p+DP_{m_0}(x_1),
    r^{k+\alpha}s+P_{m_0}(x_1),x_1)\\
&-G_{M_{ij}}(r^{k-2+\alpha}N+D^2P_{m_0}(x_2),r^{k-1+\alpha}q+DP_{m_0}(x_2),
    r^{k+\alpha}t+P_{m_0}(x_2),x_2)|\\
\leq& \tilde{K}_1\left(|M-N|+|p-q|+|s-t|+|y_1-y_2|\right),\\
&\forall (M,p,s,y_1),~(N,q,t,y_2)\in \mathcal{S}^n\times \mathbb{R}^n\times\mathbb{R}\times \bar{B}_1,
  \end{aligned}
\end{equation*}
where $\tilde{K}_1$ depends only on $k,n,\lambda, \Lambda,\alpha,\omega_0,K_1,\omega_3$ and $\omega_4$. The corresponding inequalities for $G_{p_i}$ and $G_{s}$ can be verified similarly.

Finally, we show that
\begin{equation}\label{In-e.cka-6}
\|\tilde{G}\|_{C^{k-2,\bar{\alpha}}(\bar{\mathbf{B}}_\rho\times \bar{B}_1)}\leq \tilde{\omega}_4(\rho),~\forall ~\rho>0,
\end{equation}
where $\tilde{\omega}_4$ depends only on $k,n,\lambda, \Lambda,\alpha,\omega_0,K_1,\omega_3$ and $\omega_4$. The difficulty is $r^{-(k-2+\alpha)}$ in the expression and we overcome it with the aid of \cref{In-e.tCkas-9-mu}. First, by the definition of $\tilde{G}$ and the Newton-Leibniz theorem,
\begin{equation*}
  \begin{aligned}
\tilde{G}(M,p,s,y)=&\tilde{G}(M,p,s,y)-\tilde{G}(0,0,0,y)+\tilde{G}(0,0,0,y)\\
    =&\int_{0}^{1} \Big(G_{M_{ij}}(\xi)M_{ij}+rG_{p_i}(\xi)p_i+r^2G_{s}(\xi)s\Big)d\tau\\
    &+r^{-(k-2+\alpha)}G(D^2P_{m_0}(x),DP_{m_0}(x),P_{m_0}(x),x)\\
    \coloneqq& \tilde G_1(M,p,s,y)+\tilde G_2(y),
  \end{aligned}
\end{equation*}
where
\begin{equation*}
\xi=\tau\left(r^{k-2+\alpha}M,r^{k-1+\alpha}p,r^{k+\alpha}s,x\right)
+\left(D^2P_{m_0}(x),DP_{m_0}(x),P_{m_0}(x),(1-\tau)x\right).
\end{equation*}
For $\tilde{G}_1$, there is no $r^{-(k-2+\alpha)}$ in the expression. Since $G\in C^{k-2,\bar{\alpha}}$, we have $\tilde{G}_1\in C^{k-3,\bar{\alpha}}$ and
\begin{equation}\label{e6.3}
\|\tilde{G}_1\|_{C^{k-3,\bar{\alpha}}(\bar{\mathbf{B}}_\rho\times \bar B_1)}\leq
C\|G\|_{C^{k-2,\bar{\alpha}}(\bar{\mathbf{B}}_{(\rho+C_0)}\times \bar{B}_1)}\leq C\omega_4(\rho+C_0),~\forall ~\rho>0.
\end{equation}
For $\tilde{G}_2$, by \cref{In-e.tCkas-9-mu},
\begin{equation*}
\|\tilde{G}_2\|_{L^{\infty}(B_1)}
=r^{-(k-2+\alpha)}\|G(D^2P_{m_0}(x),DP_{m_0}(x),P_{m_0}(x),x)\|_{L^{\infty}(B_{r})}
\leq \tilde{C}.
\end{equation*}
In addition, from $\|P_{m_0}\|\leq C_0$ and  $G\in C^{k-2,\bar{\alpha}}$ again,
\begin{equation*}
[\tilde{G}_2]_{C^{k-2,\bar{\alpha}}(\bar B_1)}
=r^{-\alpha}[D^{k-2}G(D^2P_{m_0}(x),DP_{m_0}(x),P_{m_0}(x),x)]_{C^{\bar{\alpha}}(\bar B_{r})}
\leq Cr^{\bar{\alpha}-\alpha}.
\end{equation*}
Thus, by interpolation in H\"{o}lder spaces,
\begin{equation}\label{e6.4}
\|\tilde{G}_2\|_{C^{k-2,\bar{\alpha}}(\bar B_1)}\leq C.
\end{equation}

By \cref{e6.3} and \cref{e6.4}, we only need to estimate the $(k-2)$-th order of $\tilde{G}_1$ in the following. From the definitions of $\tilde{G}_1$ and $\tilde{G}$, if any $(k-2)$-th derivative of $\tilde{G}_1$ involves one derivative with respect to $M,p$ or $s$, the trouble $r^{-(k-2+\alpha)}$ will be canceled and we have the desired estimate. If we take $(k-2)$-th derivatives with respect to $y$, by noting $D^{k-2}G\in C^{\bar{\alpha}}$,
\begin{equation}\label{e6.5}
  \begin{aligned}
|D^{k-2}_{y}&\tilde{G}_1(M,p,s,y)|\\
=&\left|D^{k-2}_{y}\left(\tilde{G}(M,p,s,y)-\tilde{G}(0,0,0,y)\right)\right|\\
=&r^{-\alpha}\left|D^{k-2}_{x}\Big(G(r^{k-2+\alpha}M+D^2P_{m_0}(x),
r^{k-1+\alpha}p+DP_{m_0}(x),r^{k+\alpha}s+P_{m_0}(x),x)\right.\\
&\left.-G(D^2P_{m_0}(x),DP_{m_0}(x),P_{m_0}(x),x)\Big)\right|\\
\leq& C\rho^{\bar{\alpha}}[G]_{C^{k-2,\bar{\alpha}}(\bar{\mathbf{B}}_{(\rho+C_0)}\times \bar B_r)}\\
\leq& C\rho^{\bar{\alpha}}\omega_4(\rho+C_0).
  \end{aligned}
\end{equation}
By a similar argument,
\begin{equation}\label{e6.6}
  \begin{aligned}
[D^{k-2}_{y}\tilde{G}_1]_{C^{\bar{\alpha}}(\bar{\mathbf{B}}_\rho\times \bar B_1)}
\leq C\omega_4(\rho+C_0).
  \end{aligned}
\end{equation}
Hence, by combining above arguments together, we arrive at \cref{In-e.cka-6} for some $\tilde{\omega}_4$ depending only on $k,n,\lambda, \Lambda,\alpha,\omega_0,K_1,\omega_3$ and $\omega_4$.

Choose $\delta_1$ small enough (depending only on $k,n,\lambda, \Lambda,\alpha,\omega_0,K_1,\omega_3$ and $\omega_4$) such that \Cref{In-l-Cka-mu} holds for $\tilde\omega_0,\tilde{K}_1,\tilde\omega_3,\tilde\omega_4$ and $\delta_1$. Since \cref{In-e.Ckas-F-mu} satisfies the assumptions of \Cref{In-l-Cka-mu}, there exist $\tilde{P}\in\mathcal{HP}_k$ and constants $\tilde{C}\geq C_1$ and $\eta\leq \eta_1$ depending only on $k,n,\lambda, \Lambda,\alpha,\omega_0,K_1,\omega_3$ and $\omega_4$ such that
\begin{equation*}
\begin{aligned}
    \|v-\tilde{P}\|_{L^{\infty }(B_{\eta})}&\leq \eta^{k+\alpha},
\end{aligned}
\end{equation*}
\begin{equation}\label{e.6.2}
  |\tilde G(D^2\tilde{P}(y),D\tilde{P}(y),\tilde{P}(y),y)|\leq \tilde{C}|y|^{k-2+\bar{\alpha}},
  ~~\forall ~y\in B_{1}
\end{equation}
and
\begin{equation*}
\|\tilde{P}\|\leq \tilde{C}.
\end{equation*}

Let
\begin{equation*}
P_{m_0+1}(x)=P_{m_0}(x)+r^{k+\alpha}\tilde{P}(y)=P_{m_0}(x)+r^{\alpha}\tilde{P}(x).
\end{equation*}
Then \cref{In-e.tCkas-7-mu} holds for $m_0+1$ clearly. Next, by rescaling back, \cref{e.6.2} reads
\begin{equation}\label{e6.7}
|G(D^2P_{m_0+1}(x),DP_{m_0+1}(x),P_{m_0+1}(x),x)|
\leq \tilde{C}|x|^{k-2+\bar{\alpha}}, ~~\forall ~x\in B_{r},
\end{equation}
which implies
\begin{equation*}
D^{k-2} \left(G(D^2P_{m_0+1}(x),DP_{m_0+1}(x),P_{m_0+1}(x),x)\right)\bigg |_{x=0}=0.
\end{equation*}
By combining with $G\in C^{k,\bar{\alpha}}$, \cref{e6.7} also holds for any $x\in B_1$, i.e., \cref{In-e.tCkas-9-mu} holds for $m_0+1$. By recalling \cref{In-e.tCkas-v-mu}, we have
\begin{equation*}
  \begin{aligned}
\|u-P_{m_0+1}\|_{L^{\infty}(B_{\eta^{m_0+1}})}&= \|u-P_{m_0}-r^{\alpha}\tilde{P}(x)\|_{L^{\infty}(B_{\eta r})}\\
&= \|r^{k+\alpha}v-r^{k+\alpha}\tilde{P}\|_{L^{\infty}(B_{\eta})}\\
&\leq r^{k+\alpha}\eta^{k+\alpha}=\eta^{(m_0+1)(k+\alpha)}.
  \end{aligned}
\end{equation*}
Hence, \cref{In-e.tCkas-6-mu} holds for $m=m_0+1$. By induction, the proof is completed.\qed~\\

Next, we give the~\\
\noindent\textbf{Proof of \Cref{t-Cka-i}.} As before, in the following proof,  we just make necessary normalization to satisfy the assumptions of \Cref{In-t-Ckas-mu}. Throughout this proof, $C$ always denotes a constant depending only on $k,n,\lambda, \Lambda,\alpha,\mu, b_0,c_0,\omega_0$, $\|F\|_{C^{k-2,\alpha}(0)}$, $K_1,\omega_3,\omega_4, \|f\|_{C^{k-2,\alpha}(0)}$ and $\|u\|_{L^{\infty}(B_1)}$.

For $(M,p,s,x)\in \mathcal{S}^n\times \mathbb{R}^n\times \mathbb{R}\times \bar{B}_1$, let
\begin{equation*}
F_1(M,p,s,x)=F(M,p,s,x)-P_f(x).
\end{equation*}
As usual, $P_f$ is the Taylor polynomial of $f$ at $0$. Then $u$ satisfies
\begin{equation*}
F_1(D^2u,Du,u,x)=f_1 \quad\mbox{in}~~B_1,
\end{equation*}
where $f_1(x)=f(x)-P_f(x)$. Thus,
\begin{equation*}
  |f_1(x)|\leq [f]_{C^{k-2,\alpha}(0)}|x|^{k-2+\alpha}\leq C|x|^{k-2+\alpha}, ~~\forall ~x\in B_1.
\end{equation*}

Note that $u\in C^{k-1,\alpha}(0)$ and set
\begin{equation*}
u_1(x)=u(x)-P_u(x), \quad F_2(M,p,s,x)=F_2(M+D^2P_u(x),p+DP_u(x),s+P_u(x),x).
\end{equation*}
Then $u_1$ satisfies
\begin{equation*}
F_2(D^2u_1,Du_1,u_1,x)=f_1 \quad\mbox{in}~~B_1
\end{equation*}
and
\begin{equation*}
  u_1(0)=|Du_1(0)|=\cdots=|D^{k-1}u_1(0)|=0.
\end{equation*}

Next, take $y=x/\rho$ and $u_2(y)=u_1(x)/\rho^2$, where $0<\rho<1$ is a constant to be specified later. Then $u_2$ satisfies
\begin{equation}\label{In-F6-k-mu}
F_3(D^2u_2,Du_2,u_2,y)=f_2\quad\mbox{in}~~B_1,
\end{equation}
where
\begin{equation*}
  \begin{aligned}
F_3(M,p,s,y)=F_2(M,\rho\, p,\rho^2s,x), \quad f_2(y)=f_1(x).
\end{aligned}
\end{equation*}
Finally, define fully nonlinear operators $G_1,G_2,G_{3}$ in the same way as $F_1,F_2,F_3$ (only replacing $F$ by $G$ in above definitions).

Now, we choose a proper $\rho$ such that \cref{In-F6-k-mu} satisfies the assumptions of \Cref{In-t-Ckas-mu}. Obviously,
\begin{equation*}
u_2(0)=\cdots=|D^{k-1}u_2(0)|=0.
\end{equation*}
By Combining with the $C^{k-1,\alpha}$ regularity at $0$,  we have
\begin{equation*}
  \begin{aligned}
    \|u_2\|_{L^{\infty}(B_1)}&=\rho^{-2} \|u_1\|_{L^{\infty}(B_\rho)}
    \leq [u]_{C^{k-1+\alpha}(0)}\rho^{k-3+\alpha}\leq C\rho^{k-3+\alpha}.
  \end{aligned}
\end{equation*}
In addition,
\begin{equation*}
|f_2(y)|=|f_1(x)|\leq C|x|^{k-2+\alpha}=C\rho^{k-2+\alpha}|y|^{k-2+\alpha},~\forall ~y\in B_1,
\end{equation*}

Next, it is easy to verify that $F_3$ and $G_{3}$ satisfy the structure condition \cref{SC2} with $\lambda,\Lambda,\tilde\mu,\tilde{b},\tilde{c}$ and $\tilde{\omega}_0$, where
\begin{equation*}
\tilde{\mu}= \rho^2\mu, \quad\tilde{b}= \rho b_0+C\rho\mu,\quad\tilde{c}= K_0\rho^2c_0, \quad \tilde{\omega}_0(\cdot,\cdot)=\omega_0(\cdot+C,\cdot).
\end{equation*}
Clearly, $\tilde{\omega}_0$ satisfies \cref{e.omega0-2}.

Finally, we show $F_3\in C^{k-2,\alpha}(0)$. First,
\begin{equation}\label{e5.4}
  \begin{aligned}
&|F_3(M,p,s,y)-G_{3}(M,p,s,y)|\\
&=|F(M+D^2P_u,\rho\, p+DP_u,\rho^2s+P_u,x)-G(M+D^2P_u,\rho\, p+DP_u,\rho^2s+P_u,x)|\\
&\leq \|F\|_{C^{k-2,\alpha}(0)}|x|^{k-2+\alpha}(|M|+C)\omega_3(|p|+C,|s|+C)\\
&\leq C\rho^{k-2+\alpha}|y|^{k-2+\alpha}(|M|+1)\omega_3(|p|+C,|s|+C)\\
&\coloneqq C\rho^{k-2+\alpha}|y|^{k-2+\alpha}(|M|+1)\tilde\omega_3(|p|,|s|).
\end{aligned}
\end{equation}

Furthermore, it can be verified that $G_{3}$ satisfies (i)-(iii) of \Cref{d-FP} with some $\tilde{K}_1$ and
\begin{equation*}
\|G_{3}\|_{C^{k,\bar{\alpha}}(\bar{\mathbf{B}}_r\times \bar{B}_1)}\leq \tilde{\omega}_4(r),~\forall ~r>0,
\end{equation*}
where $\tilde{K}_1$ and $\tilde{\omega}_4$ depend only on $k,n,\lambda, \Lambda,\alpha,\mu, b_0,c_0,\omega_0$, $\|F\|_{C^{k-2,\alpha}(0)},$ $K_1,\omega_3$, $\omega_4$, $\|f\|_{C^{k-2,\alpha}(0)}$ and $\|u\|_{L^{\infty}(B_1)}$.

Take $\delta_1$ small enough such that \Cref{In-t-Ckas-mu} holds for $\tilde{\omega}_0$, $\tilde{K}_1,\tilde{\omega}_3,\tilde{\omega}_4$ and $\delta_1$. From above arguments, we can choose $\rho$ small enough (depending only on $k,n,\lambda, \Lambda,\alpha,\mu, b_0$, $c_0,\omega_0$, $K_1,\omega_3,\omega_4, \|f\|_{C^{k-2,\alpha}(0)}$ and $\|u\|_{L^{\infty}(B_1)}$) such that the conditions of \Cref{In-t-Ckas-mu} are satisfied. Then $u_2$ and hence $u$ is $C^{k,\alpha}$ at $0$, and the estimates \crefrange{e.Cka-1-i}{e.Cka-2-i} hold. \qed~\\

\begin{remark}\label{r-6-3}
In fact, from above proof, we obtain that there exists $P\in \mathcal{P}_k$ such that
\begin{equation*}
  |G(D^2P(x),DP(x),P(x),x)-P_f(x)|\leq C |x|^{k-2+\bar{\alpha}}, ~~\forall ~x\in B_{1}.
\end{equation*}
Note that $F$ is $C^{k-2,\alpha}$ at $0$, \cref{e.Cka-3-i} holds for $F$.
\end{remark}
~\\

\section{Boundary \texorpdfstring{$C^{1,\alpha}$}{C1,a} regularity}\label{C1a-mu}

In the next three sections, we present the boundary pointwise regularity. In this section, we prove the boundary $C^{1,\alpha}$ regularity for the Pucci's class.
\begin{theorem}\label{t-C1a-mu}
Let $0<\alpha<\bar{\alpha}$ and $u$ be a viscosity solution of
\begin{equation*}
\left\{\begin{aligned}
&u\in S^*(\lambda,\Lambda,\mu,b,f)&& \quad\mbox{in}~~\Omega\cap B_1;\\
&u=g&& \quad\mbox{on}~~\partial \Omega\cap B_1.
\end{aligned}\right.
\end{equation*}
Suppose that
\begin{equation*}
b\in L^{p}(\Omega\cap B_1) (p=n/(1-\alpha)),\quad f\in C^{-1,\alpha}(0),\quad\partial \Omega\cap B_1\in C^{1,\alpha}(0),\quad g\in C^{1,\alpha}(0).
\end{equation*}
Then $u\in C^{1,\alpha}(0)$, i.e., there exists $P\in \mathcal{P}_1$ such that
\begin{equation}\label{e.C1a-1-mu}
  |u(x)-P(x)|\leq C |x|^{1+\alpha}, ~\forall ~x\in \Omega\cap B_{1},
\end{equation}
\begin{equation}\label{e.C1a-3-mu}
D_{x'}u(0)=D_{x'}g(0)
\end{equation}
and
\begin{equation}\label{e.C1a-2-mu}
|Du(0)|\leq C ,
\end{equation}
where $C$ depends only on $n$, $\lambda$, $\Lambda$, $\alpha$, $\mu$, $\|b\|_{L^p(\Omega\cap B_1)}$, $\|\partial \Omega\cap B_1\|_{C^{1,\alpha}(0)}$, $\|f\|_{C^{-1,\alpha}(0)}$, $\|g\|_{C^{1,\alpha}(0)}$ and $\|u\|_{L^{\infty }(\Omega\cap B_1)}$.

In particular, if $\mu=0$,
\begin{equation}\label{e.C1a-1}
  |u(x)-P(x)|\leq C |x|^{1+\alpha}\left(\|u\|_{L^{\infty }(\Omega\cap B_1)}+\|f\|_{C^{-1,\alpha}(0)}+\|g\|_{C^{1,\alpha}(0)}\right), ~\forall ~x\in \Omega\cap B_{1}
\end{equation}
and
\begin{equation}\label{e.C1a-2}
|Du(0)|\leq C \left(\|u\|_{L^{\infty }(\Omega\cap B_1)}+\|f\|_{C^{-1,\alpha}(0)}+\|g\|_{C^{1,\alpha}(0)}\right),
\end{equation}
where $C$ depends only on $n, \lambda, \Lambda,\alpha,\|b\|_{L^p(\Omega\cap B_1)}$ and $\|\partial \Omega\cap B_1\|_{C^{1,\alpha}(0)}$.
\end{theorem}

\begin{remark}\label{r-21}
Krylov \cite{MR688919} first obtained boundary $C^{1,\alpha}$ a priori estimate for fully nonlinear equations. Trudinger \cite{MR931007} proved boundary $C^{1,\alpha}$ regularity for Lipschitz continuous viscosity solutions.

In \cite{MR2853528}, Ma and Wang introduced a definition of pointwise $C^{1,\alpha}$ for the boundary, which is similar to \Cref{d-re}. Then they proved (for $\mu=b=0$) the boundary pointwise $C^{1,\alpha}$ regularity for some $\alpha<\bar{\alpha}$ since the Harnack inequality was used. Huang, Zhai and Zhou \cite{MR3933752} extended this result to equations with unbounded coefficient $b$.

Silvestre and Sirakov \cite{MR3246039} proved (for $\mu =0$, $b\in L^{\infty}$) the boundary $C^{1,\alpha}$ regularity for any  $\alpha<\bar{\alpha}$. However, it is not a pointwise regularity since it requires that the whole boundary $(\partial \Omega)_1\in C^2$ rather than the pointwise assumption $(\partial \Omega)_1\in C^{1,\alpha}(0)$. Their proof relies on the technique of flattening the boundary.

For equations with quadratic growth in the gradient, Nornberg \cite{MR3980853} obtained the boundary $C^{1,\alpha}$ regularity for some $\alpha<\bar{\alpha}$. It is not yet a pointwise regularity and $\partial \Omega\in C^{1,1}$ is needed. In addition, the proof does not apply to  the Pucci's class. Recently, this result was extended by da Silva and Nornberg \cite{MR4304555} to equations with more general nonlinear growth in the gradient.

Braga, Gomes, Moreira and Wang \cite{MR4151474} proved the boundary pointwise $C^{1,\alpha}$ regularity for any $\alpha<\bar{\alpha}$ on a flat boundary.

Wang \cite{MR1139064} gave a definition of pointwise $C^{1,\alpha}$ for the boundary for parabolic equations, which is similar to this paper as well. Then he proved (for $\mu=b=0$) the corresponding boundary pointwise $C^{1,\alpha}$ regularity for some $\alpha\leq \bar \alpha$.

The first boundary pointwise $C^{1,\alpha}$ regularity on a curved boundary for any $0<\alpha<\bar{\alpha}$ is given in \cite{MR4088470} (without lower terms). In \Cref{t-C1a-mu}, we extend this result to general equations.
\end{remark}

\begin{remark}\label{r-24}
If $u$ is a viscosity solution of $F(D^2u,Du,u,x)=f$, by the structure condition \cref{SC2},
\begin{equation*}
u\in S^*(\lambda,\Lambda,\mu,b,|f|+|F(0,0,0,\cdot)|
+c\omega_0(\|u\|_{L^{\infty}(\Omega\cap B_1)},\|u\|_{L^{\infty}(\Omega\cap B_1)})).
\end{equation*}
Hence, the boundary pointwise $C^{1,\alpha}$ regularity holds. In fact, in this case, we can obtain the boundary pointwise $C^{1,\alpha}$ regularity for any $0<\alpha<1$ if we assume that $c\in C^{-1,\alpha}(0)$ additionally (cf. \Cref{t-C2a}).
\end{remark}

\begin{remark}\label{r-2.14}
Since we treat the Pucci's class, the method of solving an auxiliary equation to approximate the solution (e.g. \cite[Lemma 3.7]{MR3980853}) is invalid since we can't solve a Dirichlet problem like
\begin{equation*}
\left\{\begin{aligned}
&h\in S^*(\lambda,\Lambda,0,0,0)&& \quad\mbox{in}~~B_1^+;\\
&h=u&& \quad\mbox{on}~~\partial B_1^+
\end{aligned}\right.
\end{equation*}
as done in \cite[Lemma 3.7]{MR3980853}.

Instead, \Cref{t-C1a-mu} can be proved by the method of compactness (e.g. \cite{MR4088470, MR3246039}) or the method of combining the Harnack inequality and a proper barrier (e.g. \cite{MR2853528} and \cite[Theorem 2.1]{MR1139064}). In fact, our proof is inspired directly by \cite{MR4088470}, which can be tracked to \cite{MR3246039} and \cite[Chapter 4]{Wang_Regularity}.
\end{remark}
~\\

Let $\gamma_1$ be as in \cref{e.C1a.beta}, which measures the oscillation of $F$ in $x$ near $x_0$. Similar to \Cref{d-2}, if $\gamma_1(\cdot, x_0)\in C^{-1,\alpha}(x_0)$ for any $x_0\in \Omega$ with the same $r_0$ and
\begin{equation*}
\|\gamma_1\|_{C^{-1,\alpha}(\bar{\Omega})}\coloneqq  \sup_{x_0\in \Omega} \|\gamma_1(\cdot, x_0)\|_{C^{-1,\alpha}(x_0)}<+\infty,
\end{equation*}
we say that $\gamma_1\in C^{-1,\alpha}(\bar{\Omega})$. We also use similar notations for $\|\beta_2\|_{C^{\alpha}(\bar{\Omega})}$, $\|\gamma_2\|_{C^{\alpha}(\bar{\Omega})}$ in the following.

By combining the interior and boundary regularity with standard covering arguments, we have the following local and global $C^{1,\alpha}$ regularity, where we use the same $r_0$ in \Cref{d-f}, \Cref{d-2} and \Cref{d-re} (similarly hereinafter). Its proof is standard and we omit it.
\begin{corollary}\label{t-C1a-global}
Let $0<\alpha<\bar{\alpha}$, $\Gamma\subset \partial \Omega $ be relatively open (may be empty) and $u$ be a viscosity solution of
\begin{equation*}
\left\{\begin{aligned}
&F(D^2u,Du,u,x)=f&& \quad\mbox{in}~~\Omega;\\
&u=g&& \quad\mbox{on}~~\Gamma.
\end{aligned}\right.
\end{equation*}
Suppose that \cref{SC2} holds and $F$ satisfies \cref{e.C1a.beta} with some $G_{x_0}$ at any $x_0\in \Omega\cup \Gamma$. Assume that
\begin{equation*}
  \begin{aligned}
&b\in L^{p}(\Omega) ~(p=n/(1-\alpha)),\quad c\in C^{-1,\alpha}(\bar{\Omega}),\quad f\in C^{-1,\alpha}(\bar{\Omega}),\quad \Gamma\in C^{1,\alpha},\\
&g\in C^{1,\alpha}(\bar{\Gamma}),\quad  \gamma_1\in C^{-1,\alpha}(\bar{\Omega}),\\
&\beta_1(x,x_0)\leq \delta_0, ~\forall ~x_0,x\in \Omega\cup \Gamma ~\mbox{ with }~|x-x_0|<r_0,
  \end{aligned}
\end{equation*}
where $0<\delta_0<1$ depends only on $n,\lambda,\Lambda$ and $\alpha$.

Then for any $\Omega'\subset\subset \Omega\cup \Gamma$, we have $u\in C^{1,\alpha}(\bar{\Omega}')$ and
\begin{equation}\label{e.C1a-1-i-mu-global}
 \|u\|_{C^{1,\alpha}(\bar{\Omega}')}\leq C,
\end{equation}
where $C$ depends only on $n,\lambda,\Lambda,\alpha,r_0,\mu,\|b\|_{L^p(\Omega)}, \|c\|_{C^{-1,\alpha}(\bar{\Omega})},\omega_0$, $\|\gamma_1\|_{C^{-1,\alpha}(\bar{\Omega})}$, $\|\partial\Omega'\cap\Gamma\|_{C^{1,\alpha}}$, $\Omega'$, $\mathrm{dist}(\Omega',\partial \Omega\backslash \Gamma)$, $ \|f\|_{C^{-1,\alpha}(\bar{\Omega})},\|g\|_{C^{1,\alpha}(\bar\Gamma)}$ and $\|u\|_{L^{\infty }(\Omega)}$.

In particular, if $F$ satisfies \cref{SC1} instead of \cref{SC2},
\begin{equation}\label{e.C1a-1-i-global}
 \|u\|_{C^{1,\alpha}(\bar{\Omega}')}\leq C \left(\|u\|_{L^{\infty }(\Omega)}+\|f\|_{C^{-1,\alpha}(\bar{\Omega})}+\|\gamma_1\|_{C^{-1,\alpha}(\bar{\Omega})}
 +\|g\|_{C^{1,\alpha}(\bar\Gamma)}\right),
\end{equation}
where $C$  depends only on $n,\lambda,\Lambda,\alpha,r_0,\|b\|_{L^p(\Omega)}, \|c\|_{C^{-1,\alpha}(\bar{\Omega})}$, $\|\partial \Omega'\cap\Gamma\|_{C^{1,\alpha}}$, $\Omega'$ and $\mathrm{dist}(\Omega',\partial \Omega\backslash \Gamma)$.
\end{corollary}

\begin{remark}\label{r-2.11}
If $F$ is uniformly continuous in $x$, we can use $F(M,0,0,x_0)$ to measure the oscillation of $F$ in $x$ in \cref{e.C1a.beta} instead of $G(M)$. Since $F$ is uniformly continuous, there exists $r_0>0$ such that
\begin{equation*}
  \beta_1(x,x_0)\leq \delta_0,~\forall ~x,x_0\in\Omega\cup \Gamma~\mbox{ with }~|x-x_0|\leq r_0.
\end{equation*}
Hence, the $C^{1,\alpha}$ regularity holds for equations with continuous coefficients.

In above theorem, if $\Gamma=\emptyset$, we will obtain the interior local $C^{1,\alpha}$ estimates analogous to \cref{e.C1a-1-i-mu-global} and \cref{e.C1a-1-i-global} on any $\bar{\Omega}'\subset\subset \Omega$; if $\Gamma=\partial \Omega$, we will obtain the global $C^{1,\alpha}$ estimates analogous to \cref{e.C1a-1-i-mu-global} and \cref{e.C1a-1-i-global} on $\bar{\Omega}$. We have similar remarks for other regularity derived in the following sections.
\end{remark}
~\\

Similar to the interior $C^{1,\alpha}$ regularity, we first prove that the solution in \Cref{t-C1a-mu} can be approximated by a linear function at some scale provided that the coefficients and the prescribed are small enough. Recall that $\Omega_1=\Omega\cap B_1$, $(\partial\Omega)_1=\partial\Omega\cap B_1$ and we use $\underset{B_r}{\mathrm{osc}}~\partial\Omega$ to characterize the oscillation of $\partial \Omega$ in $B_r$ (see \cref{e1.1}).

\begin{lemma}\label{l-C1a-mu}
For any $0<\alpha<\bar{\alpha}$, there exists $\delta>0$ depending only on $n,\lambda,\Lambda$ and $\alpha$ such that if $u$ satisfies
\begin{equation*}
\left\{\begin{aligned}
&u\in S^*(\lambda,\Lambda,\mu,b,f)&& \quad\mbox{in}~~\Omega_1;\\
&u=g&& \quad\mbox{on}~~(\partial \Omega)_1
\end{aligned}\right.
\end{equation*}
with
\begin{equation*}
  \begin{aligned}
&\|u\|_{L^{\infty}(\Omega_1)}\leq 1,\quad\max\left(\mu,\|b\|_{L^{p}(\Omega_1)},\|f\|_{L^{n}(\Omega_1)},\|g\|_{L^{\infty}((\partial \Omega)_1)},\underset{B_1}{\mathrm{osc}}~\partial\Omega\right)\leq \delta,
  \end{aligned}
\end{equation*}
then there exists a constant $a$ such that
\begin{equation*}\label{e.l4.0-mu}
  \|u-ax_n\|_{L^{\infty}(\Omega_{\eta})}\leq \eta^{1+\alpha}
\end{equation*}
and
\begin{equation*}\label{e.14.2-mu}
|a|\leq \bar{C},
\end{equation*}
where $0<\eta<1$ depends only on $n,\lambda,\Lambda$ and $\alpha$.
\end{lemma}

\proof We prove the lemma by contradiction as before. Suppose that the lemma is false. Then there exist $0<\alpha<\bar{\alpha}$ and a sequence of $(u_m,\mu_m,b_m,f_m,g_m,\Omega_m)_{m\in \mathbb{N}}$ satisfying
\begin{equation*}
\left\{\begin{aligned}
&u_m\in S^*(\lambda,\Lambda,\mu_m,b_m,f_m)&& \quad\mbox{in}~~\Omega_{m,1};\\
&u_m=g_m&& \quad\mbox{on}~~(\partial\Omega_{m})_1.
\end{aligned}\right.
\end{equation*}
In addition, $\|u_m\|_{L^{\infty}(\Omega_{m,1})}\leq 1$ and
\begin{equation*}
\max\left(\mu_m,\|b_m\|_{L^{p}(\Omega_{m,1})},
\|f_m\|_{L^{n}(\Omega_{m,1})},\|g_m\|_{L^{\infty}(\Omega_{m,1})},
\underset{B_1}{\mathrm{osc}}~\partial\Omega_m\right)\leq \frac{1}{m},
\end{equation*}
Finally, for any $|a|\leq \bar{C}$,
\begin{equation}\label{e.lC1a.1-mu}
  \|u_m-ax_n\|_{L^{\infty}(\Omega_{m,\eta})}> \eta^{1+\alpha},
\end{equation}
where $0<\eta<1$ is taken small such that
\begin{equation}\label{e.lC1a.2-mu}
\bar{C}\eta^{\bar{\alpha}-\alpha}<1/2.
\end{equation}

As before, $u_m$ are uniformly bounded and equicontinuous (by \Cref{l-3Ho}) in any $\Omega'\subset\subset B_1^+$. Hence, there exists a subsequence (denoted by $u_m$ again) and
$\bar u\colon B_1^+\rightarrow \mathbb{R}$ such that
\begin{equation*}
u_m\rightarrow \bar u~~ \mbox{ in } L^{\infty}_{\mathrm{loc}}(B_1^+).
\end{equation*}

Next, for any $B\subset\subset B_1^+$ and $\varphi\in C^2(\bar{B})$, let $\psi_m=\mathcal{M}^{+}(D^2\varphi)+\mu_m|D\varphi|^2+b_m|D\varphi|+|f_m|$ and $\psi=\mathcal{M}^{+}(D^2\varphi)$. Then,
\begin{equation*}\label{gkC1a-mu}
  \begin{aligned}
  \|\psi_m-\psi\|_{L^n(B)}&=\|\mu_m|D\varphi|^2+b_m|D\varphi|+|f_m|\|_{L^n(B)}\to 0.
  \end{aligned}
\end{equation*}
By \Cref{l-35}, $\bar u\in \underline{S}(\lambda,\Lambda,0)$ in $B_{1}^+$. Similarly, $\bar u\in \bar{S}(\lambda,\Lambda,0)$ and hence
\begin{equation*}
\bar u\in S(\lambda,\Lambda,0)\quad\mbox{ in}~B_1^+.
\end{equation*}

Next, by \Cref{l-33},
\begin{equation*}
|u_m(x)|\leq Cx_n+\frac{C}{m},~\forall ~x \in \Omega_{m,1/2}.
\end{equation*}
Let $m\rightarrow \infty$ and we have
\begin{equation*}
  |\bar u(x)|\leq Cx_n,~\forall ~x\in B^+_{1/2}.
\end{equation*}
Therefore, $u$ is continuous up to $T_{1/2}$ and $\bar u\equiv 0$ on $T_{1/2}$.

Above arguments show that
\begin{equation*}
\left\{\begin{aligned}
&\bar u\in S(\lambda,\Lambda,0)&& \quad\mbox{in}~~B_{1/2}^+;\\
&\bar u=0&& \quad\mbox{on}~~T_{1/2}.
\end{aligned}\right.
\end{equation*}
By \Cref{l-31}, there exists a constant $\bar{a}$ such that
\begin{equation*}
  |\bar u(x)-\bar{a}x_n|\leq \bar{C} |x|^{1+\bar{\alpha}}, ~~\forall ~x\in B_{1/2}^+
\end{equation*}
and
\begin{equation*}
  |\bar{a}|\leq \bar{C}.
\end{equation*}
Combining with \cref{e.lC1a.2-mu}, we have
\begin{equation}\label{e.lC1a.3-mu}
  \|\bar u-\bar{a}x_n\|_{L^{\infty}(B_{\eta}^+)}\leq \frac{1}{2}\eta^{1+\alpha}.
\end{equation}

However, from \cref{e.lC1a.1-mu},
\begin{equation*}
  \|u_m-\bar{a}x_n\|_{L^{\infty}(\Omega_{m,\eta})}> \eta^{1+\alpha}.
\end{equation*}
Let $m\rightarrow \infty$, we have
\begin{equation*}
    \|\bar u-\bar{a}x_n\|_{L^{\infty}(B_{\eta}^+)}\geq \eta^{1+\alpha},
\end{equation*}
which contradicts with \cref{e.lC1a.3-mu}.  ~\qed~\\

Now, we can prove the boundary pointwise $C^{1,\alpha}$ regularity.

\noindent\textbf{Proof of \Cref{t-C1a-mu}.} Since $(\partial \Omega)_1\in C^{1,\alpha}(0)$, we have
\begin{equation}\label{e.tC1a-1-mu}
|x_n|\leq \|(\partial \Omega)_1\|_{C^{1,\alpha}(0)}|x'|^{1+\alpha}, ~\forall ~x\in (\partial \Omega)_1.
\end{equation}
We assume that $P_g\equiv 0$. Otherwise, we may consider $v=u-P_g$. Then the regularity of $u$ follows easily from that of $v$. Since $g\in C^{1,\alpha}(0)$,
\begin{equation}\label{e.tC1a-2-mu}
 |g(x)| \leq [g]_{C^{1,\alpha}(0)} |x|^{1+\alpha},~\forall~x\in (\partial \Omega)_1.
\end{equation}

Let $\delta$ be as in \Cref{l-C1a-mu}, which depends only on $n,\lambda,\Lambda$ and $\alpha$. Without loss of generality, we assume that
\begin{equation}\label{e.tC1a-ass-mu}
  \begin{aligned}
    &\|u\|_{L^{\infty}(\Omega_1)}\leq 1,\quad\mu\leq \frac{\delta}{6C_0^2},\quad\|b\|_{L^p(\Omega_1)}\leq \frac{\delta}{3C_0},\\
    &\|f\|_{C^{-1,\alpha}(0)}\leq \frac{\delta}{3},\quad[g]_{C^{1,\alpha}(0)}\leq \frac{\delta}{2},\quad\|(\partial \Omega)_1\|_{C^{1,\alpha}(0)}\leq \frac{\delta}{2C_0},\\
  \end{aligned}
\end{equation}
where $C_0>1$ is a constant (depending only on $n,\lambda,\Lambda$ and $\alpha$) to be specified later. Otherwise, note that $\Omega$ satisfies the exterior cone condition at $0$ and we may consider for $0<\rho<1$,
\begin{equation}\label{e.c1a-v}
  \bar{u}(y)=\frac{u(x)}{\rho^{\alpha_0}},
\end{equation}
where $y=x/\rho$ and $0<\alpha_0<1$ is a H\"{o}lder exponent (depending only on $n,\lambda,\Lambda$, $\|b\|_{L^p(B_1)}$ and $\|(\partial \Omega)_1\|_{C^{1,\alpha}(0)}$) such that $u\in C^{2\alpha_0}(0)$ (by \Cref{l-3Ho}).
Then we have
\begin{equation*}
\left\{\begin{aligned}
&\bar{u}\in S^{*}(\lambda,\Lambda,\bar{\mu},\bar{b},\bar{f})&&
\quad\mbox{in}~~\tilde{\Omega}_1;\\
&\bar{u}=\bar{g}&& \quad\mbox{on}~~(\partial \tilde{\Omega})_1,
\end{aligned}\right.
\end{equation*}
where
\begin{equation*}\label{e.tC1a-n1-mu}
  \bar{\mu}=\rho^{\alpha_0}\mu,\quad\bar{b}(y)=\rho b(x),\quad\bar{f}(y)=\rho^{2-\alpha_0}f(x),
\quad\bar{g}(y)=\rho^{-\alpha_0}g(x),\quad\tilde{\Omega}=\rho^{-1}\Omega.
\end{equation*}
Hence,
\begin{equation}\label{e7.1}
  \begin{aligned}
  \|\bar{u}\|_{L^{\infty}(\tilde{\Omega}_1)}\leq& \rho^{\alpha_0}[u]_{C^{2\alpha_0}(0)}, \quad
    &&\bar{\mu}=\rho^{\alpha_0}\mu,\\
    \|\bar{b}\|_{L^{p}(\tilde{\Omega}_1)}=&\rho^{1-\frac{n}{p}}\|b\|_{L^{p}(\Omega_\rho)}\leq \rho^{\alpha}\|b\|_{L^{p}(\Omega_1)}, \quad
    &&\|\bar{f}\|_{C^{-1,\alpha}(0)}\leq\rho^{1-\alpha_0+\alpha}\|f\|_{C^{-1,\alpha}(0)},\\
    [\bar{g}]_{C^{1,\alpha}(0)}\leq& \rho^{1-\alpha_0+\alpha}[g]_{C^{1,\alpha}(0)}, \quad
    &&\|(\partial \tilde{\Omega})_1\|_{C^{1,\alpha}(0)}\leq \rho^{\alpha}\|(\partial \Omega)_1\|_{C^{1,\alpha}(0)}.
  \end{aligned}
\end{equation}
By choosing $\rho$ small enough (depending only on $n,\lambda,\Lambda,\alpha,\mu, \|b\|_{L^{p}(\Omega_1)},\|(\partial \Omega)_1\|_{C^{1,\alpha}(0)}$, $\|f\|_{C^{-1,\alpha}(0)},[g]_{C^{1,\alpha}(0)}$ and $\|u\|_{L^{\infty}(\Omega_1)}$), the assumptions \cref{e.tC1a-ass-mu} for $\bar{u}$ can be guaranteed. Hence, we can make the assumption \cref{e.tC1a-ass-mu} for $u$ without loss of generality.

Now, we prove that $u$ is $C^{1,\alpha}$ at $0$ and we only need to prove the following. There exists a sequence of constants $a_m$ ($m\geq -1$) such that for all $m\geq 0$,

\begin{equation}\label{e.tC1a-3-mu}
\|u-a_mx_n\|_{L^{\infty }(\Omega _{\eta^{m}})}\leq \eta ^{m(1+\alpha )}
\end{equation}
and
\begin{equation}\label{e.tC1a-4-mu}
|a_m-a_{m-1}|\leq \bar{C}\eta ^{m\alpha},
\end{equation}
where $\eta$ is as in \Cref{l-C1a-mu} .

We prove the above by induction. For $m=0$, by setting $a_0=a_{-1}=0$, the conclusion holds clearly. Suppose that the conclusion holds for $m\leq m_0$. We need to prove that the conclusion holds for $m=m_0+1$.

Let $r=\eta ^{m_{0}}$, $y=x/r$ and
\begin{equation}\label{e.tC1a-v1-mu}
  v(y)=\frac{u(x)-a_{m_0}x_n}{r^{1+\alpha}}.
\end{equation}
Then $v$ satisfies
\begin{equation*}
\left\{\begin{aligned}
&v\in S^*(\lambda,\Lambda,\tilde{\mu},\tilde{b},\tilde{f})&& \quad\mbox{in}~~\tilde{\Omega}_1;\\
&v=\tilde{g}&& \quad\mbox{on}~~(\partial \tilde{\Omega})_1,
\end{aligned}\right.
\end{equation*}
where
\begin{equation}\label{e.tC1a-n2-mu}
  \begin{aligned}
    &\tilde{\mu}=2r^{1+\alpha}\mu,\quad\tilde{b}(y)=rb(x),\quad
    \tilde{f}(y)=r^{1-\alpha}\left(|f(x)|+b(x)|a_{m_0}|+2\mu|a_{m_0}|^2\right),\\
    &\tilde{g}(y)=r^{-(1+\alpha)}\left(g(x)-a_{m_0}x_n\right),\quad \tilde{\Omega}=r^{-1}\Omega.
  \end{aligned}
\end{equation}

By \cref{e.tC1a-4-mu}, there exists a constant $C_0$ depending only on $n,\lambda,\Lambda$ and $\alpha$ such that
\begin{equation*}
|a_{m}|\leq C_0, ~\forall~0\leq m\leq m_0.
\end{equation*}
Then it is easy to verify that
\begin{equation}\label{e.tC1a-n22-mu}
  \begin{aligned}
\|v\|_{L^{\infty}(\tilde{\Omega}_1)}\leq& 1, ~\quad(\mathrm{by}~ \cref{e.tC1a-3-mu} ~\mathrm{and}~ \cref{e.tC1a-v1-mu})\\
\tilde\mu=& 2r^{1+\alpha}\mu\leq \delta, ~\quad(\mathrm{by}~ \cref{e.tC1a-ass-mu} ~\mathrm{and}~ \cref{e.tC1a-n2-mu})\\
\|\tilde{b}\|_{L^{p}(\tilde{\Omega}_1)}=&r^{1-\frac{n}{p}}\|b\|_{L^{p}(\Omega_r)}\leq  r^{\alpha}\|b\|_{L^{p}(\Omega_1)}\leq \delta,~\quad(\mathrm{by}~ \cref{e.tC1a-ass-mu}~\mathrm{and}~ \cref{e.tC1a-n2-mu})\\
\|\tilde{f}\|_{L^{n}(\tilde{\Omega}_1)}\leq& r^{-\alpha}\|f\|_{L^{n}(\Omega_r)}
    +C_0r^{-\alpha}r^{1-\frac{n}{p}}\|b\|_{L^{p}(\Omega_r)}+2 C_0^2r^{1-\alpha} \mu\\
\leq& \|f\|_{C^{-1,\alpha}(0)}+C_0\|b\|_{L^{p}(\Omega_1)}+2\mu C_0^2\\
\leq& \delta, ~\quad(\mathrm{by}~ \cref{e.tC1a-ass-mu}~\mathrm{and}~ \cref{e.tC1a-n2-mu})\\
\|\tilde{g}\|_{L^{\infty}((\partial \tilde{\Omega})_1)}\leq &
r^{-(1+\alpha)}\left([g]_{C^{1,\alpha}(0)}r^{1+\alpha}+C_0\|(\partial \Omega)_1\|_{C^{1,\alpha}(0)}r^{1+\alpha}\right)\\
\leq& \delta,  ~\quad(\mathrm{by}~\crefrange{e.tC1a-1-mu}{e.tC1a-ass-mu} ~\mathrm{and}~ \cref{e.tC1a-n2-mu})\\
\underset{B_1}{\mathrm{osc}}~\partial\tilde{\Omega}=&
r^{-1}\underset{B_r}{\mathrm{osc}}~\partial\Omega \leq \|(\partial \Omega)_1\|_{C^{1,\alpha}(0)} r^{\alpha} \leq \delta.  ~\quad(\mathrm{by}~ \cref{e.tC1a-1-mu}~\mathrm{and}~ \cref{e.tC1a-ass-mu})
  \end{aligned}
\end{equation}
Hence, by \Cref{l-C1a-mu}, there exists a constant $\tilde{a}$ such that
\begin{equation*}
\begin{aligned}
    \|v-\tilde{a}y_n\|_{L^{\infty }(\tilde{\Omega} _{\eta})}&\leq \eta ^{1+\alpha}
\end{aligned}
\end{equation*}
and
\begin{equation*}
|\tilde{a}|\leq \bar{C}.
\end{equation*}

Let $a_{m_0+1}=a_{m_0}+r^{\alpha}\tilde{a}$. Then \cref{e.tC1a-4-mu} holds for $m_0+1$. By recalling \cref{e.tC1a-v1-mu}, we have
\begin{equation*}
  \begin{aligned}
\|u-a_{m_0+1}x_n\|_{L^{\infty}(\Omega_{\eta^{m_0+1}})}
&= \|u-a_{m_0}x_n-r^{\alpha}\tilde{a}x_n\|_{L^{\infty}(\Omega_{\eta r})}\\
&= \|r^{1+\alpha}v-r^{1+\alpha}\tilde{a}y_n\|_{L^{\infty}(\tilde{\Omega}_{\eta})}\\
&\leq r^{1+\alpha}\eta^{1+\alpha}=\eta^{(m_0+1)(1+\alpha)}.
  \end{aligned}
\end{equation*}
Hence, \cref{e.tC1a-3-mu} holds for $m=m_0+1$. By induction, the proof is completed.

For the special case $\mu=0$, set
\begin{equation*}
K=\|u\|_{L^{\infty}(\Omega_1)}+\delta^{-1}\left(3\|f\|_{C^{-1,\alpha}(0)}
+2[g]_{C^{1,\alpha}(0)}\right)
\end{equation*}
and define for $0<\rho<1$
\begin{equation*}
  \bar{u}(y)=u(x)/K,
\end{equation*}
where $y=x/\rho$. Then by taking $\rho$ small enough (depending only on $n,\lambda,\Lambda,\alpha, \|b\|_{L^{p}(B_1)}$ and $\|(\partial \Omega)_1\|_{C^{1,\alpha}(0)}$), \cref{e.tC1a-ass-mu} can be guaranteed. Hence, for $\mu=0$, we have the explicit estimates \cref{e.C1a-1} and \cref{e.C1a-2}.\qed~\\

\begin{remark}\label{r-42-mu}
From above proof, it shows clearly that the assumption $\partial \Omega \in C^{1,\alpha}(0)$ is used to estimate $x_n$ on $\partial \Omega$ (see \cref{e.tC1a-n22-mu} for the estimate on $g$). This observation is originated from \cite{MR3780142} and pointed out in \cite{MR4088470}. It is also the key to higher regularity in following sections.
\end{remark}
~\\

\section{Boundary \texorpdfstring{$C^{2,\alpha}$}{C2,a} regularity}\label{C2a-mu}

In this section, we prove the boundary pointwise $C^{2,\alpha}$ regularity. The key observation is that if $u(0)=|Du(0)|=0$, the boundary  $C^{2,\alpha}$ regularity holds even if $\partial \Omega\in C^{1,\alpha}(0)$.

In fact, we can say a little more about above idea. Take the linear equation \cref{e.linear} for example. We have the following observation. ~\\
(i) If $u(0)=0$, the smoothness assumption on $c$ can be reduced by one order;~\\
(ii) if $u(0)=|Du(0)=0|$, the smoothness assumption on $c$ can be reduced by two order, and the assumptions on $b^i$, $\partial \Omega$ can be reduced by one order; ~\\
(iii) if $u(0)=|Du(0)|=|D^2u(0)|=0$, the smoothness assumption on $c$, the assumptions on $b^i$, $\partial \Omega$ and the assumption on $a^{ij}$ can be reduced by three order, two order and one order respectively.

For instance, if we intend to derive the $C^{3,\alpha}(0)$ regularity and we know $u(0)=|Du(0)|=|D^2u(0)|=0$ beforehand, then $c\in L^n$, $b^i\in L^{p} (p=n/(1-\alpha))$, $\partial \Omega\in C^{1,\alpha}$ and $a^{ij}\in C^{\alpha}$ are enough (see also \Cref{r-9.1} below). This observation can be used to prove the regularity of free boundary problems, which will be treated in a separate work.

A similar interior pointwise regularity has been obtained by Teixeira \cite{MR3428949}: if $a^{ij}\in C^{\varepsilon}(0)$ for some $0<\varepsilon<1$ and $D^2u(0)=0$, then $u\in C^{2,\alpha}(0)$ for any $0<\alpha<1$.

For the boundary pointwise $C^{2,\alpha}$ regularity, we have
\begin{theorem}\label{t-C2a}
Let $0<\alpha<\bar{\alpha}$ and $u$ be a viscosity solution of
\begin{equation*}
\left\{\begin{aligned}
&F(D^2u,Du,u,x)=f&& \quad\mbox{in}~~\Omega\cap B_1;\\
&u=g&& \quad\mbox{on}~~\partial \Omega\cap B_1.
\end{aligned}\right.
\end{equation*}
Suppose that $F$ satisfies \cref{SC2} and \cref{e.C2a-KF}, and $\omega_0$ satisfies \cref{e.omega0}. Assume that
\begin{equation*}
\beta_2\in C^{\alpha}(0),\quad f\in C^{\alpha}(0),\quad \partial \Omega\cap B_1\in C^{2,\alpha}(0),\quad g\in C^{2,\alpha}(0).
\end{equation*}

Then $u\in C^{2,\alpha}(0)$, i.e., there exists $P\in \mathcal{P}_2$ such that
\begin{equation}\label{e.C2a-1-mu}
  |u(x)-P(x)|\leq C |x|^{2+\alpha}, ~~\forall ~x\in \Omega\cap B_{1},
\end{equation}
\begin{equation}\label{e.C2a-3-mu}
|F(D^2P,DP(x),P(x),x)-f(0)|\leq C|x|^{\alpha}, ~~\forall ~x\in \Omega\cap B_{1},
\end{equation}
\begin{equation}\label{e.C2a-4-mu}
D^l_{x'}\left(u(x',P_{\Omega}(x'))\right)
\coloneqq D^l_{x'}\left(P(x',P_{\Omega}(x'))\right)
=D^l_{x'}\left(g(x',P_{\Omega}(x'))\right)~\mbox{at}~0,~\forall ~0\leq l\leq 2
\end{equation}
and
\begin{equation}\label{e.C2a-2-mu}
|Du(0)|+|D^2u(0)|\leq C,
\end{equation}
where $P_{\Omega}\in \mathcal{HP}_2$ denotes the polynomial in \cref{e-re} and $C$ depends only on $n, \lambda, \Lambda,\alpha,\mu$, $b_0$, $c_0$, $\omega_0$, $\|\beta_2\|_{C^{\alpha}(0)}$, $\omega_2$, $\|\partial \Omega\cap B_1\|_{C^{2,\alpha}(0)}, \|f\|_{C^{\alpha}(0)},\|g\|_{C^{2,\alpha}(0)}$ and $\|u\|_{L^{\infty }(\Omega\cap B_1)}$.

In particular, if $F$ satisfies \cref{SC1} and \cref{e.C2a-KF-0} with $\gamma_2\in C^{\alpha}(0)$,
\begin{equation}\label{e.C2a-1}
  |u(x)-P(x)|\leq \tilde{C} |x|^{2+\alpha}, ~~\forall ~x\in \Omega\cap B_{1},
\end{equation}
\begin{equation}\label{e.C2a-3}
\begin{aligned}
&|F(D^2P,DP(x),P(x),x)-f(0)|\leq \tilde{C}|x|^{\alpha},~\forall ~x\in \Omega\cap B_{1},
\end{aligned}
\end{equation}
\begin{equation}\label{e.C2a-2}
|Du(0)|+|D^2u(0)|\leq \tilde{C}
\end{equation}
and
\begin{equation*}
\tilde{C}=C\left(\|u\|_{L^{\infty }(\Omega\cap B_1)}
+\|f\|_{C^{\alpha}(0)}+\|\gamma_2\|_{C^{\alpha}(0)}+\|g\|_{C^{2,\alpha}(0)}\right)
\end{equation*}
where $C$ depends only on $n, \lambda, \Lambda,\alpha,b_0,c_0$,$\|\beta_2\|_{C^{\alpha}(0)}$ and $\|\partial \Omega\cap B_1\|_{C^{2,\alpha}(0)}$.
\end{theorem}

\begin{remark}\label{r-2.9}
Boundary $C^{2,\alpha}$ a priori estimate for fully nonlinear elliptic equations was first obtained by Krylov \cite{MR688919} (see also \cite[Theorem 9.5]{MR1351007}). Safonov \cite{MR765302} proved boundary $C^{2,\alpha}$ a priori estimate for the Bellman's equation. Silvestre and Sirakov \cite{MR3246039} first derived (for $\mu=c=0$) the boundary $C^{2,\alpha}$ regularity for viscosity solutions. The boundary \emph{pointwise} $C^{2,\alpha}$ regularity was proved (for $\mu=b=c=0$) in \cite{MR4088470}. To the best of our knowledge, \Cref{t-C2a} is the first result concerning the boundary pointwise $C^{2,\alpha}$ regularity for fully nonlinear elliptic equations with lower order terms.
\end{remark}

\begin{remark}\label{r-2.7}
We do not require that $F$ is convex in $M$, which is different from the interior $C^{2,\alpha}$ regularity. The reason is the following. The regularity in \Cref{t-C2a} is obtained as a perturbation of the regularity of the model problem, i.e., \Cref{l-32}. In \Cref{l-32}, the convexity of $F$ in $M$ is not needed. Hence, it is not necessary for \Cref{t-C2a} as well.
\end{remark}
~\\

Similar to $C^{1,\alpha}$ regularity, we have the following corollary.
\begin{corollary}\label{t-C2a-global}
Let $0<\alpha<\bar{\alpha}$, $\Gamma\subset \partial \Omega $ be relatively open and $u$ be a viscosity solution of
\begin{equation*}
\left\{\begin{aligned}
&F(D^2u,Du,u,x)=f&& \quad\mbox{in}~~\Omega;\\
&u=g&& \quad\mbox{on}~~\Gamma.
\end{aligned}\right.
\end{equation*}
Suppose that \cref{SC2} holds and $F$ satisfies \cref{e.C2a-KF} with some $G_{x_0}$ at any $x_0\in \Omega\cup \Gamma$, where $G_{x_0}$ is convex in $M$. Assume that $\omega_0$ satisfies \cref{e.omega0} and
\begin{equation*}
\beta_2\in C^{\alpha}(\bar\Omega),\quad f\in C^{\alpha}(\bar\Omega),\quad \Gamma\in C^{2,\alpha},\quad g\in C^{2,\alpha}(\bar{\Gamma}).
\end{equation*}

Then for any $\Omega'\subset\subset \Omega\cup \Gamma$, we have $u\in C^{2,\alpha}(\bar{\Omega}')$ and
\begin{equation}\label{e.C2a-1-i-mu-global}
 \|u\|_{C^{2,\alpha}(\bar{\Omega}')}\leq C,
\end{equation}
where $C$ depends only on $n, \lambda, \Lambda,\alpha,r_0,\mu,b_0,c_0,\omega_0,\|\beta_2\|_{C^{\alpha}(\bar\Omega)},\omega_2$,
$\|\partial \Omega'\cap \Gamma\|_{C^{2,\alpha}}$, $\Omega'$, $\mathrm{dist}(\Omega',\partial \Omega\backslash \Gamma)$, $\|f\|_{C^{\alpha}(\bar\Omega)}$, $\|g\|_{C^{2,\alpha}(\bar\Gamma)}$ and $\|u\|_{L^{\infty }(\Omega)}$.

In particular, if $F$ satisfies \cref{SC1} and \cref{e.C2a-KF-0} with $\gamma_2\in C^{\alpha}(\bar\Omega)$,
\begin{equation}\label{e.C2a-1-i-global}
 \|u\|_{C^{2,\alpha}(\bar{\Omega}')}\leq C \left(\|u\|_{L^{\infty }(\Omega)}+\|f\|_{C^{\alpha}(\bar\Omega)}+
 \|\gamma_2\|_{C^{\alpha}(\bar\Omega)}+\|g\|_{C^{2,\alpha}(\bar\Gamma)}\right),
\end{equation}
\sloppy where $C$ depends only on $n$, $\lambda$, $\Lambda$, $\alpha$, $r_0$, $b_0$, $c_0$, $\|\beta_2\|_{C^{\alpha}(\bar\Omega)}$, $\|\partial \Omega'\cap\Gamma\|_{C^{2,\alpha}}$, $\Omega'$ and $\mathrm{dist}(\Omega',\partial \Omega\backslash \Gamma)$.
\end{corollary}

\begin{remark}\label{re7.1}
Since $\Gamma\in C^{2,\alpha}$, this corollary can be proved by the method of flattening the boundary as well (see \cite{MR3246039}).
\end{remark}
~\\

Before proving the boundary $C^{2,\alpha}$ regularity, we first prove the following key step as before. Note that the key to proving the boundary regularity is to estimate $x_n$ on the boundary. Hence, for boundary regularity, we will frequently use the polynomials in $\mathcal{SP}_k$, i.e., the polynomials which contain at least one power of $x_n$ (see \Cref{no1.1} for details).
\begin{lemma}\label{l-C2a-mu}
Suppose that $F$ satisfies \cref{e.C2a-KF}. For any $0<\alpha<\bar{\alpha}$, there exists $\delta>0$ depending only on $n,\lambda,\Lambda,\alpha$ and $\omega_2$ such that if $u$ satisfies
\begin{equation*}
\left\{\begin{aligned}
&F(D^2u,Du,u,x)=f&& \quad\mbox{in}~~\Omega_1;\\
&u=g&& \quad\mbox{on}~~(\partial \Omega)_1
\end{aligned}\right.
\end{equation*}
with
\begin{equation*}
  \begin{aligned}
&u(0)=|Du(0)|=0,\quad \max\left(\|u\|_{L^{\infty}(\Omega_1)},\omega_0(1,1)\right)\leq 1,\\
&\max\left(\mu,b_0,c_0,\|\beta_2\|_{L^{\infty}(\Omega_1)},\|f\|_{L^{\infty}(\Omega_1)},
\|g\|_{C^{1,\alpha}(0)},\|(\partial \Omega)_1\|_{C^{1,\alpha}(0)}\right)\leq \delta,\\
  \end{aligned}
\end{equation*}
then there exists $P\in\mathcal{SP}_2$ such that
\begin{equation*}
  \|u-P\|_{L^{\infty}(\Omega_{\eta})}\leq \eta^{2+\alpha},
\end{equation*}
\begin{equation*}
G(D^2P,0,0)=0
\end{equation*}
and
\begin{equation*}
\|P\|\leq \bar{C}+1,
\end{equation*}
where $\eta$ depends only on $n,\lambda,\Lambda$ and $\alpha$.
\end{lemma}

\proof We prove the lemma by contradiction as before. Since the proof is similar to the previous proofs, we only point out the outline. Suppose that the lemma is false. Then there exist a sequence of $(F_m,u_m,f_m,g_m,\Omega_m)_{m\in \mathbb{N}}$ satisfying
\begin{equation*}
\left\{\begin{aligned}
&F_m(D^2u_m,Du_m,u_m,x)=f_m&& \quad\mbox{in}~~\Omega_{m,1};\\
&u_m=g_m&& \quad\mbox{on}~~(\partial \Omega_m)_1.
\end{aligned}\right.
\end{equation*}
But for any $P\in\mathcal{SP}_2$ satisfying $\|P\|\leq \bar{C}+1$ and
\begin{equation}\label{e.lC2a-3-mu}
G_{m}(D^2P,0,0)=0,
\end{equation}
we have
\begin{equation}\label{e.lC2a-1-mu}
  \|u_m-P\|_{L^{\infty}(\Omega_{m,\eta})}> \eta^{2+\alpha},
\end{equation}
where $0<\eta<1$ is taken small such that
\begin{equation}\label{e.lC2a-2-mu}
\bar{C}\eta^{\bar{\alpha}-\alpha}<1/2.
\end{equation}

As before, there exist $\bar u\colon B_1^+\cup T_1\rightarrow \mathbb{R}$ and $\bar G\colon \mathcal{S}^n\rightarrow \mathbb{R}$ such that
\begin{equation*}
  \begin{aligned}
&u_m\rightarrow \bar u \quad\mbox{in}~L^{\infty}_{loc}(B_1^+),\quad
G_{m}(\cdot,0,0)\rightarrow \bar G \quad\mbox{in}~L^{\infty}_{loc}(\mathcal{S}^n).
  \end{aligned}
\end{equation*}
Furthermore, $\bar u$ is a viscosity solution of
\begin{equation*}
\left\{\begin{aligned}
&\bar G(D^2\bar u)=0&& \quad\mbox{in}~~B_{1}^+;\\
&\bar u=0&& \quad\mbox{on}~~T_{1}.
\end{aligned}\right.
\end{equation*}

Note that $u_m\in S^*(\lambda,\Lambda,\mu_m,b_m,|f_m|+c_m+|F_m(0,0,0,\cdot)|)$. By the boundary $C^{1,\alpha}$ estimate for $u_m$ (see \Cref{t-C1a-mu}) and noting  $u_m(0)=|Du_m(0)|=0$, we have
\begin{equation*}
\|u_m\|_{L^{\infty }(\Omega_{m,r})}\leq Cr^{1+\alpha}, ~~~~\forall~0<r<1.
\end{equation*}
Let $m\to \infty$,
\begin{equation*}
\|\bar u\|_{L^{\infty }(B_r^+)}\leq Cr^{1+\alpha}, ~~~~\forall~0<r<1,
\end{equation*}
which implies
\begin{equation}\label{e8.2}
\bar u(0)=|D\bar u(0)|=0.
\end{equation}

By \Cref{l-32} and noting \cref{e8.2}, there exists $\bar{P}\in\mathcal{SP}_2$ such that
\begin{equation}\label{e.c2a-3-mu}
  |\bar u(x)-\bar {P}(x)|\leq \bar{C} |x|^{2+\bar{\alpha}}, ~~\forall ~x\in B_{1}^+,
\end{equation}
\begin{equation*}
  \bar G(D^2\bar {P})=0
\end{equation*}
and
\begin{equation*}
 \|\bar {P}\|\leq \bar{C}.
\end{equation*}
By combining \cref{e.lC2a-2-mu} and \cref{e.c2a-3-mu},
\begin{equation}\label{e.lC2a-4-mu}
  \|\bar u-\bar {P}\|_{L^{\infty}(B_{\eta}^+)}\leq \frac{1}{2}\eta^{2+\alpha}.
\end{equation}

Similar to the interior $C^{2,\alpha}$ regularity (see \Cref{In-l-C2a-mu} and its proof), there exist a sequence of constants $t_m\rightarrow 0$ such that
\begin{equation*}
  \begin{aligned}
G_{m}(D^2P_m,0,0)=0,
  \end{aligned}
\end{equation*}
where $P_m(x)=\bar {P}(x)+t_mx_n^2/2$.

Hence, \cref{e.lC2a-1-mu} holds for $P_m$, i.e.,
\begin{equation*}
  \|u_m-P_m\|_{L^{\infty}(\Omega_{m,\eta})}> \eta^{2+\alpha}.
\end{equation*}
Let $m\rightarrow \infty$, we have
\begin{equation*}
    \|\bar u-\bar {P}\|_{L^{\infty}(B_{\eta}^+)}\geq \eta^{2+\alpha},
\end{equation*}
which contradicts with \cref{e.lC2a-4-mu}.  ~\qed~\\

The following is a special result for the boundary pointwise $C^{2,\alpha}$ regularity, i.e., if $u(0)=|Du(0)|=0$, the $C^{2,\alpha}$ regularity holds even if $\partial \Omega\in C^{1,\alpha}(0)$. This was first observed in \cite{MR4088470}.
\begin{lemma}\label{t-C2as-mu}
Let $0<\alpha <\bar{\alpha}$ and $F$ satisfy \cref{e.C2a-KF}. Suppose that $\omega_0$ satisfies \cref{e.omega0} and $u$ is a viscosity solution of
\begin{equation*}
\left\{\begin{aligned}
&F(D^2u,Du,u,x)=f&& \quad\mbox{in}~~\Omega_1;\\
&u=g&& \quad\mbox{on}~~(\partial \Omega)_1.
\end{aligned}\right.
\end{equation*}
Assume that
\begin{equation}\label{e.C2as-be-mu}
\begin{aligned}
&\|u\|_{L^{\infty}(\Omega_1)}\leq 1,\quad u(0)=|Du(0)|=0,\quad \mu\leq \frac{\delta_1}{4C_0},\quad b_0\leq \frac{\delta_1}{2},\quad c_0\leq \frac{\delta_1}{K_0},\\
&\omega_0(1+C_0,1)\leq 1,\quad \|(\partial \Omega)_1\|_{C^{1,\alpha}(0)} \leq \frac{\delta_1}{2C_0},\\
&|\beta_2(x)|\leq \frac{\delta_1}{C_0}|x|^{\alpha},\quad|f(x)|\leq \delta_1|x|^{\alpha},\quad
|g(x)|\leq \frac{\delta_1}{2}|x|^{2+\alpha}, ~\forall ~x\in \Omega_1,\\
\end{aligned}
\end{equation}
where $C_0$ depends only on $n,\lambda,\Lambda$ and $\alpha$, and $0<\delta_1<1$ depends also on $\omega_0$ and $\omega_2$.

Then $u\in C^{2,\alpha}(0)$, i.e., there exists $P\in\mathcal{SP}_2$ such that
\begin{equation}\label{e.C2as-1-mu}
  |u(x)-P(x)|\leq C |x|^{2+\alpha}, ~~\forall ~x\in \Omega_{1},
\end{equation}
\begin{equation}\label{e.C2as-3-mu}
G(D^2P,0,0)=0
\end{equation}
and
\begin{equation}\label{e.C2as-2-mu}
\|P\| \leq C,
\end{equation}
where $C$ depends only on $n, \lambda, \Lambda$ and $\alpha$.
\end{lemma}

\proof As before, to prove that $u$ is $C^{2,\alpha}$ at $0$, we only need to prove the following. There exist a sequence of $P_m\in\mathcal{SP}_2$ ($m\geq -1$) such that for all $m\geq 0$,
\begin{equation}\label{e.C2as-4-mu}
\|u-P_m\|_{L^{\infty }(\Omega_{\eta^{m}})}\leq \eta ^{m(2+\alpha )},
\end{equation}
\begin{equation}\label{e.C2as-5-mu}
G(D^2P_m,0,0)=0
\end{equation}
and
\begin{equation}\label{e.C2as-6-mu}
\|P_m-P_{m-1}\|\leq (\bar{C}+1)\eta ^{(m-1)\alpha},
\end{equation}
where $\eta$ is as in \Cref{l-C2a-mu}.

We prove the above by induction. For $m=0$, by setting $P_0=P_{-1}\equiv 0$, \crefrange{e.C2as-4-mu}{e.C2as-6-mu} hold clearly. Suppose that the conclusion holds for $m\leq m_0$. We need to prove that the conclusion holds for $m=m_0+1$.

Let $r=\eta ^{m_{0}}$, $y=x/r$ and
\begin{equation}\label{e.C2as-v-mu}
  v(y)=\frac{u(x)-P_{m_0}(x)}{r^{2+\alpha}}.
\end{equation}
Then $v$ satisfies
\begin{equation}\label{e.C2as-F-mu}
\left\{\begin{aligned}
&\tilde{F}(D^2v,Dv,v,y)=\tilde{f}&& \quad\mbox{in}~~\tilde{\Omega}_1;\\
&v=\tilde{g}&& \quad\mbox{on}~~(\partial \tilde{\Omega})_1,
\end{aligned}\right.
\end{equation}
where for $(M,p,s,y)\in \mathcal{S}^n\times \mathbb{R}^n\times \mathbb{R}\times \bar{\tilde{\Omega}}_1$,
\begin{equation*}\label{e.c2a.b}
  \begin{aligned}
&\tilde{F}(M,p,s,y)=r^{-\alpha}
F(r^{\alpha}M+D^2P_{m_0}(x),r^{1+\alpha}p+DP_{m_0}(x),r^{2+\alpha}s
+P_{m_0}(x),x),\\
&\tilde{f}(y)=r^{-\alpha}f(x),\quad\tilde{g}(y)=r^{-(2+\alpha)}\left(g(x)-P_{m_0}(x)\right),
\quad\tilde{\Omega}=r^{-1}\Omega.\\
  \end{aligned}
\end{equation*}
In addition, define
\begin{equation*}
\tilde{G}(M,p,s)=r^{-\alpha}G(r^{\alpha}M+D^2P_{m_0},r^{1+\alpha}p
,r^{2+\alpha}s).
\end{equation*}

In the following, we show that \cref{e.C2as-F-mu} satisfies the assumptions of \Cref{l-C2a-mu}. First, it is easy to verify that
\begin{equation*}
\begin{aligned}
\|v\|_{L^{\infty}(\tilde{\Omega}_1)}\leq& 1, \quad  v(0)=|Dv(0)|=0,
~\quad(\mathrm{by}~ \cref{e.C2as-be-mu},~\cref{e.C2as-4-mu}~\mbox{and}~ \cref{e.C2as-v-mu})\\
\|\tilde{f}\|_{L^{\infty}(\tilde{\Omega}_1)}
=&r^{-\alpha}\|f\|_{L^{\infty}(\Omega_r)}\leq \delta_1,
 ~\quad(\mathrm{by}~\cref{e.C2as-be-mu})\\
\|(\partial \tilde{\Omega})_1\|_{C^{1,\alpha}(0)} \leq& r^{\alpha}\|(\partial\Omega)_1\|_{C^{1,\alpha}(0)}\leq \delta_1,
~\quad(\mathrm{by}~\cref{e.C2as-be-mu})\\
\tilde{G}(0,0,0)=&r^{-\alpha}G(D^2P_{m_0},0,0)=0.~\quad(\mathrm{by}~\cref{e.C2as-5-mu})
\end{aligned}
\end{equation*}

By \cref{e.C2as-6-mu}, there exists a constant $C_0$ depending only on $n,\lambda,\Lambda$ and $\alpha$ such that
\begin{equation*}
\|P_m\|\leq C_0,~\forall~0\leq m\leq m_0.
\end{equation*}
For any $0<\rho<1$ (note that $P_{m_0}\in\mathcal{SP}_2$),
\begin{equation*}
\|\tilde{g}\|_{L^{\infty}((\partial \tilde{\Omega})_{\rho})}
\leq \frac{1}{r^{2+\alpha}}\left(\frac{\delta _1}{2}(\rho r)^{2+\alpha}+
C_0\cdot \frac{\delta_1}{2C_0}(\rho r)^{2+\alpha}\right)\leq \delta_1 \rho^{2+\alpha}.
~~(\mathrm{by}~\cref{e.C2as-be-mu})
\end{equation*}
Hence,
\begin{equation*}
  \|\tilde{g}\|_{C^{1,\alpha}(0)}\leq \|\tilde{g}\|_{C^{2,\alpha}(0)}\leq\delta_1.
\end{equation*}

It is easy to verify that $\tilde{F}$ and $\tilde{G}$ satisfy the structure condition \cref{SC2} with $\lambda,\Lambda,\tilde{\mu},\tilde{b},\tilde{c}$ and $\tilde{\omega}_0$, where
\begin{equation*}
\tilde{\mu}= r^{2+\alpha}\mu,\quad\tilde{b}= rb_0+2 C_0 r\mu ,\quad\tilde{c}= K_0 r^{\alpha+\alpha^2}c_0,\quad\tilde{\omega}_0(\cdot,\cdot)=\omega_0(\cdot+C_0,\cdot).
\end{equation*}
Hence, from \cref{e.C2as-be-mu},
\begin{equation*}
 \tilde{\mu}\leq \delta_1,\quad
 \tilde{b}\leq \delta_1, \quad\tilde{c}\leq \delta_1,\quad\tilde{\omega}_0(1,1)\leq 1.
\end{equation*}

Similar to the interior $C^{2,\alpha}$ regularity, by combining \cref{SC2}, \cref{e.C2a-KF} and \cref{e.C2as-be-mu},
\begin{equation*}
\begin{aligned}
|\tilde{F}(&M,p,s,y)-\tilde{G}(M,p,s)|\\
\leq& r^{-\alpha}\Big(2C_0r^{1+\alpha}\mu |p||x|+ C_0^2\mu|x|^2+C_0 b_0|x|+K_0 c_0\omega_0(|s|+C_0,C_0)|x|^{2\alpha}\\
&+\beta_2(x)(|M|+C_0)\omega_2(|p|,|s|)\Big)\\
\leq& \tilde{\beta}_2(y)(|M|+1)\tilde{\omega}_2(|p|,|s|),
\end{aligned}
\end{equation*}
where
\begin{equation*}
\tilde{\beta}_2(y)\coloneqq \delta_1|y|^{\alpha},\quad\tilde{\omega}_2(|p|,|s|)\coloneqq
 \omega_2(|p|,|s|)+|p|+C_0+\omega_0(|s|+C_0,C_0).
\end{equation*}
Then $\|\tilde{\beta}_2\|_{L^{\infty}(\tilde{\Omega}_1)}\leq \delta_1$.

Choose $\delta_1$ small enough (depending only on $n,\lambda,\Lambda,\alpha,\omega_0$ and $\omega_2$) such that \Cref{l-C2a-mu} holds for $\tilde{\omega}_0,\tilde{\omega}_2$ and $\delta_1$. Since \cref{e.C2as-F-mu} satisfies the assumptions of \Cref{l-C2a-mu}, there exists $\tilde{P}(y)\in\mathcal{SP}_2$ such that
\begin{equation*}
\begin{aligned}
    \|v-\tilde{P}\|_{L^{\infty }(\tilde{\Omega} _{\eta})}&\leq \eta ^{2+\alpha},
\end{aligned}
\end{equation*}
\begin{equation*}
\tilde{G}(D^2\tilde{P},0,0)=0
\end{equation*}
and
\begin{equation*}
\|\tilde{P}\|\leq \bar{C}+1.
\end{equation*}

Let
\begin{equation*}
P_{m_0+1}(x)=P_{m_0}(x)+r^{2+\alpha}\tilde{P}(y)=P_{m_0}(x)+r^{\alpha}\tilde{P}(x).
\end{equation*}
Then \cref{e.C2as-5-mu} and \cref{e.C2as-6-mu} hold for $m_0+1$. By recalling \cref{e.C2as-v-mu}, \begin{equation*}
  \begin{aligned}
\|u-P_{m_0+1}\|_{L^{\infty}(\Omega_{\eta^{m_0+1}})}
&= \|u-P_{m_0}(x)-r^{\alpha}\tilde{P}(x)\|_{L^{\infty}(\Omega_{\eta r})}\\
&= \|r^{2+\alpha}v-r^{2+\alpha}\tilde{P}(y)\|_{L^{\infty}(\tilde{\Omega}_{\eta})}\\
&\leq r^{2+\alpha}\eta^{2+\alpha}=\eta^{(m_0+1)(2+\alpha)}.
  \end{aligned}
\end{equation*}
Hence, \cref{e.C2as-4-mu} holds for $m=m_0+1$. By induction, the proof is completed.\qed~\\

Now, we give the~\\
\noindent\textbf{Proof of \Cref{t-C2a}.} As before, the following proof is mere a normalization procedure. We prove the theorem in two cases.

\textbf{Case 1:} the general case, i.e., $F$ satisfies \cref{SC2} and \cref{e.C2a-KF}. Throughout the proof for this case, $C$ always denotes a constant depending only on $n, \lambda,\Lambda,\alpha,\mu,b_0,c_0,\omega_0$, $\|\beta_2\|_{C^{\alpha}(0)}$, $\omega_2$, $\|(\partial \Omega)_1\|_{C^{2,\alpha}(0)}, \|f\|_{C^{\alpha}(0)},
\|g\|_{C^{2,\alpha}(0)}$ and $\|u\|_{L^{\infty }(\Omega_1)}$.

Let
\begin{equation*}
  \begin{aligned}
&f_1(x)=f(x)-f(0),\quad g_1(x)=g(x)-P_{g}(x),\quad u_1(x)=u(x)-P_{g}(x),\\
&F_1(M,p,s,x)=F(M+D^2P_{g},p+DP_{g}(x),s+P_{g}(x),x)-f(0).
  \end{aligned}
\end{equation*}
Then $u_1$ satisfies
\begin{equation*}
\left\{\begin{aligned}
&F_1(D^2u_1,Du_1,u_1,x)=f_1&& \quad\mbox{in}~~\Omega_1;\\
&u_1=g_1&& \quad\mbox{on}~~(\partial \Omega)_1
\end{aligned}\right.
\end{equation*}
and
\begin{equation*}
  \begin{aligned}
  |F_1(0,0,0,x)|=&\left|F\left(D^2P_{g}(x),DP_{g}(x),P_{g},x\right)-f(0)\right|\\
  =&\big|F\left(D^2P_{g}(x),DP_{g}(x),P_{g}(x),x\right)-G(D^2P_{g}(x),DP_{g}(x),P_{g}(x))\\
  &+G(D^2P_{g}(x),DP_{g}(x),P_{g})-G(0,0,0)-f(0)\big|\\
  \leq& \beta_2(x)(C+1)\omega_2(C,C)+C+|f(0)|\leq C.
  \end{aligned}
\end{equation*}

Note that
\begin{equation*}
  u_1\in S^{*}(\lambda,\Lambda,\mu,\hat{b},|f_1|
  +c_0\omega_0(\|u\|_{L^{\infty}(\Omega_1)}+\|g\|_{C^{2,\alpha}(0)},u_1)+|F_1(0,0,0,\cdot)|),
\end{equation*}
where $\hat{b}=b_0+2\mu\|g\|_{C^{2,\alpha}(0)}$. By \Cref{t-C1a-mu}, $u_1\in C^{1,\alpha}(0)$,
\begin{equation*}
Du_1(0)=(D_{x'}u_1(0),(u_1)_n(0))=(D_{x'}g_1(0),(u_1)_n(0))=(0,...,0,(u_1)_n(0))
\end{equation*}
and
\begin{equation}\label{e.tC2a-1-mu}
  \begin{aligned}
|(u_1)_n(0)| &\leq C.
  \end{aligned}
\end{equation}

Let
\begin{equation*}
P_u(x)=(u_1)_n(0)\left(x_n-P_{\Omega}(x')\right),
\end{equation*}
where $P_{\Omega}\in \mathcal{HP}_2$ is from \cref{e-re}. Define
\begin{equation*}
u_2=u_1-P_u,\quad F_2(M,p,s,x)=F_1(M+D^2P_u,p+DP_u(x),s+P_u(x),x),
\end{equation*}
Then $u_2$ satisfies
\begin{equation*}
\left\{\begin{aligned}
&F_2(D^2u_2,Du_2,u_2,x)=f_1&& \quad\mbox{in}~~\Omega_1;\\
&u_2=g_2&& \quad\mbox{on}~~(\partial \Omega)_1,
\end{aligned}\right.
\end{equation*}
where $g_2=g_1-P_u$ and
\begin{equation*}
u_2(0)=|Du_2(0)|=0.
\end{equation*}
Since $g\in C^{2,\alpha}(0)$ and $\partial \Omega\in C^{2,\alpha}(0)$,
\begin{equation*}
  |g_2(x)|\leq |g_1(x)|+|P_u(x)|\leq |g(x)-P_g(x)|+C|x_n-P_{\Omega}(x')|\leq C|x|^{2+\alpha}, ~~\forall ~x\in (\partial \Omega)_1.
\end{equation*}
We remark here that this step shows clearly where the condition $\partial \Omega\in C^{2,\alpha}(0)$ is used. If $\partial \Omega \in C^{2,\alpha}(0)$, we can assume that $Du_2(0)=0$ and keep that $g_2$ has a decay of order $2+\alpha$ near $0$.

Next, for $\tau\in \mathbb{R}$ (to be chosen later), set
\begin{equation*}
  \begin{aligned}
&u_3(x)=u_2(x)-\tau x_n^2,\quad g_3(x)=g_2(x)-\tau x_n^2~\\
&F_3(M,p,s,x)=F_2(M+2\tau\tilde{I},p+2\tau x_n,s+\tau x_n^2,x).
  \end{aligned}
\end{equation*}
Then $u_3$ is a viscosity solution of
\begin{equation*}
\left\{\begin{aligned}
&F_3(D^2u_3,Du_3,u_3,x)=f_1&& \quad\mbox{in}~~\Omega_1;\\
&u_3=g_3&& \quad\mbox{on}~~(\partial \Omega)_1.
\end{aligned}\right.
\end{equation*}
Moreover, $u_3(0)=|Du_3(0)|=0$ and
\begin{equation*}
  |g_3(x)|\leq C|x|^{2+\alpha}, ~~\forall ~x\in (\partial \Omega)_1.
\end{equation*}

Define fully nonlinear operators $G_{1},G_{2}$ and $G_{3}$ in a similar way as $F_1,F_2$ and $F_3$. By \cref{SC2}, there exists $\tau\in \mathbb{R}$ such that $G_{3}(0,0,0)=0$ and
\begin{equation}\label{e.tC2a-2-mu}
  \begin{aligned}
|\tau|&\leq |G_{2}(0,0,0)|/\lambda\leq C.
  \end{aligned}
\end{equation}

Finally, for $0<\rho<1$, let
\begin{equation*}
  \begin{aligned}
&y=\rho^{-1}x,\quad u_4(y)=\rho^{-1}u_3(x),\quad f_2(y)=\rho f_1(x),\quad g_4(y)=\rho^{-1}g_3(x),\quad \tilde{\Omega}=\rho^{-1}\Omega,\\
&F_4(M,p,s,y)=\rho F_3\left(\rho^{-1}M, p,\rho s,x\right),\quad G_{4}(M,p,s)=\rho G_{3}\left(\rho^{-1}M, p,\rho s\right).
  \end{aligned}
\end{equation*}
Then $u_4$ satisfies
\begin{equation}\label{F5-mu}
\left\{\begin{aligned}
&F_4(D^2u_4,Du_4,u_4,y)=f_2&& \quad\mbox{in}~~\tilde{\Omega}_1;\\
&u_4=g_4&& \quad\mbox{on}~~(\partial \tilde\Omega)_1.
\end{aligned}\right.
\end{equation}

Then, it can checked as before that \cref{F5-mu} satisfies the conditions of \Cref{t-C2as-mu} by choosing a sufficiently small $\rho$. By \Cref{t-C2as-mu}, $u_4$ and hence $u$ is $C^{2,\alpha}$ at $0$, and the estimates \crefrange{e.C2a-1-mu}{e.C2a-2-mu} hold.

\textbf{Case 2:} $F$ satisfies \cref{SC1} and \cref{e.C2a-KF-0}. Let
\begin{equation*}
K=\|u\|_{L^{\infty }(\Omega_1)}+\|f\|_{C^{\alpha}(0)}+\|g\|_{C^{2,\alpha}(0)}
+\|\gamma_2\|_{C^{\alpha}(0)}, \quad u_1=u/K.
\end{equation*}
Then $u_1$ satisfies
\begin{equation}\label{e.8.1}
\left\{\begin{aligned}
&F_1(D^2u_1,Du_1,u_1,x)=f_1&& \quad\mbox{in}~~\Omega_1;\\
&u_1=g_1&& \quad\mbox{on}~~(\partial \Omega)_1,
\end{aligned}\right.
\end{equation}
where
\begin{equation*}
F_1(M,p,s,x)=F(KM,Kp,Ks,x)/K, \quad f_1=f/K,\quad g_1=g/K.
\end{equation*}

Then by applying \textbf{Case 1} to \cref{e.8.1} as before, we obtain that $u_1$ and hence $u$ is $C^{2,\alpha}$ at $0$, and the estimates \crefrange{e.C2a-1}{e.C2a-2} hold. \qed~\\

\section{Boundary \texorpdfstring{$C^{k,\alpha}$}{Ck,a} regularity}\label{Cka-mu}
In this section, we prove the following boundary pointwise $C^{k,\alpha}$ regularity for $k\geq 3$.
\begin{theorem}\label{t-Cka}
Let $k\geq 3$, $0<\alpha<\bar{\alpha}$ and $u$ be a viscosity solution of
\begin{equation*}
\left\{\begin{aligned}
&F(D^2u,Du,u,x)=f&& \quad\mbox{in}~~\Omega\cap B_1;\\
&u=g&& \quad\mbox{on}~~\partial \Omega\cap B_1.
\end{aligned}\right.
\end{equation*}
Suppose that $\omega_0$ satisfies \cref{e.omega0-2} and
\begin{equation*}
F\in C^{k-2,\alpha}(0),\quad f\in C^{k-2,\alpha}(0),\quad \partial \Omega\cap B_1\in C^{k,\alpha}(0),\quad g\in C^{k,\alpha}(0).
\end{equation*}

Then $u\in C^{k,\alpha}(0)$, i.e., there exists $P\in \mathcal{P}_k$ such that
\begin{equation}\label{e.Cka-1}
  |u(x)-P(x)|\leq C |x|^{k+\alpha}, ~~\forall ~x\in \Omega\cap B_{1},
\end{equation}
\begin{equation}\label{e.Cka-3}
|F(D^2P(x),DP(x),P(x),x)-P_f(x)|\leq C|x|^{k-2+\alpha}, ~~\forall ~x\in \Omega\cap B_{1},
\end{equation}
\begin{equation}\label{e.Cka-5}
D^l_{x'}u(x',P_{\Omega}(x'))=D^l_{x'}g(x',P_{\Omega}(x'))~\mbox{at}~0,~\forall ~0\leq l\leq k
\end{equation}
and
\begin{equation}\label{e.Cka-2}
|Du(0)|+\cdots+|D^ku(0)|\leq C,
\end{equation}
where $C$ depends only on $k,n,\lambda, \Lambda,\alpha,\mu, b_0,c_0$, $\omega_0$, $\|F\|_{C^{k-2,\alpha}(0)}$, $K_1,\omega_3,\omega_4$, $\|\partial \Omega\cap B_1\|_{C^{k,\alpha}(0)}$, $\|f\|_{C^{k-2,\alpha}(0)}$, $\|g\|_{C^{k,\alpha}(0)}$ and $\|u\|_{L^{\infty}(\Omega\cap B_1)}$.
\end{theorem}

\begin{remark}\label{r-4}
As we know, this boundary pointwise $C^{k,\alpha}$ ($k\geq 3$) regularity is new even for the linear equations.

The \cref{e.Cka-5} is not only a conclusion of the theorem, but also the key to the proof of boundary pointwise $C^{k,\alpha}$ regularity. Indeed, based on \cref{e.Cka-5}, we can construct a polynomial $P$ such that all its derivatives coincide with that of $u$ (see \cref{e8.1}) and $P$ has the desired decay on the boundary (see \cref{e9.3}).
\end{remark}
~\\

Similar to $C^{1,\alpha}$ and $C^{2,\alpha}$ regularity, we have the following corollary.
\begin{corollary}\label{t-Cka-global}
Let $k\geq 3$, $0<\alpha<\bar{\alpha}$, $\Gamma\subset \partial \Omega $ be relatively open and $u$ be a viscosity solution of
\begin{equation*}
\left\{\begin{aligned}
&F(D^2u,Du,u,x)=f&& \quad\mbox{in}~~\Omega;\\
&u=g&& \quad\mbox{on}~~\Gamma.
\end{aligned}\right.
\end{equation*}
Suppose that \cref{SC2} and \cref{e.omega0-2} hold, and
\begin{equation*}
F\in C^{k-2,\alpha}(\bar \Omega),\quad f\in C^{k-2,\alpha}(\bar\Omega),\quad \Gamma\in C^{k,\alpha},\quad g\in C^{k,\alpha}(\bar\Gamma).
\end{equation*}

Then for any $\Omega'\subset\subset \Omega\cup \Gamma$, we have $u\in C^{k,\alpha}(\bar{\Omega}')$ and
\begin{equation}\label{e.Cka-1-i-mu-global}
 \|u\|_{C^{k,\alpha}(\bar{\Omega}')}\leq C,
\end{equation}
where $C$ depends only on $k,n,\lambda, \Lambda,\alpha,r_0,\mu,b_0,c_0,\|F\|_{C^{k-2,\alpha}(\bar\Omega)},K_1,\omega_3,
\omega_4$, $\|\partial \Omega'\cap\Gamma\|_{C^{k,\alpha}}$, $\|f\|_{C^{k-2,\alpha}(\bar\Omega)}$,
$\|g\|_{C^{k,\alpha}(\bar\Gamma)}$ and $\|u\|_{L^{\infty }(\Omega)}$.
\end{corollary}

\begin{remark}\label{re8.1}
For local and global $C^{k,\alpha}$ regularity, the condition ``$F\in C^{k-2,\alpha}(\bar \Omega)$'' (see \Cref{d-FP}) can be replaced by~\\
(i) $F$ is convex in $M$;~\\
(ii) $F \in C^{k,\bar{\alpha}}(\mathcal{S}^n\times \mathbb{R}^n\times \mathbb{R}\times \bar{\Omega})$;~\\
(iii) For some $0<\varepsilon<1$,
\begin{equation*}
|F(M,p,s,x)-F(M,p,s,y)| \leq K_1|x-y|^{\varepsilon}(|M|+1)\omega(|p|,|s|).
\end{equation*}

Indeed, by \Cref{t-C2a-i} and \Cref{t-C2a}, $u\in C^{2,\alpha}(\bar{\Omega}')$. Then we can redefine $F$ outside a compact set of $\mathcal{S}^n\times \mathbb{R}^n\times \mathbb{R}\times \bar{\Omega}$ such that $F\in C^{k-2,\alpha}(\bar \Omega)$ (cf. \Cref{re1.4}).
\end{remark}
~\\

Similar to \Crefrange{l-3modin1}{l-32}, we first prove the boundary $C^{k,\alpha}$ regularity for a model problem (see \Cref{l-62}). Then the general boundary pointwise $C^{k,\alpha}$ regularity follows as a perturbation of the regularity for the model problem.

The following lemma states the local $C^{2,\alpha}$ regularity up to the boundary which is obtained by combining the interior pointwise regularity and the boundary pointwise regularity. The proof is standard and we borrow it from \cite{MR3246039}.
\begin{lemma}\label{l-61}
Let $0<\alpha<\bar{\alpha}$. Suppose that $F$ satisfies \cref{e.C2a-KF} with some $G_{x_0}$ at any $x_0\in \bar{B}_1^+$, where $G_{x_0}$ is convex in $M$. Assume that $\beta_2\in C^{\alpha}(\bar{B}_1^+)$, $\omega_0$ satisfies \cref{e.omega0} and $u$ satisfies
\begin{equation*}
\left\{\begin{aligned}
&F(D^2u,Du,u,x)=0&& \quad\mbox{in}~~B_1^+;\\
&u=0&& \quad\mbox{on}~~T_1.
\end{aligned}\right.
\end{equation*}
Then $u\in C^{2,\alpha}(\bar{B}^+_{1/2})$ and
\begin{equation*}
 \|u\|_{C^{2,\alpha}(\bar{B}^+_{1/2})}\leq C,
\end{equation*}
where $C$ depends only on $n,\lambda, \Lambda,\alpha,\mu,b_0,c_0,\omega_0,
\|\beta_2\|_{C^{\alpha}(\bar{B}_1^+)},\omega_2$ and $\|u\|_{L^{\infty }(B_1^+)}$.
\end{lemma}
\proof We only need to prove that given $x_0=(x'_0,r_0)\in \bar{B}^+_{1/2}$, there exists $P_{x_0}\in \mathcal{P}_2$ such that
\begin{equation}\label{e.l61-3}
  |u(x)-P_{x_0}(x)|\leq C|x-x_0|^{2+\alpha},~\forall ~x\in \bar{B}^+_{1/2},
\end{equation}
where $C$ depends only on $n,\lambda, \Lambda,\alpha,\mu,b_0,c_0,\omega_0,\omega_2,
\|\beta_2\|_{C^{\alpha}(\bar{B}_1^+)}$ and $\|u\|_{L^{\infty }(B_1^+)}$. Throughout this proof, $C$ always denotes a constant having the same dependence.

Let $\tilde{x}_0=(x'_0,0)$. By the boundary $C^{2,\alpha}$ regularity \Cref{t-C2a}, there exists $P_{\tilde{x}_0}\in \mathcal{P}_2$ such that
\begin{equation}\label{e9.4}
  |u(x)-P_{\tilde{x}_0}(x)|\leq C|x-\tilde{x}_0|^{2+\alpha}, ~\forall ~x\in \bar{B}_{1}^{+}.
\end{equation}
Set
\begin{equation*}
v=u-P_{\tilde{x}_0}
\end{equation*}
and then $v$ satisfies
\begin{equation}\label{e.Cka-4}
  \tilde{F}(D^2v,Dv,v,x)=0~~~~\mathrm{in}~B_{r_0}(x_0),
\end{equation}
where
\begin{equation*}
\tilde{F}(M,p,s,x)\coloneqq F(M+D^2P_{\tilde{x}_0},p+DP_{\tilde{x}_0}(x),s+P_{\tilde{x}_0}(x),x).
\end{equation*}
Then \cref{e.Cka-4} satisfies the conditions of the interior $C^{2,\alpha}$ regularity (see \Cref{t-C2a-i}). Hence, there exists $P\in \mathcal{P}_2$ such that (with the aid of \cref{e9.4} and the definition of $v$)
\begin{equation}\label{e9.5}
  \begin{aligned}
|v(x)-P(x)|\leq& C\|v\|_{L^\infty(B_{r_0}(x_0))}r_0^{-(2+\alpha)}|x-x_0|^{2+\alpha}
  \leq C |x-x_0|^{2+\alpha}~~\mathrm{in}~B_{r_0/2}(x_0),\\
|P(x_0)|=&|v(x_0)|=|u(x_0)-P_{\tilde{x}_0}(x_0)|\leq C|x_0-\tilde{x}_0|=Cr_0^{2+\alpha},\\
|D P(x_0)|\leq& Cr_0^{-1}\|v\|_{L^\infty(B_{r_0}(x_0))}\leq Cr_0^{1+\alpha},\\
|D^2P(x_0)|\leq& Cr_0^{-2}\|v\|_{L^\infty(B_{r_0}(x))}\leq Cr_0^{\alpha}.
  \end{aligned}
\end{equation}

Let
\begin{equation*}
P_{x_0}=P_{\tilde{x}_0}+P.
\end{equation*}
If $|x-x_0|<r_0/2$, by the first inequality in \cref{e9.5},
\begin{equation*}
  |u(x)-P_{x_0}(x)|=|v(x)-P(x)|\leq C|x-x_0|^{2+\alpha}.
\end{equation*}
If $|x-x_0|\geq r_0/2$, we have
\begin{equation*}
  \begin{aligned}
|u(x)-P_{x_0}(x)|&\leq |u(x)-P_{\tilde{x}_0}(x)|+|P(x)|\\
&\leq C|x-\tilde{x}_0|^{2+\alpha}
+\left|\sum_{|\sigma|\leq 2}\frac{D^{\sigma}P(x_0)}{\sigma!}(x-x_0)^{\sigma}\right|\\
&\leq C|x-\tilde{x}_0|^{2+\alpha}+C(r_0^{2+\alpha}+|x-x_0|r_0^{1+\alpha}+|x-x_0|^2r_0^{\alpha})\\
&\leq C|x-x_0|^{2+\alpha}.
  \end{aligned}
\end{equation*}
Therefore, \cref{e.l61-3} holds.~\qed \\

Now, we prove the boundary regularity for the model problem.

\begin{lemma}\label{l-62}
Suppose that $G$ satisfies (i)-(iii) in \Cref{d-FP} (with $G \in C^{k-2,\bar{\alpha}}$ instead of $G \in C^{k,\bar{\alpha}}$) and $\omega_0$ satisfies \cref{e.omega0-2}. Let $u$ be a viscosity solution of
\begin{equation*}
\left\{\begin{aligned}
&G(D^2u,Du,u,x)=0&& \quad\mbox{in}~~B_1^+;\\
&u=0&& \quad\mbox{on}~~T_1.
\end{aligned}\right.
\end{equation*}
Then $u\in C^{k,\alpha}(\bar{B}^+_{1/2})$ for any $0<\alpha<\bar{\alpha}$ and
\begin{equation*}
\|u\|_{C^{k,\alpha}(\bar{B}^+_{1/2})}\leq C_k,
\end{equation*}
where $C_k$ depends only on $k,n,\lambda, \Lambda,\alpha,\mu,b_0,c_0,\omega_0,K_1,\omega_4$ and $\|u\|_{L^{\infty }(B_1^+)}$.

In particular, $u\in C^{k,\alpha}(0)$ and there exists $P\in \mathcal{P}_k$ such that
\begin{equation}\label{e.l62-1}
  |u(x)-P(x)|\leq C_k |x|^{k+\alpha}, ~~\forall ~x\in B_{1}^+,
\end{equation}
\begin{equation}\label{e.l62-2}
  |G(D^2P(x),DP(x),P(x),x)|\leq C_k|x|^{k-2+\bar{\alpha}}, ~~\forall ~x\in B_{1}^+
\end{equation}
and
\begin{equation}\label{e.l62-3}
  \|P\|\leq C_k.
\end{equation}
Moreover, $P$ can be written as
\begin{equation*}
P=\sum_{l=1}^{k} P_l,~P_l\in\mathcal{SP}_l.
\end{equation*}
\end{lemma}

\proof The proof is similar to \Cref{In-l-62} and we only give the outline. Since $G$ is smooth and satisfies (iii) in \Cref{d-FP}, as in the proof of \Cref{In-l-62}, $G$ satisfies \cref{e.C2a-KF}. By \Cref{l-61}, $u\in C^{2,\alpha}(\bar{B}_{3/4}^+)$ for any $0<\alpha<\bar{\alpha}$ and
\begin{equation*}
\|u\|_{C^{2,\alpha}(\bar{B}_{3/4}^+)}\leq C.
\end{equation*}

Once we have $u\in C^{2,\alpha}(\bar{B}_{3/4}^+)$, the higher regularity can be derived by the standard technique of difference quotient. The only difference from the interior case is that we cannot take the difference quotient along $e_n$. Let $h>0$ be small and $1\leq l\leq n-1$. By taking the difference quotient along $e_l$ on both sides of the equation and applying the Schauder estimates for linear equations, we obtain that $u_l \in C^{2,\alpha^2}(\bar B_{1/2}^+)$ and
\begin{equation*}
\|u_l\|_{C^{2,\alpha^2}(\bar{B}_{1/2}^+)}\leq C.
\end{equation*}
It follows that $u_{il}\in C^{1,\alpha^2}$ ($i+l<2n$). By combining with $G(D^2u,Du,u,x)=0$ and the implicit function theorem, $u_{nn}\in C^{1,\alpha^2}$. Thus, $u\in C^{3,\alpha^2}$ and
\begin{equation*}
\|u\|_{C^{3,\alpha^2}(\bar B_{1/2}^+)}\leq C.
\end{equation*}
Then by a similar proof to that of \Cref{In-l-62}, we have
\begin{equation*}
\|u\|_{C^{k,\alpha}(\bar{B}^+_{1/2})}\leq C.
\end{equation*}

In particular, $u\in C^{k,\alpha}(0)$ and there exists $P\in \mathcal{P}_k$ such that \cref{e.l62-1} and \cref{e.l62-3} hold. Note that
\begin{equation*}
D^{\sigma}u(0)=0~~\mbox{ if }~\sigma_n=0.
\end{equation*}
Hence, $P$ can be written as
\begin{equation*}
P=\sum_{l=1}^{k} P_l,~P_l\in \mathcal{SP}_l.
\end{equation*}
Finally, similar to \Cref{In-l-62}, \cref{e.l62-2} holds. \qed~\\

In the following, we prove the boundary pointwise $C^{k,\alpha}(k\geq 3)$ regularity \Cref{t-Cka} by induction. For $k=3$, $F\in C^{1,\alpha}(0)$. With the aid of \cref{holder}, $F$ satisfies \cref{e.C2a-KF} (see \cref{e5.3}). Then from the boundary pointwise
$C^{2,\alpha}$ regularity \Cref{t-C2a}, $u\in C^{2,\alpha}(0)$. Hence, we can assume that the boundary $C^{k-1,\alpha}(0)$ regularity holds if $F\in C^{k-2,\alpha}(0)$ and we need to prove the boundary $C^{k,\alpha}(0)$ regularity.

The following lemma is a higher order counterpart of \Cref{l-C1a-mu} and \Cref{l-C2a-mu}.
\begin{lemma}\label{l-Cka-mu}
Let $0<\alpha<\bar{\alpha}$, $F\in C^{k-2,\alpha}(0)$ and $\omega_0$ satisfy \cref{e.omega0-2}. Then there exists $\delta>0$ depending only on $k,n,\lambda,\Lambda,\alpha,\omega_0,K_1,\omega_3$ and $\omega_4$ such that if $u$ satisfies
\begin{equation*}
\left\{\begin{aligned}
&F(D^2u,Du,u,x)=f&& \quad\mbox{in}~~\Omega_1;\\
&u=g&& \quad\mbox{on}~~(\partial \Omega)_1
\end{aligned}\right.
\end{equation*}
with
\begin{equation*}
  \begin{aligned}
&u(0)=|Du(0)|=\cdots=|D^{k-1}u(0)|=0,\quad \max\left(\|u\|_{L^{\infty}(\Omega_1)},\mu,b_0,c_0\right)\leq 1,\\
&\max\left(\|F\|_{C^{k-2,\alpha}(0)},\|f\|_{C^{k-3,\alpha}(0)}, \|g\|_{C^{k-1,\alpha}(0)},
\|(\partial \Omega)_1\|_{C^{1,\alpha}(0)}\right)\leq \delta,
  \end{aligned}
\end{equation*}
then there exists $P\in\mathcal{SP}_k$ such that
\begin{equation*}
  \|u-P\|_{L^{\infty}(\Omega_{\eta})}\leq \eta^{k+\alpha},
\end{equation*}
\begin{equation*}
  |G(D^2P(x),DP(x),P(x),x)|\leq C |x|^{k-2+\bar{\alpha}},~\forall ~x\in \Omega_1
\end{equation*}
and
\begin{equation*}\label{e.16.2-mu}
\|P\|\leq C,
\end{equation*}
where $C$ and $\eta$ depend only on $k,n,\lambda, \Lambda,\alpha,\omega_0,K_1,\omega_3$ and $\omega_4$.
\end{lemma}

\begin{remark}\label{re9.1}
The proof is similar to that of \Cref{In-l-Cka-mu} and we omit it. There is only one thing we need to be care of. In the procedure of constructing a sequence of polynomials (see \cref{e6.8}), we need $P_m\in \mathcal{SP}_k$ other than $P_m\in \mathcal{HP}_k$. We will use \Cref{le2.6} instead of \Cref{le2.2} used in the proof of \Cref{In-l-Cka-mu}.
\end{remark}

\begin{remark}\label{r-9.2}
Note that we only assume that $\partial \Omega\in C^{1,\alpha}(0)$. Clearly, the $C^{1,\alpha}(0)$ regularity holds. Since $u(0)=0$ and $Du(0)=0$, from the boundary pointwise $C^{2,\alpha}$ regularity (see last section), $u\in C^{2,\alpha}(0)$. By induction, we have $u\in C^{k-1,\alpha}(0)$.
\end{remark}
~\\

The following is a boundary pointwise $C^{k,\alpha}$ regularity. As before, the key observation is that if $u(0)=|Du(0)|\cdots=D^{k-1}u(0)=0$, the boundary pointwise $C^{k,\alpha}$ regularity holds even if $\partial \Omega\in C^{1,\alpha}(0)$.
\begin{lemma}\label{t-Ckas-mu}
Let $0<\alpha <\bar{\alpha}$, $F\in C^{k-2,\alpha}(0)$ and $\omega_0$ satisfy \cref{e.omega0-2}. Suppose that $u$ is a viscosity solution of
\begin{equation*}
\left\{\begin{aligned}
&F(D^2u,Du,u,x)=f&& \quad\mbox{in}~~\Omega_1;\\
&u=g&& \quad\mbox{on}~~(\partial \Omega)_1.
\end{aligned}\right.
\end{equation*}
Assume that
\begin{equation}\label{e.tCkas-be-mu}
\begin{aligned}
&\|u\|_{L^{\infty}(\Omega_1)}\leq 1,\quad u(0)=|Du(0)|=\cdots|D^{k-1}u(0)|=0,\\
&\mu\leq \frac{1}{4C_0},\quad b_0\leq \frac{1}{2},\quad c_0\leq 1, \quad
\|F\|_{C^{k-2+\alpha}(0)}\leq \frac{\delta_1}{K_0},\\
&|f(x)|\leq \delta_1|x|^{k-2+\alpha}, ~\forall ~x\in \Omega_1,\\
&|g(x)|\leq \frac{\delta_1}{2}|x|^{k+\alpha}, ~\forall ~x\in (\partial \Omega)_1
\quad\mbox{and}\quad\|(\partial \Omega)_1\|_{C^{1,\alpha}(0)} \leq \frac{\delta_1}{2C_0},\\
\end{aligned}
\end{equation}
where $\delta_1\leq \delta$ ($\delta$ is as in \Cref{l-Cka-mu}) and $C_0$ depend only on $k,n,\lambda, \Lambda,\alpha,\omega_0,K_1,\omega_3$ and $\omega_4$.

Then $u\in C^{k,\alpha}(0)$ and there exists $P\in\mathcal{SP}_k$ such that
\begin{equation}\label{e.tCkas-1-mu}
  |u(x)-P(x)|\leq C |x|^{k+\alpha}, ~~\forall ~x\in \Omega_{1},
\end{equation}
\begin{equation*}
  |G(D^2P(x),DP(x),P(x),x)|\leq C |x|^{k-2+\bar{\alpha}}, ~~\forall ~x\in \Omega_{1}
\end{equation*}
and
\begin{equation}\label{e.tCkas-2-mu}
\|P\|\leq C,
\end{equation}
where $C$ depends only on $k,n,\lambda, \Lambda,\alpha,\omega_0,K_1,\omega_3$ and $\omega_4$.
\end{lemma}

\proof As before, to prove that $u$ is $C^{k,\alpha}$ at $0$, we only need to prove the following. There exist a sequence of $P_m\in\mathcal{SP}_k$ ($m\geq 0$) such that for all $m\geq 1$,

\begin{equation}\label{e.tCkas-6-mu}
\|u-P_m\|_{L^{\infty }(\Omega _{\eta^{m}})}\leq \eta ^{m(k+\alpha )},
\end{equation}
\begin{equation}\label{e.tCkas-9-mu}
|G(D^2P_m(x),DP_m(x),P_m(x),x)|\leq \tilde{C}|x|^{k-2+\bar{\alpha}}, ~~\forall ~x\in \Omega_{1}
\end{equation}
and
\begin{equation}\label{e.tCkas-7-mu}
\|P_m-P_{m-1}\|\leq \tilde{C}\eta ^{(m-1)\alpha},
\end{equation}
where $\tilde{C}$ and $0<\eta<1$ depend only on $k,n,\lambda, \Lambda,\alpha,\omega_0,K_1,\omega_3$ and $\omega_4$.

We prove the above by induction. For $m=1$, by \Cref{l-Cka-mu}, there exists $P_1\in\mathcal{SP}_k$ such that \crefrange{e.tCkas-6-mu}{e.tCkas-7-mu} hold for some $C_1$ and $\eta_1$ depending only on $k,n,\lambda, \Lambda,\alpha,\omega_0,K_1,\omega_3$ and $\omega_4$ where $P_0\equiv 0$. Take $\tilde{C}\geq C_1, \eta\leq \eta_1$ and then the conclusion holds for $m=1$. Suppose that the conclusion holds for $m\leq m_0$. We need to prove that the conclusion holds for $m=m_0+1$.

Let $r=\eta ^{m_{0}}$, $y=x/r$ and
\begin{equation}\label{e.tCkas-v-mu}
  v(y)=\frac{u(x)-P_{m_0}(x)}{r^{k+\alpha}}.
\end{equation}
Then $v$ satisfies
\begin{equation}\label{e.Ckas-F-mu}
\left\{\begin{aligned}
&\tilde{F}(D^2v,Dv,v,y)=\tilde{f}&& \quad\mbox{in}~~\tilde{\Omega}_1;\\
&v=\tilde{g}&& \quad\mbox{on}~~(\partial \tilde{\Omega})_1,
\end{aligned}\right.
\end{equation}
where for $(M,p,s,y)\in \mathcal{S}^n\times \mathbb{R}^n\times \mathbb{R}\times \bar{\tilde\Omega}_1$,
\begin{equation*}
  \begin{aligned}
&\tilde{F}(M,p,s,y)=r^{-l}F(r^{l}M+D^2P_{m_0}(x),r^{l+1}p+DP_{m_0}(x),r^{l+2}s+P_{m_0}(x),x),\\
&\tilde{f}(y)=r^{-l}f(x), \quad\tilde{g}(y)=r^{-(l+2)}\left(g(x)-P_{m_0}(x)\right),
\quad\tilde{\Omega}=r^{-1}\Omega, \quad l=k-2+\alpha.\\
  \end{aligned}
\end{equation*}
In addition, define $\tilde{G}$ in a similar way to the definition of $\tilde{F}$.

In the following, we show that \cref{e.Ckas-F-mu} satisfies the assumptions of \Cref{l-Cka-mu}. First, by \cref{e.tCkas-7-mu}, there exists a constant $C_{0}$ depending only on $k,n,\lambda, \Lambda,\alpha,\omega_0,K_1,\omega_3$ and $\omega_4$ such that $\|P_m\|\leq C_{0}$ ($\forall~0\leq m\leq m_0$). Then it is easy to verify that (note that $P_{m_0}\in\mathcal{SP}_k$)
\begin{equation*}
\begin{aligned}
\|v\|_{L^{\infty}(\tilde{\Omega}_1)}\leq& 1, \quad  v(0)=\cdots=|D^{k-1}v(0)|=0,~\quad(\mathrm{by}~ \cref{e.tCkas-be-mu},~ \cref{e.tCkas-6-mu} ~\mbox{and}~\cref{e.tCkas-v-mu})\\
|\tilde{f}(y)|\leq& r^{-(k-2+\alpha)}|f(x)|\leq \delta_1|y|^{k-2+\alpha}, ~\forall ~y\in \tilde{\Omega}_1,~\quad(\mathrm{by}~ \cref{e.tCkas-be-mu})\\
|\tilde{g}(y)|\leq& \frac{1}{r^{k+\alpha}}\left(\frac{\delta_1(r|y|)^{k+\alpha}}{2}+
\frac{C_{0}\delta_1(r|y|)^{k+\alpha}}{2C_0}\right)\\
\leq& \delta_1|y|^{k+\alpha},~\forall ~y\in (\partial\tilde{\Omega})_1,~\quad(\mathrm{by}~ \cref{e.tCkas-be-mu})\\
\|(\partial \tilde{\Omega})_1\|_{C^{1,\alpha}(0)}\leq& r^{\alpha}\|(\partial\Omega)_1\|_{C^{1,\alpha}(0)}\leq \delta_1.~\quad(\mathrm{by}~ \cref{e.tCkas-be-mu})
  \end{aligned}
\end{equation*}

Next, it is easy to check that $\tilde{F}$ and $\tilde{G}$ satisfy the structure condition \cref{SC2} with $\lambda,\Lambda,\tilde{\mu},\tilde{b},\tilde{c}$ and $\tilde{\omega}_0$, where
\begin{equation*}
\begin{aligned}
&\tilde{\mu}\leq \mu\leq 1,\quad\tilde{b}\leq b_0+2C_0\mu\leq 1,\quad\tilde{c}\leq c_0\leq 1,\quad
\tilde{\omega}_0(\cdot,\cdot)=\omega_0(\cdot+C_0,\cdot).
\end{aligned}
\end{equation*}
In addition, by an argument similar to \cref{e.6.1},
\begin{equation}\label{e.9.1}
|\tilde{F}(M,p,s,y)-\tilde{G}(M,p,s,y)|\leq \delta_1|y|^{k-2+\alpha}(|M|+1)\tilde{\omega}_3(|p|,|s|),
\end{equation}
where
\begin{equation*}
 \tilde{\omega}_3(|p|,|s|)\coloneqq\omega_3(|p|+C_0,|s|+C_0).
\end{equation*}
Hence, $\|\tilde{F}\|_{C^{k-2,\alpha}(0)}\leq \delta_1$.

Finally, with the aid of \cref{e.tCkas-9-mu}, we can show that (similar to the proof of \Cref{In-t-Ckas-mu}) $\tilde{G}$ satisfies (i)-(iii) of \Cref{d-FP} with some $\tilde{K}_1$ and
\begin{equation*}
\|\tilde{G}\|_{C^{k-2,\bar{\alpha}}(\bar{\mathbf{B}}_\rho\times \bar{\tilde\Omega}_1)}\leq \tilde{\omega}_4(\rho),~\forall ~\rho>0,
\end{equation*}
where $\tilde{\omega}_4$ depends only on $k,n,\lambda, \Lambda,\alpha,\omega_0,K_1,\omega_3$ and $\omega_4$.

Choose $\delta_1$ small enough (depending only on $k,n,\lambda, \Lambda,\alpha,\omega_0,K_1,\omega_3$ and $\omega_4$) such that \Cref{l-Cka-mu} holds for $\tilde\omega_0,\tilde{K}_1,\tilde\omega_3,\tilde\omega_4$ and $\delta_1$. Since \cref{e.Ckas-F-mu} satisfies the assumptions of \Cref{l-Cka-mu}, there exist $\tilde{P}\in\mathcal{SP}_k$, constants $\tilde{C}\geq C_1$ and $\eta\leq \eta_1$ such that
\begin{equation*}
\begin{aligned}
    \|v-\tilde{P}\|_{L^{\infty }(\tilde{\Omega} _{\eta})}&\leq \eta^{k+\alpha},
\end{aligned}
\end{equation*}
\begin{equation*}
  |\tilde G(D^2\tilde{P}(y),D\tilde{P}(y),\tilde{P}(y),y)|\leq \tilde{C}|y|^{k-2+\bar{\alpha}}
  , ~~\forall ~y\in \tilde{\Omega}_{1}
\end{equation*}
and
\begin{equation*}
\|\tilde{P}\|\leq \tilde{C}.
\end{equation*}

Let
\begin{equation*}
P_{m_0+1}(x)=P_{m_0}(x)+r^{k+\alpha}\tilde{P}(y)=P_{m_0}(x)+r^{\alpha}\tilde{P}(x).
\end{equation*}
Then \cref{e.tCkas-9-mu} and \cref{e.tCkas-7-mu} hold for $m_0+1$. By recalling \cref{e.tCkas-v-mu}, we have
\begin{equation*}
  \begin{aligned}
\|u-P_{m_0+1}\|_{L^{\infty}(\Omega_{\eta^{m_0+1}})}&= \|u-P_{m_0}-r^{\alpha}\tilde{P}(x)\|_{L^{\infty}(\Omega_{\eta r})}\\
&= \|r^{k+\alpha}v-r^{k+\alpha}\tilde{P}\|_{L^{\infty}(\tilde{\Omega}_{\eta})}\\
&\leq r^{k+\alpha}\eta^{k+\alpha}=\eta^{(m_0+1)(k+\alpha)}.
  \end{aligned}
\end{equation*}
Hence, \cref{e.tCkas-6-mu} holds for $m=m_0+1$. By induction, the proof is completed.\qed~\\

\begin{remark}\label{r-9.1}
Roughly speaking, in above proof, the condition $F\in C^{k-2,\alpha}(0)$ is only used in \cref{e.9.1} (see also \cref{e.6.1}). Suppose that $\omega_3$ satisfies for some constant $K_0$ (e.g., $\omega_3(|p|,|s|)=|p|+|s|$)
\begin{equation*}
r\omega_3\left(r^{-1}|p|,r^{-1}|s|\right)\leq K_0\omega_3(|p|,|s|),~\forall ~r>0.
\end{equation*}
Then we can obtain \cref{e.9.1} by only assuming $F\in C^{\alpha}$.
\end{remark}
~\\

Now, we give the~\\
\noindent\textbf{Proof of \Cref{t-Cka}.} As before, in the following proof,  we just make necessary normalization to satisfy the conditions of \Cref{t-Ckas-mu}. Throughout this proof, $C$ always denotes a constant depending only on $k,n,\lambda, \Lambda,\alpha,\mu, b_0,c_0,\omega_0$, $\|F\|_{C^{k-2,\alpha}(0)},$ $K_1,\omega_3,\omega_4,\|(\partial \Omega)_1\|_{C^{k,\alpha}(0)}$, $\|f\|_{C^{k-2,\alpha}(0)}$, $\|g\|_{C^{k,\alpha}(0)}$ and $\|u\|_{L^{\infty}(\Omega_1)}$.

For $(M,p,s,x)\in \mathcal{S}^n\times \mathbb{R}^n\times \mathbb{R}\times \bar{\Omega}_1$, let
\begin{equation*}
F_1(M,p,s,x)=F(M,p,s,x)-P_f(x), \quad f_1=f-P_f.
\end{equation*}
Then $u$ satisfies
\begin{equation*}
\left\{\begin{aligned}
&F_1(D^2u,Du,u,x)=f_1&& \quad\mbox{in}~~\Omega_1;\\
&u=g&& \quad\mbox{on}~~(\partial \Omega)_1
\end{aligned}\right.
\end{equation*}
and
\begin{equation*}
  |f_1(x)|\leq [f]_{C^{k-2,\alpha}(0)}|x|^{k-2+\alpha}\leq C|x|^{k-2+\alpha}, ~~\forall ~x\in \Omega_1.
\end{equation*}

Next, set
\begin{equation*}
u_1=u-P_g,\quad F_2(M,p,s,x)=F_1(M+D^2P_g(x),p+DP_g(x),s+P_g(x),x).
\end{equation*}
Then $u_1$ is a viscosity solution of
\begin{equation*}
\left\{\begin{aligned}
&F_2(D^2u_1,Du_1,u_1,x)=f_1&& \quad\mbox{in}~~\Omega_1;\\
&u_1=g_1&& \quad\mbox{on}~~(\partial \Omega)_1,
\end{aligned}\right.
\end{equation*}
where $g_1=g-P_g$. Hence,
\begin{equation}\label{e.cka-1}
  |g_1(x)|\leq [g]_{C^{k,\alpha}(0)}|x|^{k+\alpha}\leq C|x|^{k+\alpha}, ~~\forall ~x\in (\partial \Omega)_1.
\end{equation}

Note that $u_1\in C^{k-1,\alpha}(0)$ and define
\begin{equation*}
P_u(x)=\left(\sum_{|\sigma|\leq k-1,\sigma_n\geq 1} \frac{1}{\sigma !} D^{\sigma} u_1(0) x^{\sigma-e_n}\right)\left(x_n-P_{\Omega}(x')\right).
\end{equation*}
By combining \cref{e.Cka-5} with \cref{e.cka-1}, we have
\begin{equation}\label{e8.1}
D^l P_u(0)=D^l u_1(0),~\forall ~0\leq l\leq k-1.
\end{equation}
Let
\begin{equation*}
u_2=u_1-P_u,\quad F_3(M,p,s,x)=F_2(M+D^2P_u(x),p+DP_u(x),s+P_u(x),x).
\end{equation*}
Then $u_2$ satisfies
\begin{equation*}
\left\{\begin{aligned}
&F_3(D^2u_2,Du_2,u_2,x)=f_1&& \quad\mbox{in}~~\Omega_1;\\
&u_2=g_2&& \quad\mbox{on}~~(\partial \Omega)_1
\end{aligned}\right.
\end{equation*}
and
\begin{equation*}
  u_2(0)=|Du_2(0)|=\cdots=|D^{k-1}u_2(0)|=0,
\end{equation*}
where $g_2=g_1-P_u$. By the boundary $C^{k-1,\alpha}$ regularity,
\begin{equation}\label{e.tCka-1-mu}
  \begin{aligned}
|Du_1(0)|+\cdots+|D^{k-1}u_1(0)| \leq C.
  \end{aligned}
\end{equation}
Hence,
\begin{equation}\label{e9.3}
  \begin{aligned}
    |g_2(x)| \leq |g_1(x)|+|P_u(x)|\leq C|x|^{k+\alpha}, ~~\forall ~x\in (\partial \Omega)_1.
  \end{aligned}
\end{equation}

Finally, take $y=x/\rho$ and $u_3(y)=u_2(x)/\rho^2$, where $0<\rho<1$ is a constant to be specified later. Then $u_3$ satisfies
\begin{equation}\label{F6-k-mu}
\left\{\begin{aligned}
&F_4(D^2u_3,Du_3,u_3,y)=f_2&& \quad\mbox{in}~~\tilde{\Omega}_1;\\
&u_3=g_3&& \quad\mbox{on}~~(\partial \tilde{\Omega})_1,
\end{aligned}\right.
\end{equation}
where
\begin{equation*}
  \begin{aligned}
&F_4(M,p,s,y)=F_3\left(M,\rho\, p,\rho^2s,\rho y\right),~~
f_2(y)=f_1(x),~~ g_3(y)=\rho^{-2}g_2(x),~~\tilde{\Omega}=\rho^{-1}\Omega.
\end{aligned}
\end{equation*}
Define fully nonlinear operators $G_1,G_2,G_3,G_{4}$ in the same way as $F_1,F_2,F_3,F_4$.

Now, we choose a proper $\rho$ such that \cref{F6-k-mu} satisfies the conditions of \Cref{t-Ckas-mu}. Let $N(x)=N_g(x)+N_u(x)$. Obviously,
\begin{equation*}
u_3(0)=\cdots=|D^{k-1}u_3(0)|=0.
\end{equation*}
By combining with the $C^{k-1,\alpha}$ regularity at $0$,  we have
\begin{equation*}
  \begin{aligned}
     \|u_3\|_{L^{\infty}(\tilde{\Omega}_1)}&= \rho^{-2}\|u_2\|_{L^{\infty}(\Omega_\rho)}
    \leq C\rho^{k-3+\alpha}.
  \end{aligned}
\end{equation*}
Next, we can deduce easily
\begin{equation*}
\begin{aligned}
|f_2(y)|=&|f_1(x)|\leq C|x|^{k-2+\alpha}=C\rho^{k-2+\alpha}|y|^{k-2+\alpha}\quad\mbox{in}~\tilde{\Omega}_1,\\
|g_3(y)|=&\rho^{-2}|g_2(x)|\leq C \rho^{-2}|x|^{k+\alpha}\leq C\rho^{k-2+\alpha}|y|^{k+\alpha}
\quad\mbox{on}~\partial \tilde{\Omega}_1,\\
\|(\partial \tilde{\Omega})_1\|_{C^{1,\alpha}(0)}\leq& \rho\|(\partial \Omega)_1\|_{C^{1,\alpha}(0)} \leq \rho\|(\partial \Omega)_1\|_{C^{k,\alpha}(0)}\leq C\rho.
\end{aligned}
\end{equation*}
Furthermore, $F_4$ and $G_{4}$ satisfy the structure condition \cref{SC2} with $\lambda,\Lambda,\tilde\mu,\tilde{b},\tilde{c}$ and $\tilde{\omega}_0$, where
\begin{equation*}
\tilde{\mu}= \rho^2\mu,\quad\tilde{b}= \rho b_0+2C\rho\mu,\quad\tilde{c}= \rho^2c_0, \quad
\tilde{\omega}_0(\cdot,\cdot)=\omega_0(\cdot+C,\cdot).
\end{equation*}

Finally, we can show that $F_4\in C^{k-2,\alpha}(0)$. First, similar to the previous argument (e.g. \cref{e5.4}),
\begin{equation*}
|F_4(M,p,s,y)-G_{4}(M,p,s,y)|\leq C\rho^{k-2+\alpha}|y|^{k-2+\alpha}(|M|+1)\tilde\omega_3(|p|,|s|),
\end{equation*}
where
\begin{equation*}
\tilde\omega_3(|p|,|s|)\coloneqq \omega_3(|p|+C,|s|+C).
\end{equation*}
In addition, it can be verified that $G_{4}$ satisfies (i)-(iii) of \Cref{d-FP} with some $\tilde{K}_1$ and
\begin{equation*}
\|G_{4}\|_{C^{k,\alpha}(\bar{\mathbf{B}}_r\times \bar{\tilde\Omega}_1)}\leq \tilde{\omega}_4(r),~\forall ~r>0,
\end{equation*}
\sloppy where $\tilde{\omega}_4$ depends only on $k$, $n$, $\lambda$, $\Lambda$, $\alpha$, $\mu, b_0,c_0,\omega_0$, $\|F\|_{C^{k-2,\alpha}(0)}$, $K_1$, $\omega_3$, $\omega_4$, $\|(\partial \Omega)_1\|_{C^{k,\alpha}(0)}$, $\|f\|_{C^{k-2,\alpha}(0)}, \|g\|_{C^{k,\alpha}(0)}$ and $\|u\|_{L^{\infty}(\Omega_1)}$.

Take $\delta_1$ small enough such that \Cref{t-Ckas-mu} holds for $\tilde{\omega}_0$, $\tilde{K}_1$, $\tilde{\omega}_3$, $\tilde{\omega}_4$ and $\delta_1$. From above arguments, we can choose $\rho$ small enough (depending only on $k,n,\lambda, \Lambda,\alpha,\mu, b_0$, $c_0,\omega_0$, $\|F\|_{C^{k-2,\alpha}(0)}$, $K_1,\omega_3,\omega_4,\|(\partial \Omega)_1\|_{C^{k,\alpha}(0)}, \|f\|_{C^{k-2,\alpha}(0)}, \|g\|_{C^{k,\alpha}(0)}$ and $\|u\|_{L^{\infty}(\Omega_1)}$) such that the assumptions of \Cref{t-Ckas-mu} are satisfied. Then $u_3$ and hence $u$ is $C^{k,\alpha}$ at $0$, and the estimates \crefrange{e.Cka-1}{e.Cka-2} hold. \qed~\\

\section{The \texorpdfstring{$C^{k}$}{Ck} regularity}\label{C1C2Ck}

Besides the classical $C^{k,\alpha}$ regularity ($0<\alpha<1$), we can also obtain other types of pointwise regularity. It is well-known that even for the Poisson equation, the Schauder estimate fails for $\alpha=0$ or $1$ (see \cite[Problem 4.9]{MR1814364}). On the other hand, if $f$ is Dini continuous, the $C^2$ regularity holds; if $f\in C^{0,1}$, one can obtain $C^{2,\mathrm{lnL}}$ regularity (the so-called ``ln-Lipschitz'' regularity). Based on the same idea to the $C^{k,\alpha}$ regularity, we can derive systematically these two kinds of regularity for fully nonlinear elliptic equations.

In this section, we list the pointwise $C^k$ regularity (precisely $C^{k,\omega}$ regularity for some $\omega$) for the reader's convenience. Since their proofs are similar to that of $C^{k,\alpha}$ regularity (with some modifications), we only prove the interior $C^{2}$ regularity and the boundary $C^{1}$ regularity for instance.

Since we consider $C^{k}$ regularity, we will encounter the modulus of continuity frequently. For convenience, we introduce the following notion. Given $0<\alpha,\eta<1$, $r_0>0$ and a Dini function $\omega$, we define the following modulus of continuity
\begin{equation}\label{e10.1}
\bar{\omega}(r)\coloneqq
r^{\alpha}+r^{\alpha}\int_{r/\rho}^{1} \frac{\omega(r_0\tau)}{\tau^{1+\alpha}}d\tau
+\int_{0}^{r/(\eta^2\rho)}\frac{\omega(r_0\tau)d\tau}{\tau}, \quad 0<r\leq \eta^2\rho
\end{equation}
and write
\begin{equation}\label{e10.2}
\bar{\omega}=W(\omega,\alpha,r_0,\eta,\rho).
\end{equation}

Now, we state the main results in this section. Recall that for pointwise regularity, we always assume that $r_0=1$ if we say that some function $f\in C^{k,\mathrm{Dini}}(0)$ (see the conventions at the beginning of \Cref{In-C1a-mu}).
\begin{theorem}\label{t-C1-i}
Let $p>n$ and $u$ be a viscosity solution of
\begin{equation*}
F(D^2u, Du, u,x)=f ~~\mbox
{in}~~ B_1.
\end{equation*}
Suppose that $F$ satisfies \cref{SC2} and \cref{e.C1a.beta}. Assume that
\begin{equation*}
b\in L^{p}(B_1),~~~ c\in C^{-1,\mathrm{Dini}}(0),~~~ \|\beta_1\|_{C^{-1,1}(0)}\leq \delta_0,~~~
\gamma_1\in C^{-1,\mathrm{Dini}}(0),~~~ f\in C^{-1,\mathrm{Dini}}(0),
\end{equation*}
where $0<\delta_0<1$ depends only on $n,\lambda,\Lambda$ and $p$.

Then $u\in C^{1,\omega_u}(0)$, i.e., there exist $P\in \mathcal{P}_1$ and a modulus of continuity $\omega_u$ such that
\begin{equation}\label{e.C1-1-i-mu}
  |u(x)-P(x)|\leq C|x|\omega_u (|x|), ~\forall ~x\in B_{\eta^2\rho}
\end{equation}
and
\begin{equation}\label{e.C1-2-i-mu}
|Du(0)|\leq C,
\end{equation}
where
\begin{equation*}
  \begin{aligned}
\omega_u=W(\tilde{\omega},\tilde{\alpha},1,\eta,\rho), \quad
\tilde{\alpha}=\min(\bar{\alpha}/2,1-n/p),\quad\tilde{\omega}=\max(\omega_c,\omega_{\gamma_1},\omega_f);
  \end{aligned}
\end{equation*}
$0<\eta<1$ depends only on $n$, $\lambda$, $\Lambda$ and $p$; $C$ and $\rho$ depend also on $\mu$, $\|b\|_{L^p(B_1)}$, $\|c\|_{C^{-1,\mathrm{Dini}}(0)}$, $\omega_0$, $\|\gamma_1\|_{C^{-1,\mathrm{Dini}}(0)}$, $\|f\|_{C^{-1,\mathrm{Dini}}(0)}$ and $\|u\|_{L^{\infty }(B_1)}$.

In particular, if $F$ satisfies \cref{SC1}, we have the following explicit estimates
\begin{equation}\label{e.C1-1-i}
  |u(x)-P(x)|\leq \tilde{C}|x|\omega_u (|x|),~\forall ~x\in B_{\eta^2\rho},
\end{equation}
\begin{equation}\label{e.C1-2-i}
|Du(0)|\leq \tilde{C}
\end{equation}
and
\begin{equation*}
\tilde{C}=C\left(\|u\|_{L^{\infty }(B_1)} +\|f\|_{C^{-1,\mathrm{Dini}}(0)}+\|\gamma_1\|_{C^{-1,\mathrm{Dini}}(0)}\right),
\end{equation*}
where $C$ and $\rho$ depend only on $n,\lambda,\Lambda,p,\|b\|_{L^p(B_1)}$ and $\|c\|_{C^{-1,\mathrm{Dini}}(0)}$.
\end{theorem}
\begin{remark}\label{re10.2}
We do not know any interior pointwise $C^1$ regularity result for fully nonlinear uniformly elliptic equations so far. Recently, interior $C^1$ regularity has been derived (see \cite{MR4477428} and \cite{MR4632592}) for viscosity solutions of some special degenerate fully nonlinear equations.
\end{remark}
~\\

For $C^2$ regularity, we assume that $\omega_0$ can be written as:
\begin{equation}\label{e.omega0-c-1}
\omega_0(K,s)=\hat\omega_0(K)\check{\omega}_{0}(s), ~\forall ~K,s>0,
\end{equation}
where $\hat{\omega}_0$ is non-decreasing and $\check{\omega}_{0}$ is a Dini function.

For example, consider the following equation:
\begin{equation}\label{e10.3}
\Delta u+\frac{1}{\ln^2|u|+1}=0 \quad\mbox{ in}~B_1.
\end{equation}
Then it satisfies the structure condition \cref{SC2} with some $\omega_0$, which can be written as in \cref{e.omega0-c-1}. Hence, by the following \Cref{t-C2-i}, we conclude that $u\in C^{2}(x_0)$ for any $x_0\in B_{1/2}$, i.e., $u\in C^2(\bar{B}_{1/2})$.

Indeed, for \cref{e10.3}, we have
\begin{equation*}
F(M,p,s,x)=tr(M)+\frac{1}{\ln^2|s|+1}.
\end{equation*}
Now, we show that $F$ satisfies the structure condition \cref{SC2} with \cref{e.omega0-c-1} holding. Let
\begin{equation*}
H(s)=\frac{1}{\ln^2|s|+1}.
\end{equation*}
Without loss of generality, we consider $t>s>0$. Set $r=t-s$ and we have
\begin{equation*}
  \begin{aligned}
H(t)-H(s)=&\int_{0}^{1} H'(s+\tau (t-s))(t-s)d\tau\\
    =&\int_{0}^{1} H'(s+\tau r)rd\tau
    =-2\int_{0}^{1} \frac{r\ln (s+\tau r)}{\left(\ln^2(s+\tau r)+1\right)^2(s+\tau r)} d\tau.
  \end{aligned}
\end{equation*}
Hence,
\begin{equation*}
|H(t)-H(s)|\leq 2\int_{0}^{1} \frac{r\big|\ln (s+\tau r)\big|}{\left(\ln^2(s+\tau r)+1\right)^2(s+\tau r)} d\tau.
\end{equation*}

Let
\begin{equation*}
h(\xi)=\frac{\big|\ln \xi\big|}{\left(\ln^2\xi+1\right)^2\xi},~\xi>0.
\end{equation*}
By analyzing the behaviour of $h$, we know that there exists $0<c_0<1/2$ such that
\begin{equation*}
h ~~\mbox{is strictly decreasing in}~~(0,c_0)~~\mbox{  and  }~~h(\xi)\leq h(c_0),~\forall ~\xi\in [c_0,+\infty).
\end{equation*}
Hence, for $r\leq c_0$,
\begin{equation*}
h(s+\tau r)\leq h(\tau r)= \frac{\big|\ln (\tau r)\big|}{\left(\ln^2(\tau r)+1\right)^2\tau r}, ~\forall ~s>0.
\end{equation*}
Thus,
\begin{equation*}
|H(t)-H(s)|\leq 2\int_{0}^{1} \frac{r\big|\ln (\tau r)\big|}{\left(\ln^2(\tau r)+1\right)^2\tau r} d\tau.
\end{equation*}

For any $K,r>0$, we take
\begin{equation*}
\omega_0(K,r)\coloneqq \hat\omega_0(K)\check{\omega}_{0}(r),
\end{equation*}
where
\begin{equation*}
\hat\omega_0(K)\coloneqq 1, \quad \check{\omega}_0(r)\coloneqq 2\int_{0}^{1} \frac{r\big|\ln (\tau r)\big|}{\left(\ln^2(\tau r)+1\right)^2\tau r} d\tau.
\end{equation*}
Then $F$ satisfies \cref{SC2} with this $\omega_0$. Note that
\begin{equation*}
  \begin{aligned}
\int_{0}^{c_0}\frac{\check{\omega}_0(r)}{r}dr
=&2\int_{0}^{c_0}\int_{0}^{1} \frac{\big|\ln (\tau r)\big|}{\left(\ln^2(\tau r)+1\right)^2\tau r}d\tau dr\\
=&-2\int_{0}^{c_0}\int_{0}^{1} \frac{\ln\tau r}{\left(\ln^2\tau r+1\right)^2\tau r}d\tau dr\\
=&\int_{0}^{c_0}\frac{dr}{\left(\ln^2r+1\right)r} dr\leq C.
  \end{aligned}
\end{equation*}
That is, $\check{\omega}_0$ is a Dini function. Therefore, \cref{e.omega0-c-1} holds.

On the other hand, there is another way to conclude that $u\in C^2(\bar{B}_{1/2})$. Let
\begin{equation*}
f(x)=\frac{-1}{\ln^2|u(x)|+1},~\forall ~x\in B_1
\end{equation*}
and we regard $f$ as the right-hand term. Since $u$ is a viscosity solution, $f\in L^{\infty}(B_1)$. By the interior $C^{1,\alpha}$ regularity, $u\in C^{1,\alpha}(\bar{B}_{3/4})$ for any $0<\alpha<1$. Then by a similar analysis as above,
\begin{equation*}
|f(x)-f(x_0)|\leq \frac{C_1}{\ln^2C_2|x-x_0|+1},~\forall ~x_0,x\in B_{3/4}.
\end{equation*}
Hence, $f\in C^{\mathrm{Dini}}(\bar{B}_{3/4})$. By the following \Cref{t-C2-i}, $u\in C^2(\bar{B}_{1/2})$.

Now, we state the interior pointwise $C^2$ regularity.
\begin{theorem}\label{t-C2-i}
Let $u$ be a viscosity solution of
\begin{equation*}
F(D^2u, Du, u,x)=f \quad\mbox{in}~~ B_1.
\end{equation*}
Suppose that $F$ satisfies \cref{SC2} and \cref{e.C2a-KF} with $G$ being convex in $M$. Assume that $\omega_0$ satisfies \cref{e.omega0-c-1}, $\beta_2\in C^{\mathrm{Dini}}(0)$ and $f\in C^{\mathrm{Dini}}(0)$.

Then $u\in C^{2,\omega_u}(0)$, i.e., there exist $P\in \mathcal{P}_2$ and $\omega_u$ such that
\begin{equation}\label{e.C2-1-i-mu}
  |u(x)-P(x)|\leq C|x|^2\omega_u(|x|), ~\forall ~x\in B_{\eta^2\rho},
\end{equation}
\begin{equation}\label{e.C2-3-i-mu}
|F(D^2P,DP(x),P(x),x)-f(0)|\leq C\omega_{\beta_2}(|x|), ~\forall ~x\in B_{\eta^2\rho}
\end{equation}
and
\begin{equation}\label{e.C2-2-i-mu}
|Du(0)|+|D^2u(0)|\leq C,
\end{equation}
where
\begin{equation*}
  \begin{aligned}
\omega_u=W(\tilde{\omega},\bar{\alpha}/2,1,\eta,\rho), \quad \tilde{\omega}=\max(\check{\omega}_{0},\omega_{\beta_2},\omega_f);
  \end{aligned}
\end{equation*}
$0<\eta<1$ depends only on $n$, $\lambda$ and $\Lambda$; $C$ and $\rho$ depend also on $\mu$, $b_0$, $c_0$, $\omega_0$,
$\|\beta_2\|_{C^{\mathrm{Dini}}(0)}$, $\omega_{\beta_2}$, $\omega_{2}$,
$\|f\|_{C^{\mathrm{Dini}}(0)}$ and $\|u\|_{L^{\infty}(B_1)}$.

In particular,  if $F$ satisfies \cref{SC1} and \cref{e.C2a-KF-0} with $\gamma_2\in C^{\mathrm{Dini}}(0)$,
\begin{equation}\label{e.C2-1-i}
|u(x)-P(x)|\leq \tilde{C}|x|^2\omega_u(|x|), ~\forall ~x\in B_{\eta^2\rho},
\end{equation}
\begin{equation}\label{e.C2-3-i}
\begin{aligned}
|F(D^2P,DP(x),P(x),x)-f(0)|\leq \tilde{C}\omega_{\beta_2}(|x|),~\forall ~x\in B_{\eta^2\rho},
\end{aligned}
\end{equation}
\begin{equation}\label{e.C2-2-i}
|Du(0)|+|D^2u(0)|\leq \tilde{C}
\end{equation}
and
\begin{equation*}
\tilde{C}=C\left(\|u\|_{L^{\infty }(B_1)}+\|f\|_{C^{\mathrm{Dini}}(0)}
+\|\gamma_2\|_{C^{\mathrm{Dini}}(0)}\right),
\end{equation*}
where $C$ and $\rho$ depend only on $n,\lambda,\Lambda,b_0,c_0,\|\beta_2\|_{C^{\mathrm{Dini}}(0)}$ and $\omega_{\beta_2}$.
\end{theorem}

\begin{remark}\label{r-2.16}
The $C^2$ regularity has been investigated much more extensively than $C^1$ regularity. For Poisson equation, it can be derived by the integral representation (see\cite[(4.47)]{MR1814364}). For linear equations, Burch \cite{MR521856} obtained the $C^2$ regularity for generalized solutions. Sperner \cite{MR650493} proved the existence of $C^2$ solutions for Dirichlet problems. Wang \cite{MR2273802} deduced explicit formula for the modulus of continuity of $D^2u$ for $C^2$ solutions by a simple method.

For fully nonlinear elliptic equations, Kovats \cite{MR1629506} first obtained the $C^2$ regularity under Dini conditions (see also \cite{MR1713596}). Wang \cite{MR2273802} also proved the $C^2$ regularity for fully nonlinear elliptic equations with additional assumption $F\in C^{1,1}$. \Cref{t-C2-i} is the first pointwise $C^2$ regularity for fully nonlinear elliptic equations with lower order terms.
\end{remark}

\begin{remark}\label{r-2.13}
In fact, $\bar{\alpha}/2$ can be replaced by any $0<\alpha_0<\bar{\alpha}$ in the expression of
$\omega_u$ (see \Cref{r-10.2}). Similar remarks can be made for other $C^k$ ($k\geq 2$) regularity.

Note that the constant $C$ depends not only $\|\beta_2\|_{C^{\mathrm{Dini}}(0)}$ but also $\omega_{\beta_2}$. The reason is that $C^2$ regularity is a critical case for $\beta_2$ and the normalization procedure have to depend on $\omega_{\beta_2}$ (see \cref{e9.1}).
\end{remark}
~\\

The following is the interior pointwise $C^k$ ($k\geq 3$) regularity.
\begin{theorem}\label{t-Ck-i}
Let $k\geq 3$ and $u$ be a viscosity solution of
\begin{equation*}
F(D^2u, Du, u,x)=f \quad\mbox{in}~~ B_1.
\end{equation*}
Suppose that $F\in C^{k-2,\mathrm{Dini}}(0)$ satisfies \cref{SC2}, $\omega_0$ satisfies \cref{e.omega0-2} and $f\in C^{k-2,\mathrm{Dini}}(0)$.

Then $u\in C^{k,\omega_u}(0)$, i.e., there exist $P\in \mathcal{P}_k$ and $\omega_u$ such that
\begin{equation}\label{e.Ck-i-1}
  |u(x)-P(x)|\leq C|x|^k\omega_u(|x|), ~\forall ~x\in B_{\eta^2\rho},
\end{equation}
\begin{equation}\label{e.Ck-i-3}
|F(D^2P(x),DP(x),P(x),x)-P_f(x)|\leq C|x|^{k-2}\omega_{F}(|x|), ~\forall ~x\in B_{\eta^2\rho}
\end{equation}
and
\begin{equation}\label{e.Ck2-i-2}
|Du(0)|+\cdots+|D^ku(0)|\leq C,
\end{equation}
where
\begin{equation*}
  \begin{aligned}
\omega_u=W(\tilde{\omega},\bar{\alpha}/2,1,\eta,\rho),\quad \tilde{\omega}=\max(\omega_{F},\omega_f);
  \end{aligned}
\end{equation*}
$\eta$ depends only on $k,n,\lambda$ and $\Lambda$; $C$ and $\rho$ depend also on $\mu,b_0,c_0,\omega_0$, $\|F\|_{C^{k-2,\mathrm{Dini}}(0)}$,
$\omega_3,\omega_4$, $\|f\|_{C^{k-2,\mathrm{Dini}}(0)}$ and $\|u\|_{L^{\infty}(B_1)}$.
\end{theorem}
\begin{remark}\label{re10.4}
Similar to pointwise $C^{k,\alpha}$ regularity, we have not seen any pointwise $C^{k}$ ($k\geq 3$) regularity.
\end{remark}
~\\

For the boundary pointwise $C^{1}$, $C^{2}$ and $C^{k}$ regularity, we have
\begin{theorem}\label{t-C1}
Let $u$ satisfy
\begin{equation}
\left\{\begin{aligned}
&u\in S^*(\lambda,\Lambda,\mu,b,f)&& \quad\mbox{in}~~\Omega\cap B_1;\\
&u=g&& \quad\mbox{on}~~\partial \Omega\cap B_1.
\end{aligned}\right.
\end{equation}
Suppose that
\begin{equation*}
b\in L^{p}(\Omega\cap B_1)(p>n),\quad f\in C^{-1,\mathrm{Dini}}(0),\quad \partial\Omega\cap B_1\in C^{1,\mathrm{Dini}}(0),\quad g\in C^{1,\mathrm{Dini}}(0).
\end{equation*}

Then $u$ is $C^{1,\omega_u}$ at $0$, i.e., there exist $P\in \mathcal{P}_1$ and $\omega_u$ such that
\begin{equation}\label{e.C1-1-mu}
  |u(x)-P(x)|\leq C|x|\omega_u(|x|), ~\forall ~x\in \Omega\cap B_{\eta^2\rho},
\end{equation}
\begin{equation}\label{e.C1-3-mu}
D_{x'}u(0)=D_{x'}g(0)
\end{equation}
and
\begin{equation}\label{e.C1-2-mu}
|Du(0)|\leq C ,
\end{equation}
where
\begin{equation*}
  \begin{aligned}
\omega_u=W(\tilde{\omega},\tilde{\alpha},1,\eta,\rho), \quad
\tilde{\alpha}=\min(\bar{\alpha}/2,1-n/p), \quad  \tilde{\omega}=\max(\omega_f, \omega_\Omega, \omega_g);
  \end{aligned}
\end{equation*}
$\eta$ depends only on $n$, $\lambda$, $\Lambda$ and $p$; $C$ and $\rho$ depend also on
$\mu,\|b\|_{L^p(\Omega\cap B_1)},\|\partial \Omega\cap B_1\|_{C^{1,\mathrm{Dini}}(0)}, \omega_{\Omega}$, $\|f\|_{C^{-1,\mathrm{Dini}}(0)}$, $\|g\|_{C^{1,\mathrm{Dini}}(0)}$ and $\|u\|_{L^{\infty }(\Omega\cap B_1)}$.

In particular, if $\mu=0$,
\begin{equation}\label{e.C1-1}
  |u(x)-P(x)|\leq \tilde{C}|x|\omega_u (|x|), ~~\forall ~x\in \Omega\cap B_{\eta^2\rho},
\end{equation}
\begin{equation}\label{e.C1-2}
|Du(0)|\leq \tilde{C}
\end{equation}
and
\begin{equation*}
\tilde{C}=C\left(\|u\|_{L^{\infty }(\Omega\cap B_1)}+\|f\|_{C^{-1,\mathrm{Dini}}(0)}+\|g\|_{C^{1,\mathrm{Dini}}(0)}\right),
\end{equation*}
where $C$ and $\rho$ depend only on $n, \lambda, \Lambda,p,\|b\|_{L^p(\Omega\cap B_1)},\|\partial \Omega\cap B_1\|_{C^{1,\mathrm{Dini}}(0)}$ and $\omega_{\Omega}$.
\end{theorem}

\begin{remark}\label{r-2.18}
Ma and Wang \cite{MR2853528} have proved ($\mu=b=0$) the boundary pointwise $C^1$ regularity for fully nonlinear equations. Huang, Zhai and Zhou \cite{MR3933752} obtained the boundary pointwise $C^1$ regularity for linear equations with $b\in C^{-1,\mathrm{Dini}}(0)$. Braga, Gomes, Moreira and Wang \cite{MR4151474} gave the boundary pointwise $C^{1}$ regularity on a flat boundary with more restrictive modulus of continuity. \Cref{t-C1} is the first boundary pointwise $C^1$ regularity for Pucci's class with lower order terms.
\end{remark}

\begin{remark}\label{r-2.19}
Note that the constant $C$ depends not only
$\|\partial \Omega\cap B_1\|_{C^{1,\mathrm{Dini}}(0)}$ but also $\omega_{\Omega}$. The reason is similar to the $C^2$ regularity for $\beta_2$ (see \cref{e9.2}).
\end{remark}
~\\

\begin{theorem}\label{t-C2}
Let $u$ be a viscosity solution of
\begin{equation*}
\left\{\begin{aligned}
&F(D^2u,Du,u,x)=f&& \quad\mbox{in}~~\Omega\cap B_1;\\
&u=g&& \quad\mbox{on}~~\partial \Omega\cap B_1.
\end{aligned}\right.
\end{equation*}
Suppose that $F$ satisfies \cref{SC2}, \cref{e.C2a-KF} and $\omega_0$ satisfies \cref{e.omega0-c-1}. Assume that
\begin{equation*}
\beta_2\in C^{\mathrm{Dini}}(0),\quad f\in C^{\mathrm{Dini}}(0),\quad
\partial\Omega\cap B_1\in C^{2,\mathrm{Dini}}(0),\quad g\in C^{2,\mathrm{Dini}}(0).
\end{equation*}

Then $u\in C^{2,\omega_u}(0)$, i.e., there exist $P\in \mathcal{P}_2$ and $\omega_u$ such that
\begin{equation}\label{e.C2-1-mu}
  |u(x)-P(x)|\leq C|x|^2\omega_u(|x|), ~~\forall ~x\in \Omega \cap B_{\eta^2\rho},
\end{equation}
\begin{equation}\label{e.C2-3-mu}
|F(D^2P,DP(x),P(x),x)-f(0)|\leq C\omega_{\beta_2}(|x|), ~~\forall ~x\in \Omega \cap B_{\eta^2\rho},
\end{equation}
\begin{equation}\label{e.C2-4-mu}
D^l_{x'}u(x',P_{\Omega}(x'))=D^l_{x'}g(x',P_{\Omega}(x'))~\mbox{at}~0,~\forall ~0\leq l\leq 2
\end{equation}
and
\begin{equation}\label{e.C2-2-mu}
|Du(0)|+|D^2u(0)|\leq C,
\end{equation}
where
\begin{equation*}
\omega_u=W(\tilde{\omega},\bar{\alpha}/2,1,\eta,\rho),\quad \tilde\omega=\max(\omega_{\beta_2}, \omega_f, \omega_\Omega, \omega_g);
\end{equation*}
$\eta$ depends only on $n$, $\lambda$ and $\Lambda$; $C$ and $\rho$ depend also on
$\mu$, $b_0$, $c_0,\omega_0$,
$\|\beta_2\|_{C^{\mathrm{Dini}}(0)}$, $\omega_{\beta_2},\omega_{2}$,
$\|\partial \Omega\cap B_1\|_{C^{2,\mathrm{Dini}}(0)}$, $\|f\|_{C^{\mathrm{Dini}}(0)}$,
$\|g\|_{C^{2,\mathrm{Dini}}(0)}$ and $\|u\|_{L^{\infty}(\Omega\cap B_1)}$.

In particular, if $F$ satisfies \cref{SC1} and \cref{e.C2a-KF-0} with $\gamma_2\in C^{\mathrm{Dini}}(0)$,
\begin{equation}\label{e.C2-1}
\begin{aligned}
  |u(x)-P(x)|\leq \tilde{C}|x|^2\omega_u(|x|), ~\forall ~x\in \Omega\cap B_{\eta^2\rho},
\end{aligned}
\end{equation}
\begin{equation}\label{e.C2-3}
\begin{aligned}
|&F(D^2P,DP(x),P(x),x)-f(0)|\leq \tilde{C}\omega_{\beta_2}(|x|), ~\forall ~x\in \Omega \cap B_{\eta^2\rho},
\end{aligned}
\end{equation}
\begin{equation}\label{e.C2-2}
|Du(0)|+|D^2u(0)|\leq \tilde{C}
\end{equation}
and
\begin{equation*}
\tilde{C}=C \left(\|u\|_{L^{\infty }(\Omega\cap B_1)}+\|f\|_{C^{\mathrm{Dini}}(0)}+\|\gamma_2\|_{C^{\mathrm{Dini}}(0)}
+\|g\|_{C^{2,\mathrm{Dini}}(0)}\right),
\end{equation*}
where $C$ and $\rho$ depend only on $n,\lambda,\Lambda,b_0,c_0,\|\beta_2\|_{C^{\mathrm{Dini}}(0)}$
$\omega_{\beta_2}$ and $\|\partial \Omega\cap B_1\|_{C^{2,\mathrm{Dini}}(0)}$.
\end{theorem}
\begin{remark}\label{re10.3}
Zou and Chen \cite{MR1924273} obtained boundary $C^2$ regularity on a flat boundary for viscosity solutions of fully nonlinear equations.
\end{remark}
~\\

The following is the boundary pointwise $C^k$ ($k\geq 3$) regularity, which is new for linear equations.
\begin{theorem}\label{t-Ck}
Let $k\geq 3$ and $u$ be a viscosity solution of
\begin{equation*}
\left\{\begin{aligned}
&F(D^2u,Du,u,x)=f&& \quad\mbox{in}~~\Omega\cap B_1;\\
&u=g&& \quad\mbox{on}~~\partial\Omega\cap B_1.
\end{aligned}\right.
\end{equation*}
Suppose that \cref{SC2} and \cref{e.omega0-2} hold, and
\begin{equation*}
F\in C^{k-2,\mathrm{Dini}}(0),\quad f\in C^{k-2,\mathrm{Dini}}(0),\quad
\partial\Omega\cap B_1\in C^{k,\mathrm{Dini}}(0),\quad g\in C^{k,\mathrm{Dini}}(0).
\end{equation*}

Then $u\in C^{k,\omega_u}(0)$, i.e., there exist $P\in \mathcal{P}_k$ and $\omega_u$ such that
\begin{equation}\label{e.Ck-1}
  |u(x)-P(x)|\leq  C|x|^k\omega_u(|x|), ~~\forall ~x\in \Omega\cap B_{\eta^2\rho},
\end{equation}
\begin{equation}\label{e.Ck-3}
|F(D^2P(x),DP(x),P(x),x)-P_f(x)|\leq C|x|^{k-2}\omega_{F}(|x|), ~~\forall ~x\in \Omega \cap B_{\eta^2\rho},
\end{equation}
\begin{equation}\label{e.Ck-5}
D^l_{x'}u(x',P_\Omega(x'))=D^l_{x'}g(x',P_\Omega(x'))~\mbox{at}~0,~\forall ~0\leq l\leq k
\end{equation}
and
\begin{equation}\label{e.Ck2-2}
|Du(0)|+\cdots+|D^ku(0)|\leq C,
\end{equation}
where
\begin{equation*}
\omega_u=W(\tilde{\omega},\bar{\alpha}/2,1,\eta,\rho),\quad \tilde\omega=\max(\omega_{F},\omega_f, \omega_\Omega, \omega_g);
\end{equation*}
$\eta$ depends only on $k,n$, $\lambda$ and $\Lambda$; $C$ and $\rho$ depend also on
$\mu$, $b_0$, $c_0$, $\omega_0$,
$\|F\|_{C^{k-2,\mathrm{Dini}}(0)}$, $\omega_3$, $\omega_4$,
$\|\partial \Omega\cap B_1\|_{C^{k,\mathrm{Dini}}(0)}$,
$\|f\|_{C^{k-2,\mathrm{Dini}}(0)},\|g\|_{C^{k,\mathrm{Dini}}(0)}$ and
$\|u\|_{L^{\infty}(\Omega\cap B_1)}$.
\end{theorem}

Similar to $C^{k,\alpha}$ regularity, by combining the interior and boundary regularity, we have the local and global $C^{k}$ regularity. For the completeness and the convention of citation, we list them as follows.
\begin{corollary}\label{t-C1-global}
Let $p>n$ and $\Gamma\subset \partial \Omega $ be relatively open (may be empty) and $u$ be a viscosity solution of
\begin{equation*}
\left\{\begin{aligned}
&F(D^2u,Du,u,x)=f&& \quad\mbox{in}~~\Omega;\\
&u=g&& \quad\mbox{on}~~\Gamma.
\end{aligned}\right.
\end{equation*}
Suppose that $F$ satisfies \cref{SC2} and \cref{e.C1a.beta} with some $G_{x_0}$ at any $x_0\in \Omega\cup \Gamma$ and
\begin{equation*}
  \begin{aligned}
&b\in L^{p}(\Omega),\quad c\in C^{-1,\mathrm{Dini}}(\bar\Omega),\quad
f\in C^{-1,\mathrm{Dini}}(\bar\Omega),\quad \Gamma\in C^{1,\mathrm{Dini}},\quad
g\in C^{1,\mathrm{Dini}}(\bar\Gamma),\\
&\gamma_1\in C^{-1,\mathrm{Dini}}(\bar{\Omega}), \quad
\beta_1(x,x_0)\leq \delta_0,~\forall ~x_0,x\in\Omega\cup \Gamma~\mbox{ with }~|x-x_0|<r_0,\quad
  \end{aligned}
\end{equation*}
where $0<\delta_0<1$ depends only on $n,\lambda,\Lambda$ and $p$.

Then for any $\Omega'\subset\subset \Omega\cup \Gamma$, we have $u\in C^{1,\omega_u}(\bar{\Omega}')$ and
\begin{equation}\label{e.C1-1-i-mu-global}
 \|u\|_{C^{1,\omega_u}(\bar{\Omega}')}\leq C,
\end{equation}
where $C$ depends only on $n$, $\lambda,\Lambda$, $r_0,\mu,\|b\|_{L^p(\Omega)}$,
$\|c\|_{C^{-1,\mathrm{Dini}}(\bar{\Omega})}$, $\omega_0,
\|\partial \Omega'\cap\Gamma\|_{C^{1,\mathrm{Dini}}}$, $\omega_{\Omega}$, $\Omega'$, $\mathrm{dist}(\Omega',\partial \Omega\backslash \Gamma)$, $\|f\|_{C^{-1,\mathrm{Dini}}(\bar{\Omega})},
\|g\|_{C^{1,\mathrm{Dini}}(\bar\Gamma)}$ and $\|u\|_{L^{\infty }(\Omega)}$; $\omega_u$ depends also on $\omega_c,\omega_{\gamma_1},\omega_f, \omega_g$.

In particular, if $F$ satisfies \cref{SC1},
\begin{equation}\label{e.C1-1-i-global}
 \|u\|_{C^{1,\omega_u}(\bar{\Omega}')}\leq C \left(\|u\|_{L^{\infty }(\Omega)}
 +\|f\|_{C^{-1,\mathrm{Dini}}(\bar{\Omega})}
 +\|\gamma_1\|_{C^{-1,\mathrm{Dini}}(\bar{\Omega})}
 +\|g\|_{C^{1,\mathrm{Dini}}(\bar\Gamma)}\right),
\end{equation}
where $C$ depends only on $n,\lambda,\Lambda,p,r_0,\|b\|_{L^p(\Omega)},
\|c\|_{C^{-1,\mathrm{Dini}}(\bar{\Omega})}$, $\|\partial \Omega'\cap\Gamma\|_{C^{1,\mathrm{Dini}}}$, $\omega_{\Omega}$, $\Omega'$ and $\mathrm{dist}(\Omega',\partial \Omega\backslash \Gamma)$.
\end{corollary}
\begin{remark}\label{re10.1}
In above theorem, we have assumed the same ``$r_0$'' used in \cref{e.C1a.beta}, $c\in C^{-1,\mathrm{Dini}}(\bar\Omega)$ etc.

Similar to the pointwise regularity, we have an explicit expression for $\omega_u(r)$ when $r$ is small:
\begin{equation*}
\begin{aligned}
\omega_u(r)=W(\tilde{\omega},\tilde{\alpha},r_0,\eta,\rho)(r),~\forall ~0<r\leq \eta^2 \rho,\\
\end{aligned}
\end{equation*}
where
\begin{equation*}
\tilde{\alpha}=\min(\bar{\alpha}/2,1-n/p), \quad  \tilde{\omega}=\max(\omega_c,\omega_{\gamma_1},\omega_f, \omega_\Omega, \omega_g);
\end{equation*}
$\eta$ depends only on $n$, $\lambda,\Lambda$ and $p$; $\rho$ depends also on
$r_0,\mu,\|b\|_{L^p(\Omega)}$,
$\|c\|_{C^{-1,\mathrm{Dini}}(\bar{\Omega})}$, $\omega_0,
\|\partial \Omega'\cap\Gamma\|_{C^{1,\mathrm{Dini}}}$, $\omega_{\Omega}$, $\Omega'$, $\mathrm{dist}(\Omega',\partial \Omega\backslash \Gamma)$, $\|f\|_{C^{-1,\mathrm{Dini}}(\bar{\Omega})},
\|g\|_{C^{1,\mathrm{Dini}}(\bar\Gamma)}$ and $\|u\|_{L^{\infty }(\Omega)}$.
\end{remark}
~\\

\begin{corollary}\label{t-C2-global}
Let $\Gamma\subset \partial \Omega $ be relatively open  and $u$ be a viscosity solution of
\begin{equation*}
\left\{\begin{aligned}
&F(D^2u,Du,u,x)=f&& \quad\mbox{in}~~\Omega;\\
&u=g&& \quad\mbox{on}~~\Gamma.
\end{aligned}\right.
\end{equation*}
Suppose that $F$ satisfies \cref{SC2} and \cref{e.C2a-KF} with some $G_{x_0}$ at any $x_0\in \Omega\cup \Gamma$ and $G_{x_0}$ is convex in $M$. Assume that $\omega_0$ satisfies \cref{e.omega0-c-1} and
\begin{equation*}
\beta_2\in C^{\mathrm{Dini}}(\bar{\Omega}),\quad f\in C^{\mathrm{Dini}}(\bar{\Omega}),\quad \Gamma\in C^{2,\mathrm{Dini}},\quad g\in C^{2,\mathrm{Dini}}(\bar\Gamma).
\end{equation*}

Then for any $\Omega'\subset\subset \Omega\cup \Gamma$, we have $u\in C^{2,\omega_u}(\bar{\Omega}')$ and
\begin{equation}\label{e.C2-1-i-mu-global}
 \|u\|_{C^{2,\omega_u}(\bar{\Omega}')}\leq C,
\end{equation}
where $C$ depends only on $n$, $\lambda, \Lambda$, $r_0,\mu,b_0,c_0,\omega_0$,
$\|\beta_2\|_{C^{\mathrm{Dini}}(\bar{\Omega})}$, $\omega_{\beta_2}$, $\omega_2$, $\|\partial\Omega'\cap\Gamma\|_{C^{2,\mathrm{Dini}}}$, $\|f\|_{C^{\mathrm{Dini}}(\bar{\Omega})}$,
$\|g\|_{C^{2,\mathrm{Dini}}(\bar\Gamma)}$ and $\|u\|_{L^{\infty }(\Omega)}$; $\omega_u$ depends also on $\omega_f, \omega_\Omega, \omega_g$.

In particular, if $F$ satisfies \cref{SC1} and \cref{e.C2a-KF-0} with $\gamma_2\in C^{\mathrm{Dini}}(\bar\Gamma)$,
\begin{equation}\label{e.C2-1-i-global}
 \|u\|_{C^{2,\omega_u}(\bar{\Omega}')}\leq C \left(\|u\|_{L^{\infty }(\Omega)}+\|f\|_{C^{\mathrm{Dini}}(\bar{\Omega})}
 +\|\gamma_2\|_{C^{\mathrm{Dini}}(\bar{\Omega})}
 +\|g\|_{C^{2,\mathrm{Dini}}(\bar\Gamma)}\right),
\end{equation}
where $C$ depends only on $n,\lambda,\Lambda,r_0,\mu,b_0,c_0$, $\|\beta_2\|_{C^{\mathrm{Dini}}(\bar{\Omega})}$, $\omega_{\beta_2},\omega_2$ and $\|\Gamma\|_{C^{2,\mathrm{Dini}}}$.
\end{corollary}

\begin{corollary}\label{t-Ck-global}
Let $k\geq 3$, $\Gamma\subset \partial \Omega $ be relatively open and $u$ be a viscosity solution of
\begin{equation*}
\left\{\begin{aligned}
&F(D^2u,Du,u,x)=f&& \quad\mbox{in}~~\Omega;\\
&u=g&& \quad\mbox{on}~~\Gamma.
\end{aligned}\right.
\end{equation*}
Suppose that \cref{SC2} and \cref{e.omega0-2} hold, and
\begin{equation*}
F\in C^{k-2,\mathrm{Dini}}(\bar \Omega),\quad f\in C^{k-2,\mathrm{Dini}}(\bar\Omega),\quad
\Gamma\in C^{k,\mathrm{Dini}},\quad g\in C^{k,\mathrm{Dini}}(\bar\Gamma).
\end{equation*}

Then for any $\Omega'\subset\subset \Omega\cup \Gamma$, we have $u\in C^{k,\omega_u}(\bar{\Omega}')$ and
\begin{equation}\label{e.Ck-1-i-mu-global}
 \|u\|_{C^{k,\omega_u}(\bar{\Omega}')}\leq C,
\end{equation}
where $C$ depends only on $k,n$, $\lambda,\Lambda$, $r_0,\mu,b_0,c_0$,
$\|F\|_{C^{k-2,\mathrm{Dini}}(\bar\Omega)}$, $\omega_3,\omega_4$,
$\|\partial \Omega'\cap\Gamma\|_{C^{k,\mathrm{Dini}}}$,  $\|f\|_{C^{k-2,\mathrm{Dini}}(\bar\Omega)}$,
$\|g\|_{C^{k,\mathrm{Dini}}(\bar\Gamma)}$ and $\|u\|_{L^{\infty }(\Omega)}$; $\omega_u$ depends also on $\omega_{F},\omega_f, \omega_\Omega, \omega_g$.
\end{corollary}
~\\

In the rest of this section, we give the outline of the proofs of the interior $C^{2}$ regularity and the boundary $C^1$ regularity. The following is the ``key step'' for interior $C^{2}$ regularity, whose proof is almost the same as that of \Cref{In-l-C2a-mu} and we omit it.
\begin{lemma}\label{In-l-C2-mu}
Suppose that $F$ satisfies \cref{e.C2a-KF} with $G$ being convex in $M$. There exists $\delta>0$ depending only on $n,\lambda,\Lambda,\alpha$ and $\omega_2$ such that if $u$ satisfies
\begin{equation*}
F(D^2u, Du, u,x)=f \quad\mbox{in}~~ B_1
\end{equation*}
with
\begin{equation*}
  \begin{aligned}
&u(0)=|Du(0)|=0,\quad \max\left(\|u\|_{L^{\infty}(B_1)},\omega_0(1,1)\right)\leq 1,\\
&\max\left(\mu,b_0,c_0,\|\beta_2\|_{L^{\infty}(B_1)},\|f\|_{L^{\infty}(B_1)}\right)\leq \delta,\\
  \end{aligned}
\end{equation*}
then there exists $P\in\mathcal{HP}_2$ such that
\begin{equation*}
  \|u-P\|_{L^{\infty}(B_{\eta})}\leq \eta^{2+\bar{\alpha}/2},
\end{equation*}
\begin{equation*}
G(D^2P,0,0)=0
\end{equation*}
and
\begin{equation*}
\|P\|\leq \bar{C}+1,
\end{equation*}
where $0<\eta<1$ depends only on $n,\lambda$ and $\Lambda$.
\end{lemma}

Now, we give the ``scaling argument'' of interior $C^2$ regularity.
\begin{lemma}\label{In-t-C2s-mu}
Suppose that $F$ satisfies \cref{e.C2a-KF} with $G$ being convex in $M$. Let $\omega_0$ satisfy \cref{e.omega0-c-1} and $u$ be a viscosity solution of
\begin{equation*}
F(D^2u, Du, u,x)=f \quad\mbox{in}~~ B_1.
\end{equation*}
Assume that
\begin{equation}\label{In-e.C2s-be-mu}
\begin{aligned}
&\|u\|_{L^{\infty}(B_1)}\leq 1,\quad u(0)=|Du(0)|=0, \\
&\mu\leq \frac{\delta_1}{4C_0},\quad b_0\leq \frac{\delta_1}{2},\quad c_0\leq \frac{\delta_1}{\check{\omega}_{0}(C_0)+1},\quad\hat\omega_0(1+C_0)\leq 1,\\
&|\beta_2(x)|\leq \frac{\delta_1}{C_0}\omega_{\beta}(|x|),\quad|f(x)|\leq \delta_1\omega_f(|x|), ~\forall ~x\in B_1,\\
&J_{\tilde\omega}
\leq 1,~\quad(\tilde\omega\coloneqq\max(\check{\omega}_{0},\omega_{\beta},\omega_f))
\end{aligned}
\end{equation}
where $C_0$ depends only on $n,\lambda$ and $\Lambda$, and $\delta_1$ depends also on $\omega_0$ and $\omega_2$.

Then $u\in C^{2,\omega_u}(0)$, i.e., there exist $P\in\mathcal{HP}_2$ and $\omega_u$ such that
\begin{equation}\label{In-e.C2s-1-mu}
  |u(x)-P(x)|\leq |x|^{2}\omega_{u}(|x|), ~\forall ~x\in B_{\eta^2},
\end{equation}
\begin{equation}\label{In-e.C2s-3-mu}
G(D^2P,0,0)=0
\end{equation}
and
\begin{equation}\label{In-e.C2s-2-mu}
\|P\| \leq C,
\end{equation}
where
\begin{equation*}
\omega_u=W(\tilde{\omega},\bar{\alpha}/2,1,\eta,1),
\end{equation*}
$\eta$ (as in \Cref{In-l-C2-mu}) and $C$ depend only on $n,\lambda$ and $\Lambda$.
\end{lemma}
\proof By \Cref{le2.3}, to prove that $u\in C^{2,\omega_u}(0)$, we only need to prove the following. There exist a sequence of $P_m\in\mathcal{HP}_2$ ($m\geq -1$, $P_{-1}\equiv 0$) such that for all $m\geq 0$,
\begin{equation}\label{In-e.C2s-4-mu}
\|u-P_m\|_{L^{\infty }(B_{\eta^{m}})}\leq \eta ^{2m}A_m,
\end{equation}
\begin{equation}\label{In-e.C2s-7-mu}
G(D^2P_m,0,0)=0
\end{equation}
and
\begin{equation}\label{In-e.C2s-6-mu}
\|P_m-P_{m-1}\|\leq (\bar{C}+1)A_{m-1},
\end{equation}
where
\begin{equation}\label{In-e.tC2-Ak}
A_{-1}=A_0=1,~A_m=\max(\tilde\omega(\eta^{m}),\eta^{\bar{\alpha}/2} A_{m-1})
(m\geq 1)
\end{equation}
and $0<\eta<1$ is as in \Cref{In-l-C2-mu}.

We prove \crefrange{In-e.C2s-4-mu}{In-e.C2s-6-mu} by induction. For $m=0$, by setting $P_0\equiv 0$, \crefrange{In-e.C2s-4-mu}{In-e.C2s-6-mu} hold clearly. Suppose that the conclusion holds for $m\leq m_0$. We need to prove that the conclusion holds for $m=m_0+1$.

Let $r=\eta ^{m_{0}}$, $y=x/r$ and
\begin{equation}\label{In-e.C2s-v-mu}
  v(y)=\frac{u(x)-P_{m_0}(x)}{r^2A_{m_0}}.
\end{equation}
Then $v$ satisfies
\begin{equation}\label{In-e.C2s-F-mu}
 \tilde{F}(D^2v,Dv,v,y)=\tilde{f} \quad\mbox{in}~~ B_1,
\end{equation}
where for $(M,p,s,y)\in \mathcal{S}^n\times \mathbb{R}^n\times \mathbb{R}\times \bar B_1$,
\begin{equation*}
  \begin{aligned}
&\tilde{F}(M,p,s,y)=A_{m_0}^{-1}F(A_{m_0}M+D^2P_{m_0},rA_{m_0}p+DP_{m_0}(x),
r^2A_{m_0}s+P_{m_0}(x),x),\\
&\tilde{f}(y)=A_{m_0}^{-1}f(x).\\
  \end{aligned}
\end{equation*}
In addition, define
\begin{equation*}
\tilde{G}(M,p,s)=A_{m_0}^{-1}G(A_{m_0}M+D^2P_{m_0},rA_{m_0}p
,r^{2}A_{m_0}s).
\end{equation*}

In the following, we show that \cref{In-e.C2s-F-mu} satisfies the assumptions of \Cref{In-l-C2-mu}. First, it is easy to verify that
\begin{equation*}
\begin{aligned}
\|v\|_{L^{\infty}(B_1)}\leq& 1, \quad  v(0)=|Dv(0)|=0,
~\quad(\mathrm{by}~\cref{In-e.C2s-be-mu},~\cref{In-e.C2s-4-mu} ~\mbox{and}~\cref{In-e.C2s-v-mu})\\
\|\tilde{f}\|_{L^{\infty}(B_1)}=&A_{m_0}^{-1}\|f\|_{L^{\infty}(B_r)}\leq \delta_1,
 ~\quad(\mathrm{by}~ \cref{In-e.C2s-be-mu})\\
\tilde{G}(0,0,0)=&A_{m_0}^{-1}G(D^2P_{m_0},0,0)=0.~\quad(\mathrm{by}~\cref{In-e.C2s-7-mu})
\end{aligned}
\end{equation*}

By \cref{In-e.C2s-6-mu}, we can choose a constant $C_0$ depending only on $n,\lambda$ and $\Lambda$ such that
\begin{equation*}
\|P_{m}\|\leq C_0,~\forall~0\leq m\leq m_0.
\end{equation*}
Then $\tilde{F}$ and $\tilde{G}$ satisfy the structure condition \cref{SC2} with $\lambda,\Lambda,\tilde{\mu},\tilde{b},\tilde{c}$ and $\tilde{\omega}_0$, where
\begin{equation*}
\begin{aligned}
&\tilde{\mu}=r^2A_{m_0}\mu,\quad \tilde{b}=rb_0+2C_0r\mu,\quad
\tilde{c}=c_0,\quad\\
&\tilde{\omega}_0(K,s)=A_{m_0}^{-1}\omega_0(K+C_0,r^2A_{m_0}s),~\forall ~K,s >0.
\end{aligned}
\end{equation*}
Hence, from \cref{In-e.C2s-be-mu},
\begin{equation*}
\begin{aligned}
 \tilde{\mu}\leq \mu\leq \delta_1,\quad\tilde{b}\leq b_0+2C_0\mu\leq \delta_1,\quad
 \tilde{c}\leq c_0\leq \delta_1
 \end{aligned}
\end{equation*}
and (note that $\omega_0$ satisfies \cref{e.omega0-c-1})
\begin{equation*}
\tilde{\omega}_0(1,1)=A_{m_0}^{-1}\omega_0(1+C_0,r^2A_{m_0})
 =\hat\omega_0(1+C_0)\check{\omega}_{0}(r^2A_{m_0})/A_{m_0}\leq 1.
\end{equation*}
Moreover, $\tilde{\omega}_0$ satisfies \cref{e.omega0-c-1} with
\begin{equation*}
\hat{\tilde{\omega}}_0(K)=\hat{\omega}_0(K+C_0),\quad
\check{\tilde{\omega}}_0(s)=A_{m_0}^{-1}\check{\omega}_{0}(r^2A_{m_0}s),~\forall ~K,s>0.
\end{equation*}

Finally, take $\eta$ small enough such that
\begin{equation*}
  \eta C_0\leq 1.
\end{equation*}
As before (cf. \cref{e5.7}), by combining \cref{SC2}, \cref{e.C2a-KF} and \cref{In-e.C2s-be-mu}, we have for
 $(M,p,s,y)\in \mathcal{S}^n\times \mathbb{R}^n\times \mathbb{R}\times \bar B_1$ (note that $r^{\bar{\alpha}/2}=\eta^{m_0\bar{\alpha}/2}\leq A_{m_0}$),
\begin{equation*}
\begin{aligned}
|\tilde{F}(M&,p,s,y)-\tilde{G}(M,p,s)|\\
\leq& A_{m_0}^{-1}\big(\beta_2(x)(|M|+C_0)\omega_2(|p|,|s|)
+2C_0 r^2A_{m_0}\mu|p|+ C_0^2r^2\mu+C_0rb_0\\
&+c_0\hat\omega_0(|s|+C_0)\check{\omega}_{0}(C_0 r^2)\big)\\
\leq&  \delta_1(|M|+1)\omega_2(|p|,|s|)+\delta_1|p|+\delta_1C_0
+\delta_1\hat\omega_0(|s|+C_0)\\
\coloneqq &\tilde{\beta}_2(y)(|M|+1)\tilde{\omega}_2(|p|,|s|),
\end{aligned}
\end{equation*}
where
\begin{equation*}
\tilde{\beta}_2(y)\equiv\delta_1,\quad \tilde{\omega}_2(|p|,|s|)=
 \omega_2(|p|,|s|)+|p|+C_0+\hat\omega_0(|s|+C_0).
\end{equation*}
Then $\|\tilde{\beta}_2\|_{L^{\infty}(B_1)}\leq \delta_1$.

Choose $\delta_1$ small enough (depending only on $n,\lambda,\Lambda,\omega_0$ and $\omega_2$) such that \Cref{In-l-C2-mu} holds for $\tilde{\omega}_0,\tilde{\omega}_2$ and $\delta_1$. Since \cref{In-e.C2s-F-mu} satisfies the assumptions of \Cref{In-l-C2-mu}, there exists $\tilde{P}(y)\in\mathcal{HP}_2$ such that
\begin{equation*}
\begin{aligned}
    \|v-\tilde{P}\|_{L^{\infty }(B_{\eta})}&\leq \eta ^{2+\bar{\alpha}/2},
\end{aligned}
\end{equation*}
\begin{equation*}
  \tilde{G}(D^2\tilde{P},0,0)=0
\end{equation*}
and
\begin{equation*}
\|\tilde{P}\|\leq \bar{C}+1.
\end{equation*}

Let
\begin{equation*}
P_{m_0+1}(x)=P_{m_0}(x)+r^2A_{m_0}\tilde{P}(y)=P_{m_0}(x)+A_{m_0}\tilde{P}(x).
\end{equation*}
Then \cref{In-e.C2s-7-mu} and \cref{In-e.C2s-6-mu} hold for $m_0+1$. By recalling \cref{In-e.C2s-v-mu}, we have (note that $\eta^{\bar{\alpha}/2}A_{m_0}\leq A_{m_0+1}$ by the definition of $A_{m_0}$)
\begin{equation*}
  \begin{aligned}
\|u-P_{m_0+1}\|_{L^{\infty}(B_{\eta^{m_0+1}})}&= \|u-P_{m_0}-A_{m_0}\tilde{P}(x)\|_{L^{\infty}(B_{\eta r})}\\
&= \|r^2A_{m_0}v-r^2A_{m_0}\tilde{P}(y)\|_{L^{\infty}(B_{\eta})}\\
&\leq r^2A_{m_0}\eta^{2+\bar{\alpha}/2}\leq\eta^{2(m_0+1)}A_{m_0+1}.
  \end{aligned}
\end{equation*}
Hence, \cref{In-e.C2s-4-mu} holds for $m=m_0+1$. By induction, the proof is completed.\qed~\\

\begin{remark}\label{r-10.2}
From above proof we know that $\bar{\alpha}/2$ can be replaced by any $0<\alpha_0<\bar{\alpha}$ in the expression of
$\omega_u$ since we can define
\begin{equation*}
A_m=\max(\tilde\omega(\eta^{m}),\eta^{\alpha_0} A_{m-1})
\end{equation*}
in \cref{In-e.tC2-Ak}.
\end{remark}
~\\

Now, we give the~\\
\noindent\textbf{Proof of \Cref{t-C2-i}.} As before, we prove the theorem in two cases.

\textbf{Case 1:} the general case, i.e., $F$ satisfies \cref{SC2} and \cref{e.C2a-KF}. Throughout the proof for this case, $C$ always denotes a constant depending only on $n, \lambda,\Lambda,\mu,b_0,c_0,\omega_0$, $\|\beta_2\|_{C^{\mathrm{Dini}}(0)}$, $\omega_{\beta_2}$, $\omega_2$, $\|f\|_{C^{\mathrm{Dini}}(0)}$ and $\|u\|_{L^{\infty }(B_1)}$.

Since $\omega_{\beta_2}$ is a Dini function, there exists $0<\rho_1<1$ depending only on $\|\beta_2\|_{C^{\mathrm{Dini}}(0)}$ and $\omega_{\beta_2}$ such that
\begin{equation*}
\|\beta_2\|_{L^{\infty }(B_{\rho_1})}\leq \delta_0,
\end{equation*}
where $\delta_0$ is as in \Cref{t-C1a-i} with $\alpha=\bar{\alpha}/2$ there. By \Cref{t-C1a-i}, $u\in C^{1,\bar{\alpha}/2}(0)$ and
\begin{equation*}
\|u\|_{ C^{1,\bar{\alpha}/2}(0)}\leq C.
\end{equation*}

Similar to the previous proof, for $\tau\in \mathbb{R}$, let
\begin{equation*}
  \begin{aligned}
&u_1(x)=u(x)-P_u(x)-\tau x_n^2,\quad f_1(x)=f(x)-f(0),~\\
&F_1(M,p,s,x)=F(M+2\tau\tilde{I},p+DP_u+2\tau x_n,s+P_u(0)+\tau x_n^2,x)-f(0).
  \end{aligned}
\end{equation*}
Recall that $\tilde{I}$ denotes the matrix whose entries are all $0$ except $\tilde{I}_{nn}=1$ (see \Cref{no1.1}). Then
\begin{equation*}
u_1(0)=|Du_1(0)|=0,\quad
|f_1(x)|\leq [f]_{C^{\mathrm{Dini}}(0)}\cdot\frac{\omega_f(|x|)}{J_{\omega_f}}, ~~\forall ~x\in B_1
\end{equation*}
and $u_1$ is a viscosity solution of
\begin{equation*}
F_1(D^2u_1,Du_1,u_1,x)=f_1\quad\mbox{in}~~B_1.
\end{equation*}

As before, define the fully nonlinear operators $G_{1}$ similarly. By the structure condition, there exists $\tau\in \mathbb{R}$ such that $G_{1}(0,0,0)=0$ and
\begin{equation*}
  \begin{aligned}
|\tau|&\leq |G(0,Du(0),u(0))-f(0)|/\lambda\leq C.
  \end{aligned}
\end{equation*}

For $0<\rho\leq \rho_1$, define $y=x/\rho$ and
\begin{equation*}
  \begin{aligned}
&u_2(y)=\rho^{-1}u_1(x),\quad F_2(M,p,s,y)=\rho F_1\left(\rho^{-1}M, p, \rho s,x\right),\\
&G_{2}(M,p,s)=\rho G_{1}\left(\rho^{-1}M, p,\rho s\right).
  \end{aligned}
\end{equation*}
Then $u_2$ satisfies
\begin{equation}\label{In-C2-F5-mu}
F_2(D^2u_2,Du_2,u_2,y)=f_2\quad\mbox{in}~~B_1,
\end{equation}
where $f_2(y)=\rho f_1(x)$.

Now, we can check that \cref{In-C2-F5-mu} satisfies the conditions of \Cref{In-t-C2s-mu} by choosing a proper $\rho$. First, it can be verified easily that
\begin{equation*}
  \begin{aligned}
    u_2(0)=|Du_2(0)|=0
  \end{aligned}
\end{equation*}
and
\begin{equation*}
|f_2(y)|= \rho|f_1(x)|\leq C\rho\cdot\frac{\omega_f(\rho|y|)}{J_{\omega_f}}
\coloneqq C\rho \omega_{f_2}(|y|),~\forall ~y\in B_1.
\end{equation*}

Next, by the interior $C^{1,\bar{\alpha}/2}$ regularity for $u$,
\begin{equation*}
  \|u_2\|_{L^{\infty}(B_1)}\leq \rho^{-1} \|u_1\|_{L^{\infty}(B_{\rho})}
  \leq \rho^{-1}\left(C\rho^{1+\bar{\alpha}/2}+C\rho^{2}\right)\leq C\rho^{\bar{\alpha}/2}.
\end{equation*}
Furthermore, $G_{2}(0,0,0)=\rho G_{1}(0,0,0)=0$ and $F_2,G_{2}$ satisfy the structure condition \cref{SC2} with $\lambda,\Lambda,\tilde{\mu},\tilde{b},\tilde{c}$ and $\tilde\omega_0$, where
\begin{equation*}
\tilde{\mu}= \rho\mu,\quad\tilde{b}= \rho b_0+C\rho\mu,\quad\tilde{c}= \rho^{1/2}
c_0,\quad \tilde{\omega}_0(K,s)=\rho^{1/2}\omega_0(K+C,\rho s), ~\forall ~K,s>0.
\end{equation*}
Thus, $\tilde{\omega}_0$ satisfies \cref{e.omega0-c-1} with
\begin{equation*}
\hat{\tilde\omega}_0(K)=\rho^{1/4}\hat{\omega}_0(K+C), \quad
\check{\tilde \omega}_{0}(s)=\rho^{1/4}\check{\omega}_{0}(\rho s), ~\forall ~K,s>0.
\end{equation*}

Finally, it can be checked as before (cf. \cref{e5.8} in the proof of \Cref{t-C2a-i}) that
\begin{equation*}
  \begin{aligned}
|F_2(&M,p,s,y)-G_{2}(M,p,s)|\\
\leq& C\omega_{\beta_2}(|x|)(|M|+1)\omega_2(|p|+C,|s|+C)+C\rho(|p|+1)|x|
+\rho\hat\omega_0(|s|+C)\check{\omega}_{0}(C|x|).
  \end{aligned}
\end{equation*}
Define
\begin{equation}\label{e9.1}
  \begin{aligned}
&C_{\rho}=\int_{0}^{3\rho}\frac{\omega_{\beta_2}(\tau)}{\tau}d\tau, \quad
\omega_{\tilde{\beta}_2}(|y|)=\frac{\omega_{\beta_2}(\rho|y|)}{2C_{\rho}}
+\rho\check{\omega}_{0}(C\rho|y|), \quad \tilde{\beta}_2(y)=CC_{\rho}\omega_{\tilde{\beta}_2}(|y|),\\
&\tilde{\omega}_2\left(|p|,|s|\right)=\omega_2\left(|p|+C,|s|+C\right)
+C|p|+C+\hat\omega_0(|s|+C).
  \end{aligned}
\end{equation}
Then
\begin{equation*}
  \begin{aligned}
|F_2(M,p,s,y)-G_{3}(M,p,s)|\leq \tilde{\beta}_2(y)(|M|+1)\tilde{\omega}_2\left(|p|,|s|\right).
  \end{aligned}
\end{equation*}

Since $\omega_{\beta_2}$ is a Dini function,
\begin{equation*}
C_{\rho}\rightarrow 0~~\mbox{ as }~~\rho\rightarrow 0
\end{equation*}
and
\begin{equation*}
\int_{0}^{1}\frac{\omega_{\tilde{\beta}_2}(\tau)}{\tau}d\tau
=\frac{1}{2C_{\rho}}\int_{0}^{\rho}\frac{\omega_{\beta_2}(\tau)}{\tau}d\tau
+\rho\int_{0}^{C\rho}\frac{\check{\omega}_{0}(\tau)}{\tau}d\tau
\leq \frac{1}{2}+\rho\int_{0}^{C\rho}\frac{\check{\omega}_{0}(\tau)}{\tau}d\tau.
\end{equation*}
Moreover,
\begin{equation*}
\frac{1}{2}+\rho\int_{0}^{3C\rho}\frac{\check{\omega}_{0}(\tau)}{\tau}d\tau
\geq\int_{0}^{3}\frac{\omega_{\tilde{\beta}_2}(\tau)}{\tau}d\tau
\geq \int_{1}^{3}\frac{\omega_{\tilde{\beta}_2}(\tau)}{\tau}d\tau
\geq \omega_{\tilde{\beta}_2}(1)\int_{1}^{3}\frac{d\tau}{\tau}
= \omega_{\tilde{\beta}_2}(1)\ln 3.
\end{equation*}
Thus,
\begin{equation*}
\omega_{\tilde{\beta}_2}(1)\leq \frac{1}{2\ln 3}
+\frac{\rho}{\ln 3}\int_{0}^{3C\rho}\frac{\check{\omega}_{0}(\tau)}{\tau}d\tau.
\end{equation*}

Take $\delta_1$ small enough such that \Cref{In-t-C2s-mu} holds with $\tilde{\omega}_0,\tilde{\omega}_2$ and $\delta_1$. From above arguments, we can choose $\rho$ small enough (depending only on $n, \lambda,\Lambda,\mu,b_0,c_0,\omega_0$, $\|\beta_2\|_{C^{\mathrm{Dini}}(0)}$, $\omega_{\beta_2}$, $\omega_2$, $\|f\|_{C^{\mathrm{Dini}}(0)}$ and $\|u\|_{L^{\infty }(B_1)}$) such that
\begin{equation*}
\begin{aligned}
&\|u_2\|_{L^{\infty }(B_1)}\leq 1,\quad \tilde\mu \leq \frac{\delta_1}{4C_0},\quad
\tilde{b}\leq \frac{\delta_1}{2},\quad
\tilde c\leq \frac{\delta_1}{\omega_{\tilde\omega_0}(C_0)+1},\quad
\hat{\tilde{\omega}}_0(1+C_0)\leq 1,\\
&|\tilde\beta_2(y)|\leq \delta_1\omega_{\tilde{\beta}_2}(|y|),\quad |f_2(y)|\leq \delta_1\omega_{f_2}(|y|), ~~\forall ~y\in B_1,\\
&J_{\tilde\omega}\leq 1,~\quad(\tilde\omega\coloneqq \max(\check{\tilde\omega}_{0},\omega_{\tilde\beta_2},\omega_{f_2}))
\end{aligned}
\end{equation*}
where $C_0$ depending only on $n,\lambda$ and $\Lambda$, is as in \Cref{In-t-C2s-mu}. Therefore, the assumptions in \Cref{In-t-C2s-mu} are satisfied for \cref{In-C2-F5-mu}. By \Cref{In-t-C2s-mu}, $u_2\in C^{2}(0)$. By rescaling back to $u$, we conclude that $u\in C^{2}(0)$ and the estimates \crefrange{e.C2-1-i-mu}{e.C2-2-i-mu} hold. We point out here that the $\rho$ in \Cref{t-C2-i} is exactly the same as the one chosen in above argument.

\textbf{Case 2:} $F$ satisfies \cref{SC1} and \cref{e.C2a-KF-0}. As before (see the proofs of \Cref{t-C2a-i} and \Cref{t-C2a}), let $K=\|u\|_{L^{\infty }(B_1)}+\|f\|_{C^{\mathrm{Dini}}(0)}+\|\gamma_2\|_{C^{\mathrm{Dini}}(0)}$ and $u_1=u/K$. Hence, $u_1$ satisfies
\begin{equation}\label{e.10.4}
F_1(D^2u_1,Du_1,u_1,x)=f_1\quad\mbox{in}~~B_1,
\end{equation}
where $F_1(M,p,s,x)=F(KM,Kp,Ks,x)/K$ and $f_1=f/K$. Then by applying \textbf{Case 1} to \cref{e.10.4}, we obtain that $u_1$ and hence $u$ is $C^2$ at $0$, and the estimates \crefrange{e.C2-1-i}{e.C2-2-i} hold. \qed
~\\

In the second half of this section, we prove the boundary $C^{1}$ regularity.

\noindent\textbf{Proof of \Cref{t-C1}.} As before, we assume that $P_g\equiv 0$. Let $\tilde\alpha=\min(\bar{\alpha}/2,1-n/p)$ and $\delta$ be as in \Cref{l-C1a-mu}, which depends only on $n,\lambda,\Lambda$ and $p$. We also assume that
\begin{equation}\label{e.tC1-ass}
  \begin{aligned}
&\|u\|_{L^{\infty}(\Omega_1)}\leq 1,\quad \mu\leq \frac{\delta}{6C_0^2},\quad \|b\|_{L^p(\Omega_1)}\leq \frac{\delta}{3C_0},\\
&\|f\|_{L^n(\Omega_r)}\leq  \frac{\delta}{3}\omega_f(r), \quad
\|g\|_{L^{\infty}((\partial\Omega)_r)}\leq \frac{\delta}{2}r\omega_g(r), ~\forall ~0<r<1,\\
&\underset{B_{r}}{\mathrm{osc}}~\partial\Omega\leq \frac{\delta}{2C_0}r\omega_{\Omega}(r), ~\forall ~0<r<1,\\
&J_{\tilde\omega}\leq 2,~\quad(\tilde\omega\coloneqq\max(\omega_f,\omega_{g},\omega_{\Omega})),
  \end{aligned}
\end{equation}
where $C_0>1$ is a constant (depending only on $n,\lambda,\Lambda$ and $p$) to be specified later.

 Otherwise, note that $\Omega$ satisfies the exterior cone condition at $0$ and we may consider for $0<\rho<1$,
\begin{equation*}
  \bar{u}(y)=\frac{u(x)}{\rho^{\alpha_0}},
\end{equation*}
where $y=x/\rho$ and $0<\alpha_0<1$ is a H\"{o}lder exponent (depending only on $n,\lambda,\Lambda,\|b\|_{L^p(B_1)}$ and $\|(\partial \Omega)_1\|_{C^{1,\mathrm{Dini}}(0)}$)
such that $u\in C^{2\alpha_0}(0)$ (by \Cref{l-3Ho}). Then we have
\begin{equation*}
\left\{\begin{aligned}
&\bar{u}\in S^{*}(\lambda,\Lambda,\bar{\mu},\bar{b},\bar{f})&&
\quad\mbox{in}~~\tilde{\Omega}\cap B_1;\\
&\bar{u}=\bar{g}&& \quad\mbox{on}~~\partial \tilde{\Omega}\cap B_1,
\end{aligned}\right.
\end{equation*}
where
\begin{equation*}\label{e.tC1-n1-mu}
  \bar{\mu}=\rho^{\alpha_0}\mu,\quad\bar{b}(y)=\rho b(x),\quad\bar{f}(y)=\rho^{2-\alpha_0}f(x),
\quad\bar{g}(y)=\rho^{-\alpha_0}g(x), \quad\tilde{\Omega}=\rho^{-1}\Omega.
\end{equation*}
Hence,
\begin{equation}\label{e9.2}
  \begin{aligned}
\|\bar{u}\|_{L^{\infty}(\tilde{\Omega}_1)}\leq& \rho^{\alpha_0}[u]_{C^{2\alpha_0}(0)},\quad \bar{\mu}=\rho^{\alpha_0}\mu,\\
\|\bar{b}\|_{L^{p}(\tilde{\Omega}_1)}=&\rho^{1-\frac{n}{p}}\|b\|_{L^{p}(\Omega_\rho)}\leq \rho^{\tilde\alpha}\|b\|_{L^{p}(\Omega_1)},\\
\|\bar{f}\|_{L^{n}(\tilde{\Omega}_r)}=&\rho^{1-\alpha_0}\|f\|_{L^{n}(\Omega_{\rho r})}
\leq \rho^{1-\alpha_0}\|f\|_{C^{-1,\mathrm{Dini}}(0)}\frac{\omega_f(\rho r)}{J_{\omega_f}}\\
\coloneqq& \rho^{1-\alpha_0}\|f\|_{C^{-1,\mathrm{Dini}}(0)}\omega_{\bar{f}}(r),\\
\|\bar{g}\|_{L^{\infty}((\partial\tilde{\Omega})_r)}
=&\rho^{-\alpha_0}\|g\|_{L^{\infty}((\partial\Omega)_{\rho r})}\leq
\rho^{1-\alpha_0}[g]_{C^{1,\mathrm{Dini}}(0)}r\frac{\omega_g(\rho r)}{J_{\omega_g}}\\
\coloneqq& \rho^{1-\alpha_0}[g]_{C^{1,\mathrm{Dini}}(0)}r\omega_{\bar{g}}(r),\\
\underset{B_{r}}{\mathrm{osc}}~\partial\tilde{\Omega}
=&\rho^{-1}\underset{B_{\rho r}}{\mathrm{osc}}~\partial\Omega
\leq Kr\omega_{\Omega}(\rho r)
=KC_{\rho}r\frac{\omega_{\Omega}(\rho r)}{C_{\rho}}\coloneqq KC_{\rho}r\omega_{\tilde\Omega}(r),
  \end{aligned}
\end{equation}
where
\begin{equation*}
K=J_{\omega_{\Omega}}^{-1}\|(\partial \Omega)_1\|_{C^{1,\mathrm{Dini}}(0)}, \quad C_{\rho}=\int_{0}^{3\rho}\frac{\omega_{\Omega}(r)}{r}dr.
\end{equation*}

Obviously,
\begin{equation*}
J_{\omega_{\bar{f}}}\leq1, \quad J_{\omega_{\bar{g}}}\leq 1.
\end{equation*}
In addition, by noting
\begin{equation*}
C_{\rho}=\int_{0}^{3\rho}\frac{\omega(r)}{r}dr\geq \int_{\rho}^{3\rho}\frac{\omega(r)}{r}dr
\geq \ln 3 \omega(\rho),
\end{equation*}
we have
\begin{equation*}
J_{\omega_{\tilde{\Omega}}}
=\int_{0}^{1}\frac{\omega_{\tilde{\Omega}}(r)}{r}dr+\omega_{\tilde{\Omega}}(1)
=\frac{1}{C_{\rho}}\int_{0}^{\rho}\frac{\omega_{\Omega}(r)}{r}dr
+\frac{\omega_{\Omega}(\rho)}{C_{\rho}}\leq 1+\frac{1}{\ln 3}\leq 2.
\end{equation*}
Then by choosing $\rho$ small enough (depending only on $n$, $\lambda$, $\Lambda$, $p$, $\mu$, $\|b\|_{L^{p}(\Omega_1)}$, $\|(\partial \Omega)_1\|_{C^{1,\mathrm{Dini}}(0)}$, $\omega_{\Omega}$, $\|f\|_{C^{-1,\mathrm{Dini}}(0)}$, $[g]_{C^{1,\mathrm{Dini}}(0)}$ and $\|u\|_{L^{\infty}(\Omega_1)}$), the assumptions \cref{e.tC1-ass} for $\bar{u}$ can be guaranteed. Hence, we can make the assumption \cref{e.tC1-ass} for $u$ without loss of generality.

Now we prove that $u$ is $C^{1}$ at $0$ and we only need to prove the following. There exists a sequence $a_m$ ($m\geq -1$) such that for all $m\geq 0$,
\begin{equation}\label{e.tC1-3}
\|u-a_mx_n\|_{L^{\infty }(\Omega _{\eta^{m}})}\leq \eta^m A_m
\end{equation}
and
\begin{equation}\label{e.tC1-4}
|a_m-a_{m-1}|\leq \bar{C} A_{m-1},
\end{equation}
where
\begin{equation}\label{e.tC1-Ak}
A_{-1}=A_0=1,A_m=\max(\tilde\omega(\eta^{m}),\eta^{\tilde\alpha} A_{m-1}) (m\geq 1)
\end{equation}
and $\eta$, depending only on $n,\lambda$ and $\Lambda$, is as in \Cref{l-C1a-mu} (take $\alpha=\tilde{\alpha}$ there).

We prove the above by induction. For $m=0$, by setting $a_0=a_{-1}=0$, the conclusion holds clearly. Suppose that the conclusion holds for $m\leq m_0$. We need to prove that the conclusion holds for $m=m_0+1$.

Let $r=\eta ^{m_{0}}$, $y=x/r$ and
\begin{equation}\label{e.tC1-v1}
  v(y)=\frac{u(x)-a_{m_0}x_n}{rA_{m_0}}.
\end{equation}
Then $v$ satisfies
\begin{equation*}
\left\{\begin{aligned}
&v\in S^*(\lambda,\Lambda,\tilde{\mu},\tilde{b},\tilde{f})&& \quad\mbox{in}~~\tilde{\Omega}\cap B_1;\\
&v=\tilde{g}&& \quad\mbox{on}~~\partial \tilde{\Omega}\cap B_1,
\end{aligned}\right.
\end{equation*}
where
\begin{equation*}
\begin{aligned}
  &\tilde{\mu}=2rA_{m_0}\mu,\quad \tilde{b}(y)=rb(x),\quad
  \tilde{f}(y)=A_{m_0}^{-1}\left(r|f(x)|+rb(x)|a_{m_0}|+2r\mu|a_{m_0}|^2\right),\\
  &\tilde{g}(y)=(rA_{m_0})^{-1}\left(g(x)-a_{m_0}x_n\right),\quad\tilde{\Omega}=r^{-1}\Omega.
\end{aligned}
\end{equation*}

From the definition of $A_{m_0}$ (see \cref{e.tC1-Ak}),
\begin{equation}\label{e.C1-4}
r^{\tilde{\alpha}}=\eta^{m_0\tilde{\alpha}}\leq A_{m_0}\leq 1.
\end{equation}
By \cref{e.tC1-4}, there exists a constant $C_0$ depending only on $n,\lambda,\Lambda$ and $p$ such that $|a_{m}|\leq C_0$ ($\forall~0\leq m\leq m_0$). Then it is easy to verify that
\begin{equation}\label{e.tC1-esti}
  \begin{aligned}
\|v\|_{L^{\infty}(\tilde{\Omega}_1)}\leq& 1, ~\quad(\mathrm{by}~ \cref{e.tC1-3} ~\mathrm{and}~ \cref{e.tC1-v1})\\
\tilde\mu\leq& 2rA_{m_0}\mu\leq \delta, ~\quad(\mathrm{by}~ \cref{e.tC1-ass})\\
\|\tilde{b}\|_{L^{p}(\tilde{\Omega}_1)}=&r^{1-\frac{n}{p}}\|b\|_{L^{p}(\Omega_r)}\leq \delta,~\quad(\mathrm{by}~ \cref{e.tC1-ass})\\
\|\tilde{f}\|_{L^{n}(\tilde{\Omega}_1)}\leq& A_{m_0}^{-1}\left(\|f\|_{L^{n}(\Omega_r)}
    +C_0r^{1-\frac{n}{p}}\|b\|_{L^{p}(\Omega_r)}+2C_0^2r\mu\right)\\
\leq& A_{m_0}^{-1}\left(\frac{\delta\tilde\omega(r)}{3}
    +\frac{\delta r^{\tilde{\alpha}}}{3}+\frac{\delta r}{3}\right)\\
\leq& \delta, ~\quad(\mathrm{by} ~\cref{e.tC1-ass}~\mathrm{and}~\cref{e.tC1-Ak})\\
\|\tilde{g}\|_{L^{\infty}(\partial \tilde{\Omega}\cap B_1)}\leq& \frac{1}{rA_{m_0}}\left(\frac{\delta r\tilde\omega(r)}{2}+\frac{C_0\delta r\tilde\omega(r)}{2C_0}\right)\leq \delta,  ~\quad(\mathrm{by}~ \cref{e.tC1-ass})\\
\underset{B_1}{\mathrm{osc}}~\partial\tilde{\Omega}=&
r^{-1}\underset{B_r}{\mathrm{osc}}~\partial\Omega \leq \delta\tilde\omega(r) \leq \delta  ~\quad(\mathrm{by}~ \cref{e.tC1-ass}).
  \end{aligned}
\end{equation}

By \Cref{l-C1a-mu} (take $\alpha=\tilde{\alpha}$ there), there exists a constant $\tilde{a}$ such that
\begin{equation*}
\begin{aligned}
    \|v-\tilde{a}y_n\|_{L^{\infty }(\tilde{\Omega} _{\eta})}&\leq \eta ^{1+\tilde{\alpha}}
\end{aligned}
\end{equation*}
and
\begin{equation*}
|\tilde{a}|\leq \bar{C}.
\end{equation*}

Let $a_{m_0+1}=a_{m_0}+A_{m_0}\tilde{a}$. Then \cref{e.tC1-4} holds for $m_0+1$. By recalling \cref{e.tC1-Ak} and \cref{e.tC1-v1}, we have
\begin{equation*}
  \begin{aligned}
\|u-a_{m_0+1}x_n\|_{L^{\infty}(\Omega_{\eta^{m_0+1}})}
&= \|u-a_{m_0}x_n-A_{m_0}\tilde{a}x_n\|_{L^{\infty}(\Omega_{\eta r})}\\
&= \|rA_{m_0}v-rA_{m_0}\tilde{a}y_n\|_{L^{\infty}(\tilde{\Omega}_{\eta})}\\
&\leq rA_{m_0}\eta^{1+\tilde{\alpha}}\leq\eta^{m_0+1}A_{m_0+1}.
  \end{aligned}
\end{equation*}
Hence, \cref{e.tC1-3} holds for $m=m_0+1$. By induction, the proof is completed.

For the special case $\mu=0$, set
\begin{equation*}
K=\|u\|_{L^{\infty}(\Omega_1)}+\delta^{-1}\left(3\|f\|_{C^{-1,\mathrm{Dini}}(0)}
+2\|g\|_{C^{1,\mathrm{Dini}}(0)}\right)
\end{equation*}
and define for $0<\rho<1$
\begin{equation*}
  \bar{u}(y)=u(x)/K,
\end{equation*}
where $y=x/\rho$. Then by taking $\rho$ small enough (depending only on $n, \lambda, \Lambda,\mu$, $p,\|b\|_{L^p(\Omega_1)}$ and $\|(\partial \Omega)_1\|_{C^{1,\mathrm{Dini}}(0)}$), \cref{e.tC1-ass} can be guaranteed. Hence, for $\mu=0$, we have the explicit estimates \cref{e.C1-1} and \cref{e.C1-2}.\qed~\\

\section{The \texorpdfstring{$C^{k,\mathrm{lnL}}$}{Ck,lnL} regularity}\label{CLlnL}
In this section, we state a series of the so-called ``ln-Lipschitz'' regularity. Similar to the $C^k$ regularity, we only give the details of proofs for the interior pointwise $C^{k,\mathrm{lnL}}$ ($k\geq 2$) regularity and the boundary pointwise $C^{1,\mathrm{lnL}}$ regularity.

The first is the $C^{0,\mathrm{lnL}}$ regularity.
\begin{theorem}\label{t-C0ln-i}
Let $p>n$ and $u$ be a viscosity solution of
\begin{equation*}
F(D^2u, Du, u,x)=f \quad\mbox{in}~~ B_1.
\end{equation*}
Suppose that $F$ satisfies \cref{SC2}, \cref{e.C1a.beta} and
\begin{equation*}
b\in L^{p}(B_1),\quad c\in L^{n}(B_1),\quad \|\beta_1\|_{C^{-1,1}(0)}\leq \delta_0,\quad
\gamma_1\in L^{n}(B_1),\quad f\in L^{n}(B_1),
\end{equation*}
where $0<\delta_0<1$ depends only on $n,\lambda,\Lambda$ and $p$.

Then $u\in C^{0,\mathrm{lnL}}(0)$, i.e.,
\begin{equation}\label{e.C0ln-1-i-mu}
  |u(x)-u(0)|\leq C|x|\big|\ln |x|\big|, ~\forall ~x\in B_{1/2},
\end{equation}
where $C$ depends only on $n,\lambda,\Lambda,p,\mu,\|b\|_{L^p(B_1)}, \|c\|_{L^{n}(B_1)},\omega_0,\|\gamma_1\|_{L^{n}(B_1)},\|f\|_{L^{n}(B_1)}$ and $\|u\|_{L^{\infty }(B_1)}$.

In particular, if $F$ satisfies \cref{SC1}, we have the following explicit estimate
\begin{equation*}\label{e.C0ln-1-i}
  |u(x)-u(0)|\leq C |x|\big|\ln |x|\big|\left(\|u\|_{L^{\infty }(B_1)}+\|f\|_{L^{n}(B_1)}+\|\gamma_1\|_{L^{n}(B_1)}\right), ~\forall ~x\in B_{1/2},
\end{equation*}
where $C$ depends only on $n, \lambda, \Lambda,p,\|b\|_{L^p(B_1)}$ and $\|c\|_{L^{n}(B_1)}$.
\end{theorem}

\begin{remark}\label{r-2.20}
Usually, one obtains the $C^{\alpha}$ regularity under the condition $f\in L^n$ (e.g. \cite[Proposition 4.10]{MR1351007}). Teixeira \cite{MR3158810} proved the interior $C^{0,\mathrm{lnL}}$ regularity for fully nonlinear elliptic equations without lower order terms. Later, this was extended by da Silva and Teixeira \cite{MR3713553} to fully nonlinear parabolic equations. Recently, da Silva and Nornberg \cite{MR4304555} derived (for $c=0$) interior $C^{0,\mathrm{lnL}}$ regularity for equations with more general nonlinear growth in the gradient.

For the Poisson equation, $f\in L^n$ implies $u\in W^{2,n}$. Since $W^{2,n}\subset W^{1,\mathrm{BMO}}\subset C^{0,\mathrm{lnL}}$, \Cref{t-C0ln-i} can be regarded as a weaker version of $W^{2,n}$ and $W^{1,\mathrm{BMO}}$ regularity when $f\in L^n$.
\end{remark}
~\\

Next, we consider the $C^{1,\mathrm{lnL}}$ regularity. Note that for $C^{1,\mathrm{lnL}}$ and higher regularity, we always assume that \cref{SC2} holds with $b\equiv b_0$ and $c\equiv c_0 $ for positive constants $b_0,c_0$.
\begin{theorem}\label{t-C1ln-i}
Let $u$ be a viscosity solution of
\begin{equation*}
F(D^2u, Du, u,x)=f \quad\mbox{in}~~ B_1.
\end{equation*}
Suppose that $F$ satisfies \cref{SC2} and \cref{e.C1a.beta} with $G$ being convex in $M$. Assume that
\begin{equation*}
\|\beta_1\|_{L^{n}(B_r)}\leq \frac{\delta_0 r}{|\ln r|},~\forall ~0<r<1/2,\quad\gamma_1\in L^{\infty}(B_1),\quad f\in L^{\infty}(B_1),
\end{equation*}
where $0<\delta_0<1$ depends only on $n,\lambda$ and $\Lambda$.

Then $u$ is $C^{1,\mathrm{lnL}}$ at $0$, i.e., there exists $P\in \mathcal{P}_1$ such that
\begin{equation}\label{e.C1ln-1-i-mu}
  |u(x)-P(x)|\leq C |x|^2\big|\ln |x|\big|, ~~\forall ~x\in B_{1/2}
\end{equation}
and
\begin{equation}\label{e.C1ln-2-i-mu}
|Du(0)|\leq C ,
\end{equation}
where $C$ depends only on $n,\lambda,\Lambda,\mu,b_0, c_0,\omega_0,
\|\gamma_1\|_{L^{\infty}(B_1)}$, $\|f\|_{L^{\infty}(B_1)}$ and
$\|u\|_{L^{\infty }(B_1)}$.

In particular, if $F$ satisfies \cref{SC1},
\begin{equation*}\label{e.C1ln-1-i}
  |u(x)-P(x)|\leq C |x|^2\big|\ln |x|\big|\left(\|u\|_{L^{\infty }(B_1)}
  +\|f\|_{L^{\infty}(B_1)}+\|\gamma_1\|_{L^{\infty}(B_1)}\right), ~~\forall ~x\in  B_{1/2}
\end{equation*}
and
\begin{equation*}\label{e.C1ln-2-i}
|Du(0)|\leq C \left(\|u\|_{L^{\infty }(B_1)}
+\|f\|_{L^{\infty}(B_1)}+\|\gamma_1\|_{L^{\infty}(B_1)}\right),
\end{equation*}
where $C$ depends only on $n,\lambda,\Lambda,b_0$ and $c_0$.
\end{theorem}
\begin{remark}\label{r-2.3}
Besides the $C^{0,\mathrm{lnL}}$ regularity mentioned in \Cref{r-2.20}, Teixeira \cite{MR3158810} and da Silva, Teixeira \cite{MR3713553} obtained the interior $C^{1,\mathrm{lnL}}$ regularity for elliptic and parabolic equations respectively as well. da Silva and Nornberg \cite{MR4304555} also derived the interior $C^{1,\mathrm{lnL}}$ regularity. In addition, Silvestre and Teixeira \cite[Theorem 1.4]{MR3494911} proved $C^{1,\mathrm{lnL}}$ regularity if the recession function has a priori $C^{2,\bar{\alpha}}$ interior estimates, which corresponds to that $G$ is convex in $M$ in \Cref{t-C1ln-i}.

It is well-known that $f\in L^{\infty}$ does not imply $u\in C^{1,1}$ (see \cite[Example P. 65]{MR1669352} for instance). \Cref{t-C1ln-i} can be regarded as a substitute for the $C^{1,1}$ regularity. Interestingly, if $f$ and $\beta$ satisfy the same conditions as in \Cref{t-C1ln-i}, Caffarelli and Huang \cite{MR1978880} proved the $W^{2,\mathrm{BMO}}$ regularity, which is stronger than $C^{1,\mathrm{lnL}}$ regularity.
\end{remark}
~\\

\begin{theorem}\label{t-Ckln-i}
Let $k\geq 2$ and $u$ be a viscosity solution of
\begin{equation*}
F(D^2u, Du, u,x)=f \quad\mbox{in}~~ B_1.
\end{equation*}
Suppose that $F\in C^{k-2,1}(0)$, $\omega_0$ satisfies \cref{e.omega0-2} and $f\in C^{k-2,1}(0)$.

Then $u\in C^{k,\mathrm{lnL}}(0)$, i.e., there exists $P\in \mathcal{P}_k$ such that
\begin{equation}\label{e.Ckln-1-i}
  |u(x)-P(x)|\leq C|x|^{k+1}\big|\ln|x|\big|, ~~\forall ~x\in B_{1/2},
\end{equation}
\begin{equation}\label{e.Ckln-3-i}
|F(D^2P(x),DP(x),P(x),x)-P_f(x)|\leq C|x|^{k-1}, ~~\forall ~x\in B_{1/2}
\end{equation}
and
\begin{equation}\label{e.Ckln-2-i}
|Du(0)|+\cdots+|D^ku(0)|\leq C,
\end{equation}
where $C$ depends only on $k,n,\lambda, \Lambda$, $\mu, b_0,c_0,\omega_0,\|F\|_{C^{k-2,1}(0)},\omega_3,\omega_4$, $\|f\|_{C^{k-2,1}(0)}$ and $\|u\|_{L^{\infty}(B_1)}$.
\end{theorem}

\begin{remark}\label{r-2.22}
In \cite{MR2273802}, Wang proved the $C^{2,\mathrm{lnL}}$ regularity for $C^{2}$ solutions and $F\in C^{1,1}$. This result has been extended to the Monge-Amp\`{e}re equation by Jian and Wang \cite{MR2338423}. Up to our knowledge, \Cref{t-Ckln-i} is the first pointwise $C^{k,\mathrm{lnL}}$ ($k\geq 2$) regularity, even for linear equations.
\end{remark}
~\\

We also obtain the boundary pointwise $C^{0,\mathrm{lnL}}$, $C^{1,\mathrm{lnL}}$ and $C^{k,\mathrm{lnL}}$ regularity. To the best of our knowledge, these regularity are new even for the Laplace equation. We list them as follows.

\begin{theorem}\label{t-C0ln-mu}
Let $p>n$ and $u$ satisfy
\begin{equation*}
\left\{\begin{aligned}
&u\in S^*(\lambda,\Lambda,\mu,b,f)&& \quad\mbox{in}~~\Omega\cap B_1;\\
&u=g&& \quad\mbox{on}~~\partial \Omega\cap B_1.
\end{aligned}\right.
\end{equation*}
Suppose that
\begin{equation*}
b\in L^{p}(\Omega\cap B_1),\quad f\in L^{n}(\Omega\cap B_1),\quad g\in C^{0,1}(0),\quad
\underset{B_{r}}{\mathrm{osc}}~\partial\Omega\leq  \frac{\delta_0 r}{|\ln r|},~\forall ~0<r<1/2,
\end{equation*}
where $0<\delta_0<1$ depends only on $n,\lambda,\Lambda$ and $p$.

Then $u\in C^{0,\mathrm{lnL}}(0)$, i.e.,
\begin{equation}\label{e.C0ln-1-mu}
  |u(x)-u(0)|\leq C |x|\big|\ln |x|\big|, ~\forall ~x\in \Omega\cap B_{1/2},
\end{equation}
where $C$ depends only on $n, \lambda, \Lambda,p,\mu$, $\|b\|_{L^p(\Omega\cap B_1)}$, $\|f\|_{L^{n}(\Omega\cap B_1)}$, $\|g\|_{C^{0,1}(0)}$ and $\|u\|_{L^{\infty }(\Omega\cap B_1)}$.

In particular, if $\mu=0$,
\begin{equation*}\label{e.C0ln-1}
  |u(x)-u(0)|\leq C |x|\big|\ln |x|\big|\left(\|u\|_{L^{\infty }(\Omega\cap B_1)}+\|f\|_{L^{n}(\Omega\cap B_1)}+\|g\|_{C^{0,1}(0)}\right), ~\forall ~x\in \Omega\cap B_{1/2},
\end{equation*}
where $C$ depends only on $n,\lambda, \Lambda,p$ and $\|b\|_{L^p(\Omega\cap B_1)}$.
\end{theorem}

\begin{remark}\label{r-2.21}
Note that the condition on $\partial \Omega$ for $C^{0,\mathrm{lnL}}$ regularity is similar to the condition on $\beta$ for $C^{1,\mathrm{lnL}}$ regularity. For $C^{k}$ regularity in last subsection, there is also a similarity between them. We make a summary here.

The $\beta$ is used to characterize the oscillation of the coefficient of second order term of the equation and then the $C^2$ regularity is a critical case for $\beta$. If we intend to obtain higher regularity than $C^2$ (e.g. $C^{2}$, $C^{k,\alpha}$ ($k\geq 2$)), we need that $\beta$ has a decay (e.g. $\beta\in C^\mathrm{Dini}$, $\beta\in C^{k-2,\alpha}$). Then by a normalization procedure, we can assume that $\beta$ is small (cf. \cref{e5.9}).

On the other hand, if we intend to obtain lower regularity than $C^2$ (e.g. $C^{0,\mathrm{lnL}}$, $C^{1,\alpha}$), we do not assume that $\beta$ has a decay since it is not necessary. Thus, we cannot make $\beta$ small by normalization and we must make the assumption that $\beta$ is small in the theorem (e.g. $\|\beta\|_{L^n(B_r)}\leq \delta_0r/|\ln r|$, $\|\beta\|_{L^n(B_r)}\leq \delta_0 r$ where $\delta_0$ is a small constant).

For $\partial \Omega$, we have a similar explanation. The key to the proof of boundary regularity is to estimate $x_n$ on $\partial \Omega$ (cf. the proof in \Cref{C1a-mu}). Hence, $C^1$ regularity is critical for $\partial \Omega$. If we intend to obtain higher regularity than $C^1$ (e.g. $C^{1}$, $C^{k,\alpha}$ ($k\geq 1$)), we need that $\mathrm{osc}~\partial \Omega$ has a decay (e.g. $\partial \Omega\in C^{1,\mathrm{Dini}}$, $\partial \Omega\in C^{k,\alpha}$). Then by a normalization procedure, we can assume that $\mathrm{osc}~\partial \Omega$ is small (cf. \cref{e7.1}).

However, if we intend to obtain lower regularity than $C^1$ (e.g. $C^{0,\mathrm{lnL}}$), that $\mathrm{osc}~\partial \Omega$ has a decay is not necessary. Then we cannot make $\mathrm{osc}~\partial \Omega$ small by normalization and we must make the assumption that $\mathrm{osc}~\partial \Omega$ is small in the theorem (e.g. $\underset{B_{r}}{\mathrm{osc}}~\partial\Omega\leq \delta_0r/|\ln r|$ where $\delta_0$ is a small constant).
\end{remark}
~\\

\begin{theorem}\label{t-C1ln}
Let $u$ be a viscosity solution of
\begin{equation*}
\left\{\begin{aligned}
&F(D^2u,Du,u,x)=f&& \quad\mbox{in}~~\Omega\cap B_1;\\
&u=g&& \quad\mbox{on}~~\partial\Omega\cap B_1.
\end{aligned}\right.
\end{equation*}
Suppose that $F$ satisfies \cref{SC2} and \cref{e.C1a.beta} with $G$ being convex in $M$. Assume that
\begin{equation}\label{e11.1}
  \begin{aligned}
&\|\beta_1\|_{L^{n}(B_r)}\leq \frac{\delta_0 r}{|\ln r|},~\forall ~0<r<1/2,\quad \gamma_1\in L^{\infty}(\Omega\cap B_1),\\
&f\in L^{\infty}(\Omega\cap B_1),\quad \partial\Omega\cap B_1\in C^{1,1} (0),\quad g\in C^{1,1}(0),
  \end{aligned}
\end{equation}
where $0<\delta_0<1$ depends only on $n,\lambda$ and $\Lambda$.

Then $u$ is $C^{1,\mathrm{lnL}}$ at $0$, i.e., there exists $P\in \mathcal{P}_1$ such that
\begin{equation}\label{e.C1ln-1-mu}
  |u(x)-P(x)|\leq C |x|^2\big|\ln |x|\big|, ~~\forall ~x\in \Omega\cap B_{1/2},
\end{equation}
\begin{equation}\label{e.C1ln-3-mu}
D_{x'}u(0)=D_{x'}g(0)
\end{equation}
and
\begin{equation}\label{e.C1ln-2-mu}
|Du(0)|\leq C ,
\end{equation}
where $C$ depends only on $n$, $\lambda$, $\Lambda$, $\mu$, $b_0$, $c_0$, $\omega_0$,
$\|\gamma_1\|_{L^{\infty}(\Omega\cap B_1)}$,
$\|\partial \Omega\cap B_1\|_{C^{1,1}(0)}$, $\|f\|_{L^{\infty}(\Omega\cap B_1)}$, $\|g\|_{C^{1,1}(0)}$ and $\|u\|_{L^{\infty }(\Omega\cap B_1)}$.

In particular, if $F$ satisfies \cref{SC1},
\begin{equation}\label{e.C1ln-1}
\begin{aligned}
  |&u(x)-P(x)|\leq \tilde C |x|^2\big|\ln |x|\big|, ~~\forall ~x\in \Omega\cap B_{1/2},
\end{aligned}
\end{equation}
\begin{equation}\label{e.C1ln-2}
|Du(0)|\leq \tilde C ,
\end{equation}
and
\begin{equation*}
\tilde{C}= C\left(\|u\|_{L^{\infty }(\Omega\cap B_1)}
  +\|f\|_{L^{\infty}(\Omega\cap B_1)}+\|\gamma_1\|_{L^{\infty}(\Omega\cap B_1)}
  +\|g\|_{C^{1,1}(0)}\right),
\end{equation*}
where $C$ depends only on $n,\lambda,\Lambda,b_0, c_0$ and
$\|\partial \Omega\cap B_1\|_{C^{1,1}(0)}$.
\end{theorem}

\begin{theorem}\label{t-Ckln}
Let $k\geq 2$ and $u$ be a viscosity solution of
\begin{equation*}
\left\{\begin{aligned}
&F(D^2u,Du,u,x)=f&& \quad\mbox{in}~~\Omega\cap B_1;\\
&u=g&& \quad\mbox{on}~~\partial \Omega\cap B_1.
\end{aligned}\right.
\end{equation*}
Suppose that \cref{SC2} and \cref{e.omega0-2} hold, and
\begin{equation*}
F\in C^{k-2,1}(0),\quad f\in C^{k-2,1}(0),\quad \partial\Omega\cap B_1\in C^{k,1} (0),\quad g\in C^{k,1}(0).
\end{equation*}

Then $u\in C^{k,\mathrm{lnL}}(0)$, i.e., there exists $P\in \mathcal{P}_k$ such that
\begin{equation}\label{e.Ckln-1}
  |u(x)-P(x)|\leq C |x|^{k+1}\big|\ln|x|\big|, ~~\forall ~x\in \Omega\cap B_{1/2},
\end{equation}
\begin{equation}\label{e.Ckln-3}
|F(D^2P(x),DP(x),P(x),x)-P_f(x)|\leq C|x|^{k-1}, ~~\forall ~x\in \Omega\cap B_{1/2},
\end{equation}
\begin{equation}\label{e.Ckln-4}
D^l_{x'}u(x',P_{\Omega}(x'))=D^l_{x'}g(x',P_{\Omega}(x'))~\mbox{at}~0,~\forall ~0\leq l\leq k
\end{equation}
and
\begin{equation}\label{e.Ckln-2}
|Du(0)|+\cdots+|D^ku(0)|\leq C,
\end{equation}
where $C$ depends only on $k,n,\lambda, \Lambda$, $\mu, b_0,c_0,\omega_0,\|F\|_{C^{k-2,1}(0)},
\omega_3,\omega_4$, $\|\partial \Omega\cap B_1\|_{C^{k,1}(0)}$, $\|f\|_{C^{k-2,1}(0)}$,
$\|g\|_{C^{k,1}(0)}$ and $\|u\|_{L^{\infty}(\Omega\cap B_1)}$.
\end{theorem}

Similar to $C^{k,\alpha}$ regularity, by combining the interior and boundary regularity together, we have the local and global $C^{k,\mathrm{lnL}}$ regularity. For the completeness and the convenience of citation, we list them as follows.

\begin{corollary}\label{t-C0ln-global}
Let $p>n$, $\Gamma\subset \partial \Omega $ be relatively open (may be empty) and $u$ be a viscosity solution of
\begin{equation*}
\left\{\begin{aligned}
&F(D^2u,Du,u,x)=f&& \quad\mbox{in}~~\Omega;\\
&u=g&& \quad\mbox{on}~~\Gamma.
\end{aligned}\right.
\end{equation*}
Suppose that \cref{SC2} holds and $F$ satisfies \cref{e.C1a.beta} with some $G_{x_0}$ at any $x_0\in \Omega\cup \Gamma$ and
\begin{equation*}
  \begin{aligned}
&b\in L^{p}(\Omega),\quad c\in L^{n}(\Omega),\quad \beta_1(x,x_0)\leq \delta_0,~\forall ~x\in B_{r_0}(x_0)\cap \Omega,x_0\in \Omega\cup \Gamma,\\
&\gamma_1\in L^{n}(\Omega),\quad f\in L^{n}(\Omega),\quad \underset{B_{r}(x_0)}{\mathrm{osc}}~\partial\Omega\leq \frac{\delta_0r}{|\ln r|},~\forall ~x_0\in \Gamma, 0<r<r_0,\quad g\in C^{0,1}(\bar\Gamma),
  \end{aligned}
\end{equation*}
where $0<\delta_0<1$ depends only on $n,\lambda,\Lambda$ and $p$.

Then for any $\Omega'\subset\subset \Omega\cup \Gamma$, we have $u\in C^{0,\mathrm{lnL}}(\bar{\Omega}')$ and
\begin{equation}\label{e.C0ln-1-i-mu-global}
 \|u\|_{C^{0,\mathrm{lnL}}(\bar{\Omega}')}\leq C,
\end{equation}
where $C$ depends only on $n,\lambda,\Lambda,p,r_0,\mu,\|b\|_{L^p(\Omega)}, \|c\|_{L^{n}(\Omega)},\omega_0,\|\gamma_1\|_{L^{n}(\Omega)}$, $\|f\|_{L^{n}(\Omega)}$, $\|g\|_{C^{0,1}(\bar\Gamma)}$ and $\|u\|_{L^{\infty }(\Omega)}$.

In particular, if $F$ satisfies \cref{SC1},
\begin{equation}\label{e.C0ln-1-i-global}
 \|u\|_{C^{0,\mathrm{lnL}}(\bar{\Omega}')}\leq C \left(\|u\|_{L^{\infty }(\Omega)}+\|f\|_{L^{n}(\Omega)}+\|g\|_{C^{0,1}(\bar\Gamma)}\right),
\end{equation}
where $C$ depends only on $n,\lambda,\Lambda,p,r_0,\|b\|_{L^p(\Omega)}$ and $\|c\|_{L^{n}(\Omega)}$.
\end{corollary}

\begin{corollary}\label{t-C1ln-global}
Let $\Gamma\subset \partial \Omega $ be relatively open and $u$ be a viscosity solution of
\begin{equation*}
\left\{\begin{aligned}
&F(D^2u,Du,u,x)=f&& \quad\mbox{in}~~\Omega;\\
&u=g&& \quad\mbox{on}~~\Gamma.
\end{aligned}\right.
\end{equation*}
Suppose that \cref{SC2} holds and $F$ satisfies \cref{e.C1a.beta} with $G_{x_0}$ at any $x_0\in \Omega\cup \Gamma$, where $G_{x_0}$ is convex in $M$. Assume that
\begin{equation*}
  \begin{aligned}
&\beta_1(x,x_0)\leq \frac{\delta_0}{|\ln|x-x_0||},~\forall ~x_0,x\in \Omega\cup \Gamma
~\mbox{ with }~|x-x_0|<r_0,\quad\gamma_1\in L^{\infty}(\Omega),\\
&f\in L^{\infty}(\Omega),\quad\Gamma\in C^{1,1},\quad g\in C^{1,1}(\bar\Gamma),
  \end{aligned}
\end{equation*}
where $0<\delta_0<1$ depends only on $n,\lambda$ and $\Lambda$.

Then for any $\Omega'\subset\subset \Omega\cup \Gamma$, we have $u\in C^{1,\mathrm{lnL}}(\bar{\Omega}')$ and
\begin{equation}\label{e.C1ln-1-i-mu-global}
 \|u\|_{C^{1,\mathrm{lnL}}(\bar{\Omega}')}\leq C,
\end{equation}
where $C$ depends only on $n,\lambda,\Lambda,r_0,\mu,b_0,c_0,\omega_0,
\|\gamma_1\|_{L^{\infty}(\Omega)}, \|\partial \Omega'\cap\Gamma\|_{C^{1,1}}, \|f\|_{L^{\infty}(\Omega)}$, $\|g\|_{C^{1,1}(\bar\Gamma)}$ and $\|u\|_{L^{\infty }(\Omega)}$.

In particular, if $F$ satisfies \cref{SC1},
\begin{equation}\label{e.C1ln-1-i-global}
 \|u\|_{C^{1,\mathrm{lnL}}(\bar{\Omega}')}\leq C \left(\|u\|_{L^{\infty }(\Omega)}+\|f\|_{L^{\infty}(\Omega)}+\|g\|_{C^{1,1}(\bar\Gamma)}\right),
\end{equation}
where $C$ depends only on $n,\lambda,\Lambda,r_0,b_0,c_0$ and $\|\partial \Omega'\cap\Gamma\|_{C^{1,1}}$.
\end{corollary}

\begin{corollary}\label{t-Ckln-global}
Let $k\geq 2$, $\Gamma\subset \partial \Omega $ be relatively open and $u$ be a viscosity solution of
\begin{equation*}
\left\{\begin{aligned}
&F(D^2u,Du,u,x)=f&& \quad\mbox{in}~~\Omega;\\
&u=g&& \quad\mbox{on}~~\Gamma.
\end{aligned}\right.
\end{equation*}
Suppose that \cref{SC1} and \cref{e.omega0-2} hold, and
\begin{equation*}
F\in C^{k-2,1}(\bar \Omega),\quad f\in C^{k-2,1}(\bar\Omega),\quad \Gamma\in C^{k,1},\quad g\in C^{k,1}(\bar\Gamma).
\end{equation*}

Then for any $\Omega'\subset\subset \Omega\cup \Gamma$, we have $u\in C^{k,\mathrm{lnL}}(\bar{\Omega}')$ and
\begin{equation}\label{e.Ckln-1-i-mu-global}
 \|u\|_{C^{k,\mathrm{lnL}}(\bar{\Omega}')}\leq C,
\end{equation}
where $C$ depends only on $n,\lambda, \Lambda,r_0,\mu, b_0,c_0,\omega_0,
\|F\|_{C^{k-2,1}(\bar\Omega)},\omega_3,\omega_4$, $\|\partial \Omega'\cap\Gamma\|_{C^{k,1}}$,
$\|f\|_{C^{k-2,1}(\bar\Omega)}$, $\|g\|_{C^{k,1}(\bar{\Gamma})}$ and $\|u\|_{L^{\infty}(\Omega)}$.
\end{corollary}
~\\

In the rest of this section, we prove the interior $C^{k,\mathrm{lnL}}$ ($k\geq 2$) regularity and the boundary $C^{1,\mathrm{lnL}}$ regularity. The following lemma is the ``key step'' for the interior $C^{k,\mathrm{lnL}}$ regularity and we omit its proof.
\begin{lemma}\label{In-l-Ckln-mu}
Suppose that $F\in C^{k-2,1}(0)$ and $\omega_0$ satisfies \cref{e.omega0-2}. Then there exists $\delta>0$ depending only on $k,n,\lambda,\Lambda,\omega_0,K_1,\omega_3$ and $\omega_4$ such that if $u$ satisfies
\begin{equation*}
F(D^2u,Du,u,x)=f \quad\mbox{in}~~B_1
\end{equation*}
with
\begin{equation*}
  \begin{aligned}
&u(0)=|Du(0)|=\cdots=|D^{k}u(0)|=0,\quad \max\left(\|u\|_{L^{\infty}(B_1)},\mu,b_0,c_0\right)\leq 1,~\\
&\max\left(\|F\|_{C^{k-2,1}(0)},\|f\|_{C^{k-2,1}(0)}\right)\leq \delta,
  \end{aligned}
\end{equation*}
then there exists $P\in\mathcal{HP}_{k+1}$ such that
\begin{equation*}
  \|u-P\|_{L^{\infty}(B_{\eta})}\leq \eta^{k+1},
\end{equation*}
\begin{equation*}
  |G(D^2P(x),DP(x),P(x),x)|\leq C |x|^{k-1+\bar{\alpha}},~\forall ~x\in B_1
\end{equation*}
and
\begin{equation*}
\|P\|\leq C,
\end{equation*}
where $C$ and $\eta$ depend only on $k,n,\lambda, \Lambda,\omega_0,K_1,\omega_3$ and $\omega_4$.
\end{lemma}

\begin{lemma}\label{In-t-Cklns-mu}
Suppose that $F\in C^{k-2,1}(0)$ and $\omega_0$ satisfies \cref{e.omega0-2}. Let $u$ satisfy
\begin{equation*}
F(D^2u,Du,u,x)=f \quad\mbox{in}~~B_1.
\end{equation*}
Assume that
\begin{equation}\label{In-e.tCklns-be-mu}
\begin{aligned}
&\|u\|_{L^{\infty}(B_1)}\leq 1,\quad u(0)=|Du(0)|=\cdots|D^ku(0)|=0,\\
&\mu\leq \frac{1}{4C_0},\quad b_0\leq\frac{1}{2},\quad c_0\leq \frac{1}{K_0},\\
&\|F\|_{C^{k-2,\alpha}(0)}\leq \frac{\delta_1}{C_0},\quad |f(x)|\leq \delta_1|x|^{k-1}, ~\forall ~x\in B_1,\\
\end{aligned}
\end{equation}
where $\delta_1\leq \delta$ ($\delta$ is as in \Cref{In-l-Ckln-mu}) and $C_0$ depend only on $k,n,\lambda, \Lambda,\omega_0,K_1,\omega_3$ and $\omega_4$.

Then $u\in C^{k,\mathrm{lnL}}(0)$ and
\begin{equation}\label{In-e.tCklns-1-mu}
  |u(x)|\leq C |x|^{k+1}|\ln|x||, ~~\forall ~x\in B_{1/2},
\end{equation}
where $C$ depends only on $k,n,\lambda, \Lambda,\omega_0,K_1,\omega_3$ and $\omega_4$.
\end{lemma}

\proof By \Cref{le2.4}, to prove that $u\in C^{k,\mathrm{lnL}}(0)$, we only need to prove the following. There exist a sequence of $P_m\in\mathcal{HP}_{k+1}$ ($m\geq 0$ and $P_0\equiv 0$) such that for all $m\geq 1$,

\begin{equation}\label{In-e.tCklns-6-mu}
\|u-P_m\|_{L^{\infty }(B_{\eta^{m}})}\leq \eta ^{m(k+1)},
\end{equation}
\begin{equation}\label{In-e.tCklns-9-mu}
|G(D^2P_m(x),DP_m(x),P_m(x),x)|\leq \tilde{C}|x|^{k-1+\bar{\alpha}},~\forall ~x\in B_1
\end{equation}
and
\begin{equation}\label{In-e.tCklns-7-mu}
\|P_m-P_{m-1}\|\leq \tilde{C},
\end{equation}
where $\tilde{C}$ and $\eta$ depends only on $n,\lambda, \Lambda,\omega_0,K_1,\omega_3$ and $\omega_4$.

We prove \crefrange{In-e.tCklns-6-mu}{In-e.tCklns-7-mu} by induction. For $m=1$, by \Cref{In-l-Ckln-mu}, there exists $P_1\in\mathcal{HP}_{k+1}$ such that \crefrange{In-e.tCklns-6-mu}{In-e.tCklns-9-mu} hold for some $C_1$ and $\eta_1$ depending only on $k,n,\lambda, \Lambda,\omega_0,\omega_3$ and $\omega_4$ where $P_0\equiv 0$. Take
$\tilde{C}\geq C_1$, $\eta\leq \eta_1$ and then the conclusion holds for $m=1$. Suppose that the conclusion holds for $m\leq m_0$. We need to prove that the conclusion holds for $m=m_0+1$.

Let $r=\eta ^{m_{0}}$, $y=x/r$ and
\begin{equation}\label{In-e.tCklns-v-mu}
  v(y)=\frac{u(x)-P_{m_0}(x)}{r^{k+1}}.
\end{equation}
Then $v$ satisfies
\begin{equation}\label{In-e.Cklns-F-mu}
\tilde{F}(D^2v,Dv,v,y)=\tilde{f} \quad\mbox{in}~~B_1,
\end{equation}
where for $(M,p,s,y)\in \mathcal{S}^n\times \mathbb{R}^n\times \mathbb{R}\times \bar B_1$,
\begin{equation*}
  \begin{aligned}
&\tilde{F}(M,p,s,y)=r^{-(k-1)}F(r^{k-1}M+D^2P_{m_0}(x),r^kp+DP_{m_0}(x),
    r^{k+1}s+P_{m_0}(x),x),\\
&\tilde{f}(y)=r^{-(k-1)}f(x),
  \end{aligned}
\end{equation*}
In addition, define $\tilde{G}$ in a similar way to the definition of $\tilde{F}$.

In the following, we show that \cref{In-e.Cklns-F-mu} satisfies the assumptions of \Cref{In-l-Ckln-mu}. First, it is easy to verify that
\begin{equation*}
\begin{aligned}
\|v\|_{L^{\infty}(B_1)}\leq& 1,\quad v(0)=\cdots=|D^kv(0)|=0,~\quad(\mathrm{by}~ \cref{In-e.tCklns-be-mu},~\cref{In-e.tCklns-6-mu}~\mbox{and}~\cref{In-e.tCklns-v-mu})\\
|\tilde{f}(y)|\leq &r^{-(k-1)}|f(x)|\leq \delta_1|y|^{k-1}, ~\forall ~y\in B_1. ~\quad(\mathrm{by}~ \cref{In-e.tCklns-be-mu})
  \end{aligned}
\end{equation*}

By \cref{In-e.tCklns-7-mu},
\begin{equation*}
\|P_m\|\leq m C_{0},~\forall ~0\leq m\leq m_0,
\end{equation*}
where $C_{0}$ depends only on $k,n,\lambda, \Lambda,\omega_0,K_1,\omega_3$ and $\omega_4$. It is easy to check that $\tilde{F}$ and $\tilde{G}$ satisfy the structure condition \cref{SC2} with $\lambda,\Lambda,\tilde{\mu},\tilde{b},\tilde{c}$ and $\tilde{\omega}_0$, where (note that $\eta^m m\leq 1$, $\forall~ m\geq 0$ and $\omega_0$ satisfies \cref{e.omega0-2})
\begin{equation*}
\tilde{\mu}=r^{k+1}\mu,\quad\tilde{b}= rb_0+2m_0C_0 r\mu ,\quad\tilde{c}= K_0r^{2}c_0, \quad \tilde{\omega}_0(\cdot,\cdot)=\omega_0(\cdot+C_0,\cdot).
\end{equation*}
Hence, $\tilde{\omega}_0$ satisfies \cref{e.omega0-2} and from \cref{In-e.tCklns-be-mu},
\begin{equation*}
\begin{aligned}
&\tilde{\mu}\leq \mu\leq 1,\quad\tilde{b}\leq b_0+2C_0\mu\leq 1,\quad\tilde{c}\leq c_0\leq 1.
\end{aligned}
\end{equation*}

In addition, by \cref{In-e.tCklns-be-mu} and a computation as before (cf. \cref{e.6.1}), for $(M,p,s,y)\in \mathcal{S}^n\times \mathbb{R}^n\times \mathbb{R}\times \bar B_1$,
\begin{equation*}
\begin{aligned}
|\tilde{F}(M,p,s,y)-\tilde{G}(M,p,s,y)|&\leq \delta_1|y|^{k-1}(|M|+1)\omega_3(|p|+C_0,|s|+C_0)\\
&\coloneqq \delta_1|y|^{k-1}(|M|+1)\tilde{\omega}_3(|p|,|s|).
\end{aligned}
\end{equation*}
Hence, $\|\tilde{F}\|_{C^{k-2,1}(0)}\leq \delta_1$.

Finally, with the aid of \cref{In-e.tCklns-9-mu}, we can show that (similar to the interior $C^{k,\alpha}$ regularity) $\tilde{G}$ satisfies (i)-(iii) of \Cref{d-FP} with some $\tilde{K}_1$ and
\begin{equation}\label{In-e.Ckln-6}
\|\tilde{G}\|_{C^{k-1,\bar{\alpha}}(\bar{\mathbf{B}}_{\rho}\times \bar{B}_1)}\leq \tilde{\omega}_4(\rho),~\forall ~\rho>0,
\end{equation}
where $\tilde{\omega}_4$ depends only on $k,n,\lambda, \Lambda,\omega_0,K_1,\omega_3$ and $\omega_4$.

Choose $\delta_1$ small enough (depending only on $k,n,\lambda, \Lambda,\omega_0,K_1,\omega_3$ and $\omega_4$) such that \Cref{In-l-Ckln-mu} holds for $\tilde\omega_0,\tilde{K}_1,\tilde\omega_3,\tilde\omega_4$ and $\delta_1$. Since \cref{In-e.Cklns-F-mu} satisfies the assumptions of \Cref{In-l-Ckln-mu}, there exist $\tilde{P}\in\mathcal{HP}_{k+1}$ and constants $\tilde{C}\geq C_1$ and $\eta\leq \eta_1$ depending only on $k,n,\lambda, \Lambda,\omega_0,K_1,\omega_3$ and $\omega_4$ such that
\begin{equation*}
\begin{aligned}
    \|v-\tilde{P}\|_{L^{\infty }(B_{\eta})}&\leq \eta^{k+1},
\end{aligned}
\end{equation*}
\begin{equation*}
  |\tilde G(D^2\tilde{P}(y),D\tilde{P}(y),\tilde{P}(y),y)|\leq \tilde{C} |y|^{k-1+\bar{\alpha}},
  ~\forall ~y\in B_1
\end{equation*}
and
\begin{equation*}
\|\tilde{P}\|\leq \tilde{C}.
\end{equation*}

Let
\begin{equation*}
P_{m_0+1}(x)=P_{m_0}(x)+r^{k+1}\tilde{P}(y)=P_{m_0}(x)+\tilde{P}(x).
\end{equation*}
Then \cref{In-e.tCklns-9-mu} and \cref{In-e.tCklns-7-mu} hold for $k_0+1$. By recalling \cref{In-e.tCklns-v-mu}, we have
\begin{equation*}
  \begin{aligned}
\|u-P_{m_0+1}\|_{L^{\infty}(B_{\eta^{m_0+1}})}&= \|u-P_{m_0}-\tilde{P}(x)\|_{L^{\infty}(B_{\eta r})}\\
&= \|r^{k+1}v-r^{k+1}\tilde{P}(y)\|_{L^{\infty}(B_{\eta})}\\
&\leq r^{k+1}\eta^{k+1}=\eta^{(m_0+1)(k+1)}.
  \end{aligned}
\end{equation*}
Hence, \cref{In-e.tCklns-6-mu} hold for $m=m_0+1$. By induction, the proof is completed.\qed~\\

Now, we give the~\\
\noindent\textbf{Proof of \Cref{t-Ckln-i}.} As before, in the following proof,  we just make necessary normalization to satisfy the conditions of \Cref{In-t-Cklns-mu}. Throughout this proof, $C$ always denotes a constant depending only on $k,n,\lambda, \Lambda,\mu, b_0,c_0,\omega_0$, $\|F\|_{C^{k-2,1}(0)}$, $K_1,\omega_3,\omega_4, \|f\|_{C^{k-2,1}(0)}$ and $\|u\|_{L^{\infty}(B_1)}$.

For $(M,p,s,x)\in \mathcal{S}^n\times \mathbb{R}^n\times \mathbb{R}\times \bar{B}_1$, let
\begin{equation*}
F_1(M,p,s,x)=F(M,p,s,x)-P_f(x), \quad f_1=f-P_f.
\end{equation*}
Then $u$ satisfies
\begin{equation*}
F_1(D^2u,Du,u,x)=f_1 \quad\mbox{in}~~B_1
\end{equation*}
and
\begin{equation*}
  |f_1(x)|\leq [f]_{C^{k-2,1}(0)}|x|^{k-1}\leq C|x|^{k-1}, ~~\forall ~x\in B_1.
\end{equation*}

Note that $u\in C^{k,\bar{\alpha}/2}(0)$ (by \Cref{t-Cka-i}). We define
\begin{equation*}
u_1=u-P_u, \quad F_2(M,p,s,x)=F_1(M+D^2P_u(x),p+DP_u(x),s+P_u(x),x).
\end{equation*}
Then $u_1$ satisfies
\begin{equation*}
F_2(D^2u_1,Du_1,u_1,x)=f_1 \quad\mbox{in}~~B_1
\end{equation*}
and
\begin{equation*}
  u_1(0)=|Du_1(0)|=\cdots |D^ku_1(0)|=0.
\end{equation*}

Next, take $y=x/\rho$ and $u_2(y)=u_1(x)/\rho^2$, where $0<\rho<1$ is a constant to be specified later. Then $u_2$ satisfies
\begin{equation}\label{In-Ckln-F6-k-mu}
F_3(D^2u_2,Du_2,u_2,y)=f_2 \quad\mbox{in}~~B_1,
\end{equation}
where
\begin{equation*}
  \begin{aligned}
&F_3(M,p,s,y)=F_2\left(M,\rho\, p,\rho^2s,\rho y\right), \quad f_2(y)=f_1(x).
\end{aligned}
\end{equation*}
Finally, define fully nonlinear operators $G_1,G_2,G_{3}$ in the same way as $F_1,F_2,F_3$.

Now, we choose a proper $\rho$ such that \cref{In-Ckln-F6-k-mu} satisfies the conditions of \Cref{In-t-Cklns-mu}. First, $u_2(0)=\cdots=|D^ku_2(0)|=0$ clearly. By combining with the $C^{k,\bar{\alpha}/2}$ regularity for $u$, we have
\begin{equation*}
  \begin{aligned}
    \|u_2\|_{L^{\infty}(B_1)}&= \rho^{-2}\|u_1\|_{L^{\infty}(B_\rho)}
    \leq C\rho^{k-2+\bar{\alpha}/2}.
  \end{aligned}
\end{equation*}
Next,
\begin{equation*}
|f_2(y)|=|f_1(x)|\leq C|x|^{k-1}=C\rho^{k-1}|y|^{k-1},~\forall ~y\in B_1,
\end{equation*}

It is easy to verify that $F_3$ and $G_{3}$ satisfy the structure condition \cref{SC2} with $\lambda,\Lambda,\tilde\mu,\tilde{b},\tilde{c}$ and $\tilde{\omega}_0$, where
\begin{equation*}
\tilde{\mu}= \rho^2\mu,\quad\tilde{b}= \rho b_0+C\rho\mu,\quad\tilde{c}= K_0\rho^2c_0,\quad
\tilde{\omega}_0(\cdot,\cdot)=\omega_0(\cdot+C,\cdot).
\end{equation*}
In addition, $\tilde{\omega}_0$ satisfies \cref{e.omega0-2}.

Finally, we show $F_3\in C^{k-2,1}(0)$. Indeed, by a calculation as before (cf. \cref{e5.4}),
\begin{equation*}
  \begin{aligned}
&|F_3(M,p,s,y)-G_{3}(M,p,s,y)|\leq C\rho^{k-1}|y|^{k-1}(|M|+1)\tilde\omega_3(|p|,|s|),
\end{aligned}
\end{equation*}
where
\begin{equation*}
\tilde\omega_3(|p|,|s|)\coloneqq\omega_3(|p|+C,|s|+C).
\end{equation*}
Moreover, it can be verified that $G_{3}$ satisfies (i)-(iii) of \Cref{d-FP} with some $\tilde{K}_1$ and
\begin{equation*}
\|G_{3}\|_{C^{k-1,\bar{\alpha}}(\bar{\mathbf{B}}_r\times \bar{B}_1)}\leq \tilde{\omega}_4(r),~\forall ~r>0,
\end{equation*}
where $\tilde{\omega}_4$ depends only on $k,n,\lambda, \Lambda,\mu, b_0,c_0,\omega_0$, $\|F\|_{C^{k-2,1}(0)}$, $K_1,\omega_3,\omega_4$, $\|f\|_{C^{k-2,1}(0)}$, and $\|u\|_{L^{\infty}(B_1)}$.

Take $\delta_1$ small enough such that \Cref{In-t-Cklns-mu} holds for $\tilde{\omega}_0,K_1,\tilde{\omega}_3,\tilde{\omega}_4$ and $\delta_1$.
From above arguments, we can choose $\rho$ small enough (depending only on $k,n,\lambda, \Lambda,\mu$, $b_0,c_0$, $\omega_0$, $\|F\|_{C^{k-2,1}(0)}$, $K_1,\omega_3,\omega_4, \|f\|_{C^{k-2,1}(0)}$ and $\|u\|_{L^{\infty}(B_1)}$) such that the conditions of \Cref{In-t-Cklns-mu} are satisfied. Then $u_2$ and hence $u$ is $C^{k,\mathrm{lnL}}$ at $0$, and the estimates \crefrange{e.Ckln-1-i}{e.Ckln-2-i} hold. \qed~\\

In the following, we give the proof of the boundary $C^{1,\mathrm{lnL}}$ regularity. The following is the corresponding ``key step'' and we omit its proof.
\begin{lemma}\label{l-C1ln-mu}
Suppose that $F$ satisfies \cref{e.C1a.beta}. There exists $\delta>0$ depending only on $n,\lambda$ and $\Lambda$ such that if $u$ satisfies
\begin{equation*}
\left\{\begin{aligned}
&F(D^2u,Du,u,x)=f&& \quad\mbox{in}~~\Omega_1;\\
&u=g&& \quad\mbox{on}~~(\partial \Omega)_1
\end{aligned}\right.
\end{equation*}
with
\begin{equation*}
  \begin{aligned}
&u(0)=|Du(0)|=0,\quad \max\left(\|u\|_{L^{\infty}(\Omega_1)},\omega_0(1,1)\right)\leq 1,\\
&\max\left(\mu,b_0,c_0,\|\beta_1\|_{L^{n}(\Omega_1)}, \|f\|_{L^{\infty}(\Omega_1)},
\|(\partial \Omega)_1\|_{C^{1,1}(0)}, \|g\|_{C^{1,1}(0)} \right)\leq \delta,
  \end{aligned}
\end{equation*}
then there exists $P\in\mathcal{SP}_{2}$ such that
\begin{equation*}
  \|u-P\|_{L^{\infty}(\Omega_{\eta})}\leq \eta^{2},
\end{equation*}
\begin{equation*}
G(D^2P,0,0)=0
\end{equation*}
and
\begin{equation*}
\|P\|\leq \bar{C}+1,
\end{equation*}
where $0<\eta<1$ depends only on $n,\lambda$ and $\Lambda$.
\end{lemma}

\begin{lemma}\label{t-C1lns}
Suppose that $F$ satisfies \cref{e.C1a.beta} and $u$ satisfies
\begin{equation*}
\left\{\begin{aligned}
&F(D^2u,Du,u,x)=f&& \quad\mbox{in}~~\Omega_1;\\
&u=g&& \quad\mbox{on}~~(\partial \Omega)_1.
\end{aligned}\right.
\end{equation*}
Assume that
\begin{equation}\label{e.C1ln-ubeg}
\begin{aligned}
&\|u\|_{L^{\infty}(\Omega_1)}\leq 1,\quad  u(0)=|Du(0)|=0,\quad \mu\leq \frac{\delta_1}{4C_0^2},\quad  b_0\leq \frac{\delta_1}{4C_0},\quad c_0\leq \frac{\delta_1}{4},\\
&\omega_0(1+C_0,C_0)\leq1,\quad\|\beta_1\|_{L^{n}(\Omega_r)}\leq \frac{\delta_1 r}{C_0|\ln r|}, ~\forall ~0<r<1/2,\quad\|\gamma_1\|_{L^{\infty}(\Omega_1)}\leq \frac{\delta_1}{4}, \\ &\|f\|_{L^{\infty}(\Omega_1)}\leq \delta_1,\quad
|g(x)|\leq \frac{\delta_1}{2}|x|^2,~\forall ~x\in(\partial \Omega)_1,\quad \|(\partial \Omega)_1\|_{C^{1,1}(0)} \leq \frac{\delta_1}{2C_0},\\
\end{aligned}
\end{equation}
where $\delta_1$ and $C_0$ depend only on $n,\lambda$ and $\Lambda$.

Then $u$ is $C^{1,\mathrm{lnL}}$ at $0$, i.e.,
\begin{equation*}
  |u(x)|\leq C |x|^2\big|\ln |x|\big|, ~~\forall ~x\in \Omega_{1/2},
\end{equation*}
where $C$ depends only on $n,\lambda$ and $\Lambda$.
\end{lemma}

\proof  To prove that $u$ is $C^{1,\mathrm{lnL}}$ at $0$, we only need to prove the following. There exist a sequence of $P_m\in\mathcal{SP}_{2}$ ($m\geq -1$) such that for all $m\geq 0$,

\begin{equation}\label{e.C1ln-u}
\|u-P_m\|_{L^{\infty }(\Omega_{\eta^{m}})}\leq \eta ^{2m},
\end{equation}
\begin{equation}\label{e.C1lns-F}
  G(D^2P_m)=0
\end{equation}
and
\begin{equation}\label{e.C1ln-bk}
\|P_m-P_{m-1}\|\leq \bar{C}+1,
\end{equation}
where $\eta$ depending only on $n,\lambda$ and $\Lambda$, is as in \Cref{l-C1ln-mu} .

Now we prove \crefrange{e.C1ln-u}{e.C1ln-bk} by induction. For $m=0$, by setting $P_0=P_{-1}\equiv 0$, the conclusion holds clearly. Suppose that the conclusion holds for $m\leq m_0$. We need to prove that the conclusion holds for $m=m_0+1$.

Let $r=\eta ^{m_{0}}$, $y=x/r$ and
\begin{equation}\label{e.C1ln-v1}
  v(y)=\frac{u(x)-P_{m_0}(x)}{r^2}.
\end{equation}
Then $v$ satisfies
\begin{equation}\label{e.C1lns-F-mu}
\left\{\begin{aligned}
&\tilde{F}(D^2v,Dv,v,y)=\tilde{f}&& \quad\mbox{in}~~\tilde{\Omega}_1;\\
&v=\tilde{g}&& \quad\mbox{on}~~(\partial \tilde{\Omega})_1,
\end{aligned}\right.
\end{equation}
where for $(M,p,s,x)\in \mathcal{S}^n\times \mathbb{R}^n\times \mathbb{R}\times \bar{\tilde{\Omega}}_1$,
\begin{equation*}
  \begin{aligned}
&\tilde{F}(M,p,s,y)=F(M+D^2P_{m_0},rp+DP_{m_0}(x),r^2s+P_{m_0}(x),x),\\
&\tilde{f}(y)=f(x),\quad\tilde{g}(y)=r^{-2}\left(g(x)-P_{m_0}(x)\right),\quad\tilde{\Omega}=r^{-1}\Omega.
  \end{aligned}
\end{equation*}
In addition, define
\begin{equation*}
\tilde{G}(M)=G(M+D^2P_{m_0}).
\end{equation*}

In the following, we show that \cref{e.C1lns-F-mu} satisfies the assumptions of \Cref{l-C1ln-mu}. First, it is easy to verify that
\begin{equation*}
\begin{aligned}
\|v\|_{L^{\infty}(\tilde{\Omega}_1)}\leq& 1,\quad v(0)=|Dv(0)|=0,
~\quad(\mathrm{by}~\cref{e.C1ln-ubeg},~ \cref{e.C1ln-u} ~\mbox{and}~ \cref{e.C1ln-v1})\\
\|\tilde{f}\|_{L^{\infty}(\tilde{\Omega}_1)}=&\|f\|_{L^{\infty}(\Omega_r)}\leq \delta_1,
 ~\quad(\mathrm{by}~\cref{e.C1ln-ubeg})\\
 \|(\partial \tilde{\Omega})_1\|_{C^{1,1}(0)} \leq& r\|(\partial\Omega)_1\|_{C^{1,1}(0)}\leq \delta_1,
~\quad(\mathrm{by}~\cref{e.C1ln-ubeg})\\
\tilde{G}(0)=&G(D^2P_{m_0})=0.~\quad(\mathrm{by}~ \cref{e.C1lns-F})
\end{aligned}
\end{equation*}

By \cref{e.C1ln-bk},
\begin{equation*}
\|P_m\|\leq m C_{0},~\forall ~0\leq m\leq m_0,
\end{equation*}
where $C_{0}$ depends only on $n,\lambda$ and $\Lambda$. Thus, by noting that $P_{m_0}\in\mathcal{SP}_{2}$ and $m\eta^m\leq 1$ for any $m\geq 1$,
\begin{equation*}
\|\tilde{g}\|_{L^{\infty}((\partial \tilde{\Omega})_{\rho})} \leq \frac{1}{r^{2}}\left(\frac{\delta _1}{2}(\rho r)^{2}+
m_0C_0\cdot \frac{\delta_1}{2C_0}(\rho r)^{3}\right)\leq \delta_1 \rho^{2},~\forall ~0<\rho<1.
~~(\mathrm{by}~\cref{e.C1ln-ubeg})
\end{equation*}
Hence,
\begin{equation*}
  \|\tilde{g}\|_{C^{1,1}(0)}\leq \delta_1.
\end{equation*}

As before, it is easy to verify that $\tilde{F}$ and $\tilde{G}$ satisfy the structure condition \cref{SC2} with $\lambda,\Lambda,\tilde{\mu},\tilde{b},\tilde{c}$ and $\tilde{\omega}_0$, where (use $m\eta^m\leq 1$ again)
\begin{equation*}
\tilde{\mu}=r^2\mu, \quad\tilde{b}= rb_0+2m_0C_0r^2\mu,\quad \tilde{c}= c_0, \quad\tilde{\omega}_0(\cdot,\cdot)=\omega_0(\cdot+C_0,\cdot).
\end{equation*}
Therefore,
\begin{equation*}
\begin{aligned}
&\tilde{\mu}\leq \mu\leq \delta_1,\quad \tilde{b}\leq b_0+2C_0\mu\leq \delta_1,\quad
\tilde{c}\leq c_0\leq \delta_1.
 \end{aligned}
\end{equation*}

Moreover, by combining the structure condition \cref{SC2}, \cref{e.C1a.beta} and \cref{e.C1ln-ubeg}, for $(M,y)\in \mathcal{S}^n\times \bar{\tilde{\Omega}}_1$ (use $m\eta^m\leq 1$ again),
\begin{equation*}
\begin{aligned}
|\tilde{F}&(M,0,0,y)-\tilde{G}(M)|\\
&= |F(M+D^2P_{m_0},DP_{m_0}(x),P_{m_0}(x),x)-G(M+D^2P_{m_0})|\\
&\leq |F(M+D^2P_{m_0},DP_{m_0}(x),P_{m_0}(x),x)-F(M+D^2P_{m_0},0,0,x)|\\
&~~~~+|F(M+D^2P_{m_0},0,0,x)-G(M+D^2P_{m_0})|\\
&\leq C_0^2\mu+C_0b_0+c_0\omega_0(C_0,C_0)+\beta_1(x)(|M|+m_0C_0)+\gamma_1(x)\\
&\leq m_0C_0\beta_1(x)(|M|+1)+\delta_1\\
&\coloneqq \tilde{\beta}_1(y)(|M|+1)+\tilde{\gamma}_1(y),
\end{aligned}
\end{equation*}
where
\begin{equation*}
\tilde{\beta}_1(y)=m_0C_0\beta_1(x),\quad\tilde{\gamma}_1(y)\equiv \delta_1.
\end{equation*}
Then
\begin{equation*}
\|\tilde{\beta}_1\|_{L^{n}(\tilde\Omega_1)}=\frac{m_0C_0}{r}\|\beta_1\|_{L^{n}(\Omega_r)}\leq \frac{m_0\delta_1}{|\ln r|}= \frac{m_0\delta_1}{m_0|\ln \eta|}\leq\delta_1.
\end{equation*}

Choose $\delta_1$ small enough (depending only on $n,\lambda$ and $\Lambda$) such that \Cref{l-C1ln-mu} holds for $\tilde\omega_0$ and $\delta_1$. Since \cref{e.C1lns-F-mu} satisfies the assumptions of \Cref{l-C1ln-mu}, there exists $\tilde{P}(y)\in\mathcal{SP}_{2}$ such that
\begin{equation*}
\begin{aligned}
    \|v-\tilde{P}\|_{L^{\infty }(\tilde{\Omega}_{\eta})}&\leq \eta^2,
\end{aligned}
\end{equation*}
\begin{equation*}
  \tilde{G}(D^2\tilde{P})=0
\end{equation*}
and
\begin{equation*}
\|\tilde{P}\|\leq \bar{C}+1.
\end{equation*}

Let
\begin{equation*}
P_{m_0+1}(x)=P_{m_0}(x)+r^2\tilde{P}(y)=P_{m_0}(x)+\tilde{P}(x).
\end{equation*}
Then \cref{e.C1lns-F} and \cref{e.C1ln-bk} hold for $m_0+1$. By recalling \cref{e.C1ln-v1}, we have
\begin{equation*}
  \begin{aligned}
\|u-P_{m_0+1}\|_{L^{\infty}(\Omega_{\eta^{m_0+1}})}
&= \|u-P_{m_0}-\tilde{P}(x)\|_{L^{\infty}(\Omega_{\eta r})}\\
&= \|r^{2}v-r^{2}\tilde{P}(y)\|_{L^{\infty}(\tilde{\Omega}_{\eta})}\\
&\leq r^{2}\eta^{2}= \eta^{2(m_0+1)}.
  \end{aligned}
\end{equation*}
Hence, \cref{e.C1ln-u} holds for $m=m_0+1$. By induction, the proof is completed.\qed~\\

Now, we give the ~\\
\noindent\textbf{Proof of \Cref{t-C1ln}.} As before, we prove the theorem in two cases.

\textbf{Case 1:} the general case, i.e., $F$ satisfies \cref{SC2} and \cref{e.C1a.beta}. Throughout the proof for this case, $C$ always denotes a constant depending only on $n,\lambda,\Lambda,\mu,b_0, c_0,\omega_0$, $\|(\partial \Omega)_1\|_{C^{1,1}(0)}$, $\|f\|_{L^{\infty}(\Omega_1)}$, $\|g\|_{C^{1,1}(0)}$ and $\|u\|_{L^{\infty }(\Omega_1)}$.

Let
\begin{equation*}
u_1=u-P_{g}, \quad g_1=g-P_{g}, \quad F_1(M,p,s,x)=F(M,p+DP_{g},s+P_{g}(x),x).
\end{equation*}
Then $u_1$ satisfies
\begin{equation*}
\left\{\begin{aligned}
&F_1(D^2u_1,Du_1,u_1,x)=f&& \quad\mbox{in}~~\Omega_1;\\
&u_1=g_1&& \quad\mbox{on}~~(\partial \Omega)_1.
\end{aligned}\right.
\end{equation*}
Hence,
\begin{equation*}
  |g_1(x)|\leq C|x|^{2}, ~~\forall ~x\in (\partial \Omega)_1
\end{equation*}
and
\begin{equation*}
  \begin{aligned}
  |F_1(0,0,0,x)|=\left|F\left(0,DP_{g},P_{g},x\right)\right|\leq C.
  \end{aligned}
\end{equation*}

Note that
\begin{equation*}
  u_1\in S^{*}(\lambda,\Lambda,\mu,\hat{b},|f|
  +c_0\omega_0(\|u\|_{L^{\infty}(\Omega_1)}+\|g\|_{C^{1,1}(0)},u_1)+|F_1(0,0,0,\cdot)|),
\end{equation*}
where $\hat{b}=b_0+2\mu\|g\|_{C^{1,1}(0)}$. By \Cref{t-C1a-mu}, $u_1\in C^{1,\alpha}(0)$ and
\begin{equation*}
  \begin{aligned}
Du_1(0)=(0,...,0,(u_1)_n(0)), \quad |(u_1)_n(0)| \leq C.
  \end{aligned}
\end{equation*}

Define
\begin{equation*}
u_2=u_1-P_{u_1}=u_1-(u_1)_n(0)x_n, \quad F_2(M,p,s,x)=F_1(M,p+DP_{u_1},s+P_{u_1}(x),x).
\end{equation*}
Then $u_2$ satisfies
\begin{equation*}
\left\{\begin{aligned}
&F_2(D^2u_2,Du_2,u_2,x)=f&& \quad\mbox{in}~~\Omega_1;\\
&u_2=g_2&& \quad\mbox{on}~~(\partial \Omega)_1,
\end{aligned}\right.
\end{equation*}
where $g_2=g_1-P_{u_1}$. Moreover, $u_2(0)=|Du_2(0)|=0$. Since $\partial \Omega$ is $C^{1,1}$ at $0$,
\begin{equation*}
  |g_2(x)|\leq |g_1(x)|+C|x_n|\leq C|x|^{2}, ~~\forall ~x\in (\partial \Omega)_1.
\end{equation*}

Finally, let
\begin{equation*}
y=x/\rho, \quad u_3(y)=u_2(x)/\rho, \quad F_3(M,p,s,y)=\rho F_2(M/\rho, p,\rho s,x).
\end{equation*}
Then $u_3$ satisfies
\begin{equation}\label{F4-Cln-mu}
\left\{\begin{aligned}
&F_3(D^2u_3,Du_3,u_3,y)=f_1&& \quad\mbox{in}~~\tilde{\Omega}_1;\\
&u_3=g_3&& \quad\mbox{on}~~(\partial \tilde\Omega)_1,
\end{aligned}\right.
\end{equation}
where
\begin{equation*}
f_1(y)=\rho f(x), \quad g_3(y)=g_2(x)/\rho, \quad \tilde{\Omega}=\Omega/\rho.
\end{equation*}
Finally, define the fully nonlinear operator
\begin{equation*}
G_{3}(M)= \rho G(\rho ^{-1}M).
\end{equation*}

Now, we can check that \cref{F4-Cln-mu} satisfies the conditions of \Cref{t-C1lns} by choosing a proper $\rho$. First, it can be checked easily that
\begin{equation*}
  \begin{aligned}
u_3(0)=&|Du_3(0)|=0,\quad\|f_1\|_{L^{\infty}(\tilde{\Omega}_1)}= \rho\|f_1\|_{L^{\infty}(\Omega_{\rho})}\leq C\rho,\\
|g_3(y)|=& \rho^{-1}|g_2(x)|\leq C\rho |y|^2,~\forall ~y\in(\partial \tilde{\Omega})_1,\\
\|(\partial \tilde{\Omega})_1\|_{C^{1,1}(0)}\leq& \rho\|(\partial \Omega)_1\|_{C^{1,1}(0)}
\leq C\rho.
  \end{aligned}
\end{equation*}

Next, by the boundary $C^{1,\bar{\alpha}/2}$ regularity for $u_1$,
\begin{equation*}
\|u_3\|_{L^{\infty}(\tilde{\Omega}_1)}= \rho^{-1}\|u_2\|_{L^{\infty}(\Omega_\rho)}
\leq C\rho^{\bar{\alpha}/2}.
\end{equation*}
It can be checked that $F_3$ satisfies the structure condition \cref{SC2} with $\lambda,\Lambda,\tilde{\mu},\tilde{b},\tilde{c}$ and $\tilde{\omega}_0$, where
\begin{equation*}
\tilde{\mu}= \rho\mu,\quad\tilde{b}= \rho b_0+C\rho\mu,\quad\tilde{c}= \rho^{1/2}c_0,
\quad \tilde{\omega}_0(\cdot,\cdot)=\rho^{1/2}\omega_0(\cdot+C,\cdot).
\end{equation*}
Finally, we check the oscillation of $F_3$ in $y$. For simplicity, introduce
\begin{equation*}
P(x)=P_g(x)+P_{u_1}(x).
\end{equation*}
We compute
\begin{equation*}
  \begin{aligned}
|F_3(&M,0,0,y)-G_{3}(M)|\\
=&|\rho F(\rho^{-1}M,DP,P(x),x)-\rho G(\rho^{-1}M)|\\
\leq & |\rho F(\rho^{-1}M,DP,P(x),x)-\rho F(\rho^{-1}M,0,0,x)|
+|\rho F(\rho^{-1}M,0,0,x)-\rho G(\rho^{-1}M)|\\
\leq& \rho(C^2\mu+Cb_0+c_0\omega_0(C,C))+\beta_1(x)|M|+\rho\gamma_1(x)\\
\coloneqq  &\tilde{\beta}_1(y)|M|+\tilde{\gamma}_1(y),
  \end{aligned}
\end{equation*}
where
\begin{equation*}
  \tilde{\beta}_1(y)=\beta_1(x),\quad
  \tilde{\gamma}_1=\rho\gamma_1(x)+\rho(C^2\mu+Cb_0+c_0\omega_0(C_0,C_0)).
\end{equation*}
By the assumption on $\beta$ (see \cref{e11.1}), for any $0<r<1/2$,
\begin{equation*}
\|\tilde{\beta}_1\|_{L^{n}(\tilde\Omega_r)}=\frac{1}{\rho}\|\beta_1\|_{L^{n}(\Omega_{\rho r})}
\leq \frac{1}{\rho}\cdot \frac{\delta_0 \rho r}{|\ln \rho r|}
= \frac{\delta_0r}{|\ln \rho|+|\ln r|}
\leq \frac{\delta_0r}{|\ln r|}.
\end{equation*}

Take $\delta_1$ small enough such that \Cref{t-C1lns} holds with $\tilde{\omega}_0$ and $\delta_1$. From above arguments, we take $\delta_0= \delta_1/C_0$ and $\rho$ small enough (depending only on $n,\lambda,\Lambda,\mu,b_0, c_0,\omega_0$, $\|(\partial \Omega)_1\|_{C^{1,1}(0)}$, $\|f\|_{L^{\infty}(\Omega_1)},\|g\|_{C^{1,1}(0)}$ and $\|u\|_{L^{\infty }(\Omega_1)}$) such that
\begin{equation*}
\begin{aligned}
&\|u_3\|_{L^{\infty}(\Omega_1)}\leq 1,\quad \tilde\mu\leq \frac{\delta_1}{4C_0^2},\quad \tilde b\leq \frac{\delta_1}{4C_0},\quad \tilde c\leq \frac{\delta_1}{4},\quad\tilde\omega_0(1+C_0,C_0)\leq1, \\ &\|\tilde{\beta}_1\|_{L^{n}(\tilde\Omega_r)}\leq \frac{\delta_1r}{C_0|\ln r|}, ~\forall ~0<r<1/2,\quad
\|\tilde\gamma_1\|_{L^{\infty}(\Omega_1)}\leq \frac{\delta_1}{4},\quad \|f_1\|_{L^{\infty}(\Omega_1)}\leq \delta_1,\\
&|g_2(x)|\leq \frac{\delta_1}{2}|x|^2, ~\forall ~x\in(\partial \Omega)_1\quad\mbox{and}\quad\|(\partial \tilde\Omega)_1\|_{C^{1,1}(0)} \leq \frac{\delta_1}{2C_0},\\
\end{aligned}
\end{equation*}
where $C_0$ depending only on $n,\lambda$ and $\Lambda$, is as in \Cref{t-C1lns}. Therefore, the assumptions in \Cref{t-C1lns} are satisfied. By \Cref{t-C1lns}, $u_3$ and hence $u$ is $C^{1,\mathrm{lnL}}$ at $0$, and the estimates \cref{e.C1ln-1} and \cref{e.C1ln-2} hold.

\textbf{Case 2:} $F$ satisfies \cref{SC1}. Let
\begin{equation*}
K=\|u\|_{L^{\infty }(\Omega_1)}+\|f\|_{L^{\infty}(\Omega_1)}+\|\gamma_1\|_{L^{\infty}(\Omega_1)}
+\|g\|_{C^{1,1}(0)}, \quad u_1=u/K.
\end{equation*}
Then $u_1$ satisfies
\begin{equation}\label{e.11.1}
\left\{\begin{aligned}
&F_1(D^2u_1,Du_1,u_1,x)=f_1&& \quad\mbox{in}~~\Omega_1;\\
&u_1=g_1&& \quad\mbox{on}~~(\partial \Omega)_1,
\end{aligned}\right.
\end{equation}
where
\begin{equation*}
F_1(M,p,s,x)=F(KM,Kp,Ks,x)/K, \quad f_1=f/K, \quad g_1=g/K.
\end{equation*}
Then by applying \textbf{Case 1} to \cref{e.11.1}, we obtain that $u_1$ and hence $u$ is $C^{1,\mathrm{lnL}}$ at $0$, and the estimates \cref{e.C1ln-1} and \cref{e.C1ln-2} hold. \qed~\\

\noindent\textbf{Acknowledgement} Part of the work was carried out when the first and the third author worked at Northwestern Polytechnical University (Xi'an) and Shanghai Jiao Tong University (Shanghai). We would like to thank the referee for many valuable comments, which improve our manuscript a lot.

%
\printbibliography
\end{document}